\documentclass[11pt]{article}

\usepackage[a4paper,truedimen,scale={.80,.85},centering]{geometry}
\usepackage{latexsym} 
\usepackage{mathrsfs,amsmath,amssymb,amsthm,amscd,ascmac} 
\usepackage[titletoc,title]{appendix}
\usepackage{verbatim}

{
\theoremstyle{plain}
\newtheorem{thm}{Theorem}[section]

\newtheorem{kmRemark}{\textbf{Remark}}[section]
} 

\newcommand{\rank}{{\rm rank}}

\newcommand{\frakC}{\mathfrak{C}}

\newcommand{\fraksp}{\mathfrak{s}\mathfrak{p}}
\newcommand{\frakham}{\mathfrak{h}\mathfrak{a}\mathfrak{m}}
\newcommand{\Pkt}[2]{\{#1,#2\}}%
\newcommand{\mR}{\ensuremath{\mathbb{R}}} 

\newcommand{\ds}{\displaystyle }

\newcommand{\mydz}{d_{}\, }
\newcommand{\mydo}{d_{0}\, }

\newcommand{\nz}[2]{{z}^{#1}_{#2}}
\newcommand{\W}[2]{{z}^{#1}_{#2}}
\newcommand{\kmqed}{\hfill \rule{1ex}{1.5ex}\par}

\newcommand{\CGF}[4]{\text{C}^{#1}_{\rm GF} ( {\frakham}_{#2}^{#3}, {\fraksp}(#2,\mR))_{#4}}
\newcommand{\HGF}[4]{\text{H}^{#1}_{\rm GF} ( {\frakham}_{#2}^{#3}, {\fraksp}(#2,\mR))_{#4}}

\newcommand{\km}{\ensuremath{ }} 
\allowdisplaybreaks

\begin{document}
\parindent=0pt

\title{\Large Another proof to 
Kotschick-Morita's Theorem of Kontsevich homomorphism}

\author{Kentaro \textsc{Mikami} 
\thanks{Supported by Grant-in-Aid for 
  Scientific Research (C) of JSPS, No.26400063, No.23540067 and No.20540059} 
}




\maketitle

\begin{abstract} In \cite{KOT:MORITA}, Kotschick and Morita showed 
 that the Gel'fand-Kalinin-Fuks class 
 in 
$\ds \HGF{7}{2}{ }{8}$  
 is decomposed as a product $\eta\wedge \omega $
  of some leaf cohomology class $\eta$ and 
  a transverse symplectic class $\omega$. 
In other words, the Kontsevich homomorphism $\ds\omega\wedge : 
\HGF{5}{2}{0}{10} \rightarrow  
\HGF{7}{2}{ }{8}$  is isomorphic.  

In this paper, we give proof for the 
Kotschick and Morita's theorem by using 
 the Gr\"obner Basis theory and  
  computer symbol calculations.  
\end{abstract}

\section{Introduction} On the symplectic space $\ds
(\mR^{2n},\omega)$, let 
$\ds \frakham_{2n}^{ }$ be the Lie algebra of the
formal Hamiltonian vector fields, and let   
$\ds \HGF{\bullet}{2n}{ }{w}$ be the relative 
Gel'fand-Fuks cohomolgy group with the weight $w$.  

When $n=1$,  
Gel'fand-Kalinin-Fuks (\cite{MR0312531}) showed that 
$\displaystyle \HGF{\bullet}{2}{ }{w} = 0$ for the weight $w=2,4,6$ and the 
$\displaystyle \HGF{7}{2}{ }{8} \cong \mR$ whose generator is called    
the Gel'fand-Kalinin-Fuks class.   
The next non-trivial result in this context is 
$\displaystyle \HGF{9}{2}{ }{14} \cong \mR$, which is proved by 
S.~Metoki (\cite{metoki:shinya}) in 1999.

D.~Kotschick and S.~Morita (\cite{KOT:MORITA}) 
studied 
$\ds \HGF{\bullet}{2n}{0}{w} $ and determined the whole space
for $n=1$ and $w\le 10$, where   
$\displaystyle \frakham_{2n}^{0}$ is the Lie subalgebra of the
formal Hamiltonian vector fields which vanish at the origin of  
$\displaystyle \mR^{2n}$.  

There is a natural homomorphism due to Kontsevich (\cite{Kont:RW})  
$$ \omega^n : 
\HGF{\bullet}{2n}{0}{w} \longrightarrow 
\HGF{\bullet+2n}{2n}{}{w-2n}$$

D.~Kotschick and S.~Morita  show in \cite{KOT:MORITA} the next
theorem. 

\begin{thm}[\cite{KOT:MORITA}]\label{thm:morita}\rm 
There is a unique element $\displaystyle \eta \in
\text{H}_{\text{GF}}^{5}(\frakham_{2}^{0})_{10}^{Sp} \cong \mR$ such that 
$$ \text{Gel'fand-Kalinin-Fuks class } = \eta \wedge \omega 
\in \HGF{7}{2}{ }{8}  $$
where $\omega$ is the cochain associated with the 
linear symplectic form of $\displaystyle \mR^2$. 
\end{thm}

About mathematical background, we refer to \cite{KOT:MORITA} or       
a draft 
``{\em An affirmative answer to a conjecture for Metoki class}''
by K.~Mikami(\cite{mikami:affirm}).  For more precise notations or notions in this
paper, we refer to 
\cite{KOT:MORITA},  \cite{M:N:K} or \cite{mikami:affirm}.  

Our aim of this draft is to give another proof of
the theorem above by using Gr\"{o}bner basis theory (cf.\ \cite{GB:J:Merker}
or \cite{mikami:affirm}).     

We use Maple Groebner Package for computing 
Groebner Basis and the normal form.  

There are several symbol
calculus softwares beside Maple, Mathematica, Risa/Asir and so on.
Risa/Asir is popular among Japanese mathematicians because it is bundled in
Math Libre Disk which is distributed at annual meetings of the Mathematical
Society of Japan.
So, 
the author presents the source code and the output about Risa/Asir 
concerning to the Theorem by D.~Kotschick and S.~Morita   
in Appendix XYZ.  
 
You can compare the results by Maple and Risa/Asir  and you will
see  that the both are the same, up to non-zero scalar multiples.

\section{Preliminaries} 
Let $x,y$ be the standard basis of $\ds\mR^2$ with the Poisson bracket is 
$\Pkt{x}{y}=1$. We denote the standard basis of $A$-homogeneous polynomials of
$x$ and $y$ as
$\ds \frac{x^a}{a!} \frac{ y^{A-a}}{(A-a)!}$ and the dual basis is written by 
$\ds \nz{a}{A}$. 

\section{About $\ds\CGF{\bullet}{2}{0}{w}$}
Using the method in 
\cite{M:N:K}, we understand    
the structures of $\displaystyle \CGF{\bullet}{2}{0}{10}$ concretely.  
We denote
$\displaystyle \CGF{\bullet}{2}{0}{10}$  
by $\displaystyle C^{\bullet}$.     
We choose our concrete bases as 
$\displaystyle \{\mathbf{q}_i\}_{i=1}^{9}$ of $C^{4}$, 
$\displaystyle \{\mathbf{w}_i\}_{i=1}^{12}$ of $C^{5}$, and 
$\displaystyle \{\mathbf{r}_i\}_{i=1}^{4}$ of $C^{6}$.   
\subsection{Basis of $\ds \CGF{4}{2}{0}{10}$} 
\begin{scriptsize}
\begin{align*}
\mathbf{q}_{1} &= 
\frac{1}{8} \W{0}{3}\km\W{1}{3}\km\W{0}{4}\km\W{8}{8}
+\frac{3}{4} \W{0}{3}\km\W{1}{3}\km\W{2}{4}\km\W{6}{8}
-\frac{1}{2} \W{0}{3}\km\W{1}{3}\km\W{3}{4}\km\W{5}{8}
-\frac{1}{2} \W{0}{3}\km\W{1}{3}\km\W{1}{4}\km\W{7}{8}
+\frac{1}{8} \W{0}{3}\km\W{1}{3}\km\W{4}{4}\km\W{4}{8}
-\frac{1}{4} \W{0}{3}\km\W{2}{3}\km\W{0}{4}\km\W{7}{8}
+\W{0}{3}\km\W{2}{3}\km\W{1}{4}\km\W{6}{8}
-\frac{3}{2} \W{0}{3}\km\W{2}{3}\km\W{2}{4}\km\W{5}{8}
+\W{0}{3}\km\W{2}{3}\km\W{3}{4}\km\W{4}{8}
\\&
+\frac{1}{8} \W{0}{3}\km\W{3}{3}\km\W{0}{4}\km\W{6}{8}
-\frac{1}{2} \W{0}{3}\km\W{3}{3}\km\W{1}{4}\km\W{5}{8}
-\frac{1}{4} \W{0}{3}\km\W{2}{3}\km\W{4}{4}\km\W{3}{8}
+\frac{3}{4} \W{0}{3}\km\W{3}{3}\km\W{2}{4}\km\W{4}{8}
-\frac{1}{2} \W{0}{3}\km\W{3}{3}\km\W{3}{4}\km\W{3}{8}
+\frac{1}{8} \W{0}{3}\km\W{3}{3}\km\W{4}{4}\km\W{2}{8}
+\frac{3}{8} \W{1}{3}\km\W{2}{3}\km\W{0}{4}\km\W{6}{8}
-\frac{3}{2} \W{1}{3}\km\W{2}{3}\km\W{1}{4}\km\W{5}{8}
-\frac{3}{2} \W{1}{3}\km\W{2}{3}\km\W{3}{4}\km\W{3}{8}
\\&
+\frac{3}{8} \W{1}{3}\km\W{2}{3}\km\W{4}{4}\km\W{2}{8}
+\frac{9}{4} \W{1}{3}\km\W{2}{3}\km\W{2}{4}\km\W{4}{8}
-\frac{1}{4} \W{1}{3}\km\W{3}{3}\km\W{0}{4}\km\W{5}{8}
+\W{1}{3}\km\W{3}{3}\km\W{1}{4}\km\W{4}{8}
-\frac{3}{2} \W{1}{3}\km\W{3}{3}\km\W{2}{4}\km\W{3}{8}
+\W{1}{3}\km\W{3}{3}\km\W{3}{4}\km\W{2}{8}
-\frac{1}{4} \W{1}{3}\km\W{3}{3}\km\W{4}{4}\km\W{1}{8}
+\frac{1}{8} \W{2}{3}\km\W{3}{3}\km\W{0}{4}\km\W{4}{8}
-\frac{1}{2} \W{2}{3}\km\W{3}{3}\km\W{1}{4}\km\W{3}{8}
\\&
+\frac{3}{4} \W{2}{3}\km\W{3}{3}\km\W{2}{4}\km\W{2}{8}
-\frac{1}{2} \W{2}{3}\km\W{3}{3}\km\W{3}{4}\km\W{1}{8}
+\frac{1}{8} \W{2}{3}\km\W{3}{3}\km\W{4}{4}\km\W{0}{8}
\\
\mathbf{q}_{2} &= 
-\frac{3}{11} \W{0}{3}\km\W{2}{3}\km\W{5}{5}\km\W{2}{7}
-\frac{14}{11} \W{1}{3}\km\W{3}{3}\km\W{2}{5}\km\W{3}{7}
+\frac{6}{11} \W{1}{3}\km\W{3}{3}\km\W{3}{5}\km\W{2}{7}
+\frac{1}{11} \W{1}{3}\km\W{3}{3}\km\W{4}{5}\km\W{1}{7}
-\frac{9}{11} \W{1}{3}\km\W{2}{3}\km\W{4}{5}\km\W{2}{7}
+\frac{1}{11} \W{0}{3}\km\W{3}{3}\km\W{0}{5}\km\W{6}{7}
-\frac{3}{11} \W{0}{3}\km\W{3}{3}\km\W{1}{5}\km\W{5}{7}
+\frac{2}{11} \W{0}{3}\km\W{1}{3}\km\W{5}{5}\km\W{3}{7}
-\frac{1}{11} \W{0}{3}\km\W{2}{3}\km\W{0}{5}\km\W{7}{7}
\\&
+\frac{1}{11} \W{0}{3}\km\W{2}{3}\km\W{1}{5}\km\W{6}{7}
+\frac{2}{11} \W{0}{3}\km\W{3}{3}\km\W{3}{5}\km\W{3}{7}
-\frac{3}{11} \W{0}{3}\km\W{3}{3}\km\W{4}{5}\km\W{2}{7}
-\frac{1}{11} \W{1}{3}\km\W{3}{3}\km\W{5}{5}\km\W{0}{7}
+\frac{2}{11} \W{2}{3}\km\W{3}{3}\km\W{0}{5}\km\W{4}{7}
-\frac{8}{11} \W{2}{3}\km\W{3}{3}\km\W{1}{5}\km\W{3}{7}
+\frac{2}{11} \W{0}{3}\km\W{3}{3}\km\W{2}{5}\km\W{4}{7}
-\frac{8}{11} \W{0}{3}\km\W{1}{3}\km\W{2}{5}\km\W{6}{7}
+\frac{12}{11} \W{0}{3}\km\W{1}{3}\km\W{3}{5}\km\W{5}{7}
\\&
-\frac{8}{11} \W{0}{3}\km\W{1}{3}\km\W{4}{5}\km\W{4}{7}
+\frac{3}{11} \W{1}{3}\km\W{2}{3}\km\W{5}{5}\km\W{1}{7}
-\frac{3}{11} \W{1}{3}\km\W{3}{3}\km\W{0}{5}\km\W{5}{7}
+\frac{6}{11} \W{1}{3}\km\W{2}{3}\km\W{2}{5}\km\W{4}{7}
+\frac{6}{11} \W{1}{3}\km\W{2}{3}\km\W{3}{5}\km\W{3}{7}
+\frac{2}{11} \W{0}{3}\km\W{1}{3}\km\W{1}{5}\km\W{7}{7}
+\W{1}{3}\km\W{3}{3}\km\W{1}{5}\km\W{4}{7}
+\frac{1}{11} \W{0}{3}\km\W{3}{3}\km\W{5}{5}\km\W{1}{7}
+\frac{3}{11} \W{1}{3}\km\W{2}{3}\km\W{0}{5}\km\W{6}{7}
\\&
-\frac{9}{11} \W{1}{3}\km\W{2}{3}\km\W{1}{5}\km\W{5}{7}
+\frac{12}{11} \W{2}{3}\km\W{3}{3}\km\W{2}{5}\km\W{2}{7}
-\frac{8}{11} \W{2}{3}\km\W{3}{3}\km\W{3}{5}\km\W{1}{7}
+\frac{2}{11} \W{2}{3}\km\W{3}{3}\km\W{4}{5}\km\W{0}{7}
+\frac{6}{11} \W{0}{3}\km\W{2}{3}\km\W{2}{5}\km\W{5}{7}
-\frac{14}{11} \W{0}{3}\km\W{2}{3}\km\W{3}{5}\km\W{4}{7}
+\W{0}{3}\km\W{2}{3}\km\W{4}{5}\km\W{3}{7}
\\
\mathbf{q}_{3} &= 
-\frac{1}{8} \W{0}{3}\km\W{0}{4}\km\W{4}{4}\km\W{5}{7}
-\frac{5}{8} \W{3}{3}\km\W{0}{4}\km\W{2}{4}\km\W{4}{7}
+\frac{1}{4} \W{3}{3}\km\W{1}{4}\km\W{4}{4}\km\W{1}{7}
+\W{1}{3}\km\W{1}{4}\km\W{3}{4}\km\W{4}{7}
+\W{3}{3}\km\W{1}{4}\km\W{2}{4}\km\W{3}{7}
-\W{3}{3}\km\W{1}{4}\km\W{3}{4}\km\W{2}{7}
-\frac{5}{8} \W{0}{3}\km\W{2}{4}\km\W{4}{4}\km\W{3}{7}
+\frac{1}{4} \W{0}{3}\km\W{3}{4}\km\W{4}{4}\km\W{2}{7}
+\W{2}{3}\km\W{1}{4}\km\W{3}{4}\km\W{3}{7}
\\&
+\frac{1}{2} \W{0}{3}\km\W{1}{4}\km\W{4}{4}\km\W{4}{7}
-\frac{3}{4} \W{1}{3}\km\W{1}{4}\km\W{4}{4}\km\W{3}{7}
-\frac{3}{2} \W{1}{3}\km\W{2}{4}\km\W{3}{4}\km\W{3}{7}
+\frac{1}{4} \W{0}{3}\km\W{0}{4}\km\W{3}{4}\km\W{6}{7}
+\frac{1}{2} \W{3}{3}\km\W{2}{4}\km\W{3}{4}\km\W{1}{7}
+\frac{1}{8} \W{1}{3}\km\W{0}{4}\km\W{4}{4}\km\W{4}{7}
-\frac{3}{2} \W{2}{3}\km\W{1}{4}\km\W{2}{4}\km\W{4}{7}
-\frac{3}{4} \W{2}{3}\km\W{0}{4}\km\W{3}{4}\km\W{4}{7}
+\W{0}{3}\km\W{2}{4}\km\W{3}{4}\km\W{4}{7}
\\&
-\frac{1}{2} \W{2}{3}\km\W{0}{4}\km\W{1}{4}\km\W{6}{7}
+\frac{9}{8} \W{2}{3}\km\W{0}{4}\km\W{2}{4}\km\W{5}{7}
+\frac{1}{8} \W{2}{3}\km\W{0}{4}\km\W{4}{4}\km\W{3}{7}
+\frac{1}{2} \W{0}{3}\km\W{1}{4}\km\W{2}{4}\km\W{6}{7}
-\W{0}{3}\km\W{1}{4}\km\W{3}{4}\km\W{5}{7}
+\frac{1}{2} \W{3}{3}\km\W{0}{4}\km\W{3}{4}\km\W{3}{7}
-\frac{1}{8} \W{3}{3}\km\W{0}{4}\km\W{4}{4}\km\W{2}{7}
-\frac{1}{8} \W{0}{3}\km\W{0}{4}\km\W{2}{4}\km\W{7}{7}
+\frac{1}{4} \W{1}{3}\km\W{0}{4}\km\W{1}{4}\km\W{7}{7}
\\&
-\frac{3}{8} \W{1}{3}\km\W{0}{4}\km\W{2}{4}\km\W{6}{7}
-\frac{1}{8} \W{3}{3}\km\W{2}{4}\km\W{4}{4}\km\W{0}{7}
-\frac{3}{8} \W{2}{3}\km\W{2}{4}\km\W{4}{4}\km\W{1}{7}
+\frac{1}{4} \W{2}{3}\km\W{3}{4}\km\W{4}{4}\km\W{0}{7}
+\frac{1}{4} \W{3}{3}\km\W{0}{4}\km\W{1}{4}\km\W{5}{7}
+\frac{9}{8} \W{1}{3}\km\W{2}{4}\km\W{4}{4}\km\W{2}{7}
-\frac{1}{2} \W{1}{3}\km\W{3}{4}\km\W{4}{4}\km\W{1}{7}
\allowdisplaybreaks
\\
\mathbf{q}_{4} &= 
-\frac{1}{6} \W{3}{3}\km\W{0}{4}\km\W{4}{5}\km\W{2}{6}
+\frac{1}{30} \W{3}{3}\km\W{0}{4}\km\W{5}{5}\km\W{1}{6}
+\frac{1}{6} \W{0}{3}\km\W{4}{4}\km\W{1}{5}\km\W{4}{6}
+\frac{1}{2} \W{1}{3}\km\W{3}{4}\km\W{4}{5}\km\W{1}{6}
-\frac{1}{2} \W{1}{3}\km\W{3}{4}\km\W{1}{5}\km\W{4}{6}
+\frac{1}{6} \W{3}{3}\km\W{1}{4}\km\W{4}{5}\km\W{1}{6}
-\frac{1}{30} \W{3}{3}\km\W{1}{4}\km\W{5}{5}\km\W{0}{6}
+\frac{1}{3} \W{0}{3}\km\W{4}{4}\km\W{3}{5}\km\W{2}{6}
-\W{2}{3}\km\W{1}{4}\km\W{3}{5}\km\W{3}{6}
\\&
+\frac{1}{2} \W{2}{3}\km\W{1}{4}\km\W{4}{5}\km\W{2}{6}
-\frac{1}{6} \W{3}{3}\km\W{1}{4}\km\W{1}{5}\km\W{4}{6}
+\frac{1}{3} \W{3}{3}\km\W{1}{4}\km\W{2}{5}\km\W{3}{6}
+\frac{1}{10} \W{2}{3}\km\W{1}{4}\km\W{0}{5}\km\W{6}{6}
-\frac{1}{3} \W{0}{3}\km\W{4}{4}\km\W{2}{5}\km\W{3}{6}
-\frac{1}{2} \W{1}{3}\km\W{2}{4}\km\W{4}{5}\km\W{2}{6}
+\frac{1}{10} \W{1}{3}\km\W{2}{4}\km\W{5}{5}\km\W{1}{6}
-\frac{1}{3} \W{3}{3}\km\W{1}{4}\km\W{3}{5}\km\W{2}{6}
-\W{1}{3}\km\W{3}{4}\km\W{3}{5}\km\W{2}{6}
\\&
+\frac{1}{10} \W{1}{3}\km\W{3}{4}\km\W{0}{5}\km\W{5}{6}
+\frac{1}{6} \W{0}{3}\km\W{3}{4}\km\W{4}{5}\km\W{2}{6}
-\frac{1}{30} \W{0}{3}\km\W{3}{4}\km\W{5}{5}\km\W{1}{6}
+\frac{1}{30} \W{3}{3}\km\W{1}{4}\km\W{0}{5}\km\W{5}{6}
-\frac{1}{2} \W{2}{3}\km\W{2}{4}\km\W{4}{5}\km\W{1}{6}
+\frac{1}{10} \W{2}{3}\km\W{2}{4}\km\W{5}{5}\km\W{0}{6}
+\W{2}{3}\km\W{2}{4}\km\W{3}{5}\km\W{2}{6}
-\frac{1}{30} \W{0}{3}\km\W{4}{4}\km\W{0}{5}\km\W{5}{6}
+\W{1}{3}\km\W{3}{4}\km\W{2}{5}\km\W{3}{6}
\\&
-\frac{1}{6} \W{0}{3}\km\W{3}{4}\km\W{1}{5}\km\W{5}{6}
+\frac{1}{3} \W{0}{3}\km\W{3}{4}\km\W{2}{5}\km\W{4}{6}
+\frac{1}{30} \W{0}{3}\km\W{4}{4}\km\W{5}{5}\km\W{0}{6}
-\frac{1}{10} \W{1}{3}\km\W{3}{4}\km\W{5}{5}\km\W{0}{6}
+\frac{1}{2} \W{2}{3}\km\W{2}{4}\km\W{1}{5}\km\W{4}{6}
-\W{2}{3}\km\W{2}{4}\km\W{2}{5}\km\W{3}{6}
-\W{1}{3}\km\W{2}{4}\km\W{2}{5}\km\W{4}{6}
+\frac{1}{30} \W{0}{3}\km\W{3}{4}\km\W{0}{5}\km\W{6}{6}
-\frac{1}{30} \W{3}{3}\km\W{0}{4}\km\W{0}{5}\km\W{6}{6}
\\&
-\frac{1}{10} \W{1}{3}\km\W{2}{4}\km\W{0}{5}\km\W{6}{6}
+\frac{1}{2} \W{1}{3}\km\W{2}{4}\km\W{1}{5}\km\W{5}{6}
+\W{1}{3}\km\W{2}{4}\km\W{3}{5}\km\W{3}{6}
+\frac{1}{6} \W{3}{3}\km\W{0}{4}\km\W{1}{5}\km\W{5}{6}
-\frac{1}{3} \W{3}{3}\km\W{0}{4}\km\W{2}{5}\km\W{4}{6}
-\frac{1}{3} \W{0}{3}\km\W{3}{4}\km\W{3}{5}\km\W{3}{6}
+\frac{1}{3} \W{3}{3}\km\W{0}{4}\km\W{3}{5}\km\W{3}{6}
-\frac{1}{10} \W{2}{3}\km\W{1}{4}\km\W{5}{5}\km\W{1}{6}
-\frac{1}{6} \W{0}{3}\km\W{4}{4}\km\W{4}{5}\km\W{1}{6}
\\&
-\frac{1}{2} \W{2}{3}\km\W{1}{4}\km\W{1}{5}\km\W{5}{6}
+\W{2}{3}\km\W{1}{4}\km\W{2}{5}\km\W{4}{6}
-\frac{1}{10} \W{2}{3}\km\W{2}{4}\km\W{0}{5}\km\W{5}{6}
\\
\mathbf{q}_{5} &= 
-\frac{3}{4} \W{1}{3}\km\W{4}{4}\km\W{0}{5}\km\W{4}{6}
-\frac{7}{6} \W{3}{3}\km\W{0}{4}\km\W{4}{5}\km\W{2}{6}
+\frac{1}{12} \W{3}{3}\km\W{0}{4}\km\W{5}{5}\km\W{1}{6}
+\frac{7}{6} \W{0}{3}\km\W{4}{4}\km\W{1}{5}\km\W{4}{6}
+\frac{11}{4} \W{1}{3}\km\W{3}{4}\km\W{4}{5}\km\W{1}{6}
+\frac{3}{4} \W{2}{3}\km\W{0}{4}\km\W{5}{5}\km\W{2}{6}
-\frac{3}{4} \W{3}{3}\km\W{2}{4}\km\W{4}{5}\km\W{0}{6}
-\frac{3}{4} \W{1}{3}\km\W{3}{4}\km\W{1}{5}\km\W{4}{6}
+\frac{23}{12} \W{3}{3}\km\W{1}{4}\km\W{4}{5}\km\W{1}{6}
\\&
-\frac{1}{12} \W{3}{3}\km\W{1}{4}\km\W{5}{5}\km\W{0}{6}
+\frac{11}{6} \W{0}{3}\km\W{4}{4}\km\W{3}{5}\km\W{2}{6}
+\frac{5}{12} \W{1}{3}\km\W{4}{4}\km\W{4}{5}\km\W{0}{6}
-\frac{4}{3} \W{1}{3}\km\W{0}{4}\km\W{5}{5}\km\W{3}{6}
+\frac{11}{6} \W{2}{3}\km\W{1}{4}\km\W{3}{5}\km\W{3}{6}
+\frac{3}{4} \W{2}{3}\km\W{1}{4}\km\W{4}{5}\km\W{2}{6}
-\frac{5}{2} \W{2}{3}\km\W{4}{4}\km\W{1}{5}\km\W{2}{6}
+\W{2}{3}\km\W{4}{4}\km\W{2}{5}\km\W{1}{6}
-\W{1}{3}\km\W{1}{4}\km\W{3}{5}\km\W{4}{6}
\\&
+\frac{19}{12} \W{3}{3}\km\W{1}{4}\km\W{1}{5}\km\W{4}{6}
-\frac{1}{6} \W{3}{3}\km\W{1}{4}\km\W{2}{5}\km\W{3}{6}
+\frac{5}{12} \W{2}{3}\km\W{1}{4}\km\W{0}{5}\km\W{6}{6}
-\frac{5}{2} \W{0}{3}\km\W{4}{4}\km\W{2}{5}\km\W{3}{6}
-3 \W{1}{3}\km\W{2}{4}\km\W{4}{5}\km\W{2}{6}
-\frac{3}{4} \W{1}{3}\km\W{2}{4}\km\W{5}{5}\km\W{1}{6}
-\frac{17}{6} \W{3}{3}\km\W{1}{4}\km\W{3}{5}\km\W{2}{6}
-\frac{3}{2} \W{1}{3}\km\W{3}{4}\km\W{3}{5}\km\W{2}{6}
+3 \W{0}{3}\km\W{2}{4}\km\W{4}{5}\km\W{3}{6}
\\&
-\frac{1}{4} \W{0}{3}\km\W{2}{4}\km\W{5}{5}\km\W{2}{6}
+\frac{7}{4} \W{1}{3}\km\W{3}{4}\km\W{0}{5}\km\W{5}{6}
-\frac{19}{12} \W{0}{3}\km\W{3}{4}\km\W{4}{5}\km\W{2}{6}
+\frac{5}{12} \W{0}{3}\km\W{3}{4}\km\W{5}{5}\km\W{1}{6}
-\frac{1}{2} \W{3}{3}\km\W{4}{4}\km\W{0}{5}\km\W{2}{6}
+\W{3}{3}\km\W{4}{4}\km\W{1}{5}\km\W{1}{6}
+\frac{3}{4} \W{0}{3}\km\W{2}{4}\km\W{1}{5}\km\W{6}{6}
-\frac{2}{3} \W{0}{3}\km\W{1}{4}\km\W{5}{5}\km\W{3}{6}
-\frac{1}{2} \W{3}{3}\km\W{4}{4}\km\W{2}{5}\km\W{0}{6}
\\&
-\frac{5}{12} \W{3}{3}\km\W{1}{4}\km\W{0}{5}\km\W{5}{6}
+\frac{3}{4} \W{2}{3}\km\W{2}{4}\km\W{4}{5}\km\W{1}{6}
+\frac{1}{2} \W{0}{3}\km\W{0}{4}\km\W{3}{5}\km\W{6}{6}
+\W{2}{3}\km\W{2}{4}\km\W{5}{5}\km\W{0}{6}
-\frac{5}{12} \W{2}{3}\km\W{0}{4}\km\W{1}{5}\km\W{6}{6}
-2 \W{2}{3}\km\W{0}{4}\km\W{3}{5}\km\W{4}{6}
-\frac{1}{3} \W{2}{3}\km\W{0}{4}\km\W{4}{5}\km\W{3}{6}
-\frac{1}{12} \W{0}{3}\km\W{4}{4}\km\W{0}{5}\km\W{5}{6}
-\frac{7}{6} \W{2}{3}\km\W{3}{4}\km\W{4}{5}\km\W{0}{6}
\\& 
-\frac{11}{6} \W{1}{3}\km\W{3}{4}\km\W{2}{5}\km\W{3}{6}
-\frac{23}{12} \W{0}{3}\km\W{3}{4}\km\W{1}{5}\km\W{5}{6}
+\frac{17}{6} \W{0}{3}\km\W{3}{4}\km\W{2}{5}\km\W{4}{6}
+\frac{5}{2} \W{1}{3}\km\W{1}{4}\km\W{5}{5}\km\W{2}{6}
-\frac{5}{2} \W{2}{3}\km\W{3}{4}\km\W{0}{5}\km\W{4}{6}
-\W{0}{3}\km\W{0}{4}\km\W{4}{5}\km\W{5}{6}
+\frac{1}{2} \W{0}{3}\km\W{0}{4}\km\W{5}{5}\km\W{4}{6}
+2 \W{1}{3}\km\W{4}{4}\km\W{2}{5}\km\W{2}{6}
-2 \W{1}{3}\km\W{4}{4}\km\W{3}{5}\km\W{1}{6}
\\&
+\frac{1}{6} \W{2}{3}\km\W{4}{4}\km\W{3}{5}\km\W{0}{6}
-\frac{5}{12} \W{1}{3}\km\W{3}{4}\km\W{5}{5}\km\W{0}{6}
+3 \W{2}{3}\km\W{2}{4}\km\W{1}{5}\km\W{4}{6}
-\frac{11}{2} \W{2}{3}\km\W{2}{4}\km\W{2}{5}\km\W{3}{6}
-3 \W{3}{3}\km\W{2}{4}\km\W{1}{5}\km\W{3}{6}
+\frac{9}{2} \W{3}{3}\km\W{2}{4}\km\W{2}{5}\km\W{2}{6}
-\W{3}{3}\km\W{2}{4}\km\W{3}{5}\km\W{1}{6}
+\frac{1}{3} \W{1}{3}\km\W{4}{4}\km\W{1}{5}\km\W{3}{6}
+\frac{4}{3} \W{2}{3}\km\W{4}{4}\km\W{0}{5}\km\W{3}{6}
\\&
-2 \W{3}{3}\km\W{3}{4}\km\W{2}{5}\km\W{1}{6}
+\frac{4}{3} \W{3}{3}\km\W{3}{4}\km\W{3}{5}\km\W{0}{6}
+\frac{7}{6} \W{1}{3}\km\W{1}{4}\km\W{1}{5}\km\W{6}{6}
+\frac{1}{12} \W{0}{3}\km\W{3}{4}\km\W{0}{5}\km\W{6}{6}
+\frac{5}{2} \W{1}{3}\km\W{0}{4}\km\W{4}{5}\km\W{4}{6}
+\frac{2}{3} \W{3}{3}\km\W{3}{4}\km\W{0}{5}\km\W{3}{6}
+\W{0}{3}\km\W{2}{4}\km\W{2}{5}\km\W{5}{6}
+\frac{8}{3} \W{2}{3}\km\W{3}{4}\km\W{1}{5}\km\W{3}{6}
-\W{1}{3}\km\W{2}{4}\km\W{0}{5}\km\W{6}{6}
\\&
-\frac{3}{4} \W{1}{3}\km\W{2}{4}\km\W{1}{5}\km\W{5}{6}
+2 \W{2}{3}\km\W{0}{4}\km\W{2}{5}\km\W{5}{6}
+\frac{11}{2} \W{1}{3}\km\W{2}{4}\km\W{3}{5}\km\W{3}{6}
+\frac{5}{12} \W{3}{3}\km\W{0}{4}\km\W{1}{5}\km\W{5}{6}
-\frac{11}{6} \W{3}{3}\km\W{0}{4}\km\W{2}{5}\km\W{4}{6}
+\frac{1}{6} \W{0}{3}\km\W{3}{4}\km\W{3}{5}\km\W{3}{6}
+\frac{5}{2} \W{3}{3}\km\W{0}{4}\km\W{3}{5}\km\W{3}{6}
-\frac{1}{6} \W{1}{3}\km\W{0}{4}\km\W{2}{5}\km\W{6}{6}
-\W{1}{3}\km\W{0}{4}\km\W{3}{5}\km\W{5}{6}
\\&
-\frac{7}{4} \W{2}{3}\km\W{1}{4}\km\W{5}{5}\km\W{1}{6}
-\frac{5}{12} \W{0}{3}\km\W{4}{4}\km\W{4}{5}\km\W{1}{6}
-\frac{11}{4} \W{2}{3}\km\W{1}{4}\km\W{1}{5}\km\W{5}{6}
+\frac{3}{2} \W{2}{3}\km\W{1}{4}\km\W{2}{5}\km\W{4}{6}
-\frac{8}{3} \W{1}{3}\km\W{1}{4}\km\W{4}{5}\km\W{3}{6}
+\frac{1}{4} \W{3}{3}\km\W{2}{4}\km\W{0}{5}\km\W{4}{6}
+\frac{3}{4} \W{2}{3}\km\W{2}{4}\km\W{0}{5}\km\W{5}{6}
+\W{2}{3}\km\W{3}{4}\km\W{2}{5}\km\W{2}{6}
-\frac{4}{3} \W{0}{3}\km\W{1}{4}\km\W{2}{5}\km\W{6}{6}
\\&
+2 \W{0}{3}\km\W{1}{4}\km\W{3}{5}\km\W{5}{6}
-\frac{9}{2} \W{0}{3}\km\W{2}{4}\km\W{3}{5}\km\W{4}{6}
\\
\mathbf{q}_{6} &= 
-\W{1}{3}\km\W{4}{4}\km\W{0}{5}\km\W{4}{6}
-\frac{25}{12} \W{3}{3}\km\W{0}{4}\km\W{4}{5}\km\W{2}{6}
+\frac{13}{60} \W{3}{3}\km\W{0}{4}\km\W{5}{5}\km\W{1}{6}
+\frac{25}{12} \W{0}{3}\km\W{4}{4}\km\W{1}{5}\km\W{4}{6}
+\frac{17}{4} \W{1}{3}\km\W{3}{4}\km\W{4}{5}\km\W{1}{6}
+\W{2}{3}\km\W{0}{4}\km\W{5}{5}\km\W{2}{6}
-\W{3}{3}\km\W{2}{4}\km\W{4}{5}\km\W{0}{6}
-\frac{5}{4} \W{1}{3}\km\W{3}{4}\km\W{1}{5}\km\W{4}{6}
+\frac{37}{12} \W{3}{3}\km\W{1}{4}\km\W{4}{5}\km\W{1}{6}
\\&
-\frac{13}{60} \W{3}{3}\km\W{1}{4}\km\W{5}{5}\km\W{0}{6}
+\frac{35}{12} \W{0}{3}\km\W{4}{4}\km\W{3}{5}\km\W{2}{6}
+\frac{3}{4} \W{1}{3}\km\W{4}{4}\km\W{4}{5}\km\W{0}{6}
-2 \W{1}{3}\km\W{0}{4}\km\W{5}{5}\km\W{3}{6}
+\frac{5}{2} \W{2}{3}\km\W{1}{4}\km\W{3}{5}\km\W{3}{6}
+\frac{5}{4} \W{2}{3}\km\W{1}{4}\km\W{4}{5}\km\W{2}{6}
-\frac{15}{4} \W{2}{3}\km\W{4}{4}\km\W{1}{5}\km\W{2}{6}
+\frac{3}{2} \W{2}{3}\km\W{4}{4}\km\W{2}{5}\km\W{1}{6}
-\W{1}{3}\km\W{1}{4}\km\W{2}{5}\km\W{5}{6}
\\&
+\frac{35}{12} \W{3}{3}\km\W{1}{4}\km\W{1}{5}\km\W{4}{6}
-\frac{5}{6} \W{3}{3}\km\W{1}{4}\km\W{2}{5}\km\W{3}{6}
+\frac{13}{20} \W{2}{3}\km\W{1}{4}\km\W{0}{5}\km\W{6}{6}
-\frac{25}{6} \W{0}{3}\km\W{4}{4}\km\W{2}{5}\km\W{3}{6}
-\frac{15}{4} \W{1}{3}\km\W{2}{4}\km\W{4}{5}\km\W{2}{6}
-\frac{27}{20} \W{1}{3}\km\W{2}{4}\km\W{5}{5}\km\W{1}{6}
-\frac{25}{6} \W{3}{3}\km\W{1}{4}\km\W{3}{5}\km\W{2}{6}
-\frac{5}{2} \W{1}{3}\km\W{3}{4}\km\W{3}{5}\km\W{2}{6}
+5 \W{0}{3}\km\W{2}{4}\km\W{4}{5}\km\W{3}{6}
\\&
-\frac{1}{2} \W{0}{3}\km\W{2}{4}\km\W{5}{5}\km\W{2}{6}
+\frac{53}{20} \W{1}{3}\km\W{3}{4}\km\W{0}{5}\km\W{5}{6}
-\frac{35}{12} \W{0}{3}\km\W{3}{4}\km\W{4}{5}\km\W{2}{6}
+\frac{47}{60} \W{0}{3}\km\W{3}{4}\km\W{5}{5}\km\W{1}{6}
-\frac{3}{4} \W{3}{3}\km\W{4}{4}\km\W{0}{5}\km\W{2}{6}
+\frac{3}{2} \W{3}{3}\km\W{4}{4}\km\W{1}{5}\km\W{1}{6}
+\W{0}{3}\km\W{2}{4}\km\W{1}{5}\km\W{6}{6}
-\W{0}{3}\km\W{1}{4}\km\W{5}{5}\km\W{3}{6}
-\frac{3}{4} \W{3}{3}\km\W{4}{4}\km\W{2}{5}\km\W{0}{6}
\\&
-\frac{47}{60} \W{3}{3}\km\W{1}{4}\km\W{0}{5}\km\W{5}{6}
+\frac{3}{4} \W{2}{3}\km\W{2}{4}\km\W{4}{5}\km\W{1}{6}
+\frac{3}{4} \W{0}{3}\km\W{0}{4}\km\W{3}{5}\km\W{6}{6}
+\frac{33}{20} \W{2}{3}\km\W{2}{4}\km\W{5}{5}\km\W{0}{6}
-\frac{3}{4} \W{2}{3}\km\W{0}{4}\km\W{1}{5}\km\W{6}{6}
-\frac{15}{4} \W{2}{3}\km\W{0}{4}\km\W{3}{5}\km\W{4}{6}
-\frac{13}{60} \W{0}{3}\km\W{4}{4}\km\W{0}{5}\km\W{5}{6}
-2 \W{2}{3}\km\W{3}{4}\km\W{4}{5}\km\W{0}{6}
-\frac{5}{2} \W{1}{3}\km\W{3}{4}\km\W{2}{5}\km\W{3}{6}
\\&
-\frac{37}{12} \W{0}{3}\km\W{3}{4}\km\W{1}{5}\km\W{5}{6}
+\frac{25}{6} \W{0}{3}\km\W{3}{4}\km\W{2}{5}\km\W{4}{6}
+4 \W{1}{3}\km\W{1}{4}\km\W{5}{5}\km\W{2}{6}
-4 \W{2}{3}\km\W{3}{4}\km\W{0}{5}\km\W{4}{6}
-\frac{3}{2} \W{0}{3}\km\W{0}{4}\km\W{4}{5}\km\W{5}{6}
+\frac{3}{4} \W{0}{3}\km\W{0}{4}\km\W{5}{5}\km\W{4}{6}
+\frac{15}{4} \W{1}{3}\km\W{4}{4}\km\W{2}{5}\km\W{2}{6}
-\frac{7}{2} \W{1}{3}\km\W{4}{4}\km\W{3}{5}\km\W{1}{6}
+\frac{1}{4} \W{2}{3}\km\W{4}{4}\km\W{3}{5}\km\W{0}{6}
\\&
-\frac{1}{30} \W{0}{3}\km\W{4}{4}\km\W{5}{5}\km\W{0}{6}
-\frac{13}{20} \W{1}{3}\km\W{3}{4}\km\W{5}{5}\km\W{0}{6}
+\frac{15}{4} \W{2}{3}\km\W{2}{4}\km\W{1}{5}\km\W{4}{6}
-\frac{15}{2} \W{2}{3}\km\W{2}{4}\km\W{2}{5}\km\W{3}{6}
-5 \W{3}{3}\km\W{2}{4}\km\W{1}{5}\km\W{3}{6}
+\frac{15}{2} \W{3}{3}\km\W{2}{4}\km\W{2}{5}\km\W{2}{6}
-2 \W{3}{3}\km\W{2}{4}\km\W{3}{5}\km\W{1}{6}
+2 \W{2}{3}\km\W{4}{4}\km\W{0}{5}\km\W{3}{6}
-3 \W{3}{3}\km\W{3}{4}\km\W{2}{5}\km\W{1}{6}
\\&
+2 \W{3}{3}\km\W{3}{4}\km\W{3}{5}\km\W{0}{6}
+2 \W{1}{3}\km\W{1}{4}\km\W{1}{5}\km\W{6}{6}
+\frac{13}{60} \W{0}{3}\km\W{3}{4}\km\W{0}{5}\km\W{6}{6}
+\frac{15}{4} \W{1}{3}\km\W{0}{4}\km\W{4}{5}\km\W{4}{6}
+\W{3}{3}\km\W{3}{4}\km\W{0}{5}\km\W{3}{6}
+\frac{1}{30} \W{3}{3}\km\W{0}{4}\km\W{0}{5}\km\W{6}{6}
+2 \W{0}{3}\km\W{2}{4}\km\W{2}{5}\km\W{5}{6}
+5 \W{2}{3}\km\W{3}{4}\km\W{1}{5}\km\W{3}{6}
-\frac{33}{20} \W{1}{3}\km\W{2}{4}\km\W{0}{5}\km\W{6}{6}
\\&
-\frac{3}{4} \W{1}{3}\km\W{2}{4}\km\W{1}{5}\km\W{5}{6}
+\frac{7}{2} \W{2}{3}\km\W{0}{4}\km\W{2}{5}\km\W{5}{6}
+\frac{15}{2} \W{1}{3}\km\W{2}{4}\km\W{3}{5}\km\W{3}{6}
+\frac{7}{12} \W{3}{3}\km\W{0}{4}\km\W{1}{5}\km\W{5}{6}
-\frac{35}{12} \W{3}{3}\km\W{0}{4}\km\W{2}{5}\km\W{4}{6}
+\frac{5}{6} \W{0}{3}\km\W{3}{4}\km\W{3}{5}\km\W{3}{6}
+\frac{25}{6} \W{3}{3}\km\W{0}{4}\km\W{3}{5}\km\W{3}{6}
-\frac{1}{4} \W{1}{3}\km\W{0}{4}\km\W{2}{5}\km\W{6}{6}
-\frac{3}{2} \W{1}{3}\km\W{0}{4}\km\W{3}{5}\km\W{5}{6}
\\&
-\frac{53}{20} \W{2}{3}\km\W{1}{4}\km\W{5}{5}\km\W{1}{6}
-\frac{7}{12} \W{0}{3}\km\W{4}{4}\km\W{4}{5}\km\W{1}{6}
-\frac{17}{4} \W{2}{3}\km\W{1}{4}\km\W{1}{5}\km\W{5}{6}
+\frac{5}{2} \W{2}{3}\km\W{1}{4}\km\W{2}{5}\km\W{4}{6}
-5 \W{1}{3}\km\W{1}{4}\km\W{4}{5}\km\W{3}{6}
+\frac{1}{2} \W{3}{3}\km\W{2}{4}\km\W{0}{5}\km\W{4}{6}
+\frac{27}{20} \W{2}{3}\km\W{2}{4}\km\W{0}{5}\km\W{5}{6}
+\W{2}{3}\km\W{3}{4}\km\W{3}{5}\km\W{1}{6}
-2 \W{0}{3}\km\W{1}{4}\km\W{2}{5}\km\W{6}{6}
\\&
+3 \W{0}{3}\km\W{1}{4}\km\W{3}{5}\km\W{5}{6}
-\frac{15}{2} \W{0}{3}\km\W{2}{4}\km\W{3}{5}\km\W{4}{6}
\allowdisplaybreaks
\\
\mathbf{q}_{7} =& 
-\frac{3}{4} \W{1}{3}\km\W{4}{4}\km\W{0}{5}\km\W{4}{6}
-\frac{3}{4} \W{3}{3}\km\W{0}{4}\km\W{4}{5}\km\W{2}{6}
+\frac{3}{4} \W{0}{3}\km\W{4}{4}\km\W{1}{5}\km\W{4}{6}
+\frac{3}{2} \W{1}{3}\km\W{3}{4}\km\W{4}{5}\km\W{1}{6}
+\frac{3}{4} \W{2}{3}\km\W{0}{4}\km\W{5}{5}\km\W{2}{6}
-\frac{3}{4} \W{3}{3}\km\W{2}{4}\km\W{4}{5}\km\W{0}{6}
-\frac{3}{2} \W{1}{3}\km\W{3}{4}\km\W{1}{5}\km\W{4}{6}
+\frac{3}{2} \W{3}{3}\km\W{1}{4}\km\W{4}{5}\km\W{1}{6}
+\frac{7}{4} \W{0}{3}\km\W{4}{4}\km\W{3}{5}\km\W{2}{6}
\\&
+\frac{1}{2} \W{1}{3}\km\W{4}{4}\km\W{4}{5}\km\W{0}{6}
-\W{1}{3}\km\W{0}{4}\km\W{5}{5}\km\W{3}{6}
+\W{2}{3}\km\W{1}{4}\km\W{3}{5}\km\W{3}{6}
+\frac{3}{2} \W{2}{3}\km\W{1}{4}\km\W{4}{5}\km\W{2}{6}
-\frac{9}{4} \W{2}{3}\km\W{4}{4}\km\W{1}{5}\km\W{2}{6}
+\frac{3}{2} \W{2}{3}\km\W{4}{4}\km\W{2}{5}\km\W{1}{6}
+\W{3}{3}\km\W{1}{4}\km\W{1}{5}\km\W{4}{6}
+\W{3}{3}\km\W{1}{4}\km\W{2}{5}\km\W{3}{6}
+\frac{1}{2} \W{2}{3}\km\W{1}{4}\km\W{0}{5}\km\W{6}{6}
\\&
-2 \W{0}{3}\km\W{4}{4}\km\W{2}{5}\km\W{3}{6}
-\frac{9}{4} \W{1}{3}\km\W{2}{4}\km\W{4}{5}\km\W{2}{6}
-3 \W{3}{3}\km\W{1}{4}\km\W{3}{5}\km\W{2}{6}
+3 \W{0}{3}\km\W{2}{4}\km\W{4}{5}\km\W{3}{6}
-\frac{3}{4} \W{0}{3}\km\W{2}{4}\km\W{5}{5}\km\W{2}{6}
+\frac{3}{2} \W{1}{3}\km\W{3}{4}\km\W{0}{5}\km\W{5}{6}
-\W{0}{3}\km\W{3}{4}\km\W{4}{5}\km\W{2}{6}
+\frac{1}{2} \W{0}{3}\km\W{3}{4}\km\W{5}{5}\km\W{1}{6}
-\frac{1}{4} \W{3}{3}\km\W{4}{4}\km\W{0}{5}\km\W{2}{6}
\\&
+\frac{1}{2} \W{3}{3}\km\W{4}{4}\km\W{1}{5}\km\W{1}{6}
+\frac{3}{4} \W{0}{3}\km\W{2}{4}\km\W{1}{5}\km\W{6}{6}
-\frac{1}{4} \W{3}{3}\km\W{4}{4}\km\W{2}{5}\km\W{0}{6}
-\frac{1}{2} \W{3}{3}\km\W{1}{4}\km\W{0}{5}\km\W{5}{6}
+\frac{1}{4} \W{0}{3}\km\W{0}{4}\km\W{3}{5}\km\W{6}{6}
+\frac{3}{4} \W{2}{3}\km\W{2}{4}\km\W{5}{5}\km\W{0}{6}
-\frac{1}{2} \W{2}{3}\km\W{0}{4}\km\W{1}{5}\km\W{6}{6}
-\frac{3}{4} \W{2}{3}\km\W{0}{4}\km\W{3}{5}\km\W{4}{6}
-\W{2}{3}\km\W{0}{4}\km\W{4}{5}\km\W{3}{6}
\\&
-\frac{1}{2} \W{2}{3}\km\W{3}{4}\km\W{4}{5}\km\W{0}{6}
-\W{1}{3}\km\W{3}{4}\km\W{2}{5}\km\W{3}{6}
-\frac{3}{2} \W{0}{3}\km\W{3}{4}\km\W{1}{5}\km\W{5}{6}
+3 \W{0}{3}\km\W{3}{4}\km\W{2}{5}\km\W{4}{6}
+\frac{3}{2} \W{1}{3}\km\W{1}{4}\km\W{5}{5}\km\W{2}{6}
-\frac{3}{2} \W{2}{3}\km\W{3}{4}\km\W{0}{5}\km\W{4}{6}
-\frac{1}{2} \W{0}{3}\km\W{0}{4}\km\W{4}{5}\km\W{5}{6}
+\frac{1}{4} \W{0}{3}\km\W{0}{4}\km\W{5}{5}\km\W{4}{6}
+\frac{3}{4} \W{1}{3}\km\W{4}{4}\km\W{2}{5}\km\W{2}{6}
\\&
-\frac{3}{2} \W{1}{3}\km\W{4}{4}\km\W{3}{5}\km\W{1}{6}
-\frac{1}{4} \W{2}{3}\km\W{4}{4}\km\W{3}{5}\km\W{0}{6}
-\frac{1}{2} \W{1}{3}\km\W{3}{4}\km\W{5}{5}\km\W{0}{6}
+\frac{9}{4} \W{2}{3}\km\W{2}{4}\km\W{1}{5}\km\W{4}{6}
-3 \W{2}{3}\km\W{2}{4}\km\W{2}{5}\km\W{3}{6}
-3 \W{3}{3}\km\W{2}{4}\km\W{1}{5}\km\W{3}{6}
+3 \W{3}{3}\km\W{2}{4}\km\W{2}{5}\km\W{2}{6}
+\W{1}{3}\km\W{4}{4}\km\W{1}{5}\km\W{3}{6}
+\W{2}{3}\km\W{4}{4}\km\W{0}{5}\km\W{3}{6}
\\&
-2 \W{3}{3}\km\W{3}{4}\km\W{2}{5}\km\W{1}{6}
+\W{3}{3}\km\W{3}{4}\km\W{3}{5}\km\W{0}{6}
+\frac{1}{2} \W{1}{3}\km\W{1}{4}\km\W{1}{5}\km\W{6}{6}
+\frac{9}{4} \W{1}{3}\km\W{0}{4}\km\W{4}{5}\km\W{4}{6}
+\W{3}{3}\km\W{3}{4}\km\W{1}{5}\km\W{2}{6}
+2 \W{2}{3}\km\W{3}{4}\km\W{1}{5}\km\W{3}{6}
-\frac{3}{4} \W{1}{3}\km\W{2}{4}\km\W{0}{5}\km\W{6}{6}
+\frac{3}{2} \W{2}{3}\km\W{0}{4}\km\W{2}{5}\km\W{5}{6}
+3 \W{1}{3}\km\W{2}{4}\km\W{3}{5}\km\W{3}{6}
\\&
+\frac{1}{2} \W{3}{3}\km\W{0}{4}\km\W{1}{5}\km\W{5}{6}
-\frac{7}{4} \W{3}{3}\km\W{0}{4}\km\W{2}{5}\km\W{4}{6}
-\W{0}{3}\km\W{3}{4}\km\W{3}{5}\km\W{3}{6}
+2 \W{3}{3}\km\W{0}{4}\km\W{3}{5}\km\W{3}{6}
+\frac{1}{4} \W{1}{3}\km\W{0}{4}\km\W{2}{5}\km\W{6}{6}
-\frac{3}{2} \W{1}{3}\km\W{0}{4}\km\W{3}{5}\km\W{5}{6}
-\W{0}{3}\km\W{1}{4}\km\W{4}{5}\km\W{4}{6}
-\frac{3}{2} \W{2}{3}\km\W{1}{4}\km\W{5}{5}\km\W{1}{6}
-\frac{1}{2} \W{0}{3}\km\W{4}{4}\km\W{4}{5}\km\W{1}{6}
\\&
-\frac{3}{2} \W{2}{3}\km\W{1}{4}\km\W{1}{5}\km\W{5}{6}
-2 \W{1}{3}\km\W{1}{4}\km\W{4}{5}\km\W{3}{6}
+\frac{3}{4} \W{3}{3}\km\W{2}{4}\km\W{0}{5}\km\W{4}{6}
-\W{0}{3}\km\W{1}{4}\km\W{2}{5}\km\W{6}{6}
+2 \W{0}{3}\km\W{1}{4}\km\W{3}{5}\km\W{5}{6}
-3 \W{0}{3}\km\W{2}{4}\km\W{3}{5}\km\W{4}{6}
\\ 
\mathbf{q}_{8} &= 
\W{3}{3}\km\W{0}{5}\km\W{2}{5}\km\W{4}{5}
-\frac{2}{5} \W{0}{3}\km\W{0}{5}\km\W{4}{5}\km\W{5}{5}
-4 \W{0}{3}\km\W{2}{5}\km\W{3}{5}\km\W{4}{5}
+\frac{3}{5} \W{1}{3}\km\W{0}{5}\km\W{3}{5}\km\W{5}{5}
+3 \W{1}{3}\km\W{1}{5}\km\W{3}{5}\km\W{4}{5}
-3 \W{1}{3}\km\W{1}{5}\km\W{2}{5}\km\W{5}{5}
+3 \W{2}{3}\km\W{1}{5}\km\W{2}{5}\km\W{4}{5}
-3 \W{2}{3}\km\W{0}{5}\km\W{3}{5}\km\W{4}{5}
+\frac{3}{5} \W{2}{3}\km\W{0}{5}\km\W{2}{5}\km\W{5}{5}
\\&
-\frac{2}{5} \W{3}{3}\km\W{0}{5}\km\W{1}{5}\km\W{5}{5}
+\W{0}{3}\km\W{1}{5}\km\W{3}{5}\km\W{5}{5}
-4 \W{3}{3}\km\W{1}{5}\km\W{2}{5}\km\W{3}{5}
\\
\mathbf{q}_{9} &= 
\W{0}{4}\km\W{2}{4}\km\W{4}{4}\km\W{3}{6}
-\frac{1}{2} \W{0}{4}\km\W{3}{4}\km\W{4}{4}\km\W{2}{6}
+2 \W{1}{4}\km\W{2}{4}\km\W{3}{4}\km\W{3}{6}
-\frac{3}{2} \W{1}{4}\km\W{2}{4}\km\W{4}{4}\km\W{2}{6}
+\W{1}{4}\km\W{3}{4}\km\W{4}{4}\km\W{1}{6}
-\frac{1}{2} \W{2}{4}\km\W{3}{4}\km\W{4}{4}\km\W{0}{6}
-\frac{1}{2} \W{0}{4}\km\W{1}{4}\km\W{2}{4}\km\W{6}{6}
+\W{0}{4}\km\W{1}{4}\km\W{3}{4}\km\W{5}{6}
-\frac{1}{2} \W{0}{4}\km\W{1}{4}\km\W{4}{4}\km\W{4}{6}
\\&
-\frac{3}{2} \W{0}{4}\km\W{2}{4}\km\W{3}{4}\km\W{4}{6}
\end{align*}

\subsection{Basis of $\ds \CGF{5}{2}{0}{10}$} 
\begin{align*}
\mathbf{w}_{1} &= 
 -\frac{1}{6} \W{0}{3}\km\W{1}{3}\km\W{2}{3}\km\W{0}{4}\km\W{7}{7}
+\frac{2}{3} \W{0}{3}\km\W{1}{3}\km\W{2}{3}\km\W{1}{4}\km\W{6}{7}
-\W{0}{3}\km\W{1}{3}\km\W{2}{3}\km\W{2}{4}\km\W{5}{7}
+\frac{2}{3} \W{0}{3}\km\W{1}{3}\km\W{2}{3}\km\W{3}{4}\km\W{4}{7}
-\frac{1}{6} \W{0}{3}\km\W{1}{3}\km\W{2}{3}\km\W{4}{4}\km\W{3}{7}
+\frac{1}{6} \W{0}{3}\km\W{1}{3}\km\W{3}{3}\km\W{0}{4}\km\W{6}{7}
-\frac{2}{3} \W{0}{3}\km\W{1}{3}\km\W{3}{3}\km\W{1}{4}\km\W{5}{7}
\\&
+\W{0}{3}\km\W{1}{3}\km\W{3}{3}\km\W{2}{4}\km\W{4}{7}
-\frac{2}{3} \W{0}{3}\km\W{1}{3}\km\W{3}{3}\km\W{3}{4}\km\W{3}{7}
+\frac{1}{6} \W{0}{3}\km\W{1}{3}\km\W{3}{3}\km\W{4}{4}\km\W{2}{7}
+\frac{1}{6} \W{1}{3}\km\W{2}{3}\km\W{3}{3}\km\W{4}{4}\km\W{0}{7}
-\frac{1}{6} \W{0}{3}\km\W{2}{3}\km\W{3}{3}\km\W{0}{4}\km\W{5}{7}
+\frac{2}{3} \W{0}{3}\km\W{2}{3}\km\W{3}{3}\km\W{1}{4}\km\W{4}{7}
-\W{0}{3}\km\W{2}{3}\km\W{3}{3}\km\W{2}{4}\km\W{3}{7}
\\&
+\frac{2}{3} \W{0}{3}\km\W{2}{3}\km\W{3}{3}\km\W{3}{4}\km\W{2}{7}
-\frac{1}{6} \W{0}{3}\km\W{2}{3}\km\W{3}{3}\km\W{4}{4}\km\W{1}{7}
+\frac{1}{6} \W{1}{3}\km\W{2}{3}\km\W{3}{3}\km\W{0}{4}\km\W{4}{7}
-\frac{2}{3} \W{1}{3}\km\W{2}{3}\km\W{3}{3}\km\W{1}{4}\km\W{3}{7}
+\W{1}{3}\km\W{2}{3}\km\W{3}{3}\km\W{2}{4}\km\W{2}{7}
-\frac{2}{3} \W{1}{3}\km\W{2}{3}\km\W{3}{3}\km\W{3}{4}\km\W{1}{7}
\\
\mathbf{w}_{2} &= 
\frac{3}{2} \W{0}{3}\km\W{1}{3}\km\W{2}{3}\km\W{1}{5}\km\W{6}{6}
-6 \W{0}{3}\km\W{1}{3}\km\W{2}{3}\km\W{2}{5}\km\W{5}{6}
+9 \W{0}{3}\km\W{1}{3}\km\W{2}{3}\km\W{3}{5}\km\W{4}{6}
-6 \W{0}{3}\km\W{1}{3}\km\W{2}{3}\km\W{4}{5}\km\W{3}{6}
+\frac{3}{2} \W{0}{3}\km\W{1}{3}\km\W{2}{3}\km\W{5}{5}\km\W{2}{6}
-\frac{1}{2} \W{0}{3}\km\W{1}{3}\km\W{3}{3}\km\W{0}{5}\km\W{6}{6}
+\W{0}{3}\km\W{1}{3}\km\W{3}{3}\km\W{1}{5}\km\W{5}{6}
\\&
+\W{0}{3}\km\W{1}{3}\km\W{3}{3}\km\W{2}{5}\km\W{4}{6}
-4 \W{0}{3}\km\W{1}{3}\km\W{3}{3}\km\W{3}{5}\km\W{3}{6}
+\frac{7}{2} \W{0}{3}\km\W{1}{3}\km\W{3}{3}\km\W{4}{5}\km\W{2}{6}
-\W{0}{3}\km\W{1}{3}\km\W{3}{3}\km\W{5}{5}\km\W{1}{6}
+\W{0}{3}\km\W{2}{3}\km\W{3}{3}\km\W{0}{5}\km\W{5}{6}
-\frac{7}{2} \W{0}{3}\km\W{2}{3}\km\W{3}{3}\km\W{1}{5}\km\W{4}{6}
+4 \W{0}{3}\km\W{2}{3}\km\W{3}{3}\km\W{2}{5}\km\W{3}{6}
\\&
-\W{0}{3}\km\W{2}{3}\km\W{3}{3}\km\W{3}{5}\km\W{2}{6}
-\W{0}{3}\km\W{2}{3}\km\W{3}{3}\km\W{4}{5}\km\W{1}{6}
+\frac{1}{2} \W{0}{3}\km\W{2}{3}\km\W{3}{3}\km\W{5}{5}\km\W{0}{6}
-\frac{3}{2} \W{1}{3}\km\W{2}{3}\km\W{3}{3}\km\W{0}{5}\km\W{4}{6}
+6 \W{1}{3}\km\W{2}{3}\km\W{3}{3}\km\W{1}{5}\km\W{3}{6}
-9 \W{1}{3}\km\W{2}{3}\km\W{3}{3}\km\W{2}{5}\km\W{2}{6}
+6 \W{1}{3}\km\W{2}{3}\km\W{3}{3}\km\W{3}{5}\km\W{1}{6}
\\&
-\frac{3}{2} \W{1}{3}\km\W{2}{3}\km\W{3}{3}\km\W{4}{5}\km\W{0}{6}
\\
\mathbf{w}_{3} &= 
\frac{1}{24} \W{0}{3}\km\W{3}{3}\km\W{0}{4}\km\W{1}{4}\km\W{6}{6}
-\frac{1}{8} \W{0}{3}\km\W{3}{3}\km\W{0}{4}\km\W{2}{4}\km\W{5}{6}
+\frac{1}{8} \W{0}{3}\km\W{3}{3}\km\W{0}{4}\km\W{3}{4}\km\W{4}{6}
-\frac{1}{24} \W{0}{3}\km\W{3}{3}\km\W{0}{4}\km\W{4}{4}\km\W{3}{6}
+\frac{1}{4} \W{0}{3}\km\W{3}{3}\km\W{1}{4}\km\W{2}{4}\km\W{4}{6}
-\frac{1}{3} \W{0}{3}\km\W{3}{3}\km\W{1}{4}\km\W{3}{4}\km\W{3}{6}
+\frac{1}{8} \W{0}{3}\km\W{3}{3}\km\W{1}{4}\km\W{4}{4}\km\W{2}{6}
\\&
+\frac{1}{4} \W{0}{3}\km\W{3}{3}\km\W{2}{4}\km\W{3}{4}\km\W{2}{6}
-\frac{1}{8} \W{0}{3}\km\W{3}{3}\km\W{2}{4}\km\W{4}{4}\km\W{1}{6}
+\frac{1}{24} \W{0}{3}\km\W{3}{3}\km\W{3}{4}\km\W{4}{4}\km\W{0}{6}
-\frac{1}{8} \W{1}{3}\km\W{2}{3}\km\W{0}{4}\km\W{1}{4}\km\W{6}{6}
+\frac{3}{8} \W{1}{3}\km\W{2}{3}\km\W{0}{4}\km\W{2}{4}\km\W{5}{6}
-\frac{3}{8} \W{1}{3}\km\W{2}{3}\km\W{0}{4}\km\W{3}{4}\km\W{4}{6}
+\frac{1}{8} \W{1}{3}\km\W{2}{3}\km\W{0}{4}\km\W{4}{4}\km\W{3}{6}
\\&
-\frac{3}{4} \W{1}{3}\km\W{2}{3}\km\W{1}{4}\km\W{2}{4}\km\W{4}{6}
-\frac{3}{8} \W{1}{3}\km\W{2}{3}\km\W{1}{4}\km\W{4}{4}\km\W{2}{6}
-\frac{3}{4} \W{1}{3}\km\W{2}{3}\km\W{2}{4}\km\W{3}{4}\km\W{2}{6}
+\frac{3}{8} \W{1}{3}\km\W{2}{3}\km\W{2}{4}\km\W{4}{4}\km\W{1}{6}
-\frac{1}{8} \W{1}{3}\km\W{2}{3}\km\W{3}{4}\km\W{4}{4}\km\W{0}{6}
+\W{1}{3}\km\W{2}{3}\km\W{1}{4}\km\W{3}{4}\km\W{3}{6}
\\
\mathbf{w}_{4} &= 
\W{0}{3}\km\W{1}{3}\km\W{1}{4}\km\W{2}{4}\km\W{6}{6}
-2 \W{0}{3}\km\W{1}{3}\km\W{1}{4}\km\W{3}{4}\km\W{5}{6}
+\W{0}{3}\km\W{1}{3}\km\W{1}{4}\km\W{4}{4}\km\W{4}{6}
+3 \W{0}{3}\km\W{1}{3}\km\W{2}{4}\km\W{3}{4}\km\W{4}{6}
-2 \W{0}{3}\km\W{1}{3}\km\W{2}{4}\km\W{4}{4}\km\W{3}{6}
+\W{0}{3}\km\W{1}{3}\km\W{3}{4}\km\W{4}{4}\km\W{2}{6}
-\frac{1}{2} \W{0}{3}\km\W{2}{3}\km\W{0}{4}\km\W{2}{4}\km\W{6}{6}
\\&
+\W{0}{3}\km\W{2}{3}\km\W{0}{4}\km\W{3}{4}\km\W{5}{6}
-\frac{1}{2} \W{0}{3}\km\W{2}{3}\km\W{0}{4}\km\W{4}{4}\km\W{4}{6}
-2 \W{0}{3}\km\W{2}{3}\km\W{2}{4}\km\W{3}{4}\km\W{3}{6}
+\frac{3}{2} \W{0}{3}\km\W{2}{3}\km\W{2}{4}\km\W{4}{4}\km\W{2}{6}
+\frac{1}{4} \W{0}{3}\km\W{3}{3}\km\W{0}{4}\km\W{2}{4}\km\W{5}{6}
-\W{0}{3}\km\W{2}{3}\km\W{3}{4}\km\W{4}{4}\km\W{1}{6}
+\frac{1}{12} \W{0}{3}\km\W{3}{3}\km\W{0}{4}\km\W{1}{4}\km\W{6}{6}
\\&
-\frac{3}{4} \W{0}{3}\km\W{3}{3}\km\W{0}{4}\km\W{3}{4}\km\W{4}{6}
+\frac{5}{12} \W{0}{3}\km\W{3}{3}\km\W{0}{4}\km\W{4}{4}\km\W{3}{6}
-\frac{1}{2} \W{0}{3}\km\W{3}{3}\km\W{1}{4}\km\W{2}{4}\km\W{4}{6}
+\frac{4}{3} \W{0}{3}\km\W{3}{3}\km\W{1}{4}\km\W{3}{4}\km\W{3}{6}
-\frac{3}{4} \W{0}{3}\km\W{3}{3}\km\W{1}{4}\km\W{4}{4}\km\W{2}{6}
-\frac{1}{2} \W{0}{3}\km\W{3}{3}\km\W{2}{4}\km\W{3}{4}\km\W{2}{6}
+\frac{1}{4} \W{0}{3}\km\W{3}{3}\km\W{2}{4}\km\W{4}{4}\km\W{1}{6}
\\&
+\frac{1}{12} \W{0}{3}\km\W{3}{3}\km\W{3}{4}\km\W{4}{4}\km\W{0}{6}
+\frac{3}{4} \W{1}{3}\km\W{2}{3}\km\W{0}{4}\km\W{1}{4}\km\W{6}{6}
-\frac{3}{4} \W{1}{3}\km\W{2}{3}\km\W{0}{4}\km\W{2}{4}\km\W{5}{6}
-\frac{3}{4} \W{1}{3}\km\W{2}{3}\km\W{0}{4}\km\W{3}{4}\km\W{4}{6}
+\frac{3}{4} \W{1}{3}\km\W{2}{3}\km\W{0}{4}\km\W{4}{4}\km\W{3}{6}
+\frac{3}{2} \W{1}{3}\km\W{2}{3}\km\W{1}{4}\km\W{2}{4}\km\W{4}{6}
-\frac{3}{4} \W{1}{3}\km\W{2}{3}\km\W{1}{4}\km\W{4}{4}\km\W{2}{6}
\\&
+\frac{3}{2} \W{1}{3}\km\W{2}{3}\km\W{2}{4}\km\W{3}{4}\km\W{2}{6}
-\frac{3}{4} \W{1}{3}\km\W{2}{3}\km\W{2}{4}\km\W{4}{4}\km\W{1}{6}
+\frac{3}{4} \W{1}{3}\km\W{2}{3}\km\W{3}{4}\km\W{4}{4}\km\W{0}{6}
-\W{1}{3}\km\W{3}{3}\km\W{0}{4}\km\W{1}{4}\km\W{5}{6}
+\frac{3}{2} \W{1}{3}\km\W{3}{3}\km\W{0}{4}\km\W{2}{4}\km\W{4}{6}
-\frac{1}{2} \W{1}{3}\km\W{3}{3}\km\W{0}{4}\km\W{4}{4}\km\W{2}{6}
-2 \W{1}{3}\km\W{3}{3}\km\W{1}{4}\km\W{2}{4}\km\W{3}{6}
\\&
+\W{1}{3}\km\W{3}{3}\km\W{1}{4}\km\W{4}{4}\km\W{1}{6}
-\frac{1}{2} \W{1}{3}\km\W{3}{3}\km\W{2}{4}\km\W{4}{4}\km\W{0}{6}
+\W{2}{3}\km\W{3}{3}\km\W{0}{4}\km\W{1}{4}\km\W{4}{6}
-2 \W{2}{3}\km\W{3}{3}\km\W{0}{4}\km\W{2}{4}\km\W{3}{6}
+\W{2}{3}\km\W{3}{3}\km\W{0}{4}\km\W{3}{4}\km\W{2}{6}
+3 \W{2}{3}\km\W{3}{3}\km\W{1}{4}\km\W{2}{4}\km\W{2}{6}
-2 \W{2}{3}\km\W{3}{3}\km\W{1}{4}\km\W{3}{4}\km\W{1}{6}
\\&
+\W{2}{3}\km\W{3}{3}\km\W{2}{4}\km\W{3}{4}\km\W{0}{6}
\allowdisplaybreaks
\\
\mathbf{w}_{5} &= 
\frac{1}{6} \W{0}{3}\km\W{1}{3}\km\W{0}{4}\km\W{3}{4}\km\W{6}{6}
-\frac{1}{6} \W{0}{3}\km\W{1}{3}\km\W{0}{4}\km\W{4}{4}\km\W{5}{6}
-\frac{1}{6} \W{0}{3}\km\W{1}{3}\km\W{1}{4}\km\W{2}{4}\km\W{6}{6}
-\frac{1}{3} \W{0}{3}\km\W{1}{3}\km\W{1}{4}\km\W{3}{4}\km\W{5}{6}
+\frac{1}{2} \W{0}{3}\km\W{1}{3}\km\W{1}{4}\km\W{4}{4}\km\W{4}{6}
+\frac{1}{2} \W{0}{3}\km\W{1}{3}\km\W{2}{4}\km\W{3}{4}\km\W{4}{6}
-\frac{2}{3} \W{0}{3}\km\W{1}{3}\km\W{2}{4}\km\W{4}{4}\km\W{3}{6}
\\&
+\frac{1}{3} \W{0}{3}\km\W{1}{3}\km\W{3}{4}\km\W{4}{4}\km\W{2}{6}
-\frac{1}{6} \W{0}{3}\km\W{2}{3}\km\W{0}{4}\km\W{2}{4}\km\W{6}{6}
+\frac{1}{6} \W{0}{3}\km\W{2}{3}\km\W{0}{4}\km\W{4}{4}\km\W{4}{6}
+\W{0}{3}\km\W{2}{3}\km\W{1}{4}\km\W{2}{4}\km\W{5}{6}
-\frac{2}{3} \W{0}{3}\km\W{2}{3}\km\W{1}{4}\km\W{3}{4}\km\W{4}{6}
-\frac{1}{3} \W{0}{3}\km\W{2}{3}\km\W{1}{4}\km\W{4}{4}\km\W{3}{6}
+\frac{1}{3} \W{0}{3}\km\W{2}{3}\km\W{2}{4}\km\W{3}{4}\km\W{3}{6}
\\&
+\frac{1}{2} \W{0}{3}\km\W{2}{3}\km\W{2}{4}\km\W{4}{4}\km\W{2}{6}
+\frac{5}{24} \W{0}{3}\km\W{3}{3}\km\W{0}{4}\km\W{2}{4}\km\W{5}{6}
-\frac{1}{3} \W{0}{3}\km\W{2}{3}\km\W{3}{4}\km\W{4}{4}\km\W{1}{6}
-\frac{1}{72} \W{0}{3}\km\W{3}{3}\km\W{0}{4}\km\W{1}{4}\km\W{6}{6}
-\frac{5}{24} \W{0}{3}\km\W{3}{3}\km\W{0}{4}\km\W{3}{4}\km\W{4}{6}
+\frac{1}{72} \W{0}{3}\km\W{3}{3}\km\W{0}{4}\km\W{4}{4}\km\W{3}{6}
-\frac{11}{12} \W{0}{3}\km\W{3}{3}\km\W{1}{4}\km\W{2}{4}\km\W{4}{6}
\\&
+\frac{10}{9} \W{0}{3}\km\W{3}{3}\km\W{1}{4}\km\W{3}{4}\km\W{3}{6}
-\frac{5}{24} \W{0}{3}\km\W{3}{3}\km\W{1}{4}\km\W{4}{4}\km\W{2}{6}
-\frac{11}{12} \W{0}{3}\km\W{3}{3}\km\W{2}{4}\km\W{3}{4}\km\W{2}{6}
+\frac{5}{24} \W{0}{3}\km\W{3}{3}\km\W{2}{4}\km\W{4}{4}\km\W{1}{6}
-\frac{1}{72} \W{0}{3}\km\W{3}{3}\km\W{3}{4}\km\W{4}{4}\km\W{0}{6}
+\frac{3}{8} \W{1}{3}\km\W{2}{3}\km\W{0}{4}\km\W{1}{4}\km\W{6}{6}
-\frac{5}{8} \W{1}{3}\km\W{2}{3}\km\W{0}{4}\km\W{2}{4}\km\W{5}{6}
\\&
+\frac{5}{8} \W{1}{3}\km\W{2}{3}\km\W{0}{4}\km\W{3}{4}\km\W{4}{6}
-\frac{3}{8} \W{1}{3}\km\W{2}{3}\km\W{0}{4}\km\W{4}{4}\km\W{3}{6}
-\frac{1}{4} \W{1}{3}\km\W{2}{3}\km\W{1}{4}\km\W{2}{4}\km\W{4}{6}
+\frac{5}{8} \W{1}{3}\km\W{2}{3}\km\W{1}{4}\km\W{4}{4}\km\W{2}{6}
-\frac{1}{4} \W{1}{3}\km\W{2}{3}\km\W{2}{4}\km\W{3}{4}\km\W{2}{6}
-\frac{5}{8} \W{1}{3}\km\W{2}{3}\km\W{2}{4}\km\W{4}{4}\km\W{1}{6}
+\frac{3}{8} \W{1}{3}\km\W{2}{3}\km\W{3}{4}\km\W{4}{4}\km\W{0}{6}
\\&
-\frac{1}{3} \W{1}{3}\km\W{3}{3}\km\W{0}{4}\km\W{1}{4}\km\W{5}{6}
+\frac{1}{2} \W{1}{3}\km\W{3}{3}\km\W{0}{4}\km\W{2}{4}\km\W{4}{6}
-\frac{1}{3} \W{1}{3}\km\W{3}{3}\km\W{0}{4}\km\W{3}{4}\km\W{3}{6}
+\frac{1}{6} \W{1}{3}\km\W{3}{3}\km\W{0}{4}\km\W{4}{4}\km\W{2}{6}
+\frac{1}{3} \W{1}{3}\km\W{3}{3}\km\W{1}{4}\km\W{2}{4}\km\W{3}{6}
-\frac{2}{3} \W{1}{3}\km\W{3}{3}\km\W{1}{4}\km\W{3}{4}\km\W{2}{6}
+\W{1}{3}\km\W{3}{3}\km\W{2}{4}\km\W{3}{4}\km\W{1}{6}
\\&
-\frac{1}{6} \W{1}{3}\km\W{3}{3}\km\W{2}{4}\km\W{4}{4}\km\W{0}{6}
+\frac{1}{3} \W{2}{3}\km\W{3}{3}\km\W{0}{4}\km\W{1}{4}\km\W{4}{6}
-\frac{2}{3} \W{2}{3}\km\W{3}{3}\km\W{0}{4}\km\W{2}{4}\km\W{3}{6}
+\frac{1}{2} \W{2}{3}\km\W{3}{3}\km\W{0}{4}\km\W{3}{4}\km\W{2}{6}
-\frac{1}{6} \W{2}{3}\km\W{3}{3}\km\W{0}{4}\km\W{4}{4}\km\W{1}{6}
+\frac{1}{2} \W{2}{3}\km\W{3}{3}\km\W{1}{4}\km\W{2}{4}\km\W{2}{6}
-\frac{1}{3} \W{2}{3}\km\W{3}{3}\km\W{1}{4}\km\W{3}{4}\km\W{1}{6}
\\&
+\frac{1}{6} \W{2}{3}\km\W{3}{3}\km\W{1}{4}\km\W{4}{4}\km\W{0}{6}
-\frac{1}{6} \W{2}{3}\km\W{3}{3}\km\W{2}{4}\km\W{3}{4}\km\W{0}{6}
\\
\mathbf{w}_{6} &= 
\frac{1}{2} \W{0}{3}\km\W{1}{3}\km\W{2}{4}\km\W{2}{5}\km\W{5}{5}
+\frac{3}{2} \W{0}{3}\km\W{1}{3}\km\W{2}{4}\km\W{3}{5}\km\W{4}{5}
-\W{0}{3}\km\W{1}{3}\km\W{3}{4}\km\W{2}{5}\km\W{4}{5}
+\frac{3}{4} \W{1}{3}\km\W{2}{3}\km\W{0}{4}\km\W{2}{5}\km\W{5}{5}
-\frac{3}{4} \W{1}{3}\km\W{2}{3}\km\W{0}{4}\km\W{3}{5}\km\W{4}{5}
-\W{2}{3}\km\W{3}{3}\km\W{1}{4}\km\W{1}{5}\km\W{3}{5}
+\frac{1}{2} \W{2}{3}\km\W{3}{3}\km\W{2}{4}\km\W{0}{5}\km\W{3}{5}
\\&
+\W{1}{3}\km\W{3}{3}\km\W{3}{4}\km\W{0}{5}\km\W{3}{5}
+\W{1}{3}\km\W{3}{3}\km\W{3}{4}\km\W{1}{5}\km\W{2}{5}
-\frac{1}{2} \W{1}{3}\km\W{3}{3}\km\W{4}{4}\km\W{0}{5}\km\W{2}{5}
+\W{0}{3}\km\W{2}{3}\km\W{1}{4}\km\W{2}{5}\km\W{5}{5}
-\W{0}{3}\km\W{1}{3}\km\W{1}{4}\km\W{3}{5}\km\W{5}{5}
+3 \W{0}{3}\km\W{3}{3}\km\W{2}{4}\km\W{2}{5}\km\W{3}{5}
+\frac{1}{6} \W{0}{3}\km\W{3}{3}\km\W{3}{4}\km\W{0}{5}\km\W{4}{5}
\\&
-\frac{4}{3} \W{0}{3}\km\W{3}{3}\km\W{3}{4}\km\W{1}{5}\km\W{3}{5}
-\frac{1}{2} \W{0}{3}\km\W{2}{3}\km\W{0}{4}\km\W{3}{5}\km\W{5}{5}
+\frac{3}{4} \W{1}{3}\km\W{2}{3}\km\W{4}{4}\km\W{0}{5}\km\W{3}{5}
-\frac{3}{4} \W{1}{3}\km\W{2}{3}\km\W{4}{4}\km\W{1}{5}\km\W{2}{5}
-\frac{1}{2} \W{1}{3}\km\W{3}{3}\km\W{0}{4}\km\W{2}{5}\km\W{4}{5}
-\frac{1}{6} \W{1}{3}\km\W{3}{3}\km\W{1}{4}\km\W{0}{5}\km\W{5}{5}
+\frac{11}{6} \W{1}{3}\km\W{3}{3}\km\W{1}{4}\km\W{1}{5}\km\W{4}{5}
\\&
-\frac{2}{3} \W{1}{3}\km\W{3}{3}\km\W{1}{4}\km\W{2}{5}\km\W{3}{5}
+\frac{3}{2} \W{2}{3}\km\W{3}{3}\km\W{2}{4}\km\W{1}{5}\km\W{2}{5}
-\W{2}{3}\km\W{3}{3}\km\W{3}{4}\km\W{0}{5}\km\W{2}{5}
+\frac{1}{2} \W{2}{3}\km\W{3}{3}\km\W{4}{4}\km\W{0}{5}\km\W{1}{5}
-\frac{1}{12} \W{0}{3}\km\W{3}{3}\km\W{4}{4}\km\W{0}{5}\km\W{3}{5}
+\frac{3}{4} \W{0}{3}\km\W{3}{3}\km\W{4}{4}\km\W{1}{5}\km\W{2}{5}
-\frac{1}{6} \W{0}{3}\km\W{2}{3}\km\W{3}{4}\km\W{0}{5}\km\W{5}{5}
\\&
+\frac{11}{6} \W{0}{3}\km\W{2}{3}\km\W{3}{4}\km\W{1}{5}\km\W{4}{5}
-\frac{2}{3} \W{0}{3}\km\W{2}{3}\km\W{3}{4}\km\W{2}{5}\km\W{3}{5}
-\frac{1}{2} \W{0}{3}\km\W{2}{3}\km\W{4}{4}\km\W{1}{5}\km\W{3}{5}
-\frac{1}{12} \W{0}{3}\km\W{3}{3}\km\W{0}{4}\km\W{2}{5}\km\W{5}{5}
+\frac{1}{2} \W{0}{3}\km\W{1}{3}\km\W{0}{4}\km\W{4}{5}\km\W{5}{5}
-\frac{3}{2} \W{1}{3}\km\W{2}{3}\km\W{3}{4}\km\W{0}{5}\km\W{4}{5}
-\frac{3}{2} \W{1}{3}\km\W{2}{3}\km\W{1}{4}\km\W{1}{5}\km\W{5}{5}
\\&
+\W{1}{3}\km\W{2}{3}\km\W{2}{4}\km\W{0}{5}\km\W{5}{5}
+\W{1}{3}\km\W{2}{3}\km\W{2}{4}\km\W{1}{5}\km\W{4}{5}
+\W{1}{3}\km\W{2}{3}\km\W{2}{4}\km\W{2}{5}\km\W{3}{5}
+\frac{1}{12} \W{2}{3}\km\W{3}{3}\km\W{0}{4}\km\W{0}{5}\km\W{5}{5}
-\frac{5}{12} \W{2}{3}\km\W{3}{3}\km\W{0}{4}\km\W{1}{5}\km\W{4}{5}
+\frac{4}{3} \W{2}{3}\km\W{3}{3}\km\W{0}{4}\km\W{2}{5}\km\W{3}{5}
-\frac{1}{2} \W{1}{3}\km\W{3}{3}\km\W{2}{4}\km\W{0}{5}\km\W{4}{5}
\\&
-2 \W{1}{3}\km\W{3}{3}\km\W{2}{4}\km\W{1}{5}\km\W{3}{5}
+\W{0}{3}\km\W{2}{3}\km\W{1}{4}\km\W{3}{5}\km\W{4}{5}
-\frac{1}{2} \W{0}{3}\km\W{2}{3}\km\W{2}{4}\km\W{1}{5}\km\W{5}{5}
-2 \W{0}{3}\km\W{2}{3}\km\W{2}{4}\km\W{2}{5}\km\W{4}{5}
+\frac{3}{4} \W{0}{3}\km\W{3}{3}\km\W{0}{4}\km\W{3}{5}\km\W{4}{5}
+\frac{1}{6} \W{0}{3}\km\W{3}{3}\km\W{1}{4}\km\W{1}{5}\km\W{5}{5}
-\frac{4}{3} \W{0}{3}\km\W{3}{3}\km\W{1}{4}\km\W{2}{5}\km\W{4}{5}
\\&
+\frac{1}{12} \W{0}{3}\km\W{1}{3}\km\W{4}{4}\km\W{0}{5}\km\W{5}{5}
-\frac{5}{12} \W{0}{3}\km\W{1}{3}\km\W{4}{4}\km\W{1}{5}\km\W{4}{5}
+\frac{4}{3} \W{0}{3}\km\W{1}{3}\km\W{4}{4}\km\W{2}{5}\km\W{3}{5}
\allowdisplaybreaks
\\
\mathbf{w}_{7} &= 
-3 \W{0}{3}\km\W{1}{3}\km\W{2}{4}\km\W{2}{5}\km\W{5}{5}
+3 \W{0}{3}\km\W{1}{3}\km\W{2}{4}\km\W{3}{5}\km\W{4}{5}
+2 \W{0}{3}\km\W{1}{3}\km\W{3}{4}\km\W{1}{5}\km\W{5}{5}
-2 \W{0}{3}\km\W{1}{3}\km\W{3}{4}\km\W{2}{5}\km\W{4}{5}
-3 \W{1}{3}\km\W{2}{3}\km\W{0}{4}\km\W{2}{5}\km\W{5}{5}
+6 \W{1}{3}\km\W{2}{3}\km\W{0}{4}\km\W{3}{5}\km\W{4}{5}
+2 \W{2}{3}\km\W{3}{3}\km\W{1}{4}\km\W{0}{5}\km\W{4}{5}
\\&
-2 \W{2}{3}\km\W{3}{3}\km\W{1}{4}\km\W{1}{5}\km\W{3}{5}
-3 \W{2}{3}\km\W{3}{3}\km\W{2}{4}\km\W{0}{5}\km\W{3}{5}
-8 \W{1}{3}\km\W{3}{3}\km\W{3}{4}\km\W{1}{5}\km\W{2}{5}
+\W{1}{3}\km\W{3}{3}\km\W{4}{4}\km\W{0}{5}\km\W{2}{5}
+2 \W{0}{3}\km\W{1}{3}\km\W{1}{4}\km\W{3}{5}\km\W{5}{5}
+\frac{3}{5} \W{0}{3}\km\W{3}{3}\km\W{2}{4}\km\W{0}{5}\km\W{5}{5}
-9 \W{0}{3}\km\W{3}{3}\km\W{2}{4}\km\W{2}{5}\km\W{3}{5}
\\&
-\W{0}{3}\km\W{3}{3}\km\W{3}{4}\km\W{0}{5}\km\W{4}{5}
+4 \W{0}{3}\km\W{3}{3}\km\W{3}{4}\km\W{1}{5}\km\W{3}{5}
+\W{0}{3}\km\W{2}{3}\km\W{0}{4}\km\W{3}{5}\km\W{5}{5}
-3 \W{1}{3}\km\W{2}{3}\km\W{4}{4}\km\W{0}{5}\km\W{3}{5}
+6 \W{1}{3}\km\W{2}{3}\km\W{4}{4}\km\W{1}{5}\km\W{2}{5}
+\W{1}{3}\km\W{3}{3}\km\W{0}{4}\km\W{1}{5}\km\W{5}{5}
-\W{1}{3}\km\W{3}{3}\km\W{0}{4}\km\W{2}{5}\km\W{4}{5}
\\&
-\frac{1}{5} \W{1}{3}\km\W{3}{3}\km\W{1}{4}\km\W{0}{5}\km\W{5}{5}
-3 \W{1}{3}\km\W{3}{3}\km\W{1}{4}\km\W{1}{5}\km\W{4}{5}
+2 \W{1}{3}\km\W{3}{3}\km\W{1}{4}\km\W{2}{5}\km\W{3}{5}
+3 \W{2}{3}\km\W{3}{3}\km\W{2}{4}\km\W{1}{5}\km\W{2}{5}
+2 \W{2}{3}\km\W{3}{3}\km\W{3}{4}\km\W{0}{5}\km\W{2}{5}
-\W{2}{3}\km\W{3}{3}\km\W{4}{4}\km\W{0}{5}\km\W{1}{5}
-\W{0}{3}\km\W{3}{3}\km\W{4}{4}\km\W{1}{5}\km\W{2}{5}
\\&
-\frac{1}{5} \W{0}{3}\km\W{2}{3}\km\W{3}{4}\km\W{0}{5}\km\W{5}{5}
-3 \W{0}{3}\km\W{2}{3}\km\W{3}{4}\km\W{1}{5}\km\W{4}{5}
+2 \W{0}{3}\km\W{2}{3}\km\W{3}{4}\km\W{2}{5}\km\W{3}{5}
+\W{0}{3}\km\W{2}{3}\km\W{4}{4}\km\W{0}{5}\km\W{4}{5}
-\W{0}{3}\km\W{2}{3}\km\W{4}{4}\km\W{1}{5}\km\W{3}{5}
-\W{0}{3}\km\W{1}{3}\km\W{0}{4}\km\W{4}{5}\km\W{5}{5}
+3 \W{1}{3}\km\W{2}{3}\km\W{3}{4}\km\W{0}{5}\km\W{4}{5}
\\&
+3 \W{1}{3}\km\W{2}{3}\km\W{1}{4}\km\W{1}{5}\km\W{5}{5}
-\frac{6}{5} \W{1}{3}\km\W{2}{3}\km\W{2}{4}\km\W{0}{5}\km\W{5}{5}
-3 \W{1}{3}\km\W{2}{3}\km\W{2}{4}\km\W{1}{5}\km\W{4}{5}
-3 \W{1}{3}\km\W{2}{3}\km\W{2}{4}\km\W{2}{5}\km\W{3}{5}
-\frac{2}{5} \W{2}{3}\km\W{3}{3}\km\W{0}{4}\km\W{0}{5}\km\W{5}{5}
+\W{2}{3}\km\W{3}{3}\km\W{0}{4}\km\W{2}{5}\km\W{3}{5}
+6 \W{1}{3}\km\W{3}{3}\km\W{2}{4}\km\W{1}{5}\km\W{3}{5}
\\&
-8 \W{0}{3}\km\W{2}{3}\km\W{1}{4}\km\W{3}{5}\km\W{4}{5}
+6 \W{0}{3}\km\W{2}{3}\km\W{2}{4}\km\W{2}{5}\km\W{4}{5}
-\W{0}{3}\km\W{3}{3}\km\W{0}{4}\km\W{3}{5}\km\W{4}{5}
-\W{0}{3}\km\W{3}{3}\km\W{1}{4}\km\W{1}{5}\km\W{5}{5}
+4 \W{0}{3}\km\W{3}{3}\km\W{1}{4}\km\W{2}{5}\km\W{4}{5}
-\frac{2}{5} \W{0}{3}\km\W{1}{3}\km\W{4}{4}\km\W{0}{5}\km\W{5}{5}
+\W{0}{3}\km\W{1}{3}\km\W{4}{4}\km\W{2}{5}\km\W{3}{5}
\\
\mathbf{w}_{8} &= 
\frac{2}{5} \W{0}{3}\km\W{2}{3}\km\W{3}{4}\km\W{0}{5}\km\W{5}{5}
-2 \W{0}{3}\km\W{2}{3}\km\W{3}{4}\km\W{1}{5}\km\W{4}{5}
+4 \W{0}{3}\km\W{2}{3}\km\W{3}{4}\km\W{2}{5}\km\W{3}{5}
-\frac{1}{5} \W{0}{3}\km\W{3}{3}\km\W{2}{4}\km\W{0}{5}\km\W{5}{5}
-2 \W{0}{3}\km\W{3}{3}\km\W{2}{4}\km\W{2}{5}\km\W{3}{5}
-\frac{3}{5} \W{1}{3}\km\W{2}{3}\km\W{2}{4}\km\W{0}{5}\km\W{5}{5}
+3 \W{1}{3}\km\W{2}{3}\km\W{2}{4}\km\W{1}{5}\km\W{4}{5}
\\&
-6 \W{1}{3}\km\W{2}{3}\km\W{2}{4}\km\W{2}{5}\km\W{3}{5}
+\frac{2}{5} \W{1}{3}\km\W{3}{3}\km\W{1}{4}\km\W{0}{5}\km\W{5}{5}
-2 \W{1}{3}\km\W{3}{3}\km\W{1}{4}\km\W{1}{5}\km\W{4}{5}
+4 \W{1}{3}\km\W{3}{3}\km\W{1}{4}\km\W{2}{5}\km\W{3}{5}
-\frac{1}{5} \W{2}{3}\km\W{3}{3}\km\W{0}{4}\km\W{0}{5}\km\W{5}{5}
+\W{2}{3}\km\W{3}{3}\km\W{0}{4}\km\W{1}{5}\km\W{4}{5}
-2 \W{2}{3}\km\W{3}{3}\km\W{0}{4}\km\W{2}{5}\km\W{3}{5}
\\&
-\frac{1}{5} \W{0}{3}\km\W{1}{3}\km\W{4}{4}\km\W{0}{5}\km\W{5}{5}
+\W{0}{3}\km\W{1}{3}\km\W{4}{4}\km\W{1}{5}\km\W{4}{5}
-2 \W{0}{3}\km\W{1}{3}\km\W{4}{4}\km\W{2}{5}\km\W{3}{5}
+\W{0}{3}\km\W{3}{3}\km\W{2}{4}\km\W{1}{5}\km\W{4}{5}
\\
\mathbf{w}_{9} &= 
\frac{1}{4} \W{1}{3}\km\W{2}{3}\km\W{0}{4}\km\W{2}{5}\km\W{5}{5}
-\frac{3}{4} \W{1}{3}\km\W{2}{3}\km\W{0}{4}\km\W{3}{5}\km\W{4}{5}
-\frac{1}{10} \W{0}{3}\km\W{3}{3}\km\W{2}{4}\km\W{0}{5}\km\W{5}{5}
+\frac{1}{2} \W{0}{3}\km\W{3}{3}\km\W{2}{4}\km\W{2}{5}\km\W{3}{5}
+\frac{1}{6} \W{0}{3}\km\W{3}{3}\km\W{3}{4}\km\W{0}{5}\km\W{4}{5}
-\frac{1}{3} \W{0}{3}\km\W{3}{3}\km\W{3}{4}\km\W{1}{5}\km\W{3}{5}
+\frac{1}{4} \W{1}{3}\km\W{2}{3}\km\W{4}{4}\km\W{0}{5}\km\W{3}{5}
\\&
-\frac{3}{4} \W{1}{3}\km\W{2}{3}\km\W{4}{4}\km\W{1}{5}\km\W{2}{5}
+\frac{1}{30} \W{1}{3}\km\W{3}{3}\km\W{1}{4}\km\W{0}{5}\km\W{5}{5}
-\frac{1}{6} \W{1}{3}\km\W{3}{3}\km\W{1}{4}\km\W{1}{5}\km\W{4}{5}
+\frac{1}{3} \W{1}{3}\km\W{3}{3}\km\W{1}{4}\km\W{2}{5}\km\W{3}{5}
-\frac{1}{12} \W{0}{3}\km\W{3}{3}\km\W{4}{4}\km\W{0}{5}\km\W{3}{5}
+\frac{1}{4} \W{0}{3}\km\W{3}{3}\km\W{4}{4}\km\W{1}{5}\km\W{2}{5}
+\frac{1}{30} \W{0}{3}\km\W{2}{3}\km\W{3}{4}\km\W{0}{5}\km\W{5}{5}
\\&
-\frac{1}{6} \W{0}{3}\km\W{2}{3}\km\W{3}{4}\km\W{1}{5}\km\W{4}{5}
+\frac{1}{3} \W{0}{3}\km\W{2}{3}\km\W{3}{4}\km\W{2}{5}\km\W{3}{5}
-\frac{1}{12} \W{0}{3}\km\W{3}{3}\km\W{0}{4}\km\W{2}{5}\km\W{5}{5}
-\frac{1}{2} \W{1}{3}\km\W{2}{3}\km\W{3}{4}\km\W{0}{5}\km\W{4}{5}
+\W{1}{3}\km\W{2}{3}\km\W{3}{4}\km\W{1}{5}\km\W{3}{5}
-\frac{1}{2} \W{1}{3}\km\W{2}{3}\km\W{1}{4}\km\W{1}{5}\km\W{5}{5}
+\W{1}{3}\km\W{2}{3}\km\W{1}{4}\km\W{2}{5}\km\W{4}{5}
\\&
+\frac{1}{5} \W{1}{3}\km\W{2}{3}\km\W{2}{4}\km\W{0}{5}\km\W{5}{5}
+\frac{1}{2} \W{1}{3}\km\W{2}{3}\km\W{2}{4}\km\W{1}{5}\km\W{4}{5}
-\frac{5}{2} \W{1}{3}\km\W{2}{3}\km\W{2}{4}\km\W{2}{5}\km\W{3}{5}
-\frac{1}{60} \W{2}{3}\km\W{3}{3}\km\W{0}{4}\km\W{0}{5}\km\W{5}{5}
+\frac{1}{12} \W{2}{3}\km\W{3}{3}\km\W{0}{4}\km\W{1}{5}\km\W{4}{5}
-\frac{1}{6} \W{2}{3}\km\W{3}{3}\km\W{0}{4}\km\W{2}{5}\km\W{3}{5}
+\frac{1}{4} \W{0}{3}\km\W{3}{3}\km\W{0}{4}\km\W{3}{5}\km\W{4}{5}
\\&
+\frac{1}{6} \W{0}{3}\km\W{3}{3}\km\W{1}{4}\km\W{1}{5}\km\W{5}{5}
-\frac{1}{3} \W{0}{3}\km\W{3}{3}\km\W{1}{4}\km\W{2}{5}\km\W{4}{5}
-\frac{1}{60} \W{0}{3}\km\W{1}{3}\km\W{4}{4}\km\W{0}{5}\km\W{5}{5}
+\frac{1}{12} \W{0}{3}\km\W{1}{3}\km\W{4}{4}\km\W{1}{5}\km\W{4}{5}
-\frac{1}{6} \W{0}{3}\km\W{1}{3}\km\W{4}{4}\km\W{2}{5}\km\W{3}{5}
\allowdisplaybreaks
\\
\mathbf{w}_{10} &= 
\frac{2}{3} \W{3}{3}\km\W{1}{4}\km\W{2}{4}\km\W{3}{4}\km\W{1}{5}
-4 \W{2}{3}\km\W{1}{4}\km\W{2}{4}\km\W{3}{4}\km\W{2}{5}
+\frac{3}{2} \W{2}{3}\km\W{1}{4}\km\W{2}{4}\km\W{4}{4}\km\W{1}{5}
+\frac{1}{3} \W{1}{3}\km\W{0}{4}\km\W{3}{4}\km\W{4}{4}\km\W{2}{5}
+4 \W{1}{3}\km\W{1}{4}\km\W{2}{4}\km\W{3}{4}\km\W{3}{5}
+\frac{1}{3} \W{3}{3}\km\W{0}{4}\km\W{2}{4}\km\W{3}{4}\km\W{2}{5}
-\frac{1}{3} \W{3}{3}\km\W{0}{4}\km\W{2}{4}\km\W{4}{4}\km\W{1}{5}
\\&
+\frac{1}{6} \W{3}{3}\km\W{0}{4}\km\W{3}{4}\km\W{4}{4}\km\W{0}{5}
-\W{1}{3}\km\W{1}{4}\km\W{2}{4}\km\W{4}{4}\km\W{2}{5}
-\frac{1}{3} \W{1}{3}\km\W{1}{4}\km\W{3}{4}\km\W{4}{4}\km\W{1}{5}
+\frac{1}{2} \W{3}{3}\km\W{0}{4}\km\W{1}{4}\km\W{2}{4}\km\W{4}{5}
-\frac{1}{2} \W{2}{3}\km\W{0}{4}\km\W{1}{4}\km\W{2}{4}\km\W{5}{5}
+\frac{1}{3} \W{2}{3}\km\W{0}{4}\km\W{1}{4}\km\W{3}{4}\km\W{4}{5}
-\frac{1}{3} \W{0}{3}\km\W{1}{4}\km\W{2}{4}\km\W{4}{4}\km\W{3}{5}
\\&
+\frac{2}{3} \W{0}{3}\km\W{1}{4}\km\W{3}{4}\km\W{4}{4}\km\W{2}{5}
-\frac{1}{6} \W{2}{3}\km\W{0}{4}\km\W{3}{4}\km\W{4}{4}\km\W{1}{5}
-\frac{2}{3} \W{3}{3}\km\W{0}{4}\km\W{1}{4}\km\W{3}{4}\km\W{3}{5}
+\frac{1}{3} \W{3}{3}\km\W{0}{4}\km\W{1}{4}\km\W{4}{4}\km\W{2}{5}
-\frac{1}{3} \W{2}{3}\km\W{0}{4}\km\W{1}{4}\km\W{4}{4}\km\W{3}{5}
+\W{2}{3}\km\W{0}{4}\km\W{2}{4}\km\W{3}{4}\km\W{3}{5}
+\frac{1}{3} \W{1}{3}\km\W{0}{4}\km\W{1}{4}\km\W{3}{4}\km\W{5}{5}
\\&
+\frac{1}{3} \W{0}{3}\km\W{0}{4}\km\W{2}{4}\km\W{4}{4}\km\W{4}{5}
-\frac{1}{2} \W{0}{3}\km\W{2}{4}\km\W{3}{4}\km\W{4}{4}\km\W{1}{5}
+\frac{1}{6} \W{1}{3}\km\W{0}{4}\km\W{1}{4}\km\W{4}{4}\km\W{4}{5}
-\frac{3}{2} \W{1}{3}\km\W{0}{4}\km\W{2}{4}\km\W{3}{4}\km\W{4}{5}
+\frac{1}{2} \W{1}{3}\km\W{2}{4}\km\W{3}{4}\km\W{4}{4}\km\W{0}{5}
-\frac{1}{3} \W{0}{3}\km\W{0}{4}\km\W{3}{4}\km\W{4}{4}\km\W{3}{5}
-\frac{2}{3} \W{0}{3}\km\W{1}{4}\km\W{2}{4}\km\W{3}{4}\km\W{4}{5}
\\&
-\frac{1}{6} \W{0}{3}\km\W{0}{4}\km\W{1}{4}\km\W{4}{4}\km\W{5}{5}
+\frac{1}{6} \W{0}{3}\km\W{0}{4}\km\W{2}{4}\km\W{3}{4}\km\W{5}{5}
-\frac{1}{6} \W{3}{3}\km\W{1}{4}\km\W{2}{4}\km\W{4}{4}\km\W{0}{5}
-\frac{1}{3} \W{2}{3}\km\W{1}{4}\km\W{3}{4}\km\W{4}{4}\km\W{0}{5}
\\
\mathbf{w}_{11} &= 
\frac{4}{3} \W{3}{3}\km\W{1}{4}\km\W{2}{4}\km\W{3}{4}\km\W{1}{5}
-\W{2}{3}\km\W{1}{4}\km\W{2}{4}\km\W{4}{4}\km\W{1}{5}
+\frac{1}{3} \W{1}{3}\km\W{0}{4}\km\W{3}{4}\km\W{4}{4}\km\W{2}{5}
-\frac{4}{3} \W{3}{3}\km\W{0}{4}\km\W{2}{4}\km\W{3}{4}\km\W{2}{5}
+\frac{1}{3} \W{3}{3}\km\W{0}{4}\km\W{2}{4}\km\W{4}{4}\km\W{1}{5}
+\frac{2}{3} \W{1}{3}\km\W{1}{4}\km\W{3}{4}\km\W{4}{4}\km\W{1}{5}
-\W{3}{3}\km\W{0}{4}\km\W{1}{4}\km\W{2}{4}\km\W{4}{5}
\\&
+\W{2}{3}\km\W{0}{4}\km\W{1}{4}\km\W{2}{4}\km\W{5}{5}
-\frac{2}{3} \W{2}{3}\km\W{0}{4}\km\W{1}{4}\km\W{3}{4}\km\W{4}{5}
+\frac{4}{3} \W{0}{3}\km\W{1}{4}\km\W{2}{4}\km\W{4}{4}\km\W{3}{5}
-\frac{4}{3} \W{0}{3}\km\W{1}{4}\km\W{3}{4}\km\W{4}{4}\km\W{2}{5}
-\frac{2}{3} \W{2}{3}\km\W{0}{4}\km\W{3}{4}\km\W{4}{4}\km\W{1}{5}
+\frac{4}{3} \W{3}{3}\km\W{0}{4}\km\W{1}{4}\km\W{3}{4}\km\W{3}{5}
-\frac{1}{3} \W{3}{3}\km\W{0}{4}\km\W{1}{4}\km\W{4}{4}\km\W{2}{5}
\\&
-\frac{1}{3} \W{2}{3}\km\W{0}{4}\km\W{1}{4}\km\W{4}{4}\km\W{3}{5}
+\W{2}{3}\km\W{0}{4}\km\W{2}{4}\km\W{4}{4}\km\W{2}{5}
-\frac{2}{3} \W{1}{3}\km\W{0}{4}\km\W{1}{4}\km\W{3}{4}\km\W{5}{5}
-\frac{1}{3} \W{0}{3}\km\W{0}{4}\km\W{2}{4}\km\W{4}{4}\km\W{4}{5}
+\W{0}{3}\km\W{2}{4}\km\W{3}{4}\km\W{4}{4}\km\W{1}{5}
+\frac{2}{3} \W{1}{3}\km\W{0}{4}\km\W{1}{4}\km\W{4}{4}\km\W{4}{5}
+\W{1}{3}\km\W{0}{4}\km\W{2}{4}\km\W{3}{4}\km\W{4}{5}
\\&
-\W{1}{3}\km\W{0}{4}\km\W{2}{4}\km\W{4}{4}\km\W{3}{5}
-\W{1}{3}\km\W{2}{4}\km\W{3}{4}\km\W{4}{4}\km\W{0}{5}
+\frac{1}{3} \W{0}{3}\km\W{0}{4}\km\W{3}{4}\km\W{4}{4}\km\W{3}{5}
-\frac{4}{3} \W{0}{3}\km\W{1}{4}\km\W{2}{4}\km\W{3}{4}\km\W{4}{5}
+\frac{1}{3} \W{0}{3}\km\W{0}{4}\km\W{2}{4}\km\W{3}{4}\km\W{5}{5}
-\frac{1}{3} \W{3}{3}\km\W{1}{4}\km\W{2}{4}\km\W{4}{4}\km\W{0}{5}
+\frac{2}{3} \W{2}{3}\km\W{1}{4}\km\W{3}{4}\km\W{4}{4}\km\W{0}{5}
\\
\mathbf{w}_{12} &= 
\W{0}{4}\km\W{1}{4}\km\W{2}{4}\km\W{3}{4}\km\W{4}{4}
\end{align*}

\subsection{Basis of $\ds \CGF{6}{2}{0}{10}$} 
\begin{align*}
\mathbf{r}_{1} &=
-\frac{1}{5} \W{0}{3}\km\W{1}{3}\km\W{2}{3}\km\W{3}{3}\km\W{0}{5}\km\W{5}{5}
+\W{0}{3}\km\W{1}{3}\km\W{2}{3}\km\W{3}{3}\km\W{1}{5}\km\W{4}{5}
-2 \W{0}{3}\km\W{1}{3}\km\W{2}{3}\km\W{3}{3}\km\W{2}{5}\km\W{3}{5}
\\ 
\mathbf{r}_{2} &=
-3 \W{0}{3}\km\W{1}{3}\km\W{2}{3}\km\W{1}{4}\km\W{2}{4}\km\W{5}{5}
+\frac{1}{2} \W{0}{3}\km\W{1}{3}\km\W{2}{3}\km\W{0}{4}\km\W{3}{4}\km\W{5}{5}
-\frac{1}{2} \W{0}{3}\km\W{1}{3}\km\W{2}{3}\km\W{0}{4}\km\W{4}{4}\km\W{4}{5}
+4 \W{0}{3}\km\W{1}{3}\km\W{2}{3}\km\W{1}{4}\km\W{3}{4}\km\W{4}{5}
-\W{0}{3}\km\W{1}{3}\km\W{2}{3}\km\W{1}{4}\km\W{4}{4}\km\W{3}{5}
-6 \W{0}{3}\km\W{1}{3}\km\W{2}{3}\km\W{2}{4}\km\W{3}{4}\km\W{3}{5}
\\&
+3 \W{0}{3}\km\W{1}{3}\km\W{2}{3}\km\W{2}{4}\km\W{4}{4}\km\W{2}{5}
-\frac{3}{2} \W{0}{3}\km\W{1}{3}\km\W{2}{3}\km\W{3}{4}\km\W{4}{4}\km\W{1}{5}
+\frac{1}{2} \W{0}{3}\km\W{1}{3}\km\W{3}{3}\km\W{0}{4}\km\W{2}{4}\km\W{5}{5}
-\frac{3}{2} \W{0}{3}\km\W{1}{3}\km\W{3}{3}\km\W{0}{4}\km\W{3}{4}\km\W{4}{5}
+\W{0}{3}\km\W{1}{3}\km\W{3}{3}\km\W{0}{4}\km\W{4}{4}\km\W{3}{5}
+\W{0}{3}\km\W{1}{3}\km\W{3}{3}\km\W{1}{4}\km\W{2}{4}\km\W{4}{5}
\\&
-\W{0}{3}\km\W{1}{3}\km\W{3}{3}\km\W{1}{4}\km\W{4}{4}\km\W{2}{5}
+2 \W{0}{3}\km\W{1}{3}\km\W{3}{3}\km\W{2}{4}\km\W{3}{4}\km\W{2}{5}
-\frac{1}{2} \W{0}{3}\km\W{1}{3}\km\W{3}{3}\km\W{2}{4}\km\W{4}{4}\km\W{1}{5}
+\frac{1}{2} \W{0}{3}\km\W{1}{3}\km\W{3}{3}\km\W{3}{4}\km\W{4}{4}\km\W{0}{5}
-\frac{1}{2} \W{0}{3}\km\W{2}{3}\km\W{3}{3}\km\W{0}{4}\km\W{1}{4}\km\W{5}{5}
+\frac{1}{2} \W{0}{3}\km\W{2}{3}\km\W{3}{3}\km\W{0}{4}\km\W{2}{4}\km\W{4}{5}
\\&
+\W{0}{3}\km\W{2}{3}\km\W{3}{3}\km\W{0}{4}\km\W{3}{4}\km\W{3}{5}
-\W{0}{3}\km\W{2}{3}\km\W{3}{3}\km\W{0}{4}\km\W{4}{4}\km\W{2}{5}
-2 \W{0}{3}\km\W{2}{3}\km\W{3}{3}\km\W{1}{4}\km\W{2}{4}\km\W{3}{5}
+\frac{3}{2} \W{0}{3}\km\W{2}{3}\km\W{3}{3}\km\W{1}{4}\km\W{4}{4}\km\W{1}{5}
-\W{0}{3}\km\W{2}{3}\km\W{3}{3}\km\W{2}{4}\km\W{3}{4}\km\W{1}{5}
-\frac{1}{2} \W{0}{3}\km\W{2}{3}\km\W{3}{3}\km\W{2}{4}\km\W{4}{4}\km\W{0}{5}
\\&
+\frac{3}{2} \W{1}{3}\km\W{2}{3}\km\W{3}{3}\km\W{0}{4}\km\W{1}{4}\km\W{4}{5}
-3 \W{1}{3}\km\W{2}{3}\km\W{3}{3}\km\W{0}{4}\km\W{2}{4}\km\W{3}{5}
+\W{1}{3}\km\W{2}{3}\km\W{3}{3}\km\W{0}{4}\km\W{3}{4}\km\W{2}{5}
+\frac{1}{2} \W{1}{3}\km\W{2}{3}\km\W{3}{3}\km\W{0}{4}\km\W{4}{4}\km\W{1}{5}
+6 \W{1}{3}\km\W{2}{3}\km\W{3}{3}\km\W{1}{4}\km\W{2}{4}\km\W{2}{5}
-4 \W{1}{3}\km\W{2}{3}\km\W{3}{3}\km\W{1}{4}\km\W{3}{4}\km\W{1}{5}
\\&
-\frac{1}{2} \W{1}{3}\km\W{2}{3}\km\W{3}{3}\km\W{1}{4}\km\W{4}{4}\km\W{0}{5}
+3 \W{1}{3}\km\W{2}{3}\km\W{3}{3}\km\W{2}{4}\km\W{3}{4}\km\W{0}{5}
\\ 
\mathbf{r}_{3} &=
\frac{1}{4} \W{0}{3}\km\W{1}{3}\km\W{2}{3}\km\W{0}{4}\km\W{3}{4}\km\W{5}{5}
-\frac{1}{4} \W{0}{3}\km\W{1}{3}\km\W{2}{3}\km\W{0}{4}\km\W{4}{4}\km\W{4}{5}
-\W{0}{3}\km\W{1}{3}\km\W{2}{3}\km\W{1}{4}\km\W{3}{4}\km\W{4}{5}
+\W{0}{3}\km\W{1}{3}\km\W{2}{3}\km\W{1}{4}\km\W{4}{4}\km\W{3}{5}
+\frac{3}{2} \W{0}{3}\km\W{1}{3}\km\W{2}{3}\km\W{2}{4}\km\W{3}{4}\km\W{3}{5}
-\frac{3}{2} \W{0}{3}\km\W{1}{3}\km\W{2}{3}\km\W{2}{4}\km\W{4}{4}\km\W{2}{5}
\\&
+\frac{3}{4} \W{0}{3}\km\W{1}{3}\km\W{2}{3}\km\W{3}{4}\km\W{4}{4}\km\W{1}{5}
-\frac{1}{4} \W{0}{3}\km\W{1}{3}\km\W{3}{3}\km\W{0}{4}\km\W{2}{4}\km\W{5}{5}
+\frac{1}{4} \W{0}{3}\km\W{1}{3}\km\W{3}{3}\km\W{0}{4}\km\W{3}{4}\km\W{4}{5}
+\W{0}{3}\km\W{1}{3}\km\W{3}{3}\km\W{1}{4}\km\W{2}{4}\km\W{4}{5}
-\W{0}{3}\km\W{1}{3}\km\W{3}{3}\km\W{1}{4}\km\W{3}{4}\km\W{3}{5}
+\frac{1}{2} \W{0}{3}\km\W{1}{3}\km\W{3}{3}\km\W{2}{4}\km\W{3}{4}\km\W{2}{5}
\\&
+\frac{1}{4} \W{0}{3}\km\W{1}{3}\km\W{3}{3}\km\W{2}{4}\km\W{4}{4}\km\W{1}{5}
-\frac{1}{4} \W{0}{3}\km\W{1}{3}\km\W{3}{3}\km\W{3}{4}\km\W{4}{4}\km\W{0}{5}
+\frac{1}{4} \W{0}{3}\km\W{2}{3}\km\W{3}{3}\km\W{0}{4}\km\W{1}{4}\km\W{5}{5}
-\frac{1}{4} \W{0}{3}\km\W{2}{3}\km\W{3}{3}\km\W{0}{4}\km\W{2}{4}\km\W{4}{5}
-\frac{1}{2} \W{0}{3}\km\W{2}{3}\km\W{3}{3}\km\W{1}{4}\km\W{2}{4}\km\W{3}{5}
+\W{0}{3}\km\W{2}{3}\km\W{3}{3}\km\W{1}{4}\km\W{3}{4}\km\W{2}{5}
\\&
-\frac{1}{4} \W{0}{3}\km\W{2}{3}\km\W{3}{3}\km\W{1}{4}\km\W{4}{4}\km\W{1}{5}
-\W{0}{3}\km\W{2}{3}\km\W{3}{3}\km\W{2}{4}\km\W{3}{4}\km\W{1}{5}
+\frac{1}{4} \W{0}{3}\km\W{2}{3}\km\W{3}{3}\km\W{2}{4}\km\W{4}{4}\km\W{0}{5}
-\frac{3}{4} \W{1}{3}\km\W{2}{3}\km\W{3}{3}\km\W{0}{4}\km\W{1}{4}\km\W{4}{5}
+\frac{3}{2} \W{1}{3}\km\W{2}{3}\km\W{3}{3}\km\W{0}{4}\km\W{2}{4}\km\W{3}{5}
-\W{1}{3}\km\W{2}{3}\km\W{3}{3}\km\W{0}{4}\km\W{3}{4}\km\W{2}{5}
\\&
+\frac{1}{4} \W{1}{3}\km\W{2}{3}\km\W{3}{3}\km\W{0}{4}\km\W{4}{4}\km\W{1}{5}
-\frac{3}{2} \W{1}{3}\km\W{2}{3}\km\W{3}{3}\km\W{1}{4}\km\W{2}{4}\km\W{2}{5}
+\W{1}{3}\km\W{2}{3}\km\W{3}{3}\km\W{1}{4}\km\W{3}{4}\km\W{1}{5}
-\frac{1}{4} \W{1}{3}\km\W{2}{3}\km\W{3}{3}\km\W{1}{4}\km\W{4}{4}\km\W{0}{5}
\\
\mathbf{r}_{4} &=
2 \W{0}{3}\km\W{1}{3}\km\W{1}{4}\km\W{2}{4}\km\W{3}{4}\km\W{4}{4}
-\W{0}{3}\km\W{2}{3}\km\W{0}{4}\km\W{2}{4}\km\W{3}{4}\km\W{4}{4}
+\frac{1}{3} \W{0}{3}\km\W{3}{3}\km\W{0}{4}\km\W{1}{4}\km\W{3}{4}\km\W{4}{4}
-\W{1}{3}\km\W{3}{3}\km\W{0}{4}\km\W{1}{4}\km\W{2}{4}\km\W{4}{4}
+2 \W{2}{3}\km\W{3}{3}\km\W{0}{4}\km\W{1}{4}\km\W{2}{4}\km\W{3}{4}
+\W{1}{3}\km\W{2}{3}\km\W{0}{4}\km\W{1}{4}\km\W{3}{4}\km\W{4}{4}
\end{align*}
\end{scriptsize}
\subsection{Matrix representation of $\mydo$}
We denote
$\displaystyle \CGF{\bullet}{2}{0}{10}$  
by $\displaystyle C^{\bullet}$.     
We choose our concrete bases as 
$\displaystyle \{\mathbf{q}_i\}_{i=1}^{9}$ of $C^{4}$, 
$\displaystyle \{\mathbf{w}_i\}_{i=1}^{12}$ of $C^{5}$, and 
$\displaystyle \{\mathbf{r}_i\}_{i=1}^{4}$ of $C^{6}$.   
Then 
the matrix representations of linear maps 
$\displaystyle \mydo : C^{4} \rightarrow C^{5}$ and 
$\displaystyle \mydo : C^{5} \rightarrow C^{6}$ are 
given as
$$ 
[\mydo (\mathbf{q}_1),\ldots, \mydo(\mathbf{q}_9)] =
[ \mathbf{w}_1,\ldots, \mathbf{w}_{12}] M $$ and
$$
[\mydo(\mathbf{w}_1),\ldots, 
\mydo(\mathbf{w}_{12}) ] = [\mathbf{r}_1,\mathbf{r}_2,\mathbf{r}_3, \mathbf{r}_{4}] N$$ 
where 
\begingroup
\begin{equation}\label{d:one:4and5}
\renewcommand{\arraystretch}{1.7} 
\setcounter{MaxMatrixCols}{13}
{}^t M = 
\begin{bmatrix} 
-\dfrac{135}{4} & 0 & -60 &\dfrac{15}{2} & -45 & -15 &\dfrac{5}{4} & -\dfrac{45}{4} &\dfrac{75}{2} & 0 & 0 & 0 \\ 
\dfrac{108}{11} &\dfrac{18}{11} & 0 & 0 & 0 &\dfrac{60}{11} &\dfrac{46}{11} & -\dfrac{90}{11} &\dfrac{156}{11} & 0 & 0 & 0 \\ 
\dfrac{27}{4} & 0 & 12 & -\dfrac{9}{2} & -9 & 0 & 0 & 0 & 0 &\dfrac{27}{4} & 18 & 0 \\ 
0 & 0 & -10 &\dfrac{2}{3} & -2 & 2 & 1 & 0 & 6 & 4 & -1 & 0 \\
0 &\dfrac{5}{2} & 29 &\dfrac{47}{3} & -23 & 43 &\dfrac{13}{2} &\dfrac{9}{2} & 25 & 16 & -\dfrac{71}{2} & 0 \\ 
0 & 5 & 45 &\dfrac{155}{6} & -40 & 65 & 10 & 0 & 50 & 20 & -\dfrac{115}{2} & 0 \\ 
0 &\dfrac{3}{2} & 18 &\dfrac{23}{2} & -3 & 30 &\dfrac{11}{2} &\dfrac{9}{2} & 9 & 6 & -33 & 0 \\ 
0 & 0 & 0 & 0 & 0 & 6 & 7 & 0 & -6 & 0 & 0 & 0 \\ 
0 & 0 & -6 & -3 & 0 & 0 & 0 & 0 & 0 & 3 & -6 & 70 
\end{bmatrix} 
\end{equation}
\endgroup
and 
\begingroup
\renewcommand{\arraystretch}{1.7} 
\begin{equation}
N = 
\begin{bmatrix} 
 0 & 140 & 0 & 0 & 0 & -15 & 15 & 30 &
      \dfrac{5}{2} & 0 & 0 & 0 \\ 
 -5 & -4 &\dfrac{1}{4} & -\dfrac{11}{2} & \dfrac{31}{12} 
 &\dfrac{31}{6} & -3 & -2 &\dfrac{5}{3} & -1 & 2 & 0 \\ 
  -16 & 32 & -2 & -12 &\dfrac{22}{3} &\dfrac{58}{3} & -18 & -12 &
 -\dfrac{5}{3} & 0 & 8 & 0 \\ 
 0 & 0 & 0 & 42 & 7 & 0 & 0 & 0 & 0 & 0 & 14 & 3 
\end{bmatrix} 
\end{equation}
\endgroup
Since $\rank M = 7$ and $\rank N = 4$, we see the dimensions of 
$\mydo ( C^{4})$ and $\ker( \mydo : C^{5} \rightarrow C^{6})$ , and so
on. The precise data of 
the structures of $\displaystyle \CGF{\bullet}{2}{0}{10}$ and
$\displaystyle \HGF{\bullet}{2}{0}{10}$  is in the table below, 
where    
$\dim$ and $\rank$ mean the dimension of  
$\displaystyle C^{\bullet}$ and the rank of $ 
\displaystyle \mydo: C^{\bullet} \rightarrow  C^{\bullet+1}$, and Betti
num is the Betti number, which is the dimension of the cohomology group 
$\displaystyle \HGF{\bullet}{2}{0}{10}$.

\begin{center}
\begin{tabular}{|c|*{12}{c}c|}
\hline
$\displaystyle \frakham^{0}_{2}$, w=10
 & $\displaystyle \mathbf{0}$ & $\rightarrow$ &
$\displaystyle C^{2} $& $ \rightarrow $ &
$\displaystyle C^{3} $& $ \rightarrow $ &
$\displaystyle C^{4} $& $ \rightarrow $ &
$\displaystyle C^{5} $& $ \rightarrow $ &
$\displaystyle C^{6} $& $ \rightarrow $ &
$\mathbf{0}$ \\
\hline
$\dim$  &   && 1 && 3 && 9 && 12 && 4&& \\
rank & & 0 && 1 && 2 && 7 && 4 && 0 &\\
Betti num &   && 0 && 0 && 0 && 1 && 0 &&\\
\hline
\end{tabular}
\end{center}

\section{About $\ds\CGF{\bullet}{2}{}{w}$}
We also know the structures of $\displaystyle
\CGF{\bullet}{2}{}{8}$ well.  
\begin{scriptsize}
\subsection{Basis of $\ds\CGF{6}{2}{}{8}$}

\begin{align*}
\widehat{\mathbf{q}}_{1} &= 
+\frac{1}{8} \W{0}{1}\km\W{1}{1}\km\W{0}{3}\km\W{1}{3}\km\W{0}{4}\km\W{8}{8}
-\frac{1}{2} \W{0}{1}\km\W{1}{1}\km\W{0}{3}\km\W{1}{3}\km\W{1}{4}\km\W{7}{8}
+\frac{3}{4} \W{0}{1}\km\W{1}{1}\km\W{0}{3}\km\W{1}{3}\km\W{2}{4}\km\W{6}{8}
+\frac{1}{8} \W{0}{1}\km\W{1}{1}\km\W{0}{3}\km\W{1}{3}\km\W{4}{4}\km\W{4}{8}
-\frac{1}{2} \W{0}{1}\km\W{1}{1}\km\W{0}{3}\km\W{1}{3}\km\W{3}{4}\km\W{5}{8}
-\frac{1}{4} \W{0}{1}\km\W{1}{1}\km\W{0}{3}\km\W{2}{3}\km\W{0}{4}\km\W{7}{8}
\\&
+\W{0}{1}\km\W{1}{1}\km\W{0}{3}\km\W{2}{3}\km\W{1}{4}\km\W{6}{8}
-\frac{3}{2} \W{0}{1}\km\W{1}{1}\km\W{0}{3}\km\W{2}{3}\km\W{2}{4}\km\W{5}{8}
+\W{0}{1}\km\W{1}{1}\km\W{0}{3}\km\W{2}{3}\km\W{3}{4}\km\W{4}{8}
-\frac{1}{4} \W{0}{1}\km\W{1}{1}\km\W{0}{3}\km\W{2}{3}\km\W{4}{4}\km\W{3}{8}
+\frac{1}{8} \W{0}{1}\km\W{1}{1}\km\W{0}{3}\km\W{3}{3}\km\W{0}{4}\km\W{6}{8}
-\frac{1}{2} \W{0}{1}\km\W{1}{1}\km\W{0}{3}\km\W{3}{3}\km\W{1}{4}\km\W{5}{8}
\\&
+\frac{3}{4} \W{0}{1}\km\W{1}{1}\km\W{0}{3}\km\W{3}{3}\km\W{2}{4}\km\W{4}{8}
-\frac{1}{2} \W{0}{1}\km\W{1}{1}\km\W{0}{3}\km\W{3}{3}\km\W{3}{4}\km\W{3}{8}
+\frac{1}{8} \W{0}{1}\km\W{1}{1}\km\W{0}{3}\km\W{3}{3}\km\W{4}{4}\km\W{2}{8}
+\frac{3}{8} \W{0}{1}\km\W{1}{1}\km\W{1}{3}\km\W{2}{3}\km\W{0}{4}\km\W{6}{8}
-\frac{3}{2} \W{0}{1}\km\W{1}{1}\km\W{1}{3}\km\W{2}{3}\km\W{1}{4}\km\W{5}{8}
+\frac{9}{4} \W{0}{1}\km\W{1}{1}\km\W{1}{3}\km\W{2}{3}\km\W{2}{4}\km\W{4}{8}
\\&
+\frac{3}{8} \W{0}{1}\km\W{1}{1}\km\W{1}{3}\km\W{2}{3}\km\W{4}{4}\km\W{2}{8}
-\frac{3}{2} \W{0}{1}\km\W{1}{1}\km\W{1}{3}\km\W{2}{3}\km\W{3}{4}\km\W{3}{8}
-\frac{1}{4} \W{0}{1}\km\W{1}{1}\km\W{1}{3}\km\W{3}{3}\km\W{0}{4}\km\W{5}{8}
-\frac{3}{2} \W{0}{1}\km\W{1}{1}\km\W{1}{3}\km\W{3}{3}\km\W{2}{4}\km\W{3}{8}
+\W{0}{1}\km\W{1}{1}\km\W{1}{3}\km\W{3}{3}\km\W{1}{4}\km\W{4}{8}
+\W{0}{1}\km\W{1}{1}\km\W{1}{3}\km\W{3}{3}\km\W{3}{4}\km\W{2}{8}
\\&
-\frac{1}{4} \W{0}{1}\km\W{1}{1}\km\W{1}{3}\km\W{3}{3}\km\W{4}{4}\km\W{1}{8} 
+\frac{1}{8} \W{0}{1}\km\W{1}{1}\km\W{2}{3}\km\W{3}{3}\km\W{0}{4}\km\W{4}{8}
-\frac{1}{2} \W{0}{1}\km\W{1}{1}\km\W{2}{3}\km\W{3}{3}\km\W{1}{4}\km\W{3}{8}
+\frac{3}{4} \W{0}{1}\km\W{1}{1}\km\W{2}{3}\km\W{3}{3}\km\W{2}{4}\km\W{2}{8}
-\frac{1}{2} \W{0}{1}\km\W{1}{1}\km\W{2}{3}\km\W{3}{3}\km\W{3}{4}\km\W{1}{8}
+\frac{1}{8} \W{0}{1}\km\W{1}{1}\km\W{2}{3}\km\W{3}{3}\km\W{4}{4}\km\W{0}{8} 
\\
\widehat{\mathbf{q}}_{2}  &= 
+\frac{1}{3} \W{0}{1}\km\W{1}{1}\km\W{0}{3}\km\W{1}{3}\km\W{1}{5}\km\W{7}{7}
+2 \W{0}{1}\km\W{1}{1}\km\W{0}{3}\km\W{1}{3}\km\W{3}{5}\km\W{5}{7}
-\frac{4}{3} \W{0}{1}\km\W{1}{1}\km\W{0}{3}\km\W{1}{3}\km\W{2}{5}\km\W{6}{7}
-\frac{4}{3} \W{0}{1}\km\W{1}{1}\km\W{0}{3}\km\W{1}{3}\km\W{4}{5}\km\W{4}{7}
+\frac{1}{3} \W{0}{1}\km\W{1}{1}\km\W{0}{3}\km\W{1}{3}\km\W{5}{5}\km\W{3}{7}
-\frac{1}{6} \W{0}{1}\km\W{1}{1}\km\W{0}{3}\km\W{2}{3}\km\W{0}{5}\km\W{7}{7}
\\&
+\frac{1}{6} \W{0}{1}\km\W{1}{1}\km\W{0}{3}\km\W{2}{3}\km\W{1}{5}\km\W{6}{7}
+\W{0}{1}\km\W{1}{1}\km\W{0}{3}\km\W{2}{3}\km\W{2}{5}\km\W{5}{7}
-\frac{7}{3} \W{0}{1}\km\W{1}{1}\km\W{0}{3}\km\W{2}{3}\km\W{3}{5}\km\W{4}{7}
+\frac{11}{6} \W{0}{1}\km\W{1}{1}\km\W{0}{3}\km\W{2}{3}\km\W{4}{5}\km\W{3}{7}
-\frac{1}{2} \W{0}{1}\km\W{1}{1}\km\W{0}{3}\km\W{2}{3}\km\W{5}{5}\km\W{2}{7}
+\frac{1}{6} \W{0}{1}\km\W{1}{1}\km\W{0}{3}\km\W{3}{3}\km\W{0}{5}\km\W{6}{7}
\\&
-\frac{1}{2} \W{0}{1}\km\W{1}{1}\km\W{0}{3}\km\W{3}{3}\km\W{1}{5}\km\W{5}{7}
+\frac{1}{3} \W{0}{1}\km\W{1}{1}\km\W{0}{3}\km\W{3}{3}\km\W{2}{5}\km\W{4}{7}
+\frac{1}{3} \W{0}{1}\km\W{1}{1}\km\W{0}{3}\km\W{3}{3}\km\W{3}{5}\km\W{3}{7}
-\frac{1}{2} \W{0}{1}\km\W{1}{1}\km\W{0}{3}\km\W{3}{3}\km\W{4}{5}\km\W{2}{7}
+\frac{1}{6} \W{0}{1}\km\W{1}{1}\km\W{0}{3}\km\W{3}{3}\km\W{5}{5}\km\W{1}{7}
+\frac{1}{2} \W{0}{1}\km\W{1}{1}\km\W{1}{3}\km\W{2}{3}\km\W{0}{5}\km\W{6}{7}
\\&
-\frac{3}{2} \W{0}{1}\km\W{1}{1}\km\W{1}{3}\km\W{2}{3}\km\W{1}{5}\km\W{5}{7}
+\W{0}{1}\km\W{1}{1}\km\W{1}{3}\km\W{2}{3}\km\W{3}{5}\km\W{3}{7}
+\W{0}{1}\km\W{1}{1}\km\W{1}{3}\km\W{2}{3}\km\W{2}{5}\km\W{4}{7}
-\frac{3}{2} \W{0}{1}\km\W{1}{1}\km\W{1}{3}\km\W{2}{3}\km\W{4}{5}\km\W{2}{7}
+\frac{1}{2} \W{0}{1}\km\W{1}{1}\km\W{1}{3}\km\W{2}{3}\km\W{5}{5}\km\W{1}{7}
-\frac{1}{2} \W{0}{1}\km\W{1}{1}\km\W{1}{3}\km\W{3}{3}\km\W{0}{5}\km\W{5}{7}
\\&
+\frac{11}{6} \W{0}{1}\km\W{1}{1}\km\W{1}{3}\km\W{3}{3}\km\W{1}{5}\km\W{4}{7}
-\frac{7}{3} \W{0}{1}\km\W{1}{1}\km\W{1}{3}\km\W{3}{3}\km\W{2}{5}\km\W{3}{7}
+\W{0}{1}\km\W{1}{1}\km\W{1}{3}\km\W{3}{3}\km\W{3}{5}\km\W{2}{7}
+\frac{1}{6} \W{0}{1}\km\W{1}{1}\km\W{1}{3}\km\W{3}{3}\km\W{4}{5}\km\W{1}{7}
-\frac{1}{6} \W{0}{1}\km\W{1}{1}\km\W{1}{3}\km\W{3}{3}\km\W{5}{5}\km\W{0}{7}
+\frac{1}{3} \W{0}{1}\km\W{1}{1}\km\W{2}{3}\km\W{3}{3}\km\W{0}{5}\km\W{4}{7}
\\&
-\frac{4}{3} \W{0}{1}\km\W{1}{1}\km\W{2}{3}\km\W{3}{3}\km\W{1}{5}\km\W{3}{7}
+2 \W{0}{1}\km\W{1}{1}\km\W{2}{3}\km\W{3}{3}\km\W{2}{5}\km\W{2}{7}
-\frac{4}{3} \W{0}{1}\km\W{1}{1}\km\W{2}{3}\km\W{3}{3}\km\W{3}{5}\km\W{1}{7}
+\frac{1}{3} \W{0}{1}\km\W{1}{1}\km\W{2}{3}\km\W{3}{3}\km\W{4}{5}\km\W{0}{7}
\\
\widehat{\mathbf{q}}_{3}  &= 
-\frac{1}{8} \W{0}{1}\km\W{1}{1}\km\W{0}{3}\km\W{0}{4}\km\W{2}{4}\km\W{7}{7}
+\frac{1}{4} \W{0}{1}\km\W{1}{1}\km\W{0}{3}\km\W{0}{4}\km\W{3}{4}\km\W{6}{7}
-\frac{1}{8} \W{0}{1}\km\W{1}{1}\km\W{0}{3}\km\W{0}{4}\km\W{4}{4}\km\W{5}{7}
+\frac{1}{2} \W{0}{1}\km\W{1}{1}\km\W{0}{3}\km\W{1}{4}\km\W{2}{4}\km\W{6}{7}
-\W{0}{1}\km\W{1}{1}\km\W{0}{3}\km\W{1}{4}\km\W{3}{4}\km\W{5}{7}
+\frac{1}{2} \W{0}{1}\km\W{1}{1}\km\W{0}{3}\km\W{1}{4}\km\W{4}{4}\km\W{4}{7} 
\\& 
+\W{0}{1}\km\W{1}{1}\km\W{0}{3}\km\W{2}{4}\km\W{3}{4}\km\W{4}{7}
-\frac{5}{8} \W{0}{1}\km\W{1}{1}\km\W{0}{3}\km\W{2}{4}\km\W{4}{4}\km\W{3}{7}
+\frac{1}{4} \W{0}{1}\km\W{1}{1}\km\W{0}{3}\km\W{3}{4}\km\W{4}{4}\km\W{2}{7}
+\frac{1}{4} \W{0}{1}\km\W{1}{1}\km\W{1}{3}\km\W{0}{4}\km\W{1}{4}\km\W{7}{7}
-\frac{3}{8} \W{0}{1}\km\W{1}{1}\km\W{1}{3}\km\W{0}{4}\km\W{2}{4}\km\W{6}{7} 
+\frac{1}{8} \W{0}{1}\km\W{1}{1}\km\W{1}{3}\km\W{0}{4}\km\W{4}{4}\km\W{4}{7}
\\&
+\W{0}{1}\km\W{1}{1}\km\W{1}{3}\km\W{1}{4}\km\W{3}{4}\km\W{4}{7}
-\frac{3}{4} \W{0}{1}\km\W{1}{1}\km\W{1}{3}\km\W{1}{4}\km\W{4}{4}\km\W{3}{7}
-\frac{3}{2} \W{0}{1}\km\W{1}{1}\km\W{1}{3}\km\W{2}{4}\km\W{3}{4}\km\W{3}{7}
+\frac{9}{8} \W{0}{1}\km\W{1}{1}\km\W{1}{3}\km\W{2}{4}\km\W{4}{4}\km\W{2}{7} 
-\frac{1}{2} \W{0}{1}\km\W{1}{1}\km\W{1}{3}\km\W{3}{4}\km\W{4}{4}\km\W{1}{7}
-\frac{1}{2} \W{0}{1}\km\W{1}{1}\km\W{2}{3}\km\W{0}{4}\km\W{1}{4}\km\W{6}{7}
\\&
+\frac{9}{8} \W{0}{1}\km\W{1}{1}\km\W{2}{3}\km\W{0}{4}\km\W{2}{4}\km\W{5}{7}
-\frac{3}{4} \W{0}{1}\km\W{1}{1}\km\W{2}{3}\km\W{0}{4}\km\W{3}{4}\km\W{4}{7}
+\frac{1}{8} \W{0}{1}\km\W{1}{1}\km\W{2}{3}\km\W{0}{4}\km\W{4}{4}\km\W{3}{7} 
-\frac{3}{2} \W{0}{1}\km\W{1}{1}\km\W{2}{3}\km\W{1}{4}\km\W{2}{4}\km\W{4}{7}
+\W{0}{1}\km\W{1}{1}\km\W{2}{3}\km\W{1}{4}\km\W{3}{4}\km\W{3}{7}
-\frac{3}{8} \W{0}{1}\km\W{1}{1}\km\W{2}{3}\km\W{2}{4}\km\W{4}{4}\km\W{1}{7}
\\&
+\frac{1}{4} \W{0}{1}\km\W{1}{1}\km\W{2}{3}\km\W{3}{4}\km\W{4}{4}\km\W{0}{7}
+\frac{1}{4} \W{0}{1}\km\W{1}{1}\km\W{3}{3}\km\W{0}{4}\km\W{1}{4}\km\W{5}{7} 
-\frac{5}{8} \W{0}{1}\km\W{1}{1}\km\W{3}{3}\km\W{0}{4}\km\W{2}{4}\km\W{4}{7}
+\frac{1}{2} \W{0}{1}\km\W{1}{1}\km\W{3}{3}\km\W{0}{4}\km\W{3}{4}\km\W{3}{7}
-\frac{1}{8} \W{0}{1}\km\W{1}{1}\km\W{3}{3}\km\W{0}{4}\km\W{4}{4}\km\W{2}{7}
+\W{0}{1}\km\W{1}{1}\km\W{3}{3}\km\W{1}{4}\km\W{2}{4}\km\W{3}{7}
\\&
-\W{0}{1}\km\W{1}{1}\km\W{3}{3}\km\W{1}{4}\km\W{3}{4}\km\W{2}{7} 
+\frac{1}{4} \W{0}{1}\km\W{1}{1}\km\W{3}{3}\km\W{1}{4}\km\W{4}{4}\km\W{1}{7}
+\frac{1}{2} \W{0}{1}\km\W{1}{1}\km\W{3}{3}\km\W{2}{4}\km\W{3}{4}\km\W{1}{7}
-\frac{1}{8} \W{0}{1}\km\W{1}{1}\km\W{3}{3}\km\W{2}{4}\km\W{4}{4}\km\W{0}{7}
\\
\widehat{\mathbf{q}}_{4} &= 
+\frac{1}{30} \W{0}{1}\km\W{1}{1}\km\W{0}{3}\km\W{3}{4}\km\W{0}{5}\km\W{6}{6}
-\frac{1}{6} \W{0}{1}\km\W{1}{1}\km\W{0}{3}\km\W{3}{4}\km\W{1}{5}\km\W{5}{6}
+\frac{1}{3} \W{0}{1}\km\W{1}{1}\km\W{0}{3}\km\W{3}{4}\km\W{2}{5}\km\W{4}{6}
-\frac{1}{3} \W{0}{1}\km\W{1}{1}\km\W{0}{3}\km\W{3}{4}\km\W{3}{5}\km\W{3}{6}
+\frac{1}{6} \W{0}{1}\km\W{1}{1}\km\W{0}{3}\km\W{3}{4}\km\W{4}{5}\km\W{2}{6}
-\frac{1}{30} \W{0}{1}\km\W{1}{1}\km\W{0}{3}\km\W{4}{4}\km\W{0}{5}\km\W{5}{6}
\\&
-\frac{1}{30} \W{0}{1}\km\W{1}{1}\km\W{0}{3}\km\W{3}{4}\km\W{5}{5}\km\W{1}{6}
+\frac{1}{6} \W{0}{1}\km\W{1}{1}\km\W{0}{3}\km\W{4}{4}\km\W{1}{5}\km\W{4}{6}
-\frac{1}{3} \W{0}{1}\km\W{1}{1}\km\W{0}{3}\km\W{4}{4}\km\W{2}{5}\km\W{3}{6}
-\frac{1}{6} \W{0}{1}\km\W{1}{1}\km\W{0}{3}\km\W{4}{4}\km\W{4}{5}\km\W{1}{6}
+\frac{1}{3} \W{0}{1}\km\W{1}{1}\km\W{0}{3}\km\W{4}{4}\km\W{3}{5}\km\W{2}{6}
+\frac{1}{30} \W{0}{1}\km\W{1}{1}\km\W{0}{3}\km\W{4}{4}\km\W{5}{5}\km\W{0}{6}
\\&
-\frac{1}{10} \W{0}{1}\km\W{1}{1}\km\W{1}{3}\km\W{2}{4}\km\W{0}{5}\km\W{6}{6}
+\frac{1}{2} \W{0}{1}\km\W{1}{1}\km\W{1}{3}\km\W{2}{4}\km\W{1}{5}\km\W{5}{6}
-\W{0}{1}\km\W{1}{1}\km\W{1}{3}\km\W{2}{4}\km\W{2}{5}\km\W{4}{6}
+\W{0}{1}\km\W{1}{1}\km\W{1}{3}\km\W{2}{4}\km\W{3}{5}\km\W{3}{6}
-\frac{1}{2} \W{0}{1}\km\W{1}{1}\km\W{1}{3}\km\W{2}{4}\km\W{4}{5}\km\W{2}{6}
+\frac{1}{10} \W{0}{1}\km\W{1}{1}\km\W{1}{3}\km\W{2}{4}\km\W{5}{5}\km\W{1}{6}
\\&
+\frac{1}{10} \W{0}{1}\km\W{1}{1}\km\W{1}{3}\km\W{3}{4}\km\W{0}{5}\km\W{5}{6}
+\W{0}{1}\km\W{1}{1}\km\W{1}{3}\km\W{3}{4}\km\W{2}{5}\km\W{3}{6}
-\frac{1}{2} \W{0}{1}\km\W{1}{1}\km\W{1}{3}\km\W{3}{4}\km\W{1}{5}\km\W{4}{6}
-\W{0}{1}\km\W{1}{1}\km\W{1}{3}\km\W{3}{4}\km\W{3}{5}\km\W{2}{6}
+\frac{1}{2} \W{0}{1}\km\W{1}{1}\km\W{1}{3}\km\W{3}{4}\km\W{4}{5}\km\W{1}{6}
-\frac{1}{10} \W{0}{1}\km\W{1}{1}\km\W{1}{3}\km\W{3}{4}\km\W{5}{5}\km\W{0}{6}
\\&
+\frac{1}{10} \W{0}{1}\km\W{1}{1}\km\W{2}{3}\km\W{1}{4}\km\W{0}{5}\km\W{6}{6}
-\frac{1}{2} \W{0}{1}\km\W{1}{1}\km\W{2}{3}\km\W{1}{4}\km\W{1}{5}\km\W{5}{6}
+\W{0}{1}\km\W{1}{1}\km\W{2}{3}\km\W{1}{4}\km\W{2}{5}\km\W{4}{6}
-\W{0}{1}\km\W{1}{1}\km\W{2}{3}\km\W{1}{4}\km\W{3}{5}\km\W{3}{6}
+\frac{1}{2} \W{0}{1}\km\W{1}{1}\km\W{2}{3}\km\W{1}{4}\km\W{4}{5}\km\W{2}{6}
-\frac{1}{10} \W{0}{1}\km\W{1}{1}\km\W{2}{3}\km\W{1}{4}\km\W{5}{5}\km\W{1}{6}
\\&
-\frac{1}{10} \W{0}{1}\km\W{1}{1}\km\W{2}{3}\km\W{2}{4}\km\W{0}{5}\km\W{5}{6}
+\frac{1}{2} \W{0}{1}\km\W{1}{1}\km\W{2}{3}\km\W{2}{4}\km\W{1}{5}\km\W{4}{6}
-\W{0}{1}\km\W{1}{1}\km\W{2}{3}\km\W{2}{4}\km\W{2}{5}\km\W{3}{6}
+\W{0}{1}\km\W{1}{1}\km\W{2}{3}\km\W{2}{4}\km\W{3}{5}\km\W{2}{6}
-\frac{1}{2} \W{0}{1}\km\W{1}{1}\km\W{2}{3}\km\W{2}{4}\km\W{4}{5}\km\W{1}{6}
+\frac{1}{10} \W{0}{1}\km\W{1}{1}\km\W{2}{3}\km\W{2}{4}\km\W{5}{5}\km\W{0}{6}
\\&
-\frac{1}{30} \W{0}{1}\km\W{1}{1}\km\W{3}{3}\km\W{0}{4}\km\W{0}{5}\km\W{6}{6}
+\frac{1}{6} \W{0}{1}\km\W{1}{1}\km\W{3}{3}\km\W{0}{4}\km\W{1}{5}\km\W{5}{6}
-\frac{1}{3} \W{0}{1}\km\W{1}{1}\km\W{3}{3}\km\W{0}{4}\km\W{2}{5}\km\W{4}{6}
+\frac{1}{3} \W{0}{1}\km\W{1}{1}\km\W{3}{3}\km\W{0}{4}\km\W{3}{5}\km\W{3}{6}
+\frac{1}{30} \W{0}{1}\km\W{1}{1}\km\W{3}{3}\km\W{0}{4}\km\W{5}{5}\km\W{1}{6}
-\frac{1}{6} \W{0}{1}\km\W{1}{1}\km\W{3}{3}\km\W{0}{4}\km\W{4}{5}\km\W{2}{6}
\\&
+\frac{1}{30} \W{0}{1}\km\W{1}{1}\km\W{3}{3}\km\W{1}{4}\km\W{0}{5}\km\W{5}{6}
-\frac{1}{6} \W{0}{1}\km\W{1}{1}\km\W{3}{3}\km\W{1}{4}\km\W{1}{5}\km\W{4}{6}
+\frac{1}{3} \W{0}{1}\km\W{1}{1}\km\W{3}{3}\km\W{1}{4}\km\W{2}{5}\km\W{3}{6}
-\frac{1}{3} \W{0}{1}\km\W{1}{1}\km\W{3}{3}\km\W{1}{4}\km\W{3}{5}\km\W{2}{6}
+\frac{1}{6} \W{0}{1}\km\W{1}{1}\km\W{3}{3}\km\W{1}{4}\km\W{4}{5}\km\W{1}{6}
-\frac{1}{30} \W{0}{1}\km\W{1}{1}\km\W{3}{3}\km\W{1}{4}\km\W{5}{5}\km\W{0}{6}
\\
\widehat{\mathbf{q}}_{5}&= 
+\frac{3}{16} \W{0}{1}\km\W{1}{1}\km\W{0}{3}\km\W{0}{4}\km\W{3}{5}\km\W{6}{6}
+\frac{3}{16} \W{0}{1}\km\W{1}{1}\km\W{0}{3}\km\W{0}{4}\km\W{5}{5}\km\W{4}{6}
-\frac{3}{8} \W{0}{1}\km\W{1}{1}\km\W{0}{3}\km\W{0}{4}\km\W{4}{5}\km\W{5}{6}
-\frac{1}{4} \W{0}{1}\km\W{1}{1}\km\W{0}{3}\km\W{1}{4}\km\W{2}{5}\km\W{6}{6}
+\frac{3}{4} \W{0}{1}\km\W{1}{1}\km\W{0}{3}\km\W{1}{4}\km\W{4}{5}\km\W{4}{6}
-\frac{1}{2} \W{0}{1}\km\W{1}{1}\km\W{0}{3}\km\W{1}{4}\km\W{5}{5}\km\W{3}{6}
\\&
-\frac{1}{16} \W{0}{1}\km\W{1}{1}\km\W{0}{3}\km\W{2}{4}\km\W{1}{5}\km\W{6}{6}
+\W{0}{1}\km\W{1}{1}\km\W{0}{3}\km\W{2}{4}\km\W{2}{5}\km\W{5}{6}
-\frac{3}{2} \W{0}{1}\km\W{1}{1}\km\W{0}{3}\km\W{2}{4}\km\W{3}{5}\km\W{4}{6}
+\frac{1}{4} \W{0}{1}\km\W{1}{1}\km\W{0}{3}\km\W{2}{4}\km\W{4}{5}\km\W{3}{6}
+\frac{5}{16} \W{0}{1}\km\W{1}{1}\km\W{0}{3}\km\W{2}{4}\km\W{5}{5}\km\W{2}{6}
+\frac{13}{120} \W{0}{1}\km\W{1}{1}\km\W{0}{3}\km\W{3}{4}\km\W{0}{5}\km\W{6}{6}
\\&
-\frac{5}{12} \W{0}{1}\km\W{1}{1}\km\W{0}{3}\km\W{3}{4}\km\W{1}{5}\km\W{5}{6}
-\frac{1}{6} \W{0}{1}\km\W{1}{1}\km\W{0}{3}\km\W{3}{4}\km\W{2}{5}\km\W{4}{6}
+\frac{7}{6} \W{0}{1}\km\W{1}{1}\km\W{0}{3}\km\W{3}{4}\km\W{3}{5}\km\W{3}{6}
-\frac{17}{24} \W{0}{1}\km\W{1}{1}\km\W{0}{3}\km\W{3}{4}\km\W{4}{5}\km\W{2}{6}
-\frac{13}{120} \W{0}{1}\km\W{1}{1}\km\W{0}{3}\km\W{4}{4}\km\W{0}{5}\km\W{5}{6}
+\frac{1}{60} \W{0}{1}\km\W{1}{1}\km\W{0}{3}\km\W{3}{4}\km\W{5}{5}\km\W{1}{6}
\\&
+\frac{23}{48} \W{0}{1}\km\W{1}{1}\km\W{0}{3}\km\W{4}{4}\km\W{1}{5}\km\W{4}{6}
-\frac{7}{12} \W{0}{1}\km\W{1}{1}\km\W{0}{3}\km\W{4}{4}\km\W{2}{5}\km\W{3}{6}
+\frac{1}{12} \W{0}{1}\km\W{1}{1}\km\W{0}{3}\km\W{4}{4}\km\W{4}{5}\km\W{1}{6}
+\frac{7}{48} \W{0}{1}\km\W{1}{1}\km\W{0}{3}\km\W{4}{4}\km\W{3}{5}\km\W{2}{6}
-\frac{1}{60} \W{0}{1}\km\W{1}{1}\km\W{0}{3}\km\W{4}{4}\km\W{5}{5}\km\W{0}{6}
-\frac{5}{16} \W{0}{1}\km\W{1}{1}\km\W{1}{3}\km\W{0}{4}\km\W{2}{5}\km\W{6}{6}
\\&
+\frac{3}{8} \W{0}{1}\km\W{1}{1}\km\W{1}{3}\km\W{0}{4}\km\W{3}{5}\km\W{5}{6}
+\frac{3}{16} \W{0}{1}\km\W{1}{1}\km\W{1}{3}\km\W{0}{4}\km\W{4}{5}\km\W{4}{6}
-\frac{1}{4} \W{0}{1}\km\W{1}{1}\km\W{1}{3}\km\W{0}{4}\km\W{5}{5}\km\W{3}{6}
+\frac{5}{8} \W{0}{1}\km\W{1}{1}\km\W{1}{3}\km\W{1}{4}\km\W{1}{5}\km\W{6}{6}
-\frac{1}{2} \W{0}{1}\km\W{1}{1}\km\W{1}{3}\km\W{1}{4}\km\W{2}{5}\km\W{5}{6}
-\W{0}{1}\km\W{1}{1}\km\W{1}{3}\km\W{1}{4}\km\W{4}{5}\km\W{3}{6}
\\&
-\frac{21}{80} \W{0}{1}\km\W{1}{1}\km\W{1}{3}\km\W{2}{4}\km\W{0}{5}\km\W{6}{6}
+\frac{7}{8} \W{0}{1}\km\W{1}{1}\km\W{1}{3}\km\W{1}{4}\km\W{5}{5}\km\W{2}{6}
-\frac{3}{8} \W{0}{1}\km\W{1}{1}\km\W{1}{3}\km\W{2}{4}\km\W{1}{5}\km\W{5}{6}
+\frac{3}{2} \W{0}{1}\km\W{1}{1}\km\W{1}{3}\km\W{2}{4}\km\W{3}{5}\km\W{3}{6}
-\frac{3}{16} \W{0}{1}\km\W{1}{1}\km\W{1}{3}\km\W{2}{4}\km\W{4}{5}\km\W{2}{6}
-\frac{27}{40} \W{0}{1}\km\W{1}{1}\km\W{1}{3}\km\W{2}{4}\km\W{5}{5}\km\W{1}{6}
\\&
+\frac{1}{5} \W{0}{1}\km\W{1}{1}\km\W{1}{3}\km\W{3}{4}\km\W{0}{5}\km\W{5}{6}
-\frac{1}{2} \W{0}{1}\km\W{1}{1}\km\W{1}{3}\km\W{3}{4}\km\W{2}{5}\km\W{3}{6}
+\frac{1}{2} \W{0}{1}\km\W{1}{1}\km\W{1}{3}\km\W{3}{4}\km\W{1}{5}\km\W{4}{6}
-\frac{5}{4} \W{0}{1}\km\W{1}{1}\km\W{1}{3}\km\W{3}{4}\km\W{3}{5}\km\W{2}{6}
+\W{0}{1}\km\W{1}{1}\km\W{1}{3}\km\W{3}{4}\km\W{4}{5}\km\W{1}{6}
+\frac{1}{20} \W{0}{1}\km\W{1}{1}\km\W{1}{3}\km\W{3}{4}\km\W{5}{5}\km\W{0}{6}
\\&
+\frac{1}{16} \W{0}{1}\km\W{1}{1}\km\W{1}{3}\km\W{4}{4}\km\W{0}{5}\km\W{4}{6}
-\frac{3}{4} \W{0}{1}\km\W{1}{1}\km\W{1}{3}\km\W{4}{4}\km\W{1}{5}\km\W{3}{6}
+\frac{21}{16} \W{0}{1}\km\W{1}{1}\km\W{1}{3}\km\W{4}{4}\km\W{2}{5}\km\W{2}{6}
-\frac{5}{8} \W{0}{1}\km\W{1}{1}\km\W{1}{3}\km\W{4}{4}\km\W{3}{5}\km\W{1}{6}
+\frac{5}{8} \W{0}{1}\km\W{1}{1}\km\W{2}{3}\km\W{0}{4}\km\W{2}{5}\km\W{5}{6}
-\frac{21}{16} \W{0}{1}\km\W{1}{1}\km\W{2}{3}\km\W{0}{4}\km\W{3}{5}\km\W{4}{6}
\\&
+\frac{3}{4} \W{0}{1}\km\W{1}{1}\km\W{2}{3}\km\W{0}{4}\km\W{4}{5}\km\W{3}{6}
-\frac{1}{16} \W{0}{1}\km\W{1}{1}\km\W{2}{3}\km\W{0}{4}\km\W{5}{5}\km\W{2}{6}
-\frac{1}{20} \W{0}{1}\km\W{1}{1}\km\W{2}{3}\km\W{1}{4}\km\W{0}{5}\km\W{6}{6}
-\W{0}{1}\km\W{1}{1}\km\W{2}{3}\km\W{1}{4}\km\W{1}{5}\km\W{5}{6}
+\frac{5}{4} \W{0}{1}\km\W{1}{1}\km\W{2}{3}\km\W{1}{4}\km\W{2}{5}\km\W{4}{6}
+\frac{1}{2} \W{0}{1}\km\W{1}{1}\km\W{2}{3}\km\W{1}{4}\km\W{3}{5}\km\W{3}{6}
\\&
-\frac{1}{2} \W{0}{1}\km\W{1}{1}\km\W{2}{3}\km\W{1}{4}\km\W{4}{5}\km\W{2}{6}
-\frac{1}{5} \W{0}{1}\km\W{1}{1}\km\W{2}{3}\km\W{1}{4}\km\W{5}{5}\km\W{1}{6}
+\frac{27}{40} \W{0}{1}\km\W{1}{1}\km\W{2}{3}\km\W{2}{4}\km\W{0}{5}\km\W{5}{6}
+\frac{3}{16} \W{0}{1}\km\W{1}{1}\km\W{2}{3}\km\W{2}{4}\km\W{1}{5}\km\W{4}{6}
-\frac{3}{2} \W{0}{1}\km\W{1}{1}\km\W{2}{3}\km\W{2}{4}\km\W{2}{5}\km\W{3}{6}
+\frac{3}{8} \W{0}{1}\km\W{1}{1}\km\W{2}{3}\km\W{2}{4}\km\W{4}{5}\km\W{1}{6}
\\&
+\frac{21}{80} \W{0}{1}\km\W{1}{1}\km\W{2}{3}\km\W{2}{4}\km\W{5}{5}\km\W{0}{6}
-\frac{7}{8} \W{0}{1}\km\W{1}{1}\km\W{2}{3}\km\W{3}{4}\km\W{0}{5}\km\W{4}{6}
+\W{0}{1}\km\W{1}{1}\km\W{2}{3}\km\W{3}{4}\km\W{1}{5}\km\W{3}{6}
+\frac{1}{2} \W{0}{1}\km\W{1}{1}\km\W{2}{3}\km\W{3}{4}\km\W{3}{5}\km\W{1}{6}
-\frac{5}{8} \W{0}{1}\km\W{1}{1}\km\W{2}{3}\km\W{3}{4}\km\W{4}{5}\km\W{0}{6}
+\frac{1}{4} \W{0}{1}\km\W{1}{1}\km\W{2}{3}\km\W{4}{4}\km\W{0}{5}\km\W{3}{6}
\\&
-\frac{3}{16} \W{0}{1}\km\W{1}{1}\km\W{2}{3}\km\W{4}{4}\km\W{1}{5}\km\W{2}{6}
-\frac{3}{8} \W{0}{1}\km\W{1}{1}\km\W{2}{3}\km\W{4}{4}\km\W{2}{5}\km\W{1}{6}
+\frac{5}{16} \W{0}{1}\km\W{1}{1}\km\W{2}{3}\km\W{4}{4}\km\W{3}{5}\km\W{0}{6}
+\frac{1}{60} \W{0}{1}\km\W{1}{1}\km\W{3}{3}\km\W{0}{4}\km\W{0}{5}\km\W{6}{6}
-\frac{1}{12} \W{0}{1}\km\W{1}{1}\km\W{3}{3}\km\W{0}{4}\km\W{1}{5}\km\W{5}{6}
-\frac{7}{48} \W{0}{1}\km\W{1}{1}\km\W{3}{3}\km\W{0}{4}\km\W{2}{5}\km\W{4}{6}
\\&
+\frac{7}{12} \W{0}{1}\km\W{1}{1}\km\W{3}{3}\km\W{0}{4}\km\W{3}{5}\km\W{3}{6}
+\frac{13}{120} \W{0}{1}\km\W{1}{1}\km\W{3}{3}\km\W{0}{4}\km\W{5}{5}\km\W{1}{6}
-\frac{23}{48} \W{0}{1}\km\W{1}{1}\km\W{3}{3}\km\W{0}{4}\km\W{4}{5}\km\W{2}{6}
-\frac{1}{60} \W{0}{1}\km\W{1}{1}\km\W{3}{3}\km\W{1}{4}\km\W{0}{5}\km\W{5}{6}
+\frac{17}{24} \W{0}{1}\km\W{1}{1}\km\W{3}{3}\km\W{1}{4}\km\W{1}{5}\km\W{4}{6}
-\frac{7}{6} \W{0}{1}\km\W{1}{1}\km\W{3}{3}\km\W{1}{4}\km\W{2}{5}\km\W{3}{6}
\\&
+\frac{1}{6} \W{0}{1}\km\W{1}{1}\km\W{3}{3}\km\W{1}{4}\km\W{3}{5}\km\W{2}{6}
+\frac{5}{12} \W{0}{1}\km\W{1}{1}\km\W{3}{3}\km\W{1}{4}\km\W{4}{5}\km\W{1}{6}
-\frac{13}{120} \W{0}{1}\km\W{1}{1}\km\W{3}{3}\km\W{1}{4}\km\W{5}{5}\km\W{0}{6}
-\frac{5}{16} \W{0}{1}\km\W{1}{1}\km\W{3}{3}\km\W{2}{4}\km\W{0}{5}\km\W{4}{6}
-\frac{1}{4} \W{0}{1}\km\W{1}{1}\km\W{3}{3}\km\W{2}{4}\km\W{1}{5}\km\W{3}{6}
+\frac{3}{2} \W{0}{1}\km\W{1}{1}\km\W{3}{3}\km\W{2}{4}\km\W{2}{5}\km\W{2}{6}
\\&
-\W{0}{1}\km\W{1}{1}\km\W{3}{3}\km\W{2}{4}\km\W{3}{5}\km\W{1}{6}
+\frac{1}{16} \W{0}{1}\km\W{1}{1}\km\W{3}{3}\km\W{2}{4}\km\W{4}{5}\km\W{0}{6}
-\frac{3}{4} \W{0}{1}\km\W{1}{1}\km\W{3}{3}\km\W{3}{4}\km\W{1}{5}\km\W{2}{6}
+\frac{1}{2} \W{0}{1}\km\W{1}{1}\km\W{3}{3}\km\W{3}{4}\km\W{0}{5}\km\W{3}{6}
+\frac{1}{4} \W{0}{1}\km\W{1}{1}\km\W{3}{3}\km\W{3}{4}\km\W{3}{5}\km\W{0}{6}
-\frac{3}{16} \W{0}{1}\km\W{1}{1}\km\W{3}{3}\km\W{4}{4}\km\W{2}{5}\km\W{0}{6}
\\&
-\frac{3}{16} \W{0}{1}\km\W{1}{1}\km\W{3}{3}\km\W{4}{4}\km\W{0}{5}\km\W{2}{6}
+\frac{3}{8} \W{0}{1}\km\W{1}{1}\km\W{3}{3}\km\W{4}{4}\km\W{1}{5}\km\W{1}{6}
\\
\widehat{\mathbf{q}}_{6} &= 
+\frac{1}{8} \W{0}{1}\km\W{1}{1}\km\W{0}{3}\km\W{0}{4}\km\W{3}{5}\km\W{6}{6}
+\frac{1}{8} \W{0}{1}\km\W{1}{1}\km\W{0}{3}\km\W{0}{4}\km\W{5}{5}\km\W{4}{6}
-\frac{1}{4} \W{0}{1}\km\W{1}{1}\km\W{0}{3}\km\W{0}{4}\km\W{4}{5}\km\W{5}{6}
-\frac{1}{6} \W{0}{1}\km\W{1}{1}\km\W{0}{3}\km\W{1}{4}\km\W{2}{5}\km\W{6}{6}
+\frac{1}{2} \W{0}{1}\km\W{1}{1}\km\W{0}{3}\km\W{1}{4}\km\W{4}{5}\km\W{4}{6}
-\frac{1}{3} \W{0}{1}\km\W{1}{1}\km\W{0}{3}\km\W{1}{4}\km\W{5}{5}\km\W{3}{6}
\\&
+\frac{1}{24} \W{0}{1}\km\W{1}{1}\km\W{0}{3}\km\W{2}{4}\km\W{1}{5}\km\W{6}{6}
+\frac{1}{3} \W{0}{1}\km\W{1}{1}\km\W{0}{3}\km\W{2}{4}\km\W{2}{5}\km\W{5}{6}
-\frac{1}{2} \W{0}{1}\km\W{1}{1}\km\W{0}{3}\km\W{2}{4}\km\W{3}{5}\km\W{4}{6}
-\frac{1}{6} \W{0}{1}\km\W{1}{1}\km\W{0}{3}\km\W{2}{4}\km\W{4}{5}\km\W{3}{6}
+\frac{7}{24} \W{0}{1}\km\W{1}{1}\km\W{0}{3}\km\W{2}{4}\km\W{5}{5}\km\W{2}{6}
+\frac{1}{90} \W{0}{1}\km\W{1}{1}\km\W{0}{3}\km\W{3}{4}\km\W{0}{5}\km\W{6}{6}
\\&
-\frac{5}{36} \W{0}{1}\km\W{1}{1}\km\W{0}{3}\km\W{3}{4}\km\W{1}{5}\km\W{5}{6}
-\frac{1}{18} \W{0}{1}\km\W{1}{1}\km\W{0}{3}\km\W{3}{4}\km\W{2}{5}\km\W{4}{6}
+\frac{7}{18} \W{0}{1}\km\W{1}{1}\km\W{0}{3}\km\W{3}{4}\km\W{3}{5}\km\W{3}{6}
-\frac{1}{9} \W{0}{1}\km\W{1}{1}\km\W{0}{3}\km\W{3}{4}\km\W{4}{5}\km\W{2}{6}
-\frac{1}{90} \W{0}{1}\km\W{1}{1}\km\W{0}{3}\km\W{4}{4}\km\W{0}{5}\km\W{5}{6}
-\frac{17}{180} \W{0}{1}\km\W{1}{1}\km\W{0}{3}\km\W{3}{4}\km\W{5}{5}\km\W{1}{6}
\\&
+\frac{7}{72} \W{0}{1}\km\W{1}{1}\km\W{0}{3}\km\W{4}{4}\km\W{1}{5}\km\W{4}{6}
-\frac{1}{9} \W{0}{1}\km\W{1}{1}\km\W{0}{3}\km\W{4}{4}\km\W{2}{5}\km\W{3}{6}
+\frac{1}{36} \W{0}{1}\km\W{1}{1}\km\W{0}{3}\km\W{4}{4}\km\W{4}{5}\km\W{1}{6}
-\frac{1}{72} \W{0}{1}\km\W{1}{1}\km\W{0}{3}\km\W{4}{4}\km\W{3}{5}\km\W{2}{6}
+\frac{1}{90} \W{0}{1}\km\W{1}{1}\km\W{0}{3}\km\W{4}{4}\km\W{5}{5}\km\W{0}{6}
-\frac{5}{24} \W{0}{1}\km\W{1}{1}\km\W{1}{3}\km\W{0}{4}\km\W{2}{5}\km\W{6}{6}
\\&
+\frac{1}{4} \W{0}{1}\km\W{1}{1}\km\W{1}{3}\km\W{0}{4}\km\W{3}{5}\km\W{5}{6}
+\frac{1}{8} \W{0}{1}\km\W{1}{1}\km\W{1}{3}\km\W{0}{4}\km\W{4}{5}\km\W{4}{6}
-\frac{1}{6} \W{0}{1}\km\W{1}{1}\km\W{1}{3}\km\W{0}{4}\km\W{5}{5}\km\W{3}{6}
+\frac{1}{4} \W{0}{1}\km\W{1}{1}\km\W{1}{3}\km\W{1}{4}\km\W{1}{5}\km\W{6}{6}
+\frac{1}{3} \W{0}{1}\km\W{1}{1}\km\W{1}{3}\km\W{1}{4}\km\W{2}{5}\km\W{5}{6}
-\W{0}{1}\km\W{1}{1}\km\W{1}{3}\km\W{1}{4}\km\W{3}{5}\km\W{4}{6}
\\&
-\frac{3}{40} \W{0}{1}\km\W{1}{1}\km\W{1}{3}\km\W{2}{4}\km\W{0}{5}\km\W{6}{6}
+\frac{5}{12} \W{0}{1}\km\W{1}{1}\km\W{1}{3}\km\W{1}{4}\km\W{5}{5}\km\W{2}{6}
-\frac{1}{2} \W{0}{1}\km\W{1}{1}\km\W{1}{3}\km\W{2}{4}\km\W{1}{5}\km\W{5}{6}
+\frac{3}{2} \W{0}{1}\km\W{1}{1}\km\W{1}{3}\km\W{2}{4}\km\W{3}{5}\km\W{3}{6}
-\frac{5}{8} \W{0}{1}\km\W{1}{1}\km\W{1}{3}\km\W{2}{4}\km\W{4}{5}\km\W{2}{6}
-\frac{3}{10} \W{0}{1}\km\W{1}{1}\km\W{1}{3}\km\W{2}{4}\km\W{5}{5}\km\W{1}{6}
\\&
+\frac{7}{60} \W{0}{1}\km\W{1}{1}\km\W{1}{3}\km\W{3}{4}\km\W{0}{5}\km\W{5}{6}
-\frac{1}{2} \W{0}{1}\km\W{1}{1}\km\W{1}{3}\km\W{3}{4}\km\W{2}{5}\km\W{3}{6}
+\frac{5}{12} \W{0}{1}\km\W{1}{1}\km\W{1}{3}\km\W{3}{4}\km\W{1}{5}\km\W{4}{6}
-\frac{2}{3} \W{0}{1}\km\W{1}{1}\km\W{1}{3}\km\W{3}{4}\km\W{3}{5}\km\W{2}{6}
+\frac{7}{12} \W{0}{1}\km\W{1}{1}\km\W{1}{3}\km\W{3}{4}\km\W{4}{5}\km\W{1}{6}
+\frac{1}{20} \W{0}{1}\km\W{1}{1}\km\W{1}{3}\km\W{3}{4}\km\W{5}{5}\km\W{0}{6}
\\&
-\frac{1}{24} \W{0}{1}\km\W{1}{1}\km\W{1}{3}\km\W{4}{4}\km\W{0}{5}\km\W{4}{6}
-\frac{1}{6} \W{0}{1}\km\W{1}{1}\km\W{1}{3}\km\W{4}{4}\km\W{1}{5}\km\W{3}{6}
+\frac{3}{8} \W{0}{1}\km\W{1}{1}\km\W{1}{3}\km\W{4}{4}\km\W{2}{5}\km\W{2}{6}
-\frac{1}{12} \W{0}{1}\km\W{1}{1}\km\W{1}{3}\km\W{4}{4}\km\W{3}{5}\km\W{1}{6}
-\frac{1}{12} \W{0}{1}\km\W{1}{1}\km\W{1}{3}\km\W{4}{4}\km\W{4}{5}\km\W{0}{6}
+\frac{1}{12} \W{0}{1}\km\W{1}{1}\km\W{2}{3}\km\W{0}{4}\km\W{1}{5}\km\W{6}{6}
\\&
+\frac{1}{12} \W{0}{1}\km\W{1}{1}\km\W{2}{3}\km\W{0}{4}\km\W{2}{5}\km\W{5}{6}
-\frac{3}{8} \W{0}{1}\km\W{1}{1}\km\W{2}{3}\km\W{0}{4}\km\W{3}{5}\km\W{4}{6}
+\frac{1}{6} \W{0}{1}\km\W{1}{1}\km\W{2}{3}\km\W{0}{4}\km\W{4}{5}\km\W{3}{6}
+\frac{1}{24} \W{0}{1}\km\W{1}{1}\km\W{2}{3}\km\W{0}{4}\km\W{5}{5}\km\W{2}{6}
-\frac{1}{20} \W{0}{1}\km\W{1}{1}\km\W{2}{3}\km\W{1}{4}\km\W{0}{5}\km\W{6}{6}
-\frac{7}{12} \W{0}{1}\km\W{1}{1}\km\W{2}{3}\km\W{1}{4}\km\W{1}{5}\km\W{5}{6}
\\&
+\frac{2}{3} \W{0}{1}\km\W{1}{1}\km\W{2}{3}\km\W{1}{4}\km\W{2}{5}\km\W{4}{6}
+\frac{1}{2} \W{0}{1}\km\W{1}{1}\km\W{2}{3}\km\W{1}{4}\km\W{3}{5}\km\W{3}{6}
-\frac{5}{12} \W{0}{1}\km\W{1}{1}\km\W{2}{3}\km\W{1}{4}\km\W{4}{5}\km\W{2}{6}
-\frac{7}{60} \W{0}{1}\km\W{1}{1}\km\W{2}{3}\km\W{1}{4}\km\W{5}{5}\km\W{1}{6}
+\frac{3}{10} \W{0}{1}\km\W{1}{1}\km\W{2}{3}\km\W{2}{4}\km\W{0}{5}\km\W{5}{6}
+\frac{5}{8} \W{0}{1}\km\W{1}{1}\km\W{2}{3}\km\W{2}{4}\km\W{1}{5}\km\W{4}{6}
\\&
-\frac{3}{2} \W{0}{1}\km\W{1}{1}\km\W{2}{3}\km\W{2}{4}\km\W{2}{5}\km\W{3}{6}
+\frac{1}{2} \W{0}{1}\km\W{1}{1}\km\W{2}{3}\km\W{2}{4}\km\W{4}{5}\km\W{1}{6}
+\frac{3}{40} \W{0}{1}\km\W{1}{1}\km\W{2}{3}\km\W{2}{4}\km\W{5}{5}\km\W{0}{6}
-\frac{5}{12} \W{0}{1}\km\W{1}{1}\km\W{2}{3}\km\W{3}{4}\km\W{0}{5}\km\W{4}{6}
+\W{0}{1}\km\W{1}{1}\km\W{2}{3}\km\W{3}{4}\km\W{2}{5}\km\W{2}{6}
-\frac{1}{3} \W{0}{1}\km\W{1}{1}\km\W{2}{3}\km\W{3}{4}\km\W{3}{5}\km\W{1}{6}
\\&
-\frac{1}{4} \W{0}{1}\km\W{1}{1}\km\W{2}{3}\km\W{3}{4}\km\W{4}{5}\km\W{0}{6}
+\frac{1}{6} \W{0}{1}\km\W{1}{1}\km\W{2}{3}\km\W{4}{4}\km\W{0}{5}\km\W{3}{6}
-\frac{1}{8} \W{0}{1}\km\W{1}{1}\km\W{2}{3}\km\W{4}{4}\km\W{1}{5}\km\W{2}{6}
-\frac{1}{4} \W{0}{1}\km\W{1}{1}\km\W{2}{3}\km\W{4}{4}\km\W{2}{5}\km\W{1}{6}
+\frac{5}{24} \W{0}{1}\km\W{1}{1}\km\W{2}{3}\km\W{4}{4}\km\W{3}{5}\km\W{0}{6}
-\frac{1}{90} \W{0}{1}\km\W{1}{1}\km\W{3}{3}\km\W{0}{4}\km\W{0}{5}\km\W{6}{6}
\\&
-\frac{1}{36} \W{0}{1}\km\W{1}{1}\km\W{3}{3}\km\W{0}{4}\km\W{1}{5}\km\W{5}{6}
+\frac{1}{72} \W{0}{1}\km\W{1}{1}\km\W{3}{3}\km\W{0}{4}\km\W{2}{5}\km\W{4}{6}
+\frac{1}{9} \W{0}{1}\km\W{1}{1}\km\W{3}{3}\km\W{0}{4}\km\W{3}{5}\km\W{3}{6}
+\frac{1}{90} \W{0}{1}\km\W{1}{1}\km\W{3}{3}\km\W{0}{4}\km\W{5}{5}\km\W{1}{6}
-\frac{7}{72} \W{0}{1}\km\W{1}{1}\km\W{3}{3}\km\W{0}{4}\km\W{4}{5}\km\W{2}{6}
+\frac{17}{180} \W{0}{1}\km\W{1}{1}\km\W{3}{3}\km\W{1}{4}\km\W{0}{5}\km\W{5}{6}
\\&
+\frac{1}{9} \W{0}{1}\km\W{1}{1}\km\W{3}{3}\km\W{1}{4}\km\W{1}{5}\km\W{4}{6}
-\frac{7}{18} \W{0}{1}\km\W{1}{1}\km\W{3}{3}\km\W{1}{4}\km\W{2}{5}\km\W{3}{6}
+\frac{1}{18} \W{0}{1}\km\W{1}{1}\km\W{3}{3}\km\W{1}{4}\km\W{3}{5}\km\W{2}{6}
+\frac{5}{36} \W{0}{1}\km\W{1}{1}\km\W{3}{3}\km\W{1}{4}\km\W{4}{5}\km\W{1}{6}
-\frac{1}{90} \W{0}{1}\km\W{1}{1}\km\W{3}{3}\km\W{1}{4}\km\W{5}{5}\km\W{0}{6}
-\frac{7}{24} \W{0}{1}\km\W{1}{1}\km\W{3}{3}\km\W{2}{4}\km\W{0}{5}\km\W{4}{6}
\\&
+\frac{1}{6} \W{0}{1}\km\W{1}{1}\km\W{3}{3}\km\W{2}{4}\km\W{1}{5}\km\W{3}{6}
+\frac{1}{2} \W{0}{1}\km\W{1}{1}\km\W{3}{3}\km\W{2}{4}\km\W{2}{5}\km\W{2}{6}
-\frac{1}{3} \W{0}{1}\km\W{1}{1}\km\W{3}{3}\km\W{2}{4}\km\W{3}{5}\km\W{1}{6}
-\frac{1}{24} \W{0}{1}\km\W{1}{1}\km\W{3}{3}\km\W{2}{4}\km\W{4}{5}\km\W{0}{6}
-\frac{1}{2} \W{0}{1}\km\W{1}{1}\km\W{3}{3}\km\W{3}{4}\km\W{1}{5}\km\W{2}{6}
+\frac{1}{3} \W{0}{1}\km\W{1}{1}\km\W{3}{3}\km\W{3}{4}\km\W{0}{5}\km\W{3}{6}
\\&
+\frac{1}{6} \W{0}{1}\km\W{1}{1}\km\W{3}{3}\km\W{3}{4}\km\W{3}{5}\km\W{0}{6}
-\frac{1}{8} \W{0}{1}\km\W{1}{1}\km\W{3}{3}\km\W{4}{4}\km\W{2}{5}\km\W{0}{6}
-\frac{1}{8} \W{0}{1}\km\W{1}{1}\km\W{3}{3}\km\W{4}{4}\km\W{0}{5}\km\W{2}{6}
+\frac{1}{4} \W{0}{1}\km\W{1}{1}\km\W{3}{3}\km\W{4}{4}\km\W{1}{5}\km\W{1}{6}
\\
\widehat{\mathbf{q}}_{7} &= 
+\frac{1}{16} \W{0}{1}\km\W{1}{1}\km\W{0}{3}\km\W{0}{4}\km\W{3}{5}\km\W{6}{6}
+\frac{1}{16} \W{0}{1}\km\W{1}{1}\km\W{0}{3}\km\W{0}{4}\km\W{5}{5}\km\W{4}{6}
-\frac{1}{8} \W{0}{1}\km\W{1}{1}\km\W{0}{3}\km\W{0}{4}\km\W{4}{5}\km\W{5}{6}
+\frac{1}{4} \W{0}{1}\km\W{1}{1}\km\W{0}{3}\km\W{1}{4}\km\W{2}{5}\km\W{6}{6}
-\W{0}{1}\km\W{1}{1}\km\W{0}{3}\km\W{1}{4}\km\W{3}{5}\km\W{5}{6}
+\frac{5}{4} \W{0}{1}\km\W{1}{1}\km\W{0}{3}\km\W{1}{4}\km\W{4}{5}\km\W{4}{6}
\\&
-\frac{1}{2} \W{0}{1}\km\W{1}{1}\km\W{0}{3}\km\W{1}{4}\km\W{5}{5}\km\W{3}{6}
-\frac{7}{16} \W{0}{1}\km\W{1}{1}\km\W{0}{3}\km\W{2}{4}\km\W{1}{5}\km\W{6}{6}
+\W{0}{1}\km\W{1}{1}\km\W{0}{3}\km\W{2}{4}\km\W{2}{5}\km\W{5}{6}
-\frac{5}{4} \W{0}{1}\km\W{1}{1}\km\W{0}{3}\km\W{2}{4}\km\W{4}{5}\km\W{3}{6}
+\frac{11}{16} \W{0}{1}\km\W{1}{1}\km\W{0}{3}\km\W{2}{4}\km\W{5}{5}\km\W{2}{6}
+\frac{13}{120} \W{0}{1}\km\W{1}{1}\km\W{0}{3}\km\W{3}{4}\km\W{0}{5}\km\W{6}{6}
\\&
+\frac{1}{3} \W{0}{1}\km\W{1}{1}\km\W{0}{3}\km\W{3}{4}\km\W{1}{5}\km\W{5}{6}
-\frac{5}{3} \W{0}{1}\km\W{1}{1}\km\W{0}{3}\km\W{3}{4}\km\W{2}{5}\km\W{4}{6}
+\frac{5}{3} \W{0}{1}\km\W{1}{1}\km\W{0}{3}\km\W{3}{4}\km\W{3}{5}\km\W{3}{6}
-\frac{5}{24} \W{0}{1}\km\W{1}{1}\km\W{0}{3}\km\W{3}{4}\km\W{4}{5}\km\W{2}{6}
-\frac{13}{120} \W{0}{1}\km\W{1}{1}\km\W{0}{3}\km\W{4}{4}\km\W{0}{5}\km\W{5}{6}
-\frac{7}{30} \W{0}{1}\km\W{1}{1}\km\W{0}{3}\km\W{3}{4}\km\W{5}{5}\km\W{1}{6}
\\&
+\frac{5}{48} \W{0}{1}\km\W{1}{1}\km\W{0}{3}\km\W{4}{4}\km\W{1}{5}\km\W{4}{6}
+\frac{5}{12} \W{0}{1}\km\W{1}{1}\km\W{0}{3}\km\W{4}{4}\km\W{2}{5}\km\W{3}{6}
+\frac{1}{3} \W{0}{1}\km\W{1}{1}\km\W{0}{3}\km\W{4}{4}\km\W{4}{5}\km\W{1}{6}
-\frac{35}{48} \W{0}{1}\km\W{1}{1}\km\W{0}{3}\km\W{4}{4}\km\W{3}{5}\km\W{2}{6}
-\frac{1}{60} \W{0}{1}\km\W{1}{1}\km\W{0}{3}\km\W{4}{4}\km\W{5}{5}\km\W{0}{6}
-\frac{7}{16} \W{0}{1}\km\W{1}{1}\km\W{1}{3}\km\W{0}{4}\km\W{2}{5}\km\W{6}{6}
\\&
+\frac{9}{8} \W{0}{1}\km\W{1}{1}\km\W{1}{3}\km\W{0}{4}\km\W{3}{5}\km\W{5}{6}
-\frac{15}{16} \W{0}{1}\km\W{1}{1}\km\W{1}{3}\km\W{0}{4}\km\W{4}{5}\km\W{4}{6}
+\frac{1}{4} \W{0}{1}\km\W{1}{1}\km\W{1}{3}\km\W{0}{4}\km\W{5}{5}\km\W{3}{6}
+\frac{3}{8} \W{0}{1}\km\W{1}{1}\km\W{1}{3}\km\W{1}{4}\km\W{1}{5}\km\W{6}{6}
-\frac{1}{2} \W{0}{1}\km\W{1}{1}\km\W{1}{3}\km\W{1}{4}\km\W{2}{5}\km\W{5}{6}
+\frac{9}{80} \W{0}{1}\km\W{1}{1}\km\W{1}{3}\km\W{2}{4}\km\W{0}{5}\km\W{6}{6}
\\&
+\frac{1}{8} \W{0}{1}\km\W{1}{1}\km\W{1}{3}\km\W{1}{4}\km\W{5}{5}\km\W{2}{6}
-\frac{3}{8} \W{0}{1}\km\W{1}{1}\km\W{1}{3}\km\W{2}{4}\km\W{1}{5}\km\W{5}{6}
+\frac{15}{16} \W{0}{1}\km\W{1}{1}\km\W{1}{3}\km\W{2}{4}\km\W{4}{5}\km\W{2}{6}
-\frac{27}{40} \W{0}{1}\km\W{1}{1}\km\W{1}{3}\km\W{2}{4}\km\W{5}{5}\km\W{1}{6}
-\frac{11}{20} \W{0}{1}\km\W{1}{1}\km\W{1}{3}\km\W{3}{4}\km\W{0}{5}\km\W{5}{6}
+\frac{5}{4} \W{0}{1}\km\W{1}{1}\km\W{1}{3}\km\W{3}{4}\km\W{1}{5}\km\W{4}{6}
\\&
-\frac{5}{4} \W{0}{1}\km\W{1}{1}\km\W{1}{3}\km\W{3}{4}\km\W{3}{5}\km\W{2}{6}
+\frac{1}{4} \W{0}{1}\km\W{1}{1}\km\W{1}{3}\km\W{3}{4}\km\W{4}{5}\km\W{1}{6}
+\frac{3}{10} \W{0}{1}\km\W{1}{1}\km\W{1}{3}\km\W{3}{4}\km\W{5}{5}\km\W{0}{6}
+\frac{7}{16} \W{0}{1}\km\W{1}{1}\km\W{1}{3}\km\W{4}{4}\km\W{0}{5}\km\W{4}{6}
-\frac{5}{4} \W{0}{1}\km\W{1}{1}\km\W{1}{3}\km\W{4}{4}\km\W{1}{5}\km\W{3}{6}
+\frac{15}{16} \W{0}{1}\km\W{1}{1}\km\W{1}{3}\km\W{4}{4}\km\W{2}{5}\km\W{2}{6}
\\&
+\frac{1}{8} \W{0}{1}\km\W{1}{1}\km\W{1}{3}\km\W{4}{4}\km\W{3}{5}\km\W{1}{6}
-\frac{1}{4} \W{0}{1}\km\W{1}{1}\km\W{1}{3}\km\W{4}{4}\km\W{4}{5}\km\W{0}{6}
+\frac{1}{4} \W{0}{1}\km\W{1}{1}\km\W{2}{3}\km\W{0}{4}\km\W{1}{5}\km\W{6}{6}
-\frac{1}{8} \W{0}{1}\km\W{1}{1}\km\W{2}{3}\km\W{0}{4}\km\W{2}{5}\km\W{5}{6}
-\frac{15}{16} \W{0}{1}\km\W{1}{1}\km\W{2}{3}\km\W{0}{4}\km\W{3}{5}\km\W{4}{6}
+\frac{5}{4} \W{0}{1}\km\W{1}{1}\km\W{2}{3}\km\W{0}{4}\km\W{4}{5}\km\W{3}{6}
\\&
-\frac{7}{16} \W{0}{1}\km\W{1}{1}\km\W{2}{3}\km\W{0}{4}\km\W{5}{5}\km\W{2}{6}
-\frac{3}{10} \W{0}{1}\km\W{1}{1}\km\W{2}{3}\km\W{1}{4}\km\W{0}{5}\km\W{6}{6}
-\frac{1}{4} \W{0}{1}\km\W{1}{1}\km\W{2}{3}\km\W{1}{4}\km\W{1}{5}\km\W{5}{6}
+\frac{5}{4} \W{0}{1}\km\W{1}{1}\km\W{2}{3}\km\W{1}{4}\km\W{2}{5}\km\W{4}{6}
-\frac{5}{4} \W{0}{1}\km\W{1}{1}\km\W{2}{3}\km\W{1}{4}\km\W{4}{5}\km\W{2}{6}
+\frac{11}{20} \W{0}{1}\km\W{1}{1}\km\W{2}{3}\km\W{1}{4}\km\W{5}{5}\km\W{1}{6}
\\&
+\frac{27}{40} \W{0}{1}\km\W{1}{1}\km\W{2}{3}\km\W{2}{4}\km\W{0}{5}\km\W{5}{6}
-\frac{15}{16} \W{0}{1}\km\W{1}{1}\km\W{2}{3}\km\W{2}{4}\km\W{1}{5}\km\W{4}{6}
+\frac{3}{8} \W{0}{1}\km\W{1}{1}\km\W{2}{3}\km\W{2}{4}\km\W{4}{5}\km\W{1}{6}
-\frac{9}{80} \W{0}{1}\km\W{1}{1}\km\W{2}{3}\km\W{2}{4}\km\W{5}{5}\km\W{0}{6}
-\frac{1}{8} \W{0}{1}\km\W{1}{1}\km\W{2}{3}\km\W{3}{4}\km\W{0}{5}\km\W{4}{6}
+\frac{1}{2} \W{0}{1}\km\W{1}{1}\km\W{2}{3}\km\W{3}{4}\km\W{3}{5}\km\W{1}{6}
\\&
-\frac{3}{8} \W{0}{1}\km\W{1}{1}\km\W{2}{3}\km\W{3}{4}\km\W{4}{5}\km\W{0}{6}
-\frac{1}{4} \W{0}{1}\km\W{1}{1}\km\W{2}{3}\km\W{4}{4}\km\W{0}{5}\km\W{3}{6}
+\frac{15}{16} \W{0}{1}\km\W{1}{1}\km\W{2}{3}\km\W{4}{4}\km\W{1}{5}\km\W{2}{6}
-\frac{9}{8} \W{0}{1}\km\W{1}{1}\km\W{2}{3}\km\W{4}{4}\km\W{2}{5}\km\W{1}{6}
+\frac{7}{16} \W{0}{1}\km\W{1}{1}\km\W{2}{3}\km\W{4}{4}\km\W{3}{5}\km\W{0}{6}
+\frac{1}{60} \W{0}{1}\km\W{1}{1}\km\W{3}{3}\km\W{0}{4}\km\W{0}{5}\km\W{6}{6}
\\&
-\frac{1}{3} \W{0}{1}\km\W{1}{1}\km\W{3}{3}\km\W{0}{4}\km\W{1}{5}\km\W{5}{6}
+\frac{35}{48} \W{0}{1}\km\W{1}{1}\km\W{3}{3}\km\W{0}{4}\km\W{2}{5}\km\W{4}{6}
-\frac{5}{12} \W{0}{1}\km\W{1}{1}\km\W{3}{3}\km\W{0}{4}\km\W{3}{5}\km\W{3}{6}
+\frac{13}{120} \W{0}{1}\km\W{1}{1}\km\W{3}{3}\km\W{0}{4}\km\W{5}{5}\km\W{1}{6}
-\frac{5}{48} \W{0}{1}\km\W{1}{1}\km\W{3}{3}\km\W{0}{4}\km\W{4}{5}\km\W{2}{6}
+\frac{7}{30} \W{0}{1}\km\W{1}{1}\km\W{3}{3}\km\W{1}{4}\km\W{0}{5}\km\W{5}{6}
\\&
+\frac{5}{24} \W{0}{1}\km\W{1}{1}\km\W{3}{3}\km\W{1}{4}\km\W{1}{5}\km\W{4}{6}
-\frac{5}{3} \W{0}{1}\km\W{1}{1}\km\W{3}{3}\km\W{1}{4}\km\W{2}{5}\km\W{3}{6}
+\frac{5}{3} \W{0}{1}\km\W{1}{1}\km\W{3}{3}\km\W{1}{4}\km\W{3}{5}\km\W{2}{6}
-\frac{1}{3} \W{0}{1}\km\W{1}{1}\km\W{3}{3}\km\W{1}{4}\km\W{4}{5}\km\W{1}{6}
-\frac{13}{120} \W{0}{1}\km\W{1}{1}\km\W{3}{3}\km\W{1}{4}\km\W{5}{5}\km\W{0}{6}
-\frac{11}{16} \W{0}{1}\km\W{1}{1}\km\W{3}{3}\km\W{2}{4}\km\W{0}{5}\km\W{4}{6}
\\&
+\frac{5}{4} \W{0}{1}\km\W{1}{1}\km\W{3}{3}\km\W{2}{4}\km\W{1}{5}\km\W{3}{6}
-\W{0}{1}\km\W{1}{1}\km\W{3}{3}\km\W{2}{4}\km\W{3}{5}\km\W{1}{6}
+\frac{7}{16} \W{0}{1}\km\W{1}{1}\km\W{3}{3}\km\W{2}{4}\km\W{4}{5}\km\W{0}{6}
-\frac{5}{4} \W{0}{1}\km\W{1}{1}\km\W{3}{3}\km\W{3}{4}\km\W{1}{5}\km\W{2}{6}
+\frac{1}{2} \W{0}{1}\km\W{1}{1}\km\W{3}{3}\km\W{3}{4}\km\W{0}{5}\km\W{3}{6}
+\W{0}{1}\km\W{1}{1}\km\W{3}{3}\km\W{3}{4}\km\W{2}{5}\km\W{1}{6}
\\&
-\frac{1}{4} \W{0}{1}\km\W{1}{1}\km\W{3}{3}\km\W{3}{4}\km\W{3}{5}\km\W{0}{6}
-\frac{1}{16} \W{0}{1}\km\W{1}{1}\km\W{3}{3}\km\W{4}{4}\km\W{2}{5}\km\W{0}{6}
-\frac{1}{16} \W{0}{1}\km\W{1}{1}\km\W{3}{3}\km\W{4}{4}\km\W{0}{5}\km\W{2}{6}
+\frac{1}{8} \W{0}{1}\km\W{1}{1}\km\W{3}{3}\km\W{4}{4}\km\W{1}{5}\km\W{1}{6} 
\\
\widehat{\mathbf{q}}_{8}  &= 
+ \W{0}{1}\km\W{1}{1}\km\W{3}{3}\km\W{0}{5}\km\W{2}{5}\km\W{4}{5}
+\W{0}{1}\km\W{1}{1}\km\W{0}{3}\km\W{1}{5}\km\W{3}{5}\km\W{5}{5}
-\frac{2}{5} \W{0}{1}\km\W{1}{1}\km\W{0}{3}\km\W{0}{5}\km\W{4}{5}\km\W{5}{5}
-4 \W{0}{1}\km\W{1}{1}\km\W{0}{3}\km\W{2}{5}\km\W{3}{5}\km\W{4}{5}
+\frac{3}{5} \W{0}{1}\km\W{1}{1}\km\W{1}{3}\km\W{0}{5}\km\W{3}{5}\km\W{5}{5}
-3 \W{0}{1}\km\W{1}{1}\km\W{1}{3}\km\W{1}{5}\km\W{2}{5}\km\W{5}{5}
\\&
+3 \W{0}{1}\km\W{1}{1}\km\W{1}{3}\km\W{1}{5}\km\W{3}{5}\km\W{4}{5}
+\frac{3}{5} \W{0}{1}\km\W{1}{1}\km\W{2}{3}\km\W{0}{5}\km\W{2}{5}\km\W{5}{5}
-3 \W{0}{1}\km\W{1}{1}\km\W{2}{3}\km\W{0}{5}\km\W{3}{5}\km\W{4}{5}
+3 \W{0}{1}\km\W{1}{1}\km\W{2}{3}\km\W{1}{5}\km\W{2}{5}\km\W{4}{5}
-\frac{2}{5} \W{0}{1}\km\W{1}{1}\km\W{3}{3}\km\W{0}{5}\km\W{1}{5}\km\W{5}{5}
-4 \W{0}{1}\km\W{1}{1}\km\W{3}{3}\km\W{1}{5}\km\W{2}{5}\km\W{3}{5}
\\
\widehat{\mathbf{q}}_{9} &= 
+\W{0}{1}\km\W{1}{1}\km\W{0}{4}\km\W{2}{4}\km\W{4}{4}\km\W{3}{6}
+\W{0}{1}\km\W{1}{1}\km\W{0}{4}\km\W{1}{4}\km\W{3}{4}\km\W{5}{6}
+\W{0}{1}\km\W{1}{1}\km\W{1}{4}\km\W{3}{4}\km\W{4}{4}\km\W{1}{6}
-\frac{1}{2} \W{0}{1}\km\W{1}{1}\km\W{0}{4}\km\W{1}{4}\km\W{2}{4}\km\W{6}{6}
-\frac{1}{2} \W{0}{1}\km\W{1}{1}\km\W{0}{4}\km\W{1}{4}\km\W{4}{4}\km\W{4}{6}
-\frac{3}{2} \W{0}{1}\km\W{1}{1}\km\W{0}{4}\km\W{2}{4}\km\W{3}{4}\km\W{4}{6}
\\&
-\frac{1}{2} \W{0}{1}\km\W{1}{1}\km\W{0}{4}\km\W{3}{4}\km\W{4}{4}\km\W{2}{6}
+2 \W{0}{1}\km\W{1}{1}\km\W{1}{4}\km\W{2}{4}\km\W{3}{4}\km\W{3}{6}
-\frac{3}{2} \W{0}{1}\km\W{1}{1}\km\W{1}{4}\km\W{2}{4}\km\W{4}{4}\km\W{2}{6}
-\frac{1}{2} \W{0}{1}\km\W{1}{1}\km\W{2}{4}\km\W{3}{4}\km\W{4}{4}\km\W{0}{6}
\\
\widehat{\mathbf{q}}_{10}  &= 
-\frac{2}{3} \W{0}{1}\km\W{0}{3}\km\W{1}{3}\km\W{2}{3}\km\W{1}{4}\km\W{6}{6}
+\frac{2}{3} \W{0}{1}\km\W{0}{3}\km\W{1}{3}\km\W{2}{3}\km\W{4}{4}\km\W{3}{6}
+2 \W{0}{1}\km\W{0}{3}\km\W{1}{3}\km\W{2}{3}\km\W{2}{4}\km\W{5}{6}
-2 \W{0}{1}\km\W{0}{3}\km\W{1}{3}\km\W{2}{3}\km\W{3}{4}\km\W{4}{6}
+\frac{1}{4} \W{0}{1}\km\W{0}{3}\km\W{1}{3}\km\W{3}{3}\km\W{0}{4}\km\W{6}{6}
-\frac{1}{3} \W{0}{1}\km\W{0}{3}\km\W{1}{3}\km\W{3}{3}\km\W{1}{4}\km\W{5}{6}
\\&
-\frac{1}{2} \W{0}{1}\km\W{0}{3}\km\W{1}{3}\km\W{3}{3}\km\W{2}{4}\km\W{4}{6}
+\W{0}{1}\km\W{0}{3}\km\W{1}{3}\km\W{3}{3}\km\W{3}{4}\km\W{3}{6}
-\frac{5}{12} \W{0}{1}\km\W{0}{3}\km\W{1}{3}\km\W{3}{3}\km\W{4}{4}\km\W{2}{6}
-\frac{1}{2} \W{0}{1}\km\W{0}{3}\km\W{2}{3}\km\W{3}{3}\km\W{0}{4}\km\W{5}{6}
+\frac{4}{3} \W{0}{1}\km\W{0}{3}\km\W{2}{3}\km\W{3}{3}\km\W{1}{4}\km\W{4}{6}
-\W{0}{1}\km\W{0}{3}\km\W{2}{3}\km\W{3}{3}\km\W{2}{4}\km\W{3}{6}
\\&
+\frac{1}{6} \W{0}{1}\km\W{0}{3}\km\W{2}{3}\km\W{3}{3}\km\W{4}{4}\km\W{1}{6}
+\frac{3}{4} \W{0}{1}\km\W{1}{3}\km\W{2}{3}\km\W{3}{3}\km\W{0}{4}\km\W{4}{6}
-\frac{7}{3} \W{0}{1}\km\W{1}{3}\km\W{2}{3}\km\W{3}{3}\km\W{1}{4}\km\W{3}{6}
+\frac{5}{2} \W{0}{1}\km\W{1}{3}\km\W{2}{3}\km\W{3}{3}\km\W{2}{4}\km\W{2}{6}
-\W{0}{1}\km\W{1}{3}\km\W{2}{3}\km\W{3}{3}\km\W{3}{4}\km\W{1}{6}
+\frac{1}{12} \W{0}{1}\km\W{1}{3}\km\W{2}{3}\km\W{3}{3}\km\W{4}{4}\km\W{0}{6}
\\&
-\frac{1}{12} \W{1}{1}\km\W{0}{3}\km\W{1}{3}\km\W{2}{3}\km\W{0}{4}\km\W{6}{6}
+\W{1}{1}\km\W{0}{3}\km\W{1}{3}\km\W{2}{3}\km\W{1}{4}\km\W{5}{6}
-\frac{5}{2} \W{1}{1}\km\W{0}{3}\km\W{1}{3}\km\W{2}{3}\km\W{2}{4}\km\W{4}{6}
-\frac{3}{4} \W{1}{1}\km\W{0}{3}\km\W{1}{3}\km\W{2}{3}\km\W{4}{4}\km\W{2}{6}
+\frac{7}{3} \W{1}{1}\km\W{0}{3}\km\W{1}{3}\km\W{2}{3}\km\W{3}{4}\km\W{3}{6}
-\frac{1}{6} \W{1}{1}\km\W{0}{3}\km\W{1}{3}\km\W{3}{3}\km\W{0}{4}\km\W{5}{6}
\\&
-\frac{4}{3} \W{1}{1}\km\W{0}{3}\km\W{1}{3}\km\W{3}{3}\km\W{3}{4}\km\W{2}{6}
+\W{1}{1}\km\W{0}{3}\km\W{1}{3}\km\W{3}{3}\km\W{2}{4}\km\W{3}{6}
+\frac{5}{12} \W{1}{1}\km\W{0}{3}\km\W{2}{3}\km\W{3}{3}\km\W{0}{4}\km\W{4}{6}
+\frac{1}{2} \W{1}{1}\km\W{0}{3}\km\W{1}{3}\km\W{3}{3}\km\W{4}{4}\km\W{1}{6}
-\W{1}{1}\km\W{0}{3}\km\W{2}{3}\km\W{3}{3}\km\W{1}{4}\km\W{3}{6}
+\frac{1}{2} \W{1}{1}\km\W{0}{3}\km\W{2}{3}\km\W{3}{3}\km\W{2}{4}\km\W{2}{6}
\\&
+\frac{1}{3} \W{1}{1}\km\W{0}{3}\km\W{2}{3}\km\W{3}{3}\km\W{3}{4}\km\W{1}{6}
-\frac{1}{4} \W{1}{1}\km\W{0}{3}\km\W{2}{3}\km\W{3}{3}\km\W{4}{4}\km\W{0}{6}
-\frac{2}{3} \W{1}{1}\km\W{1}{3}\km\W{2}{3}\km\W{3}{3}\km\W{0}{4}\km\W{3}{6}
+2 \W{1}{1}\km\W{1}{3}\km\W{2}{3}\km\W{3}{3}\km\W{1}{4}\km\W{2}{6}
-2 \W{1}{1}\km\W{1}{3}\km\W{2}{3}\km\W{3}{3}\km\W{2}{4}\km\W{1}{6}
+\frac{2}{3} \W{1}{1}\km\W{1}{3}\km\W{2}{3}\km\W{3}{3}\km\W{3}{4}\km\W{0}{6}
\\
\widehat{\mathbf{q}}_{11}&= 
-\frac{1}{3} \W{0}{1}\km\W{0}{3}\km\W{1}{3}\km\W{2}{3}\km\W{1}{4}\km\W{6}{6}
+\frac{1}{3} \W{0}{1}\km\W{0}{3}\km\W{1}{3}\km\W{2}{3}\km\W{4}{4}\km\W{3}{6}
+\W{0}{1}\km\W{0}{3}\km\W{1}{3}\km\W{2}{3}\km\W{2}{4}\km\W{5}{6}
-\W{0}{1}\km\W{0}{3}\km\W{1}{3}\km\W{2}{3}\km\W{3}{4}\km\W{4}{6}
+\frac{1}{4} \W{0}{1}\km\W{0}{3}\km\W{1}{3}\km\W{3}{3}\km\W{0}{4}\km\W{6}{6}
-\frac{2}{3} \W{0}{1}\km\W{0}{3}\km\W{1}{3}\km\W{3}{3}\km\W{1}{4}\km\W{5}{6}
\\&
+\frac{1}{2} \W{0}{1}\km\W{0}{3}\km\W{1}{3}\km\W{3}{3}\km\W{2}{4}\km\W{4}{6}
-\frac{1}{12} \W{0}{1}\km\W{0}{3}\km\W{1}{3}\km\W{3}{3}\km\W{4}{4}\km\W{2}{6}
-\frac{1}{2} \W{0}{1}\km\W{0}{3}\km\W{2}{3}\km\W{3}{3}\km\W{0}{4}\km\W{5}{6}
+\frac{5}{3} \W{0}{1}\km\W{0}{3}\km\W{2}{3}\km\W{3}{3}\km\W{1}{4}\km\W{4}{6}
-2 \W{0}{1}\km\W{0}{3}\km\W{2}{3}\km\W{3}{3}\km\W{2}{4}\km\W{3}{6}
+\W{0}{1}\km\W{0}{3}\km\W{2}{3}\km\W{3}{3}\km\W{3}{4}\km\W{2}{6}
\\&
-\frac{1}{6} \W{0}{1}\km\W{0}{3}\km\W{2}{3}\km\W{3}{3}\km\W{4}{4}\km\W{1}{6}
+\frac{3}{4} \W{0}{1}\km\W{1}{3}\km\W{2}{3}\km\W{3}{3}\km\W{0}{4}\km\W{4}{6}
-\frac{8}{3} \W{0}{1}\km\W{1}{3}\km\W{2}{3}\km\W{3}{3}\km\W{1}{4}\km\W{3}{6}
+\frac{7}{2} \W{0}{1}\km\W{1}{3}\km\W{2}{3}\km\W{3}{3}\km\W{2}{4}\km\W{2}{6}
-2 \W{0}{1}\km\W{1}{3}\km\W{2}{3}\km\W{3}{3}\km\W{3}{4}\km\W{1}{6}
+\frac{5}{12} \W{0}{1}\km\W{1}{3}\km\W{2}{3}\km\W{3}{3}\km\W{4}{4}\km\W{0}{6}
\\&
-\frac{5}{12} \W{1}{1}\km\W{0}{3}\km\W{1}{3}\km\W{2}{3}\km\W{0}{4}\km\W{6}{6}
+2 \W{1}{1}\km\W{0}{3}\km\W{1}{3}\km\W{2}{3}\km\W{1}{4}\km\W{5}{6}
-\frac{7}{2} \W{1}{1}\km\W{0}{3}\km\W{1}{3}\km\W{2}{3}\km\W{2}{4}\km\W{4}{6}
-\frac{3}{4} \W{1}{1}\km\W{0}{3}\km\W{1}{3}\km\W{2}{3}\km\W{4}{4}\km\W{2}{6}
+\frac{8}{3} \W{1}{1}\km\W{0}{3}\km\W{1}{3}\km\W{2}{3}\km\W{3}{4}\km\W{3}{6}
+\frac{1}{6} \W{1}{1}\km\W{0}{3}\km\W{1}{3}\km\W{3}{3}\km\W{0}{4}\km\W{5}{6}
\\&
-\W{1}{1}\km\W{0}{3}\km\W{1}{3}\km\W{3}{3}\km\W{1}{4}\km\W{4}{6}
-\frac{5}{3} \W{1}{1}\km\W{0}{3}\km\W{1}{3}\km\W{3}{3}\km\W{3}{4}\km\W{2}{6}
+2 \W{1}{1}\km\W{0}{3}\km\W{1}{3}\km\W{3}{3}\km\W{2}{4}\km\W{3}{6}
+\frac{1}{12} \W{1}{1}\km\W{0}{3}\km\W{2}{3}\km\W{3}{3}\km\W{0}{4}\km\W{4}{6}
+\frac{1}{2} \W{1}{1}\km\W{0}{3}\km\W{1}{3}\km\W{3}{3}\km\W{4}{4}\km\W{1}{6}
-\frac{1}{2} \W{1}{1}\km\W{0}{3}\km\W{2}{3}\km\W{3}{3}\km\W{2}{4}\km\W{2}{6}
\\&
+\frac{2}{3} \W{1}{1}\km\W{0}{3}\km\W{2}{3}\km\W{3}{3}\km\W{3}{4}\km\W{1}{6}
-\frac{1}{4} \W{1}{1}\km\W{0}{3}\km\W{2}{3}\km\W{3}{3}\km\W{4}{4}\km\W{0}{6}
-\frac{1}{3} \W{1}{1}\km\W{1}{3}\km\W{2}{3}\km\W{3}{3}\km\W{0}{4}\km\W{3}{6}
+\W{1}{1}\km\W{1}{3}\km\W{2}{3}\km\W{3}{3}\km\W{1}{4}\km\W{2}{6}
-\W{1}{1}\km\W{1}{3}\km\W{2}{3}\km\W{3}{3}\km\W{2}{4}\km\W{1}{6}
+\frac{1}{3} \W{1}{1}\km\W{1}{3}\km\W{2}{3}\km\W{3}{3}\km\W{3}{4}\km\W{0}{6}]
\\
\widehat{\mathbf{q}}_{12} &= 
+\W{1}{1}\km\W{0}{3}\km\W{1}{3}\km\W{3}{3}\km\W{1}{5}\km\W{4}{5}
+\W{0}{1}\km\W{0}{3}\km\W{2}{3}\km\W{3}{3}\km\W{0}{5}\km\W{5}{5}
+\W{0}{1}\km\W{0}{3}\km\W{2}{3}\km\W{3}{3}\km\W{1}{5}\km\W{4}{5}
+\W{1}{1}\km\W{0}{3}\km\W{1}{3}\km\W{3}{3}\km\W{0}{5}\km\W{5}{5}
+6 \W{0}{1}\km\W{0}{3}\km\W{1}{3}\km\W{2}{3}\km\W{2}{5}\km\W{5}{5}
-18 \W{0}{1}\km\W{0}{3}\km\W{1}{3}\km\W{2}{3}\km\W{3}{5}\km\W{4}{5}
\\&
-3 \W{0}{1}\km\W{0}{3}\km\W{1}{3}\km\W{3}{3}\km\W{1}{5}\km\W{5}{5}
+6 \W{0}{1}\km\W{0}{3}\km\W{1}{3}\km\W{3}{3}\km\W{2}{5}\km\W{4}{5}
-8 \W{0}{1}\km\W{0}{3}\km\W{2}{3}\km\W{3}{3}\km\W{2}{5}\km\W{3}{5}
-3 \W{0}{1}\km\W{1}{3}\km\W{2}{3}\km\W{3}{3}\km\W{0}{5}\km\W{4}{5}
+6 \W{0}{1}\km\W{1}{3}\km\W{2}{3}\km\W{3}{3}\km\W{1}{5}\km\W{3}{5}
-3 \W{1}{1}\km\W{0}{3}\km\W{1}{3}\km\W{2}{3}\km\W{1}{5}\km\W{5}{5}
\\&
+6 \W{1}{1}\km\W{0}{3}\km\W{1}{3}\km\W{2}{3}\km\W{2}{5}\km\W{4}{5}
-8 \W{1}{1}\km\W{0}{3}\km\W{1}{3}\km\W{3}{3}\km\W{2}{5}\km\W{3}{5}
-3 \W{1}{1}\km\W{0}{3}\km\W{2}{3}\km\W{3}{3}\km\W{0}{5}\km\W{4}{5}
+6 \W{1}{1}\km\W{0}{3}\km\W{2}{3}\km\W{3}{3}\km\W{1}{5}\km\W{3}{5}
+6 \W{1}{1}\km\W{1}{3}\km\W{2}{3}\km\W{3}{3}\km\W{0}{5}\km\W{3}{5}
-18 \W{1}{1}\km\W{1}{3}\km\W{2}{3}\km\W{3}{3}\km\W{1}{5}\km\W{2}{5}
\\
\widehat{\mathbf{q}}_{13} &= 
+\frac{1}{2} \W{0}{1}\km\W{0}{3}\km\W{1}{3}\km\W{0}{4}\km\W{4}{4}\km\W{5}{5}
-2 \W{0}{1}\km\W{0}{3}\km\W{1}{3}\km\W{1}{4}\km\W{4}{4}\km\W{4}{5}
+3 \W{0}{1}\km\W{0}{3}\km\W{1}{3}\km\W{2}{4}\km\W{4}{4}\km\W{3}{5}
-2 \W{0}{1}\km\W{0}{3}\km\W{1}{3}\km\W{3}{4}\km\W{4}{4}\km\W{2}{5}
-\W{0}{1}\km\W{0}{3}\km\W{2}{3}\km\W{0}{4}\km\W{3}{4}\km\W{5}{5}
+4 \W{0}{1}\km\W{0}{3}\km\W{2}{3}\km\W{1}{4}\km\W{3}{4}\km\W{4}{5}
\\&
-6 \W{0}{1}\km\W{0}{3}\km\W{2}{3}\km\W{2}{4}\km\W{3}{4}\km\W{3}{5}
+\W{0}{1}\km\W{0}{3}\km\W{2}{3}\km\W{3}{4}\km\W{4}{4}\km\W{1}{5}
+\frac{1}{2} \W{0}{1}\km\W{0}{3}\km\W{3}{3}\km\W{0}{4}\km\W{2}{4}\km\W{5}{5}
-2 \W{0}{1}\km\W{0}{3}\km\W{3}{3}\km\W{1}{4}\km\W{2}{4}\km\W{4}{5}
-\frac{1}{2} \W{0}{1}\km\W{0}{3}\km\W{3}{3}\km\W{2}{4}\km\W{4}{4}\km\W{1}{5}
+2 \W{0}{1}\km\W{0}{3}\km\W{3}{3}\km\W{2}{4}\km\W{3}{4}\km\W{2}{5}
\\&
+\frac{3}{2} \W{0}{1}\km\W{1}{3}\km\W{2}{3}\km\W{0}{4}\km\W{2}{4}\km\W{5}{5}
-6 \W{0}{1}\km\W{1}{3}\km\W{2}{3}\km\W{1}{4}\km\W{2}{4}\km\W{4}{5}
+6 \W{0}{1}\km\W{1}{3}\km\W{2}{3}\km\W{2}{4}\km\W{3}{4}\km\W{2}{5}
-\frac{3}{2} \W{0}{1}\km\W{1}{3}\km\W{2}{3}\km\W{2}{4}\km\W{4}{4}\km\W{1}{5}
-\W{0}{1}\km\W{1}{3}\km\W{3}{3}\km\W{0}{4}\km\W{1}{4}\km\W{5}{5}
+6 \W{0}{1}\km\W{1}{3}\km\W{3}{3}\km\W{1}{4}\km\W{2}{4}\km\W{3}{5}
\\&
-4 \W{0}{1}\km\W{1}{3}\km\W{3}{3}\km\W{1}{4}\km\W{3}{4}\km\W{2}{5}
+\W{0}{1}\km\W{1}{3}\km\W{3}{3}\km\W{1}{4}\km\W{4}{4}\km\W{1}{5}
+2 \W{0}{1}\km\W{2}{3}\km\W{3}{3}\km\W{0}{4}\km\W{1}{4}\km\W{4}{5}
-3 \W{0}{1}\km\W{2}{3}\km\W{3}{3}\km\W{0}{4}\km\W{2}{4}\km\W{3}{5}
+2 \W{0}{1}\km\W{2}{3}\km\W{3}{3}\km\W{0}{4}\km\W{3}{4}\km\W{2}{5}
-\frac{1}{2} \W{0}{1}\km\W{2}{3}\km\W{3}{3}\km\W{0}{4}\km\W{4}{4}\km\W{1}{5}
\\&
-\frac{1}{2} \W{1}{1}\km\W{0}{3}\km\W{1}{3}\km\W{0}{4}\km\W{4}{4}\km\W{4}{5}
+2 \W{1}{1}\km\W{0}{3}\km\W{1}{3}\km\W{1}{4}\km\W{4}{4}\km\W{3}{5}
-3 \W{1}{1}\km\W{0}{3}\km\W{1}{3}\km\W{2}{4}\km\W{4}{4}\km\W{2}{5}
+2 \W{1}{1}\km\W{0}{3}\km\W{1}{3}\km\W{3}{4}\km\W{4}{4}\km\W{1}{5}
+\W{1}{1}\km\W{0}{3}\km\W{2}{3}\km\W{0}{4}\km\W{3}{4}\km\W{4}{5}
-4 \W{1}{1}\km\W{0}{3}\km\W{2}{3}\km\W{1}{4}\km\W{3}{4}\km\W{3}{5}
\\&
+6 \W{1}{1}\km\W{0}{3}\km\W{2}{3}\km\W{2}{4}\km\W{3}{4}\km\W{2}{5}
-\W{1}{1}\km\W{0}{3}\km\W{2}{3}\km\W{3}{4}\km\W{4}{4}\km\W{0}{5}
-\frac{1}{2} \W{1}{1}\km\W{0}{3}\km\W{3}{3}\km\W{0}{4}\km\W{2}{4}\km\W{4}{5}
+2 \W{1}{1}\km\W{0}{3}\km\W{3}{3}\km\W{1}{4}\km\W{2}{4}\km\W{3}{5}
-2 \W{1}{1}\km\W{0}{3}\km\W{3}{3}\km\W{2}{4}\km\W{3}{4}\km\W{1}{5}
+\frac{1}{2} \W{1}{1}\km\W{0}{3}\km\W{3}{3}\km\W{2}{4}\km\W{4}{4}\km\W{0}{5}
\\&
-\frac{3}{2} \W{1}{1}\km\W{1}{3}\km\W{2}{3}\km\W{0}{4}\km\W{2}{4}\km\W{4}{5}
+6 \W{1}{1}\km\W{1}{3}\km\W{2}{3}\km\W{1}{4}\km\W{2}{4}\km\W{3}{5}
-6 \W{1}{1}\km\W{1}{3}\km\W{2}{3}\km\W{2}{4}\km\W{3}{4}\km\W{1}{5}
+\frac{3}{2} \W{1}{1}\km\W{1}{3}\km\W{2}{3}\km\W{2}{4}\km\W{4}{4}\km\W{0}{5}
+\W{1}{1}\km\W{1}{3}\km\W{3}{3}\km\W{0}{4}\km\W{1}{4}\km\W{4}{5}
-6 \W{1}{1}\km\W{1}{3}\km\W{3}{3}\km\W{1}{4}\km\W{2}{4}\km\W{2}{5}
\\&
+4 \W{1}{1}\km\W{1}{3}\km\W{3}{3}\km\W{1}{4}\km\W{3}{4}\km\W{1}{5}
-\W{1}{1}\km\W{1}{3}\km\W{3}{3}\km\W{1}{4}\km\W{4}{4}\km\W{0}{5}
-2 \W{1}{1}\km\W{2}{3}\km\W{3}{3}\km\W{0}{4}\km\W{1}{4}\km\W{3}{5}
+3 \W{1}{1}\km\W{2}{3}\km\W{3}{3}\km\W{0}{4}\km\W{2}{4}\km\W{2}{5}
-2 \W{1}{1}\km\W{2}{3}\km\W{3}{3}\km\W{0}{4}\km\W{3}{4}\km\W{1}{5}
+\frac{1}{2} \W{1}{1}\km\W{2}{3}\km\W{3}{3}\km\W{0}{4}\km\W{4}{4}\km\W{0}{5}
\\
\widehat{\mathbf{q}}_{14}  &= 
-\frac{1}{6} \W{0}{1}\km\W{0}{3}\km\W{1}{3}\km\W{0}{4}\km\W{4}{4}\km\W{5}{5}
+\frac{2}{3} \W{0}{1}\km\W{0}{3}\km\W{1}{3}\km\W{1}{4}\km\W{3}{4}\km\W{5}{5}
-2 \W{0}{1}\km\W{0}{3}\km\W{1}{3}\km\W{2}{4}\km\W{3}{4}\km\W{4}{5}
+\W{0}{1}\km\W{0}{3}\km\W{1}{3}\km\W{2}{4}\km\W{4}{4}\km\W{3}{5}
-\frac{2}{3} \W{0}{1}\km\W{0}{3}\km\W{1}{3}\km\W{3}{4}\km\W{4}{4}\km\W{2}{5}
+\frac{1}{3} \W{0}{1}\km\W{0}{3}\km\W{2}{3}\km\W{0}{4}\km\W{4}{4}\km\W{4}{5}
\\&
-\W{0}{1}\km\W{0}{3}\km\W{2}{3}\km\W{1}{4}\km\W{2}{4}\km\W{5}{5}
+\frac{2}{3} \W{0}{1}\km\W{0}{3}\km\W{2}{3}\km\W{1}{4}\km\W{3}{4}\km\W{4}{5}
-\W{0}{1}\km\W{0}{3}\km\W{2}{3}\km\W{1}{4}\km\W{4}{4}\km\W{3}{5}
+\W{0}{1}\km\W{0}{3}\km\W{2}{3}\km\W{2}{4}\km\W{3}{4}\km\W{3}{5}
+\frac{1}{3} \W{0}{1}\km\W{0}{3}\km\W{2}{3}\km\W{3}{4}\km\W{4}{4}\km\W{1}{5}
+\frac{1}{6} \W{0}{1}\km\W{0}{3}\km\W{3}{3}\km\W{0}{4}\km\W{2}{4}\km\W{5}{5}
\\&
-\frac{1}{3} \W{0}{1}\km\W{0}{3}\km\W{3}{3}\km\W{0}{4}\km\W{3}{4}\km\W{4}{5}
+\frac{1}{3} \W{0}{1}\km\W{0}{3}\km\W{3}{3}\km\W{1}{4}\km\W{2}{4}\km\W{4}{5}
+\frac{1}{3} \W{0}{1}\km\W{0}{3}\km\W{3}{3}\km\W{1}{4}\km\W{4}{4}\km\W{2}{5}
-\frac{1}{6} \W{0}{1}\km\W{0}{3}\km\W{3}{3}\km\W{2}{4}\km\W{4}{4}\km\W{1}{5}
-\frac{1}{3} \W{0}{1}\km\W{0}{3}\km\W{3}{3}\km\W{2}{4}\km\W{3}{4}\km\W{2}{5}
-\W{0}{1}\km\W{1}{3}\km\W{2}{3}\km\W{0}{4}\km\W{3}{4}\km\W{4}{5}
\\&
+\frac{1}{2} \W{0}{1}\km\W{1}{3}\km\W{2}{3}\km\W{0}{4}\km\W{2}{4}\km\W{5}{5}
+\W{0}{1}\km\W{1}{3}\km\W{2}{3}\km\W{1}{4}\km\W{2}{4}\km\W{4}{5}
+\W{0}{1}\km\W{1}{3}\km\W{2}{3}\km\W{1}{4}\km\W{4}{4}\km\W{2}{5}
-\W{0}{1}\km\W{1}{3}\km\W{2}{3}\km\W{2}{4}\km\W{3}{4}\km\W{2}{5}
-\frac{1}{2} \W{0}{1}\km\W{1}{3}\km\W{2}{3}\km\W{2}{4}\km\W{4}{4}\km\W{1}{5}
-\frac{1}{3} \W{0}{1}\km\W{1}{3}\km\W{3}{3}\km\W{0}{4}\km\W{1}{4}\km\W{5}{5}
\\&
+\W{0}{1}\km\W{1}{3}\km\W{3}{3}\km\W{0}{4}\km\W{3}{4}\km\W{3}{5}
-\frac{1}{3} \W{0}{1}\km\W{1}{3}\km\W{3}{3}\km\W{0}{4}\km\W{4}{4}\km\W{2}{5}
-\W{0}{1}\km\W{1}{3}\km\W{3}{3}\km\W{1}{4}\km\W{2}{4}\km\W{3}{5}
-\frac{2}{3} \W{0}{1}\km\W{1}{3}\km\W{3}{3}\km\W{1}{4}\km\W{3}{4}\km\W{2}{5}
+\W{0}{1}\km\W{1}{3}\km\W{3}{3}\km\W{2}{4}\km\W{3}{4}\km\W{1}{5}
+\frac{2}{3} \W{0}{1}\km\W{2}{3}\km\W{3}{3}\km\W{0}{4}\km\W{1}{4}\km\W{4}{5}
\\&
-\W{0}{1}\km\W{2}{3}\km\W{3}{3}\km\W{0}{4}\km\W{2}{4}\km\W{3}{5}
+\frac{1}{6} \W{0}{1}\km\W{2}{3}\km\W{3}{3}\km\W{0}{4}\km\W{4}{4}\km\W{1}{5}
+2 \W{0}{1}\km\W{2}{3}\km\W{3}{3}\km\W{1}{4}\km\W{2}{4}\km\W{2}{5}
-\frac{2}{3} \W{0}{1}\km\W{2}{3}\km\W{3}{3}\km\W{1}{4}\km\W{3}{4}\km\W{1}{5}
+\frac{1}{6} \W{1}{1}\km\W{0}{3}\km\W{1}{3}\km\W{0}{4}\km\W{4}{4}\km\W{4}{5}
-\frac{2}{3} \W{1}{1}\km\W{0}{3}\km\W{1}{3}\km\W{1}{4}\km\W{3}{4}\km\W{4}{5}
\\&
-\W{1}{1}\km\W{0}{3}\km\W{1}{3}\km\W{2}{4}\km\W{4}{4}\km\W{2}{5}
+2 \W{1}{1}\km\W{0}{3}\km\W{1}{3}\km\W{2}{4}\km\W{3}{4}\km\W{3}{5}
+\frac{2}{3} \W{1}{1}\km\W{0}{3}\km\W{1}{3}\km\W{3}{4}\km\W{4}{4}\km\W{1}{5}
-\frac{1}{3} \W{1}{1}\km\W{0}{3}\km\W{2}{3}\km\W{0}{4}\km\W{4}{4}\km\W{3}{5}
+\W{1}{1}\km\W{0}{3}\km\W{2}{3}\km\W{1}{4}\km\W{2}{4}\km\W{4}{5}
-\frac{2}{3} \W{1}{1}\km\W{0}{3}\km\W{2}{3}\km\W{1}{4}\km\W{3}{4}\km\W{3}{5}
\\&
+\W{1}{1}\km\W{0}{3}\km\W{2}{3}\km\W{1}{4}\km\W{4}{4}\km\W{2}{5}
-\W{1}{1}\km\W{0}{3}\km\W{2}{3}\km\W{2}{4}\km\W{3}{4}\km\W{2}{5}
-\frac{1}{3} \W{1}{1}\km\W{0}{3}\km\W{2}{3}\km\W{3}{4}\km\W{4}{4}\km\W{0}{5}
-\frac{1}{6} \W{1}{1}\km\W{0}{3}\km\W{3}{3}\km\W{0}{4}\km\W{2}{4}\km\W{4}{5}
+\frac{1}{3} \W{1}{1}\km\W{0}{3}\km\W{3}{3}\km\W{0}{4}\km\W{3}{4}\km\W{3}{5}
-\frac{1}{3} \W{1}{1}\km\W{0}{3}\km\W{3}{3}\km\W{1}{4}\km\W{2}{4}\km\W{3}{5}
\\&
-\frac{1}{3} \W{1}{1}\km\W{0}{3}\km\W{3}{3}\km\W{1}{4}\km\W{4}{4}\km\W{1}{5}
+\frac{1}{3} \W{1}{1}\km\W{0}{3}\km\W{3}{3}\km\W{2}{4}\km\W{3}{4}\km\W{1}{5}
+\frac{1}{6} \W{1}{1}\km\W{0}{3}\km\W{3}{3}\km\W{2}{4}\km\W{4}{4}\km\W{0}{5}
-\frac{1}{2} \W{1}{1}\km\W{1}{3}\km\W{2}{3}\km\W{0}{4}\km\W{2}{4}\km\W{4}{5}
+\W{1}{1}\km\W{1}{3}\km\W{2}{3}\km\W{0}{4}\km\W{3}{4}\km\W{3}{5}
-\W{1}{1}\km\W{1}{3}\km\W{2}{3}\km\W{1}{4}\km\W{2}{4}\km\W{3}{5}
\\&
-\W{1}{1}\km\W{1}{3}\km\W{2}{3}\km\W{1}{4}\km\W{4}{4}\km\W{1}{5}
+\W{1}{1}\km\W{1}{3}\km\W{2}{3}\km\W{2}{4}\km\W{3}{4}\km\W{1}{5}
+\frac{1}{2} \W{1}{1}\km\W{1}{3}\km\W{2}{3}\km\W{2}{4}\km\W{4}{4}\km\W{0}{5}
-\W{1}{1}\km\W{1}{3}\km\W{3}{3}\km\W{0}{4}\km\W{3}{4}\km\W{2}{5}
+\frac{1}{3} \W{1}{1}\km\W{1}{3}\km\W{3}{3}\km\W{0}{4}\km\W{1}{4}\km\W{4}{5}
+\frac{1}{3} \W{1}{1}\km\W{1}{3}\km\W{3}{3}\km\W{0}{4}\km\W{4}{4}\km\W{1}{5}
\\&
+\W{1}{1}\km\W{1}{3}\km\W{3}{3}\km\W{1}{4}\km\W{2}{4}\km\W{2}{5}
-\W{1}{1}\km\W{1}{3}\km\W{3}{3}\km\W{2}{4}\km\W{3}{4}\km\W{0}{5}
+\frac{2}{3} \W{1}{1}\km\W{1}{3}\km\W{3}{3}\km\W{1}{4}\km\W{3}{4}\km\W{1}{5}
-\frac{2}{3} \W{1}{1}\km\W{2}{3}\km\W{3}{3}\km\W{0}{4}\km\W{1}{4}\km\W{3}{5}
+\W{1}{1}\km\W{2}{3}\km\W{3}{3}\km\W{0}{4}\km\W{2}{4}\km\W{2}{5}
-\frac{1}{6} \W{1}{1}\km\W{2}{3}\km\W{3}{3}\km\W{0}{4}\km\W{4}{4}\km\W{0}{5}
\\&
-2 \W{1}{1}\km\W{2}{3}\km\W{3}{3}\km\W{1}{4}\km\W{2}{4}\km\W{1}{5}
+\frac{2}{3} \W{1}{1}\km\W{2}{3}\km\W{3}{3}\km\W{1}{4}\km\W{3}{4}\km\W{0}{5}
\\
\widehat{\mathbf{q}}_{15}  &= 
-\frac{7}{6} \W{0}{1}\km\W{0}{3}\km\W{1}{3}\km\W{0}{4}\km\W{4}{4}\km\W{5}{5}
-\frac{4}{3} \W{0}{1}\km\W{0}{3}\km\W{1}{3}\km\W{1}{4}\km\W{3}{4}\km\W{5}{5}
+6 \W{0}{1}\km\W{0}{3}\km\W{1}{3}\km\W{1}{4}\km\W{4}{4}\km\W{4}{5}
+4 \W{0}{1}\km\W{0}{3}\km\W{1}{3}\km\W{2}{4}\km\W{3}{4}\km\W{4}{5}
-11 \W{0}{1}\km\W{0}{3}\km\W{1}{3}\km\W{2}{4}\km\W{4}{4}\km\W{3}{5}
+\frac{22}{3} \W{0}{1}\km\W{0}{3}\km\W{1}{3}\km\W{3}{4}\km\W{4}{4}\km\W{2}{5}
\\&
+\frac{10}{3} \W{0}{1}\km\W{0}{3}\km\W{2}{3}\km\W{0}{4}\km\W{3}{4}\km\W{5}{5}
-\W{0}{1}\km\W{0}{3}\km\W{2}{3}\km\W{0}{4}\km\W{4}{4}\km\W{4}{5}
-\frac{32}{3} \W{0}{1}\km\W{0}{3}\km\W{2}{3}\km\W{1}{4}\km\W{3}{4}\km\W{4}{5}
+\frac{4}{3} \W{0}{1}\km\W{0}{3}\km\W{2}{3}\km\W{1}{4}\km\W{4}{4}\km\W{3}{5}
+12 \W{0}{1}\km\W{0}{3}\km\W{2}{3}\km\W{2}{4}\km\W{3}{4}\km\W{3}{5}
+2 \W{0}{1}\km\W{0}{3}\km\W{2}{3}\km\W{2}{4}\km\W{4}{4}\km\W{2}{5}
\\&
-\frac{14}{3} \W{0}{1}\km\W{0}{3}\km\W{2}{3}\km\W{3}{4}\km\W{4}{4}\km\W{1}{5}
-\frac{3}{2} \W{0}{1}\km\W{0}{3}\km\W{3}{3}\km\W{0}{4}\km\W{2}{4}\km\W{5}{5}
-\frac{1}{3} \W{0}{1}\km\W{0}{3}\km\W{3}{3}\km\W{0}{4}\km\W{3}{4}\km\W{4}{5}
+\frac{2}{3} \W{0}{1}\km\W{0}{3}\km\W{3}{3}\km\W{0}{4}\km\W{4}{4}\km\W{3}{5}
+6 \W{0}{1}\km\W{0}{3}\km\W{3}{3}\km\W{1}{4}\km\W{2}{4}\km\W{4}{5}
-\frac{4}{3} \W{0}{1}\km\W{0}{3}\km\W{3}{3}\km\W{1}{4}\km\W{4}{4}\km\W{2}{5}
\\&
+\frac{3}{2} \W{0}{1}\km\W{0}{3}\km\W{3}{3}\km\W{2}{4}\km\W{4}{4}\km\W{1}{5}
-4 \W{0}{1}\km\W{0}{3}\km\W{3}{3}\km\W{2}{4}\km\W{3}{4}\km\W{2}{5}
+\frac{1}{3} \W{0}{1}\km\W{0}{3}\km\W{3}{3}\km\W{3}{4}\km\W{4}{4}\km\W{0}{5}
-\W{0}{1}\km\W{1}{3}\km\W{2}{3}\km\W{0}{4}\km\W{3}{4}\km\W{4}{5}
+2 \W{0}{1}\km\W{1}{3}\km\W{2}{3}\km\W{0}{4}\km\W{4}{4}\km\W{3}{5}
-\frac{9}{2} \W{0}{1}\km\W{1}{3}\km\W{2}{3}\km\W{0}{4}\km\W{2}{4}\km\W{5}{5}
\\&
+18 \W{0}{1}\km\W{1}{3}\km\W{2}{3}\km\W{1}{4}\km\W{2}{4}\km\W{4}{5}
-4 \W{0}{1}\km\W{1}{3}\km\W{2}{3}\km\W{1}{4}\km\W{4}{4}\km\W{2}{5}
-12 \W{0}{1}\km\W{1}{3}\km\W{2}{3}\km\W{2}{4}\km\W{3}{4}\km\W{2}{5}
+\frac{9}{2} \W{0}{1}\km\W{1}{3}\km\W{2}{3}\km\W{2}{4}\km\W{4}{4}\km\W{1}{5}
+\W{0}{1}\km\W{1}{3}\km\W{2}{3}\km\W{3}{4}\km\W{4}{4}\km\W{0}{5}
+\frac{8}{3} \W{0}{1}\km\W{1}{3}\km\W{3}{3}\km\W{0}{4}\km\W{1}{4}\km\W{5}{5}
\\&
+\W{0}{1}\km\W{1}{3}\km\W{3}{3}\km\W{0}{4}\km\W{2}{4}\km\W{4}{5}
-\frac{4}{3} \W{0}{1}\km\W{1}{3}\km\W{3}{3}\km\W{0}{4}\km\W{4}{4}\km\W{2}{5}
-20 \W{0}{1}\km\W{1}{3}\km\W{3}{3}\km\W{1}{4}\km\W{2}{4}\km\W{3}{5}
+\frac{40}{3} \W{0}{1}\km\W{1}{3}\km\W{3}{3}\km\W{1}{4}\km\W{3}{4}\km\W{2}{5}
-4 \W{0}{1}\km\W{1}{3}\km\W{3}{3}\km\W{2}{4}\km\W{3}{4}\km\W{1}{5}
-\W{0}{1}\km\W{1}{3}\km\W{3}{3}\km\W{2}{4}\km\W{4}{4}\km\W{0}{5}
\\&
-\frac{16}{3} \W{0}{1}\km\W{2}{3}\km\W{3}{3}\km\W{0}{4}\km\W{1}{4}\km\W{4}{5}
+7 \W{0}{1}\km\W{2}{3}\km\W{3}{3}\km\W{0}{4}\km\W{2}{4}\km\W{3}{5}
-\frac{14}{3} \W{0}{1}\km\W{2}{3}\km\W{3}{3}\km\W{0}{4}\km\W{3}{4}\km\W{2}{5}
+\frac{11}{6} \W{0}{1}\km\W{2}{3}\km\W{3}{3}\km\W{0}{4}\km\W{4}{4}\km\W{1}{5}
+4 \W{0}{1}\km\W{2}{3}\km\W{3}{3}\km\W{1}{4}\km\W{2}{4}\km\W{2}{5}
+4 \W{0}{1}\km\W{2}{3}\km\W{3}{3}\km\W{2}{4}\km\W{3}{4}\km\W{0}{5}
\\&
-4 \W{0}{1}\km\W{2}{3}\km\W{3}{3}\km\W{1}{4}\km\W{3}{4}\km\W{1}{5}
-\frac{2}{3} \W{0}{1}\km\W{2}{3}\km\W{3}{3}\km\W{1}{4}\km\W{4}{4}\km\W{0}{5}
-\frac{2}{3} \W{1}{1}\km\W{0}{3}\km\W{1}{3}\km\W{0}{4}\km\W{3}{4}\km\W{5}{5}
+\frac{11}{6} \W{1}{1}\km\W{0}{3}\km\W{1}{3}\km\W{0}{4}\km\W{4}{4}\km\W{4}{5}
+4 \W{1}{1}\km\W{0}{3}\km\W{1}{3}\km\W{1}{4}\km\W{2}{4}\km\W{5}{5}
-4 \W{1}{1}\km\W{0}{3}\km\W{1}{3}\km\W{1}{4}\km\W{3}{4}\km\W{4}{5}
\\&
-\frac{14}{3} \W{1}{1}\km\W{0}{3}\km\W{1}{3}\km\W{1}{4}\km\W{4}{4}\km\W{3}{5}
+7 \W{1}{1}\km\W{0}{3}\km\W{1}{3}\km\W{2}{4}\km\W{4}{4}\km\W{2}{5}
+4 \W{1}{1}\km\W{0}{3}\km\W{1}{3}\km\W{2}{4}\km\W{3}{4}\km\W{3}{5}
-\frac{16}{3} \W{1}{1}\km\W{0}{3}\km\W{1}{3}\km\W{3}{4}\km\W{4}{4}\km\W{1}{5}
-\W{1}{1}\km\W{0}{3}\km\W{2}{3}\km\W{0}{4}\km\W{2}{4}\km\W{5}{5}
-\frac{4}{3} \W{1}{1}\km\W{0}{3}\km\W{2}{3}\km\W{0}{4}\km\W{4}{4}\km\W{3}{5}
\\&
-4 \W{1}{1}\km\W{0}{3}\km\W{2}{3}\km\W{1}{4}\km\W{2}{4}\km\W{4}{5}
+\frac{40}{3} \W{1}{1}\km\W{0}{3}\km\W{2}{3}\km\W{1}{4}\km\W{3}{4}\km\W{3}{5}
-20 \W{1}{1}\km\W{0}{3}\km\W{2}{3}\km\W{2}{4}\km\W{3}{4}\km\W{2}{5}
+\W{1}{1}\km\W{0}{3}\km\W{2}{3}\km\W{2}{4}\km\W{4}{4}\km\W{1}{5}
+\frac{8}{3} \W{1}{1}\km\W{0}{3}\km\W{2}{3}\km\W{3}{4}\km\W{4}{4}\km\W{0}{5}
+\frac{1}{3} \W{1}{1}\km\W{0}{3}\km\W{3}{3}\km\W{0}{4}\km\W{1}{4}\km\W{5}{5}
\\&
+\frac{3}{2} \W{1}{1}\km\W{0}{3}\km\W{3}{3}\km\W{0}{4}\km\W{2}{4}\km\W{4}{5}
-\frac{4}{3} \W{1}{1}\km\W{0}{3}\km\W{3}{3}\km\W{0}{4}\km\W{3}{4}\km\W{3}{5}
+\frac{2}{3} \W{1}{1}\km\W{0}{3}\km\W{3}{3}\km\W{0}{4}\km\W{4}{4}\km\W{2}{5}
-4 \W{1}{1}\km\W{0}{3}\km\W{3}{3}\km\W{1}{4}\km\W{2}{4}\km\W{3}{5}
-\frac{1}{3} \W{1}{1}\km\W{0}{3}\km\W{3}{3}\km\W{1}{4}\km\W{4}{4}\km\W{1}{5}
+6 \W{1}{1}\km\W{0}{3}\km\W{3}{3}\km\W{2}{4}\km\W{3}{4}\km\W{1}{5}
\\&
-\frac{3}{2} \W{1}{1}\km\W{0}{3}\km\W{3}{3}\km\W{2}{4}\km\W{4}{4}\km\W{0}{5}
+\W{1}{1}\km\W{1}{3}\km\W{2}{3}\km\W{0}{4}\km\W{1}{4}\km\W{5}{5}
+\frac{9}{2} \W{1}{1}\km\W{1}{3}\km\W{2}{3}\km\W{0}{4}\km\W{2}{4}\km\W{4}{5}
-4 \W{1}{1}\km\W{1}{3}\km\W{2}{3}\km\W{0}{4}\km\W{3}{4}\km\W{3}{5}
+2 \W{1}{1}\km\W{1}{3}\km\W{2}{3}\km\W{0}{4}\km\W{4}{4}\km\W{2}{5}
-12 \W{1}{1}\km\W{1}{3}\km\W{2}{3}\km\W{1}{4}\km\W{2}{4}\km\W{3}{5}
\\&
-\W{1}{1}\km\W{1}{3}\km\W{2}{3}\km\W{1}{4}\km\W{4}{4}\km\W{1}{5}
+18 \W{1}{1}\km\W{1}{3}\km\W{2}{3}\km\W{2}{4}\km\W{3}{4}\km\W{1}{5}
-\frac{9}{2} \W{1}{1}\km\W{1}{3}\km\W{2}{3}\km\W{2}{4}\km\W{4}{4}\km\W{0}{5}
+\frac{4}{3} \W{1}{1}\km\W{1}{3}\km\W{3}{3}\km\W{0}{4}\km\W{3}{4}\km\W{2}{5}
-\frac{14}{3} \W{1}{1}\km\W{1}{3}\km\W{3}{3}\km\W{0}{4}\km\W{1}{4}\km\W{4}{5}
+2 \W{1}{1}\km\W{1}{3}\km\W{3}{3}\km\W{0}{4}\km\W{2}{4}\km\W{3}{5}
\\&
-\W{1}{1}\km\W{1}{3}\km\W{3}{3}\km\W{0}{4}\km\W{4}{4}\km\W{1}{5}
+12 \W{1}{1}\km\W{1}{3}\km\W{3}{3}\km\W{1}{4}\km\W{2}{4}\km\W{2}{5}
-\frac{32}{3} \W{1}{1}\km\W{1}{3}\km\W{3}{3}\km\W{1}{4}\km\W{3}{4}\km\W{1}{5}
+\frac{10}{3} \W{1}{1}\km\W{1}{3}\km\W{3}{3}\km\W{1}{4}\km\W{4}{4}\km\W{0}{5}
+\frac{22}{3} \W{1}{1}\km\W{2}{3}\km\W{3}{3}\km\W{0}{4}\km\W{1}{4}\km\W{3}{5}
-11 \W{1}{1}\km\W{2}{3}\km\W{3}{3}\km\W{0}{4}\km\W{2}{4}\km\W{2}{5}
\\&
+6 \W{1}{1}\km\W{2}{3}\km\W{3}{3}\km\W{0}{4}\km\W{3}{4}\km\W{1}{5}
-\frac{7}{6} \W{1}{1}\km\W{2}{3}\km\W{3}{3}\km\W{0}{4}\km\W{4}{4}\km\W{0}{5}
+4 \W{1}{1}\km\W{2}{3}\km\W{3}{3}\km\W{1}{4}\km\W{2}{4}\km\W{1}{5}
-\frac{4}{3} \W{1}{1}\km\W{2}{3}\km\W{3}{3}\km\W{1}{4}\km\W{3}{4}\km\W{0}{5}
\\
\widehat{\mathbf{q}}_{16}& =
+\frac{1}{4} \W{0}{1}\km\W{0}{3}\km\W{1}{3}\km\W{0}{4}\km\W{4}{4}\km\W{5}{5}
+\frac{1}{2} \W{0}{1}\km\W{0}{3}\km\W{1}{3}\km\W{1}{4}\km\W{3}{4}\km\W{5}{5}
-\frac{3}{2} \W{0}{1}\km\W{0}{3}\km\W{1}{3}\km\W{1}{4}\km\W{4}{4}\km\W{4}{5}
-\frac{3}{2} \W{0}{1}\km\W{0}{3}\km\W{1}{3}\km\W{2}{4}\km\W{3}{4}\km\W{4}{5}
+3 \W{0}{1}\km\W{0}{3}\km\W{1}{3}\km\W{2}{4}\km\W{4}{4}\km\W{3}{5}
-2 \W{0}{1}\km\W{0}{3}\km\W{1}{3}\km\W{3}{4}\km\W{4}{4}\km\W{2}{5}
\\&
-\W{0}{1}\km\W{0}{3}\km\W{2}{3}\km\W{0}{4}\km\W{3}{4}\km\W{5}{5}
+\frac{1}{2} \W{0}{1}\km\W{0}{3}\km\W{2}{3}\km\W{0}{4}\km\W{4}{4}\km\W{4}{5}
+3 \W{0}{1}\km\W{0}{3}\km\W{2}{3}\km\W{1}{4}\km\W{3}{4}\km\W{4}{5}
-\W{0}{1}\km\W{0}{3}\km\W{2}{3}\km\W{1}{4}\km\W{4}{4}\km\W{3}{5}
-3 \W{0}{1}\km\W{0}{3}\km\W{2}{3}\km\W{2}{4}\km\W{3}{4}\km\W{3}{5}
+\W{0}{1}\km\W{0}{3}\km\W{2}{3}\km\W{3}{4}\km\W{4}{4}\km\W{1}{5}
\\&
+\frac{9}{16} \W{0}{1}\km\W{0}{3}\km\W{3}{3}\km\W{0}{4}\km\W{2}{4}\km\W{5}{5}
-\frac{1}{8} \W{0}{1}\km\W{0}{3}\km\W{3}{3}\km\W{0}{4}\km\W{3}{4}\km\W{4}{5}
-\frac{3}{16} \W{0}{1}\km\W{0}{3}\km\W{3}{3}\km\W{0}{4}\km\W{4}{4}\km\W{3}{5}
-\frac{9}{4} \W{0}{1}\km\W{0}{3}\km\W{3}{3}\km\W{1}{4}\km\W{2}{4}\km\W{4}{5}
+\W{0}{1}\km\W{0}{3}\km\W{3}{3}\km\W{1}{4}\km\W{3}{4}\km\W{3}{5}
+\frac{1}{4} \W{0}{1}\km\W{0}{3}\km\W{3}{3}\km\W{1}{4}\km\W{4}{4}\km\W{2}{5}
\\&
-\frac{3}{16} \W{0}{1}\km\W{0}{3}\km\W{3}{3}\km\W{2}{4}\km\W{4}{4}\km\W{1}{5}
-\frac{1}{8} \W{0}{1}\km\W{0}{3}\km\W{3}{3}\km\W{3}{4}\km\W{4}{4}\km\W{0}{5}
+\frac{3}{8} \W{0}{1}\km\W{1}{3}\km\W{2}{3}\km\W{0}{4}\km\W{3}{4}\km\W{4}{5}
-\frac{15}{16} \W{0}{1}\km\W{1}{3}\km\W{2}{3}\km\W{0}{4}\km\W{4}{4}\km\W{3}{5}
+\frac{21}{16} \W{0}{1}\km\W{1}{3}\km\W{2}{3}\km\W{0}{4}\km\W{2}{4}\km\W{5}{5}
-\frac{21}{4} \W{0}{1}\km\W{1}{3}\km\W{2}{3}\km\W{1}{4}\km\W{2}{4}\km\W{4}{5}
\\&
+\frac{9}{4} \W{0}{1}\km\W{1}{3}\km\W{2}{3}\km\W{1}{4}\km\W{4}{4}\km\W{2}{5}
+3 \W{0}{1}\km\W{1}{3}\km\W{2}{3}\km\W{2}{4}\km\W{3}{4}\km\W{2}{5}
-\frac{39}{16} \W{0}{1}\km\W{1}{3}\km\W{2}{3}\km\W{2}{4}\km\W{4}{4}\km\W{1}{5}
+\frac{3}{8} \W{0}{1}\km\W{1}{3}\km\W{2}{3}\km\W{3}{4}\km\W{4}{4}\km\W{0}{5}
-\W{0}{1}\km\W{1}{3}\km\W{3}{3}\km\W{0}{4}\km\W{1}{4}\km\W{5}{5}
+\frac{1}{2} \W{0}{1}\km\W{1}{3}\km\W{3}{3}\km\W{0}{4}\km\W{4}{4}\km\W{2}{5}
\\&
+6 \W{0}{1}\km\W{1}{3}\km\W{3}{3}\km\W{1}{4}\km\W{2}{4}\km\W{3}{5}
-5 \W{0}{1}\km\W{1}{3}\km\W{3}{3}\km\W{1}{4}\km\W{3}{4}\km\W{2}{5}
+3 \W{0}{1}\km\W{1}{3}\km\W{3}{3}\km\W{2}{4}\km\W{3}{4}\km\W{1}{5}
+2 \W{0}{1}\km\W{2}{3}\km\W{3}{3}\km\W{0}{4}\km\W{1}{4}\km\W{4}{5}
-3 \W{0}{1}\km\W{2}{3}\km\W{3}{3}\km\W{0}{4}\km\W{2}{4}\km\W{3}{5}
+2 \W{0}{1}\km\W{2}{3}\km\W{3}{3}\km\W{0}{4}\km\W{3}{4}\km\W{2}{5}
\\&
-\frac{3}{4} \W{0}{1}\km\W{2}{3}\km\W{3}{3}\km\W{0}{4}\km\W{4}{4}\km\W{1}{5}
-\frac{3}{2} \W{0}{1}\km\W{2}{3}\km\W{3}{3}\km\W{2}{4}\km\W{3}{4}\km\W{0}{5}
+\frac{1}{2} \W{0}{1}\km\W{2}{3}\km\W{3}{3}\km\W{1}{4}\km\W{3}{4}\km\W{1}{5}
+\frac{1}{2} \W{0}{1}\km\W{2}{3}\km\W{3}{3}\km\W{1}{4}\km\W{4}{4}\km\W{0}{5}
+\frac{1}{2} \W{1}{1}\km\W{0}{3}\km\W{1}{3}\km\W{0}{4}\km\W{3}{4}\km\W{5}{5}
-\frac{3}{4} \W{1}{1}\km\W{0}{3}\km\W{1}{3}\km\W{0}{4}\km\W{4}{4}\km\W{4}{5}
\\&
-\frac{3}{2} \W{1}{1}\km\W{0}{3}\km\W{1}{3}\km\W{1}{4}\km\W{2}{4}\km\W{5}{5}
+\frac{1}{2} \W{1}{1}\km\W{0}{3}\km\W{1}{3}\km\W{1}{4}\km\W{3}{4}\km\W{4}{5}
+2 \W{1}{1}\km\W{0}{3}\km\W{1}{3}\km\W{1}{4}\km\W{4}{4}\km\W{3}{5}
-3 \W{1}{1}\km\W{0}{3}\km\W{1}{3}\km\W{2}{4}\km\W{4}{4}\km\W{2}{5}
+2 \W{1}{1}\km\W{0}{3}\km\W{1}{3}\km\W{3}{4}\km\W{4}{4}\km\W{1}{5}
+\frac{1}{2} \W{1}{1}\km\W{0}{3}\km\W{2}{3}\km\W{0}{4}\km\W{4}{4}\km\W{3}{5}
\\&
+3 \W{1}{1}\km\W{0}{3}\km\W{2}{3}\km\W{1}{4}\km\W{2}{4}\km\W{4}{5}
-5 \W{1}{1}\km\W{0}{3}\km\W{2}{3}\km\W{1}{4}\km\W{3}{4}\km\W{3}{5}
+6 \W{1}{1}\km\W{0}{3}\km\W{2}{3}\km\W{2}{4}\km\W{3}{4}\km\W{2}{5}
-\W{1}{1}\km\W{0}{3}\km\W{2}{3}\km\W{3}{4}\km\W{4}{4}\km\W{0}{5}
-\frac{1}{8} \W{1}{1}\km\W{0}{3}\km\W{3}{3}\km\W{0}{4}\km\W{1}{4}\km\W{5}{5}
-\frac{3}{16} \W{1}{1}\km\W{0}{3}\km\W{3}{3}\km\W{0}{4}\km\W{2}{4}\km\W{4}{5}
\\&
+\frac{1}{4} \W{1}{1}\km\W{0}{3}\km\W{3}{3}\km\W{0}{4}\km\W{3}{4}\km\W{3}{5}
-\frac{3}{16} \W{1}{1}\km\W{0}{3}\km\W{3}{3}\km\W{0}{4}\km\W{4}{4}\km\W{2}{5}
+\W{1}{1}\km\W{0}{3}\km\W{3}{3}\km\W{1}{4}\km\W{3}{4}\km\W{2}{5}
-\frac{1}{8} \W{1}{1}\km\W{0}{3}\km\W{3}{3}\km\W{1}{4}\km\W{4}{4}\km\W{1}{5}
-\frac{9}{4} \W{1}{1}\km\W{0}{3}\km\W{3}{3}\km\W{2}{4}\km\W{3}{4}\km\W{1}{5}
+\frac{9}{16} \W{1}{1}\km\W{0}{3}\km\W{3}{3}\km\W{2}{4}\km\W{4}{4}\km\W{0}{5}
\\&
+\frac{3}{8} \W{1}{1}\km\W{1}{3}\km\W{2}{3}\km\W{0}{4}\km\W{1}{4}\km\W{5}{5}
-\frac{39}{16} \W{1}{1}\km\W{1}{3}\km\W{2}{3}\km\W{0}{4}\km\W{2}{4}\km\W{4}{5}
+\frac{9}{4} \W{1}{1}\km\W{1}{3}\km\W{2}{3}\km\W{0}{4}\km\W{3}{4}\km\W{3}{5}
-\frac{15}{16} \W{1}{1}\km\W{1}{3}\km\W{2}{3}\km\W{0}{4}\km\W{4}{4}\km\W{2}{5}
+3 \W{1}{1}\km\W{1}{3}\km\W{2}{3}\km\W{1}{4}\km\W{2}{4}\km\W{3}{5}
+\frac{3}{8} \W{1}{1}\km\W{1}{3}\km\W{2}{3}\km\W{1}{4}\km\W{4}{4}\km\W{1}{5}
\\&
-\frac{21}{4} \W{1}{1}\km\W{1}{3}\km\W{2}{3}\km\W{2}{4}\km\W{3}{4}\km\W{1}{5}
+\frac{21}{16} \W{1}{1}\km\W{1}{3}\km\W{2}{3}\km\W{2}{4}\km\W{4}{4}\km\W{0}{5}
-\W{1}{1}\km\W{1}{3}\km\W{3}{3}\km\W{0}{4}\km\W{3}{4}\km\W{2}{5}
+\W{1}{1}\km\W{1}{3}\km\W{3}{3}\km\W{0}{4}\km\W{1}{4}\km\W{4}{5}
+\frac{1}{2} \W{1}{1}\km\W{1}{3}\km\W{3}{3}\km\W{0}{4}\km\W{4}{4}\km\W{1}{5}
-3 \W{1}{1}\km\W{1}{3}\km\W{3}{3}\km\W{1}{4}\km\W{2}{4}\km\W{2}{5}
\\&
+3 \W{1}{1}\km\W{1}{3}\km\W{3}{3}\km\W{1}{4}\km\W{3}{4}\km\W{1}{5}
-\W{1}{1}\km\W{1}{3}\km\W{3}{3}\km\W{1}{4}\km\W{4}{4}\km\W{0}{5}
-2 \W{1}{1}\km\W{2}{3}\km\W{3}{3}\km\W{0}{4}\km\W{1}{4}\km\W{3}{5}
+3 \W{1}{1}\km\W{2}{3}\km\W{3}{3}\km\W{0}{4}\km\W{2}{4}\km\W{2}{5} 
-\frac{3}{2} \W{1}{1}\km\W{2}{3}\km\W{3}{3}\km\W{0}{4}\km\W{3}{4}\km\W{1}{5} 
+\frac{1}{4} \W{1}{1}\km\W{2}{3}\km\W{3}{3}\km\W{0}{4}\km\W{4}{4}\km\W{0}{5}
\\&
-\frac{3}{2} \W{1}{1}\km\W{2}{3}\km\W{3}{3}\km\W{1}{4}\km\W{2}{4}\km\W{1}{5}
+\frac{1}{2} \W{1}{1}\km\W{2}{3}\km\W{3}{3}\km\W{1}{4}\km\W{3}{4}\km\W{0}{5}
\\
\widehat{\mathbf{q}}_{17} &= 
+\frac{1}{12} \W{0}{1}\km\W{0}{3}\km\W{1}{3}\km\W{0}{4}\km\W{4}{4}\km\W{5}{5}
+\frac{1}{6} \W{0}{1}\km\W{0}{3}\km\W{1}{3}\km\W{1}{4}\km\W{3}{4}\km\W{5}{5}
-\frac{1}{2} \W{0}{1}\km\W{0}{3}\km\W{1}{3}\km\W{1}{4}\km\W{4}{4}\km\W{4}{5}
-\frac{1}{2} \W{0}{1}\km\W{0}{3}\km\W{1}{3}\km\W{2}{4}\km\W{3}{4}\km\W{4}{5}
+\W{0}{1}\km\W{0}{3}\km\W{1}{3}\km\W{2}{4}\km\W{4}{4}\km\W{3}{5}
-\frac{2}{3} \W{0}{1}\km\W{0}{3}\km\W{1}{3}\km\W{3}{4}\km\W{4}{4}\km\W{2}{5}
\\&
-\frac{1}{3} \W{0}{1}\km\W{0}{3}\km\W{2}{3}\km\W{0}{4}\km\W{3}{4}\km\W{5}{5}
+\frac{1}{6} \W{0}{1}\km\W{0}{3}\km\W{2}{3}\km\W{0}{4}\km\W{4}{4}\km\W{4}{5}
+\W{0}{1}\km\W{0}{3}\km\W{2}{3}\km\W{1}{4}\km\W{3}{4}\km\W{4}{5}
-\frac{1}{3} \W{0}{1}\km\W{0}{3}\km\W{2}{3}\km\W{1}{4}\km\W{4}{4}\km\W{3}{5}
-\W{0}{1}\km\W{0}{3}\km\W{2}{3}\km\W{2}{4}\km\W{3}{4}\km\W{3}{5}
+\frac{1}{3} \W{0}{1}\km\W{0}{3}\km\W{2}{3}\km\W{3}{4}\km\W{4}{4}\km\W{1}{5}
\\&
+\frac{7}{48} \W{0}{1}\km\W{0}{3}\km\W{3}{3}\km\W{0}{4}\km\W{2}{4}\km\W{5}{5}
+\frac{1}{24} \W{0}{1}\km\W{0}{3}\km\W{3}{3}\km\W{0}{4}\km\W{3}{4}\km\W{4}{5}
-\frac{5}{48} \W{0}{1}\km\W{0}{3}\km\W{3}{3}\km\W{0}{4}\km\W{4}{4}\km\W{3}{5}
-\frac{7}{12} \W{0}{1}\km\W{0}{3}\km\W{3}{3}\km\W{1}{4}\km\W{2}{4}\km\W{4}{5}
+\frac{1}{4} \W{0}{1}\km\W{0}{3}\km\W{3}{3}\km\W{1}{4}\km\W{4}{4}\km\W{2}{5}
-\frac{13}{48} \W{0}{1}\km\W{0}{3}\km\W{3}{3}\km\W{2}{4}\km\W{4}{4}\km\W{1}{5}
\\&
+\frac{1}{3} \W{0}{1}\km\W{0}{3}\km\W{3}{3}\km\W{2}{4}\km\W{3}{4}\km\W{2}{5}
+\frac{1}{24} \W{0}{1}\km\W{0}{3}\km\W{3}{3}\km\W{3}{4}\km\W{4}{4}\km\W{0}{5}
-\frac{1}{8} \W{0}{1}\km\W{1}{3}\km\W{2}{3}\km\W{0}{4}\km\W{3}{4}\km\W{4}{5}
-\frac{3}{16} \W{0}{1}\km\W{1}{3}\km\W{2}{3}\km\W{0}{4}\km\W{4}{4}\km\W{3}{5}
+\frac{9}{16} \W{0}{1}\km\W{1}{3}\km\W{2}{3}\km\W{0}{4}\km\W{2}{4}\km\W{5}{5}
-\frac{9}{4} \W{0}{1}\km\W{1}{3}\km\W{2}{3}\km\W{1}{4}\km\W{2}{4}\km\W{4}{5}
\\&
+\W{0}{1}\km\W{1}{3}\km\W{2}{3}\km\W{1}{4}\km\W{3}{4}\km\W{3}{5}
+\frac{1}{4} \W{0}{1}\km\W{1}{3}\km\W{2}{3}\km\W{1}{4}\km\W{4}{4}\km\W{2}{5} 
-\frac{3}{16} \W{0}{1}\km\W{1}{3}\km\W{2}{3}\km\W{2}{4}\km\W{4}{4}\km\W{1}{5}
-\frac{1}{8} \W{0}{1}\km\W{1}{3}\km\W{2}{3}\km\W{3}{4}\km\W{4}{4}\km\W{0}{5}
-\frac{1}{3} \W{0}{1}\km\W{1}{3}\km\W{3}{3}\km\W{0}{4}\km\W{1}{4}\km\W{5}{5}
+\frac{1}{6} \W{0}{1}\km\W{1}{3}\km\W{3}{3}\km\W{0}{4}\km\W{4}{4}\km\W{2}{5}
\\&
+2 \W{0}{1}\km\W{1}{3}\km\W{3}{3}\km\W{1}{4}\km\W{2}{4}\km\W{3}{5}
-\frac{5}{3} \W{0}{1}\km\W{1}{3}\km\W{3}{3}\km\W{1}{4}\km\W{3}{4}\km\W{2}{5}
+\W{0}{1}\km\W{1}{3}\km\W{3}{3}\km\W{2}{4}\km\W{3}{4}\km\W{1}{5}
+\frac{2}{3} \W{0}{1}\km\W{2}{3}\km\W{3}{3}\km\W{0}{4}\km\W{1}{4}\km\W{4}{5}
-\W{0}{1}\km\W{2}{3}\km\W{3}{3}\km\W{0}{4}\km\W{2}{4}\km\W{3}{5}
+\frac{2}{3} \W{0}{1}\km\W{2}{3}\km\W{3}{3}\km\W{0}{4}\km\W{3}{4}\km\W{2}{5}
\\&
-\frac{1}{4} \W{0}{1}\km\W{2}{3}\km\W{3}{3}\km\W{0}{4}\km\W{4}{4}\km\W{1}{5}
-\frac{1}{2} \W{0}{1}\km\W{2}{3}\km\W{3}{3}\km\W{2}{4}\km\W{3}{4}\km\W{0}{5}
+\frac{1}{6} \W{0}{1}\km\W{2}{3}\km\W{3}{3}\km\W{1}{4}\km\W{3}{4}\km\W{1}{5}
+\frac{1}{6} \W{0}{1}\km\W{2}{3}\km\W{3}{3}\km\W{1}{4}\km\W{4}{4}\km\W{0}{5}
+\frac{1}{6} \W{1}{1}\km\W{0}{3}\km\W{1}{3}\km\W{0}{4}\km\W{3}{4}\km\W{5}{5}
-\frac{1}{4} \W{1}{1}\km\W{0}{3}\km\W{1}{3}\km\W{0}{4}\km\W{4}{4}\km\W{4}{5}
\\&
-\frac{1}{2} \W{1}{1}\km\W{0}{3}\km\W{1}{3}\km\W{1}{4}\km\W{2}{4}\km\W{5}{5}
+\frac{1}{6} \W{1}{1}\km\W{0}{3}\km\W{1}{3}\km\W{1}{4}\km\W{3}{4}\km\W{4}{5}
+\frac{2}{3} \W{1}{1}\km\W{0}{3}\km\W{1}{3}\km\W{1}{4}\km\W{4}{4}\km\W{3}{5}
-\W{1}{1}\km\W{0}{3}\km\W{1}{3}\km\W{2}{4}\km\W{4}{4}\km\W{2}{5}
+\frac{2}{3} \W{1}{1}\km\W{0}{3}\km\W{1}{3}\km\W{3}{4}\km\W{4}{4}\km\W{1}{5}
+\frac{1}{6} \W{1}{1}\km\W{0}{3}\km\W{2}{3}\km\W{0}{4}\km\W{4}{4}\km\W{3}{5}
\\&
+\W{1}{1}\km\W{0}{3}\km\W{2}{3}\km\W{1}{4}\km\W{2}{4}\km\W{4}{5}
-\frac{5}{3} \W{1}{1}\km\W{0}{3}\km\W{2}{3}\km\W{1}{4}\km\W{3}{4}\km\W{3}{5}
+2 \W{1}{1}\km\W{0}{3}\km\W{2}{3}\km\W{2}{4}\km\W{3}{4}\km\W{2}{5}
-\frac{1}{3} \W{1}{1}\km\W{0}{3}\km\W{2}{3}\km\W{3}{4}\km\W{4}{4}\km\W{0}{5}
+\frac{1}{24} \W{1}{1}\km\W{0}{3}\km\W{3}{3}\km\W{0}{4}\km\W{1}{4}\km\W{5}{5}
-\frac{13}{48} \W{1}{1}\km\W{0}{3}\km\W{3}{3}\km\W{0}{4}\km\W{2}{4}\km\W{4}{5}
\\&
+\frac{1}{4} \W{1}{1}\km\W{0}{3}\km\W{3}{3}\km\W{0}{4}\km\W{3}{4}\km\W{3}{5}
-\frac{5}{48} \W{1}{1}\km\W{0}{3}\km\W{3}{3}\km\W{0}{4}\km\W{4}{4}\km\W{2}{5}
+\frac{1}{3} \W{1}{1}\km\W{0}{3}\km\W{3}{3}\km\W{1}{4}\km\W{2}{4}\km\W{3}{5}
+\frac{1}{24} \W{1}{1}\km\W{0}{3}\km\W{3}{3}\km\W{1}{4}\km\W{4}{4}\km\W{1}{5}
-\frac{7}{12} \W{1}{1}\km\W{0}{3}\km\W{3}{3}\km\W{2}{4}\km\W{3}{4}\km\W{1}{5}
+\frac{7}{48} \W{1}{1}\km\W{0}{3}\km\W{3}{3}\km\W{2}{4}\km\W{4}{4}\km\W{0}{5}
\\&
-\frac{1}{8} \W{1}{1}\km\W{1}{3}\km\W{2}{3}\km\W{0}{4}\km\W{1}{4}\km\W{5}{5}
-\frac{3}{16} \W{1}{1}\km\W{1}{3}\km\W{2}{3}\km\W{0}{4}\km\W{2}{4}\km\W{4}{5}
+\frac{1}{4} \W{1}{1}\km\W{1}{3}\km\W{2}{3}\km\W{0}{4}\km\W{3}{4}\km\W{3}{5}
-\frac{3}{16} \W{1}{1}\km\W{1}{3}\km\W{2}{3}\km\W{0}{4}\km\W{4}{4}\km\W{2}{5}
+\W{1}{1}\km\W{1}{3}\km\W{2}{3}\km\W{1}{4}\km\W{3}{4}\km\W{2}{5}
-\frac{1}{8} \W{1}{1}\km\W{1}{3}\km\W{2}{3}\km\W{1}{4}\km\W{4}{4}\km\W{1}{5}
\\&
-\frac{9}{4} \W{1}{1}\km\W{1}{3}\km\W{2}{3}\km\W{2}{4}\km\W{3}{4}\km\W{1}{5}
+\frac{9}{16} \W{1}{1}\km\W{1}{3}\km\W{2}{3}\km\W{2}{4}\km\W{4}{4}\km\W{0}{5}
-\frac{1}{3} \W{1}{1}\km\W{1}{3}\km\W{3}{3}\km\W{0}{4}\km\W{3}{4}\km\W{2}{5}
+\frac{1}{3} \W{1}{1}\km\W{1}{3}\km\W{3}{3}\km\W{0}{4}\km\W{1}{4}\km\W{4}{5}
+\frac{1}{6} \W{1}{1}\km\W{1}{3}\km\W{3}{3}\km\W{0}{4}\km\W{4}{4}\km\W{1}{5}
-\W{1}{1}\km\W{1}{3}\km\W{3}{3}\km\W{1}{4}\km\W{2}{4}\km\W{2}{5}
\\&
+\W{1}{1}\km\W{1}{3}\km\W{3}{3}\km\W{1}{4}\km\W{3}{4}\km\W{1}{5}
-\frac{1}{3} \W{1}{1}\km\W{1}{3}\km\W{3}{3}\km\W{1}{4}\km\W{4}{4}\km\W{0}{5}
-\frac{2}{3} \W{1}{1}\km\W{2}{3}\km\W{3}{3}\km\W{0}{4}\km\W{1}{4}\km\W{3}{5}
+\W{1}{1}\km\W{2}{3}\km\W{3}{3}\km\W{0}{4}\km\W{2}{4}\km\W{2}{5}
-\frac{1}{2} \W{1}{1}\km\W{2}{3}\km\W{3}{3}\km\W{0}{4}\km\W{3}{4}\km\W{1}{5}
+\frac{1}{12} \W{1}{1}\km\W{2}{3}\km\W{3}{3}\km\W{0}{4}\km\W{4}{4}\km\W{0}{5}
\\&
-\frac{1}{2} \W{1}{1}\km\W{2}{3}\km\W{3}{3}\km\W{1}{4}\km\W{2}{4}\km\W{1}{5}
+\frac{1}{6} \W{1}{1}\km\W{2}{3}\km\W{3}{3}\km\W{1}{4}\km\W{3}{4}\km\W{0}{5}
\end{align*}

\subsection{Basis of $\ds\CGF{7}{2}{}{8}$}
\begin{align*}
\widehat{\mathbf{w}}_{1} &= 
\W{0}{1}\km\W{1}{1}\km\W{0}{3}\km\W{2}{3}\km\W{3}{3}\km\W{2}{4}\km\W{3}{7}
+\W{0}{1}\km\W{1}{1}\km\W{0}{3}\km\W{1}{3}\km\W{2}{3}\km\W{2}{4}\km\W{5}{7}
+\frac{1}{6} \W{0}{1}\km\W{1}{1}\km\W{0}{3}\km\W{1}{3}\km\W{2}{3}\km\W{0}{4}\km\W{7}{7}
-\frac{2}{3} \W{0}{1}\km\W{1}{1}\km\W{0}{3}\km\W{1}{3}\km\W{2}{3}\km\W{1}{4}\km\W{6}{7}
-\frac{2}{3} \W{0}{1}\km\W{1}{1}\km\W{0}{3}\km\W{1}{3}\km\W{2}{3}\km\W{3}{4}\km\W{4}{7}
+\frac{1}{6} \W{0}{1}\km\W{1}{1}\km\W{0}{3}\km\W{1}{3}\km\W{2}{3}\km\W{4}{4}\km\W{3}{7}
\\&
-\frac{1}{6} \W{0}{1}\km\W{1}{1}\km\W{0}{3}\km\W{1}{3}\km\W{3}{3}\km\W{0}{4}\km\W{6}{7}
+\frac{2}{3} \W{0}{1}\km\W{1}{1}\km\W{0}{3}\km\W{1}{3}\km\W{3}{3}\km\W{1}{4}\km\W{5}{7}
-\W{0}{1}\km\W{1}{1}\km\W{0}{3}\km\W{1}{3}\km\W{3}{3}\km\W{2}{4}\km\W{4}{7}
+\frac{2}{3} \W{0}{1}\km\W{1}{1}\km\W{0}{3}\km\W{1}{3}\km\W{3}{3}\km\W{3}{4}\km\W{3}{7}
-\frac{1}{6} \W{0}{1}\km\W{1}{1}\km\W{0}{3}\km\W{1}{3}\km\W{3}{3}\km\W{4}{4}\km\W{2}{7}
+\frac{1}{6} \W{0}{1}\km\W{1}{1}\km\W{0}{3}\km\W{2}{3}\km\W{3}{3}\km\W{0}{4}\km\W{5}{7}
\\&
-\frac{2}{3} \W{0}{1}\km\W{1}{1}\km\W{0}{3}\km\W{2}{3}\km\W{3}{3}\km\W{1}{4}\km\W{4}{7}
-\frac{2}{3} \W{0}{1}\km\W{1}{1}\km\W{0}{3}\km\W{2}{3}\km\W{3}{3}\km\W{3}{4}\km\W{2}{7}
+\frac{1}{6} \W{0}{1}\km\W{1}{1}\km\W{0}{3}\km\W{2}{3}\km\W{3}{3}\km\W{4}{4}\km\W{1}{7}
-\frac{1}{6} \W{0}{1}\km\W{1}{1}\km\W{1}{3}\km\W{2}{3}\km\W{3}{3}\km\W{0}{4}\km\W{4}{7}
+\frac{2}{3} \W{0}{1}\km\W{1}{1}\km\W{1}{3}\km\W{2}{3}\km\W{3}{3}\km\W{1}{4}\km\W{3}{7}
-\W{0}{1}\km\W{1}{1}\km\W{1}{3}\km\W{2}{3}\km\W{3}{3}\km\W{2}{4}\km\W{2}{7}
\\&
+\frac{2}{3} \W{0}{1}\km\W{1}{1}\km\W{1}{3}\km\W{2}{3}\km\W{3}{3}\km\W{3}{4}\km\W{1}{7}
-\frac{1}{6} \W{0}{1}\km\W{1}{1}\km\W{1}{3}\km\W{2}{3}\km\W{3}{3}\km\W{4}{4}\km\W{0}{7}
\\ 
\widehat{\mathbf{w}}_{2} &= 
\W{0}{1}\km\W{1}{1}\km\W{0}{3}\km\W{1}{3}\km\W{3}{3}\km\W{3}{5}\km\W{3}{6}
-\frac{3}{8} \W{0}{1}\km\W{1}{1}\km\W{0}{3}\km\W{1}{3}\km\W{2}{3}\km\W{1}{5}\km\W{6}{6}
+\frac{3}{2} \W{0}{1}\km\W{1}{1}\km\W{0}{3}\km\W{1}{3}\km\W{2}{3}\km\W{2}{5}\km\W{5}{6}
-\frac{9}{4} \W{0}{1}\km\W{1}{1}\km\W{0}{3}\km\W{1}{3}\km\W{2}{3}\km\W{3}{5}\km\W{4}{6}
+\frac{3}{2} \W{0}{1}\km\W{1}{1}\km\W{0}{3}\km\W{1}{3}\km\W{2}{3}\km\W{4}{5}\km\W{3}{6}
-\frac{3}{8} \W{0}{1}\km\W{1}{1}\km\W{0}{3}\km\W{1}{3}\km\W{2}{3}\km\W{5}{5}\km\W{2}{6}
\\&
+\frac{1}{8} \W{0}{1}\km\W{1}{1}\km\W{0}{3}\km\W{1}{3}\km\W{3}{3}\km\W{0}{5}\km\W{6}{6}
-\frac{1}{4} \W{0}{1}\km\W{1}{1}\km\W{0}{3}\km\W{1}{3}\km\W{3}{3}\km\W{1}{5}\km\W{5}{6}
-\frac{1}{4} \W{0}{1}\km\W{1}{1}\km\W{0}{3}\km\W{1}{3}\km\W{3}{3}\km\W{2}{5}\km\W{4}{6}
-\frac{7}{8} \W{0}{1}\km\W{1}{1}\km\W{0}{3}\km\W{1}{3}\km\W{3}{3}\km\W{4}{5}\km\W{2}{6}
+\frac{1}{4} \W{0}{1}\km\W{1}{1}\km\W{0}{3}\km\W{1}{3}\km\W{3}{3}\km\W{5}{5}\km\W{1}{6}
-\frac{1}{4} \W{0}{1}\km\W{1}{1}\km\W{0}{3}\km\W{2}{3}\km\W{3}{3}\km\W{0}{5}\km\W{5}{6}
\\&
+\frac{7}{8} \W{0}{1}\km\W{1}{1}\km\W{0}{3}\km\W{2}{3}\km\W{3}{3}\km\W{1}{5}\km\W{4}{6}
-\W{0}{1}\km\W{1}{1}\km\W{0}{3}\km\W{2}{3}\km\W{3}{3}\km\W{2}{5}\km\W{3}{6}
+\frac{1}{4} \W{0}{1}\km\W{1}{1}\km\W{0}{3}\km\W{2}{3}\km\W{3}{3}\km\W{3}{5}\km\W{2}{6}
+\frac{1}{4} \W{0}{1}\km\W{1}{1}\km\W{0}{3}\km\W{2}{3}\km\W{3}{3}\km\W{4}{5}\km\W{1}{6}
-\frac{1}{8} \W{0}{1}\km\W{1}{1}\km\W{0}{3}\km\W{2}{3}\km\W{3}{3}\km\W{5}{5}\km\W{0}{6}
+\frac{3}{8} \W{0}{1}\km\W{1}{1}\km\W{1}{3}\km\W{2}{3}\km\W{3}{3}\km\W{0}{5}\km\W{4}{6}
\\&
-\frac{3}{2} \W{0}{1}\km\W{1}{1}\km\W{1}{3}\km\W{2}{3}\km\W{3}{3}\km\W{1}{5}\km\W{3}{6}
+\frac{9}{4} \W{0}{1}\km\W{1}{1}\km\W{1}{3}\km\W{2}{3}\km\W{3}{3}\km\W{2}{5}\km\W{2}{6}
-\frac{3}{2} \W{0}{1}\km\W{1}{1}\km\W{1}{3}\km\W{2}{3}\km\W{3}{3}\km\W{3}{5}\km\W{1}{6}
+\frac{3}{8} \W{0}{1}\km\W{1}{1}\km\W{1}{3}\km\W{2}{3}\km\W{3}{3}\km\W{4}{5}\km\W{0}{6}
\\
\widehat{\mathbf{w}}_{3} &=
\frac{3}{4} \W{0}{1}\km\W{1}{1}\km\W{0}{3}\km\W{1}{3}\km\W{1}{4}\km\W{2}{4}\km\W{6}{6}
-\frac{3}{2} \W{0}{1}\km\W{1}{1}\km\W{0}{3}\km\W{1}{3}\km\W{1}{4}\km\W{3}{4}\km\W{5}{6}
+\frac{3}{4} \W{0}{1}\km\W{1}{1}\km\W{0}{3}\km\W{1}{3}\km\W{1}{4}\km\W{4}{4}\km\W{4}{6}
-\frac{3}{2} \W{0}{1}\km\W{1}{1}\km\W{0}{3}\km\W{1}{3}\km\W{2}{4}\km\W{4}{4}\km\W{3}{6}
+\frac{3}{4} \W{0}{1}\km\W{1}{1}\km\W{0}{3}\km\W{1}{3}\km\W{3}{4}\km\W{4}{4}\km\W{2}{6}
+\frac{9}{4} \W{0}{1}\km\W{1}{1}\km\W{0}{3}\km\W{1}{3}\km\W{2}{4}\km\W{3}{4}\km\W{4}{6}
\\&
+\frac{3}{4} \W{0}{1}\km\W{1}{1}\km\W{0}{3}\km\W{2}{3}\km\W{0}{4}\km\W{3}{4}\km\W{5}{6}
-\frac{3}{8} \W{0}{1}\km\W{1}{1}\km\W{0}{3}\km\W{2}{3}\km\W{0}{4}\km\W{2}{4}\km\W{6}{6}
-\frac{3}{8} \W{0}{1}\km\W{1}{1}\km\W{0}{3}\km\W{2}{3}\km\W{0}{4}\km\W{4}{4}\km\W{4}{6}
-\frac{3}{2} \W{0}{1}\km\W{1}{1}\km\W{0}{3}\km\W{2}{3}\km\W{2}{4}\km\W{3}{4}\km\W{3}{6}
+\frac{9}{8} \W{0}{1}\km\W{1}{1}\km\W{0}{3}\km\W{2}{3}\km\W{2}{4}\km\W{4}{4}\km\W{2}{6}
-\frac{3}{4} \W{0}{1}\km\W{1}{1}\km\W{0}{3}\km\W{2}{3}\km\W{3}{4}\km\W{4}{4}\km\W{1}{6}
\\&
+\frac{1}{16} \W{0}{1}\km\W{1}{1}\km\W{0}{3}\km\W{3}{3}\km\W{0}{4}\km\W{1}{4}\km\W{6}{6}
+\frac{3}{16} \W{0}{1}\km\W{1}{1}\km\W{0}{3}\km\W{3}{3}\km\W{0}{4}\km\W{2}{4}\km\W{5}{6}
+\frac{5}{16} \W{0}{1}\km\W{1}{1}\km\W{0}{3}\km\W{3}{3}\km\W{0}{4}\km\W{4}{4}\km\W{3}{6}
-\frac{9}{16} \W{0}{1}\km\W{1}{1}\km\W{0}{3}\km\W{3}{3}\km\W{0}{4}\km\W{3}{4}\km\W{4}{6}
-\frac{3}{8} \W{0}{1}\km\W{1}{1}\km\W{0}{3}\km\W{3}{3}\km\W{1}{4}\km\W{2}{4}\km\W{4}{6}
+\W{0}{1}\km\W{1}{1}\km\W{0}{3}\km\W{3}{3}\km\W{1}{4}\km\W{3}{4}\km\W{3}{6}
\\&
-\frac{9}{16} \W{0}{1}\km\W{1}{1}\km\W{0}{3}\km\W{3}{3}\km\W{1}{4}\km\W{4}{4}\km\W{2}{6}
-\frac{3}{8} \W{0}{1}\km\W{1}{1}\km\W{0}{3}\km\W{3}{3}\km\W{2}{4}\km\W{3}{4}\km\W{2}{6}
+\frac{3}{16} \W{0}{1}\km\W{1}{1}\km\W{0}{3}\km\W{3}{3}\km\W{2}{4}\km\W{4}{4}\km\W{1}{6}
+\frac{1}{16} \W{0}{1}\km\W{1}{1}\km\W{0}{3}\km\W{3}{3}\km\W{3}{4}\km\W{4}{4}\km\W{0}{6}
+\frac{9}{16} \W{0}{1}\km\W{1}{1}\km\W{1}{3}\km\W{2}{3}\km\W{0}{4}\km\W{1}{4}\km\W{6}{6}
-\frac{9}{16} \W{0}{1}\km\W{1}{1}\km\W{1}{3}\km\W{2}{3}\km\W{0}{4}\km\W{2}{4}\km\W{5}{6}
\\&
-\frac{9}{16} \W{0}{1}\km\W{1}{1}\km\W{1}{3}\km\W{2}{3}\km\W{0}{4}\km\W{3}{4}\km\W{4}{6}
+\frac{9}{16} \W{0}{1}\km\W{1}{1}\km\W{1}{3}\km\W{2}{3}\km\W{0}{4}\km\W{4}{4}\km\W{3}{6}
+\frac{9}{8} \W{0}{1}\km\W{1}{1}\km\W{1}{3}\km\W{2}{3}\km\W{1}{4}\km\W{2}{4}\km\W{4}{6}
-\frac{9}{16} \W{0}{1}\km\W{1}{1}\km\W{1}{3}\km\W{2}{3}\km\W{1}{4}\km\W{4}{4}\km\W{2}{6}
+\frac{9}{8} \W{0}{1}\km\W{1}{1}\km\W{1}{3}\km\W{2}{3}\km\W{2}{4}\km\W{3}{4}\km\W{2}{6}
-\frac{9}{16} \W{0}{1}\km\W{1}{1}\km\W{1}{3}\km\W{2}{3}\km\W{2}{4}\km\W{4}{4}\km\W{1}{6}
\\&
+\frac{9}{16} \W{0}{1}\km\W{1}{1}\km\W{1}{3}\km\W{2}{3}\km\W{3}{4}\km\W{4}{4}\km\W{0}{6}
-\frac{3}{4} \W{0}{1}\km\W{1}{1}\km\W{1}{3}\km\W{3}{3}\km\W{0}{4}\km\W{1}{4}\km\W{5}{6}
+\frac{9}{8} \W{0}{1}\km\W{1}{1}\km\W{1}{3}\km\W{3}{3}\km\W{0}{4}\km\W{2}{4}\km\W{4}{6}
-\frac{3}{8} \W{0}{1}\km\W{1}{1}\km\W{1}{3}\km\W{3}{3}\km\W{0}{4}\km\W{4}{4}\km\W{2}{6}
-\frac{3}{2} \W{0}{1}\km\W{1}{1}\km\W{1}{3}\km\W{3}{3}\km\W{1}{4}\km\W{2}{4}\km\W{3}{6}
+\frac{3}{4} \W{0}{1}\km\W{1}{1}\km\W{1}{3}\km\W{3}{3}\km\W{1}{4}\km\W{4}{4}\km\W{1}{6}
\\&
-\frac{3}{8} \W{0}{1}\km\W{1}{1}\km\W{1}{3}\km\W{3}{3}\km\W{2}{4}\km\W{4}{4}\km\W{0}{6}
+\frac{3}{4} \W{0}{1}\km\W{1}{1}\km\W{2}{3}\km\W{3}{3}\km\W{0}{4}\km\W{1}{4}\km\W{4}{6}
-\frac{3}{2} \W{0}{1}\km\W{1}{1}\km\W{2}{3}\km\W{3}{3}\km\W{0}{4}\km\W{2}{4}\km\W{3}{6}
+\frac{3}{4} \W{0}{1}\km\W{1}{1}\km\W{2}{3}\km\W{3}{3}\km\W{0}{4}\km\W{3}{4}\km\W{2}{6}
+\frac{9}{4} \W{0}{1}\km\W{1}{1}\km\W{2}{3}\km\W{3}{3}\km\W{1}{4}\km\W{2}{4}\km\W{2}{6}
-\frac{3}{2} \W{0}{1}\km\W{1}{1}\km\W{2}{3}\km\W{3}{3}\km\W{1}{4}\km\W{3}{4}\km\W{1}{6}
\\&
+\frac{3}{4} \W{0}{1}\km\W{1}{1}\km\W{2}{3}\km\W{3}{3}\km\W{2}{4}\km\W{3}{4}\km\W{0}{6}
\\
\widehat{\mathbf{w}}_{4} &=
\frac{1}{4} \W{0}{1}\km\W{1}{1}\km\W{0}{3}\km\W{1}{3}\km\W{1}{4}\km\W{2}{4}\km\W{6}{6}
-\frac{1}{2} \W{0}{1}\km\W{1}{1}\km\W{0}{3}\km\W{1}{3}\km\W{1}{4}\km\W{3}{4}\km\W{5}{6}
+\frac{1}{4} \W{0}{1}\km\W{1}{1}\km\W{0}{3}\km\W{1}{3}\km\W{1}{4}\km\W{4}{4}\km\W{4}{6}
-\frac{1}{2} \W{0}{1}\km\W{1}{1}\km\W{0}{3}\km\W{1}{3}\km\W{2}{4}\km\W{4}{4}\km\W{3}{6}
+\frac{1}{4} \W{0}{1}\km\W{1}{1}\km\W{0}{3}\km\W{1}{3}\km\W{3}{4}\km\W{4}{4}\km\W{2}{6}
+\frac{3}{4} \W{0}{1}\km\W{1}{1}\km\W{0}{3}\km\W{1}{3}\km\W{2}{4}\km\W{3}{4}\km\W{4}{6}
\\&
+\frac{1}{4} \W{0}{1}\km\W{1}{1}\km\W{0}{3}\km\W{2}{3}\km\W{0}{4}\km\W{3}{4}\km\W{5}{6}
-\frac{1}{8} \W{0}{1}\km\W{1}{1}\km\W{0}{3}\km\W{2}{3}\km\W{0}{4}\km\W{2}{4}\km\W{6}{6}
-\frac{1}{8} \W{0}{1}\km\W{1}{1}\km\W{0}{3}\km\W{2}{3}\km\W{0}{4}\km\W{4}{4}\km\W{4}{6}
-\frac{1}{2} \W{0}{1}\km\W{1}{1}\km\W{0}{3}\km\W{2}{3}\km\W{2}{4}\km\W{3}{4}\km\W{3}{6}
+\frac{3}{8} \W{0}{1}\km\W{1}{1}\km\W{0}{3}\km\W{2}{3}\km\W{2}{4}\km\W{4}{4}\km\W{2}{6}
-\frac{1}{4} \W{0}{1}\km\W{1}{1}\km\W{0}{3}\km\W{2}{3}\km\W{3}{4}\km\W{4}{4}\km\W{1}{6}
\\&
+\frac{1}{16} \W{0}{1}\km\W{1}{1}\km\W{0}{3}\km\W{3}{3}\km\W{0}{4}\km\W{1}{4}\km\W{6}{6}
-\frac{1}{16} \W{0}{1}\km\W{1}{1}\km\W{0}{3}\km\W{3}{3}\km\W{0}{4}\km\W{2}{4}\km\W{5}{6}
+\frac{1}{16} \W{0}{1}\km\W{1}{1}\km\W{0}{3}\km\W{3}{3}\km\W{0}{4}\km\W{4}{4}\km\W{3}{6}
-\frac{1}{16} \W{0}{1}\km\W{1}{1}\km\W{0}{3}\km\W{3}{3}\km\W{0}{4}\km\W{3}{4}\km\W{4}{6}
+\frac{1}{8} \W{0}{1}\km\W{1}{1}\km\W{0}{3}\km\W{3}{3}\km\W{1}{4}\km\W{2}{4}\km\W{4}{6}
-\frac{1}{16} \W{0}{1}\km\W{1}{1}\km\W{0}{3}\km\W{3}{3}\km\W{1}{4}\km\W{4}{4}\km\W{2}{6}
\\&
+\frac{1}{8} \W{0}{1}\km\W{1}{1}\km\W{0}{3}\km\W{3}{3}\km\W{2}{4}\km\W{3}{4}\km\W{2}{6}
-\frac{1}{16} \W{0}{1}\km\W{1}{1}\km\W{0}{3}\km\W{3}{3}\km\W{2}{4}\km\W{4}{4}\km\W{1}{6}
+\frac{1}{16} \W{0}{1}\km\W{1}{1}\km\W{0}{3}\km\W{3}{3}\km\W{3}{4}\km\W{4}{4}\km\W{0}{6}
+\frac{1}{16} \W{0}{1}\km\W{1}{1}\km\W{1}{3}\km\W{2}{3}\km\W{0}{4}\km\W{1}{4}\km\W{6}{6}
+\frac{3}{16} \W{0}{1}\km\W{1}{1}\km\W{1}{3}\km\W{2}{3}\km\W{0}{4}\km\W{2}{4}\km\W{5}{6}
-\frac{9}{16} \W{0}{1}\km\W{1}{1}\km\W{1}{3}\km\W{2}{3}\km\W{0}{4}\km\W{3}{4}\km\W{4}{6}
\\&
+\frac{5}{16} \W{0}{1}\km\W{1}{1}\km\W{1}{3}\km\W{2}{3}\km\W{0}{4}\km\W{4}{4}\km\W{3}{6}
-\frac{3}{8} \W{0}{1}\km\W{1}{1}\km\W{1}{3}\km\W{2}{3}\km\W{1}{4}\km\W{2}{4}\km\W{4}{6}
+\W{0}{1}\km\W{1}{1}\km\W{1}{3}\km\W{2}{3}\km\W{1}{4}\km\W{3}{4}\km\W{3}{6}
-\frac{9}{16} \W{0}{1}\km\W{1}{1}\km\W{1}{3}\km\W{2}{3}\km\W{1}{4}\km\W{4}{4}\km\W{2}{6}
-\frac{3}{8} \W{0}{1}\km\W{1}{1}\km\W{1}{3}\km\W{2}{3}\km\W{2}{4}\km\W{3}{4}\km\W{2}{6}
+\frac{3}{16} \W{0}{1}\km\W{1}{1}\km\W{1}{3}\km\W{2}{3}\km\W{2}{4}\km\W{4}{4}\km\W{1}{6}
\\&
+\frac{1}{16} \W{0}{1}\km\W{1}{1}\km\W{1}{3}\km\W{2}{3}\km\W{3}{4}\km\W{4}{4}\km\W{0}{6}
-\frac{1}{4} \W{0}{1}\km\W{1}{1}\km\W{1}{3}\km\W{3}{3}\km\W{0}{4}\km\W{1}{4}\km\W{5}{6}
+\frac{3}{8} \W{0}{1}\km\W{1}{1}\km\W{1}{3}\km\W{3}{3}\km\W{0}{4}\km\W{2}{4}\km\W{4}{6}
-\frac{1}{8} \W{0}{1}\km\W{1}{1}\km\W{1}{3}\km\W{3}{3}\km\W{0}{4}\km\W{4}{4}\km\W{2}{6}
-\frac{1}{2} \W{0}{1}\km\W{1}{1}\km\W{1}{3}\km\W{3}{3}\km\W{1}{4}\km\W{2}{4}\km\W{3}{6}
+\frac{1}{4} \W{0}{1}\km\W{1}{1}\km\W{1}{3}\km\W{3}{3}\km\W{1}{4}\km\W{4}{4}\km\W{1}{6}
\\&
-\frac{1}{8} \W{0}{1}\km\W{1}{1}\km\W{1}{3}\km\W{3}{3}\km\W{2}{4}\km\W{4}{4}\km\W{0}{6}
+\frac{1}{4} \W{0}{1}\km\W{1}{1}\km\W{2}{3}\km\W{3}{3}\km\W{0}{4}\km\W{1}{4}\km\W{4}{6}
-\frac{1}{2} \W{0}{1}\km\W{1}{1}\km\W{2}{3}\km\W{3}{3}\km\W{0}{4}\km\W{2}{4}\km\W{3}{6}
+\frac{1}{4} \W{0}{1}\km\W{1}{1}\km\W{2}{3}\km\W{3}{3}\km\W{0}{4}\km\W{3}{4}\km\W{2}{6}
+\frac{3}{4} \W{0}{1}\km\W{1}{1}\km\W{2}{3}\km\W{3}{3}\km\W{1}{4}\km\W{2}{4}\km\W{2}{6}
-\frac{1}{2} \W{0}{1}\km\W{1}{1}\km\W{2}{3}\km\W{3}{3}\km\W{1}{4}\km\W{3}{4}\km\W{1}{6}
\\&
+\frac{1}{4} \W{0}{1}\km\W{1}{1}\km\W{2}{3}\km\W{3}{3}\km\W{2}{4}\km\W{3}{4}\km\W{0}{6}
\\
\widehat{\mathbf{w}}_{5} &=
\frac{1}{2} \W{0}{1}\km\W{1}{1}\km\W{0}{3}\km\W{1}{3}\km\W{0}{4}\km\W{4}{4}\km\W{5}{6}
-\frac{1}{2} \W{0}{1}\km\W{1}{1}\km\W{0}{3}\km\W{1}{3}\km\W{0}{4}\km\W{3}{4}\km\W{6}{6}
+3 \W{0}{1}\km\W{1}{1}\km\W{0}{3}\km\W{1}{3}\km\W{1}{4}\km\W{2}{4}\km\W{6}{6}
-4 \W{0}{1}\km\W{1}{1}\km\W{0}{3}\km\W{1}{3}\km\W{1}{4}\km\W{3}{4}\km\W{5}{6}
+\W{0}{1}\km\W{1}{1}\km\W{0}{3}\km\W{1}{3}\km\W{1}{4}\km\W{4}{4}\km\W{4}{6}
-3 \W{0}{1}\km\W{1}{1}\km\W{0}{3}\km\W{1}{3}\km\W{2}{4}\km\W{4}{4}\km\W{3}{6}
\\&
+\frac{3}{2} \W{0}{1}\km\W{1}{1}\km\W{0}{3}\km\W{1}{3}\km\W{3}{4}\km\W{4}{4}\km\W{2}{6}
+6 \W{0}{1}\km\W{1}{1}\km\W{0}{3}\km\W{1}{3}\km\W{2}{4}\km\W{3}{4}\km\W{4}{6}
+\frac{5}{2} \W{0}{1}\km\W{1}{1}\km\W{0}{3}\km\W{2}{3}\km\W{0}{4}\km\W{3}{4}\km\W{5}{6}
-\frac{3}{4} \W{0}{1}\km\W{1}{1}\km\W{0}{3}\km\W{2}{3}\km\W{0}{4}\km\W{2}{4}\km\W{6}{6}
-\frac{7}{4} \W{0}{1}\km\W{1}{1}\km\W{0}{3}\km\W{2}{3}\km\W{0}{4}\km\W{4}{4}\km\W{4}{6}
-3 \W{0}{1}\km\W{1}{1}\km\W{0}{3}\km\W{2}{3}\km\W{1}{4}\km\W{2}{4}\km\W{5}{6}
\\& 
+2 \W{0}{1}\km\W{1}{1}\km\W{0}{3}\km\W{2}{3}\km\W{1}{4}\km\W{3}{4}\km\W{4}{6}
+\W{0}{1}\km\W{1}{1}\km\W{0}{3}\km\W{2}{3}\km\W{1}{4}\km\W{4}{4}\km\W{3}{6}
-6 \W{0}{1}\km\W{1}{1}\km\W{0}{3}\km\W{2}{3}\km\W{2}{4}\km\W{3}{4}\km\W{3}{6}
+\frac{9}{4} \W{0}{1}\km\W{1}{1}\km\W{0}{3}\km\W{2}{3}\km\W{2}{4}\km\W{4}{4}\km\W{2}{6}
-\frac{3}{2} \W{0}{1}\km\W{1}{1}\km\W{0}{3}\km\W{2}{3}\km\W{3}{4}\km\W{4}{4}\km\W{1}{6}
+\frac{1}{4} \W{0}{1}\km\W{1}{1}\km\W{0}{3}\km\W{3}{3}\km\W{0}{4}\km\W{1}{4}\km\W{6}{6}
\\&
+\W{0}{1}\km\W{1}{1}\km\W{0}{3}\km\W{3}{3}\km\W{0}{4}\km\W{4}{4}\km\W{3}{6}
-\frac{5}{4} \W{0}{1}\km\W{1}{1}\km\W{0}{3}\km\W{3}{3}\km\W{0}{4}\km\W{3}{4}\km\W{4}{6}
+\frac{3}{2} \W{0}{1}\km\W{1}{1}\km\W{0}{3}\km\W{3}{3}\km\W{1}{4}\km\W{2}{4}\km\W{4}{6}
-\frac{5}{4} \W{0}{1}\km\W{1}{1}\km\W{0}{3}\km\W{3}{3}\km\W{1}{4}\km\W{4}{4}\km\W{2}{6}
+\frac{3}{2} \W{0}{1}\km\W{1}{1}\km\W{0}{3}\km\W{3}{3}\km\W{2}{4}\km\W{3}{4}\km\W{2}{6}
+\frac{1}{4} \W{0}{1}\km\W{1}{1}\km\W{0}{3}\km\W{3}{3}\km\W{3}{4}\km\W{4}{4}\km\W{0}{6}
\\&
+\frac{3}{4} \W{0}{1}\km\W{1}{1}\km\W{1}{3}\km\W{2}{3}\km\W{0}{4}\km\W{1}{4}\km\W{6}{6}
-\frac{15}{4} \W{0}{1}\km\W{1}{1}\km\W{1}{3}\km\W{2}{3}\km\W{0}{4}\km\W{3}{4}\km\W{4}{6}
+3 \W{0}{1}\km\W{1}{1}\km\W{1}{3}\km\W{2}{3}\km\W{0}{4}\km\W{4}{4}\km\W{3}{6}
+\frac{9}{2} \W{0}{1}\km\W{1}{1}\km\W{1}{3}\km\W{2}{3}\km\W{1}{4}\km\W{2}{4}\km\W{4}{6}
-\frac{15}{4} \W{0}{1}\km\W{1}{1}\km\W{1}{3}\km\W{2}{3}\km\W{1}{4}\km\W{4}{4}\km\W{2}{6}
+\frac{9}{2} \W{0}{1}\km\W{1}{1}\km\W{1}{3}\km\W{2}{3}\km\W{2}{4}\km\W{3}{4}\km\W{2}{6}
\\&
+\frac{3}{4} \W{0}{1}\km\W{1}{1}\km\W{1}{3}\km\W{2}{3}\km\W{3}{4}\km\W{4}{4}\km\W{0}{6}
-\frac{3}{2} \W{0}{1}\km\W{1}{1}\km\W{1}{3}\km\W{3}{3}\km\W{0}{4}\km\W{1}{4}\km\W{5}{6}
+\frac{9}{4} \W{0}{1}\km\W{1}{1}\km\W{1}{3}\km\W{3}{3}\km\W{0}{4}\km\W{2}{4}\km\W{4}{6}
+\W{0}{1}\km\W{1}{1}\km\W{1}{3}\km\W{3}{3}\km\W{0}{4}\km\W{3}{4}\km\W{3}{6}
-\frac{7}{4} \W{0}{1}\km\W{1}{1}\km\W{1}{3}\km\W{3}{3}\km\W{0}{4}\km\W{4}{4}\km\W{2}{6}
-6 \W{0}{1}\km\W{1}{1}\km\W{1}{3}\km\W{3}{3}\km\W{1}{4}\km\W{2}{4}\km\W{3}{6}
\\&
+2 \W{0}{1}\km\W{1}{1}\km\W{1}{3}\km\W{3}{3}\km\W{1}{4}\km\W{3}{4}\km\W{2}{6}
+\frac{5}{2} \W{0}{1}\km\W{1}{1}\km\W{1}{3}\km\W{3}{3}\km\W{1}{4}\km\W{4}{4}\km\W{1}{6}
-3 \W{0}{1}\km\W{1}{1}\km\W{1}{3}\km\W{3}{3}\km\W{2}{4}\km\W{3}{4}\km\W{1}{6}
-\frac{3}{4} \W{0}{1}\km\W{1}{1}\km\W{1}{3}\km\W{3}{3}\km\W{2}{4}\km\W{4}{4}\km\W{0}{6}
+\frac{3}{2} \W{0}{1}\km\W{1}{1}\km\W{2}{3}\km\W{3}{3}\km\W{0}{4}\km\W{1}{4}\km\W{4}{6}
-3 \W{0}{1}\km\W{1}{1}\km\W{2}{3}\km\W{3}{3}\km\W{0}{4}\km\W{2}{4}\km\W{3}{6}
\\&
+\W{0}{1}\km\W{1}{1}\km\W{2}{3}\km\W{3}{3}\km\W{0}{4}\km\W{3}{4}\km\W{2}{6}
+\frac{1}{2} \W{0}{1}\km\W{1}{1}\km\W{2}{3}\km\W{3}{3}\km\W{0}{4}\km\W{4}{4}\km\W{1}{6}
+6 \W{0}{1}\km\W{1}{1}\km\W{2}{3}\km\W{3}{3}\km\W{1}{4}\km\W{2}{4}\km\W{2}{6}
-4 \W{0}{1}\km\W{1}{1}\km\W{2}{3}\km\W{3}{3}\km\W{1}{4}\km\W{3}{4}\km\W{1}{6}
-\frac{1}{2} \W{0}{1}\km\W{1}{1}\km\W{2}{3}\km\W{3}{3}\km\W{1}{4}\km\W{4}{4}\km\W{0}{6}
+3 \W{0}{1}\km\W{1}{1}\km\W{2}{3}\km\W{3}{3}\km\W{2}{4}\km\W{3}{4}\km\W{0}{6} 
\\
\widehat{\mathbf{w}}_{6} &= 
-\W{0}{1}\km\W{1}{1}\km\W{0}{3}\km\W{1}{3}\km\W{0}{4}\km\W{4}{5}\km\W{5}{5}
+2 \W{0}{1}\km\W{1}{1}\km\W{0}{3}\km\W{1}{3}\km\W{1}{4}\km\W{3}{5}\km\W{5}{5}
-3 \W{0}{1}\km\W{1}{1}\km\W{0}{3}\km\W{1}{3}\km\W{2}{4}\km\W{2}{5}\km\W{5}{5}
+3 \W{0}{1}\km\W{1}{1}\km\W{0}{3}\km\W{1}{3}\km\W{2}{4}\km\W{3}{5}\km\W{4}{5}
-2 \W{0}{1}\km\W{1}{1}\km\W{0}{3}\km\W{1}{3}\km\W{3}{4}\km\W{2}{5}\km\W{4}{5}
-\frac{1}{10} \W{0}{1}\km\W{1}{1}\km\W{0}{3}\km\W{1}{3}\km\W{4}{4}\km\W{0}{5}\km\W{5}{5}
\\&
+2 \W{0}{1}\km\W{1}{1}\km\W{0}{3}\km\W{1}{3}\km\W{3}{4}\km\W{1}{5}\km\W{5}{5}
-\frac{3}{2} \W{0}{1}\km\W{1}{1}\km\W{0}{3}\km\W{1}{3}\km\W{4}{4}\km\W{1}{5}\km\W{4}{5}
+4 \W{0}{1}\km\W{1}{1}\km\W{0}{3}\km\W{1}{3}\km\W{4}{4}\km\W{2}{5}\km\W{3}{5}
+\W{0}{1}\km\W{1}{1}\km\W{0}{3}\km\W{2}{3}\km\W{0}{4}\km\W{3}{5}\km\W{5}{5}
-8 \W{0}{1}\km\W{1}{1}\km\W{0}{3}\km\W{2}{3}\km\W{1}{4}\km\W{3}{5}\km\W{4}{5}
+6 \W{0}{1}\km\W{1}{1}\km\W{0}{3}\km\W{2}{3}\km\W{2}{4}\km\W{2}{5}\km\W{4}{5}
\\&
-\frac{4}{5} \W{0}{1}\km\W{1}{1}\km\W{0}{3}\km\W{2}{3}\km\W{3}{4}\km\W{0}{5}\km\W{5}{5}
-4 \W{0}{1}\km\W{1}{1}\km\W{0}{3}\km\W{2}{3}\km\W{3}{4}\km\W{2}{5}\km\W{3}{5}
+\W{0}{1}\km\W{1}{1}\km\W{0}{3}\km\W{2}{3}\km\W{4}{4}\km\W{0}{5}\km\W{4}{5}
-\W{0}{1}\km\W{1}{1}\km\W{0}{3}\km\W{2}{3}\km\W{4}{4}\km\W{1}{5}\km\W{3}{5}
+\frac{13}{2} \W{0}{1}\km\W{1}{1}\km\W{0}{3}\km\W{3}{3}\km\W{0}{4}\km\W{3}{5}\km\W{4}{5}
-\frac{5}{2} \W{0}{1}\km\W{1}{1}\km\W{0}{3}\km\W{3}{3}\km\W{0}{4}\km\W{2}{5}\km\W{5}{5}
\\&
+4 \W{0}{1}\km\W{1}{1}\km\W{0}{3}\km\W{3}{3}\km\W{1}{4}\km\W{1}{5}\km\W{5}{5}
-6 \W{0}{1}\km\W{1}{1}\km\W{0}{3}\km\W{3}{3}\km\W{1}{4}\km\W{2}{5}\km\W{4}{5}
-\frac{8}{5} \W{0}{1}\km\W{1}{1}\km\W{0}{3}\km\W{3}{3}\km\W{2}{4}\km\W{0}{5}\km\W{5}{5}
-4 \W{0}{1}\km\W{1}{1}\km\W{0}{3}\km\W{3}{3}\km\W{2}{4}\km\W{1}{5}\km\W{4}{5}
+14 \W{0}{1}\km\W{1}{1}\km\W{0}{3}\km\W{3}{3}\km\W{2}{4}\km\W{2}{5}\km\W{3}{5}
+4 \W{0}{1}\km\W{1}{1}\km\W{0}{3}\km\W{3}{3}\km\W{3}{4}\km\W{0}{5}\km\W{4}{5}
\\&
-6 \W{0}{1}\km\W{1}{1}\km\W{0}{3}\km\W{3}{3}\km\W{3}{4}\km\W{1}{5}\km\W{3}{5}
-\frac{5}{2} \W{0}{1}\km\W{1}{1}\km\W{0}{3}\km\W{3}{3}\km\W{4}{4}\km\W{0}{5}\km\W{3}{5}
+\frac{13}{2} \W{0}{1}\km\W{1}{1}\km\W{0}{3}\km\W{3}{3}\km\W{4}{4}\km\W{1}{5}\km\W{2}{5}
+\frac{9}{2} \W{0}{1}\km\W{1}{1}\km\W{1}{3}\km\W{2}{3}\km\W{0}{4}\km\W{2}{5}\km\W{5}{5}
-\frac{33}{2} \W{0}{1}\km\W{1}{1}\km\W{1}{3}\km\W{2}{3}\km\W{0}{4}\km\W{3}{5}\km\W{4}{5}
-12 \W{0}{1}\km\W{1}{1}\km\W{1}{3}\km\W{2}{3}\km\W{1}{4}\km\W{1}{5}\km\W{5}{5}
\\&
+30 \W{0}{1}\km\W{1}{1}\km\W{1}{3}\km\W{2}{3}\km\W{1}{4}\km\W{2}{5}\km\W{4}{5}
+\frac{36}{5} \W{0}{1}\km\W{1}{1}\km\W{1}{3}\km\W{2}{3}\km\W{2}{4}\km\W{0}{5}\km\W{5}{5}
-54 \W{0}{1}\km\W{1}{1}\km\W{1}{3}\km\W{2}{3}\km\W{2}{4}\km\W{2}{5}\km\W{3}{5}
-12 \W{0}{1}\km\W{1}{1}\km\W{1}{3}\km\W{2}{3}\km\W{3}{4}\km\W{0}{5}\km\W{4}{5}
+30 \W{0}{1}\km\W{1}{1}\km\W{1}{3}\km\W{2}{3}\km\W{3}{4}\km\W{1}{5}\km\W{3}{5}
+\frac{9}{2} \W{0}{1}\km\W{1}{1}\km\W{1}{3}\km\W{2}{3}\km\W{4}{4}\km\W{0}{5}\km\W{3}{5}
\\&
-\frac{33}{2} \W{0}{1}\km\W{1}{1}\km\W{1}{3}\km\W{2}{3}\km\W{4}{4}\km\W{1}{5}\km\W{2}{5}
+\W{0}{1}\km\W{1}{1}\km\W{1}{3}\km\W{3}{3}\km\W{0}{4}\km\W{1}{5}\km\W{5}{5}
-\W{0}{1}\km\W{1}{1}\km\W{1}{3}\km\W{3}{3}\km\W{0}{4}\km\W{2}{5}\km\W{4}{5}
-\frac{4}{5} \W{0}{1}\km\W{1}{1}\km\W{1}{3}\km\W{3}{3}\km\W{1}{4}\km\W{0}{5}\km\W{5}{5}
-4 \W{0}{1}\km\W{1}{1}\km\W{1}{3}\km\W{3}{3}\km\W{1}{4}\km\W{2}{5}\km\W{3}{5}
+6 \W{0}{1}\km\W{1}{1}\km\W{1}{3}\km\W{3}{3}\km\W{2}{4}\km\W{1}{5}\km\W{3}{5}
\\&
-8 \W{0}{1}\km\W{1}{1}\km\W{1}{3}\km\W{3}{3}\km\W{3}{4}\km\W{1}{5}\km\W{2}{5}
+\W{0}{1}\km\W{1}{1}\km\W{1}{3}\km\W{3}{3}\km\W{4}{4}\km\W{0}{5}\km\W{2}{5}
-\frac{1}{10} \W{0}{1}\km\W{1}{1}\km\W{2}{3}\km\W{3}{3}\km\W{0}{4}\km\W{0}{5}\km\W{5}{5}
+4 \W{0}{1}\km\W{1}{1}\km\W{2}{3}\km\W{3}{3}\km\W{0}{4}\km\W{2}{5}\km\W{3}{5}
+2 \W{0}{1}\km\W{1}{1}\km\W{2}{3}\km\W{3}{3}\km\W{1}{4}\km\W{0}{5}\km\W{4}{5}
-\frac{3}{2} \W{0}{1}\km\W{1}{1}\km\W{2}{3}\km\W{3}{3}\km\W{0}{4}\km\W{1}{5}\km\W{4}{5}
\\&
-2 \W{0}{1}\km\W{1}{1}\km\W{2}{3}\km\W{3}{3}\km\W{1}{4}\km\W{1}{5}\km\W{3}{5}
-3 \W{0}{1}\km\W{1}{1}\km\W{2}{3}\km\W{3}{3}\km\W{2}{4}\km\W{0}{5}\km\W{3}{5}
+3 \W{0}{1}\km\W{1}{1}\km\W{2}{3}\km\W{3}{3}\km\W{2}{4}\km\W{1}{5}\km\W{2}{5}
+2 \W{0}{1}\km\W{1}{1}\km\W{2}{3}\km\W{3}{3}\km\W{3}{4}\km\W{0}{5}\km\W{2}{5}
\\&
-\W{0}{1}\km\W{1}{1}\km\W{2}{3}\km\W{3}{3}\km\W{4}{4}\km\W{0}{5}\km\W{1}{5}
\\
\widehat{\mathbf{w}}_{7} &= 
4 \W{0}{1}\km\W{1}{1}\km\W{1}{3}\km\W{2}{3}\km\W{3}{4}\km\W{1}{5}\km\W{3}{5}
+\W{0}{1}\km\W{1}{1}\km\W{1}{3}\km\W{2}{3}\km\W{4}{4}\km\W{0}{5}\km\W{3}{5}
-3 \W{0}{1}\km\W{1}{1}\km\W{1}{3}\km\W{2}{3}\km\W{4}{4}\km\W{1}{5}\km\W{2}{5}
+\frac{2}{3} \W{0}{1}\km\W{1}{1}\km\W{0}{3}\km\W{3}{3}\km\W{3}{4}\km\W{0}{5}\km\W{4}{5}
-\frac{4}{3} \W{0}{1}\km\W{1}{1}\km\W{0}{3}\km\W{3}{3}\km\W{3}{4}\km\W{1}{5}\km\W{3}{5}
+\W{0}{1}\km\W{1}{1}\km\W{1}{3}\km\W{2}{3}\km\W{2}{4}\km\W{1}{5}\km\W{4}{5}
\\&
-\frac{1}{3} \W{0}{1}\km\W{1}{1}\km\W{0}{3}\km\W{3}{3}\km\W{0}{4}\km\W{2}{5}\km\W{5}{5}
+\W{0}{1}\km\W{1}{1}\km\W{0}{3}\km\W{3}{3}\km\W{0}{4}\km\W{3}{5}\km\W{4}{5}
+\frac{2}{3} \W{0}{1}\km\W{1}{1}\km\W{0}{3}\km\W{3}{3}\km\W{1}{4}\km\W{1}{5}\km\W{5}{5}
-\frac{4}{3} \W{0}{1}\km\W{1}{1}\km\W{0}{3}\km\W{3}{3}\km\W{1}{4}\km\W{2}{5}\km\W{4}{5}
-\frac{1}{3} \W{0}{1}\km\W{1}{1}\km\W{0}{3}\km\W{3}{3}\km\W{2}{4}\km\W{0}{5}\km\W{5}{5}
-\frac{1}{3} \W{0}{1}\km\W{1}{1}\km\W{0}{3}\km\W{3}{3}\km\W{2}{4}\km\W{1}{5}\km\W{4}{5}
\\&
+\frac{8}{3} \W{0}{1}\km\W{1}{1}\km\W{0}{3}\km\W{3}{3}\km\W{2}{4}\km\W{2}{5}\km\W{3}{5}
-2 \W{0}{1}\km\W{1}{1}\km\W{1}{3}\km\W{2}{3}\km\W{3}{4}\km\W{0}{5}\km\W{4}{5}
+\W{0}{1}\km\W{1}{1}\km\W{1}{3}\km\W{2}{3}\km\W{0}{4}\km\W{2}{5}\km\W{5}{5}
-3 \W{0}{1}\km\W{1}{1}\km\W{1}{3}\km\W{2}{3}\km\W{0}{4}\km\W{3}{5}\km\W{4}{5}
-2 \W{0}{1}\km\W{1}{1}\km\W{1}{3}\km\W{2}{3}\km\W{1}{4}\km\W{1}{5}\km\W{5}{5}
+4 \W{0}{1}\km\W{1}{1}\km\W{1}{3}\km\W{2}{3}\km\W{1}{4}\km\W{2}{5}\km\W{4}{5}
\\&
+\W{0}{1}\km\W{1}{1}\km\W{1}{3}\km\W{2}{3}\km\W{2}{4}\km\W{0}{5}\km\W{5}{5}
-8 \W{0}{1}\km\W{1}{1}\km\W{1}{3}\km\W{2}{3}\km\W{2}{4}\km\W{2}{5}\km\W{3}{5}
-\frac{1}{3} \W{0}{1}\km\W{1}{1}\km\W{0}{3}\km\W{3}{3}\km\W{4}{4}\km\W{0}{5}\km\W{3}{5}
+\W{0}{1}\km\W{1}{1}\km\W{0}{3}\km\W{3}{3}\km\W{4}{4}\km\W{1}{5}\km\W{2}{5}
\\
\widehat{\mathbf{w}}_{8} &= 
\frac{1}{10} \W{0}{1}\km\W{1}{1}\km\W{0}{3}\km\W{1}{3}\km\W{4}{4}\km\W{0}{5}\km\W{5}{5}
-\frac{1}{2} \W{0}{1}\km\W{1}{1}\km\W{0}{3}\km\W{1}{3}\km\W{4}{4}\km\W{1}{5}\km\W{4}{5}
+\W{0}{1}\km\W{1}{1}\km\W{0}{3}\km\W{1}{3}\km\W{4}{4}\km\W{2}{5}\km\W{3}{5}
-\frac{1}{5} \W{0}{1}\km\W{1}{1}\km\W{0}{3}\km\W{2}{3}\km\W{3}{4}\km\W{0}{5}\km\W{5}{5}
+\W{0}{1}\km\W{1}{1}\km\W{0}{3}\km\W{2}{3}\km\W{3}{4}\km\W{1}{5}\km\W{4}{5}
-2 \W{0}{1}\km\W{1}{1}\km\W{0}{3}\km\W{2}{3}\km\W{3}{4}\km\W{2}{5}\km\W{3}{5}
\\&
+\frac{3}{2} \W{0}{1}\km\W{1}{1}\km\W{0}{3}\km\W{3}{3}\km\W{0}{4}\km\W{3}{5}\km\W{4}{5}
-\frac{1}{2} \W{0}{1}\km\W{1}{1}\km\W{0}{3}\km\W{3}{3}\km\W{0}{4}\km\W{2}{5}\km\W{5}{5}
+\W{0}{1}\km\W{1}{1}\km\W{0}{3}\km\W{3}{3}\km\W{1}{4}\km\W{1}{5}\km\W{5}{5}
-2 \W{0}{1}\km\W{1}{1}\km\W{0}{3}\km\W{3}{3}\km\W{1}{4}\km\W{2}{5}\km\W{4}{5}
-\frac{2}{5} \W{0}{1}\km\W{1}{1}\km\W{0}{3}\km\W{3}{3}\km\W{2}{4}\km\W{0}{5}\km\W{5}{5}
-\W{0}{1}\km\W{1}{1}\km\W{0}{3}\km\W{3}{3}\km\W{2}{4}\km\W{1}{5}\km\W{4}{5}
\\&
+5 \W{0}{1}\km\W{1}{1}\km\W{0}{3}\km\W{3}{3}\km\W{2}{4}\km\W{2}{5}\km\W{3}{5}
+\W{0}{1}\km\W{1}{1}\km\W{0}{3}\km\W{3}{3}\km\W{3}{4}\km\W{0}{5}\km\W{4}{5}
-2 \W{0}{1}\km\W{1}{1}\km\W{0}{3}\km\W{3}{3}\km\W{3}{4}\km\W{1}{5}\km\W{3}{5}
-\frac{1}{2} \W{0}{1}\km\W{1}{1}\km\W{0}{3}\km\W{3}{3}\km\W{4}{4}\km\W{0}{5}\km\W{3}{5}
+\frac{3}{2} \W{0}{1}\km\W{1}{1}\km\W{0}{3}\km\W{3}{3}\km\W{4}{4}\km\W{1}{5}\km\W{2}{5}
+\frac{3}{2} \W{0}{1}\km\W{1}{1}\km\W{1}{3}\km\W{2}{3}\km\W{0}{4}\km\W{2}{5}\km\W{5}{5}
\\&
-\frac{9}{2} \W{0}{1}\km\W{1}{1}\km\W{1}{3}\km\W{2}{3}\km\W{0}{4}\km\W{3}{5}\km\W{4}{5}
-3 \W{0}{1}\km\W{1}{1}\km\W{1}{3}\km\W{2}{3}\km\W{1}{4}\km\W{1}{5}\km\W{5}{5}
+6 \W{0}{1}\km\W{1}{1}\km\W{1}{3}\km\W{2}{3}\km\W{1}{4}\km\W{2}{5}\km\W{4}{5}
+\frac{9}{5} \W{0}{1}\km\W{1}{1}\km\W{1}{3}\km\W{2}{3}\km\W{2}{4}\km\W{0}{5}\km\W{5}{5}
-9 \W{0}{1}\km\W{1}{1}\km\W{1}{3}\km\W{2}{3}\km\W{2}{4}\km\W{2}{5}\km\W{3}{5}
-3 \W{0}{1}\km\W{1}{1}\km\W{1}{3}\km\W{2}{3}\km\W{3}{4}\km\W{0}{5}\km\W{4}{5}
\\&
+6 \W{0}{1}\km\W{1}{1}\km\W{1}{3}\km\W{2}{3}\km\W{3}{4}\km\W{1}{5}\km\W{3}{5}
+\frac{3}{2} \W{0}{1}\km\W{1}{1}\km\W{1}{3}\km\W{2}{3}\km\W{4}{4}\km\W{0}{5}\km\W{3}{5}
-\frac{9}{2} \W{0}{1}\km\W{1}{1}\km\W{1}{3}\km\W{2}{3}\km\W{4}{4}\km\W{1}{5}\km\W{2}{5}
-\frac{1}{5} \W{0}{1}\km\W{1}{1}\km\W{1}{3}\km\W{3}{3}\km\W{1}{4}\km\W{0}{5}\km\W{5}{5}
+\W{0}{1}\km\W{1}{1}\km\W{1}{3}\km\W{3}{3}\km\W{1}{4}\km\W{1}{5}\km\W{4}{5}
-2 \W{0}{1}\km\W{1}{1}\km\W{1}{3}\km\W{3}{3}\km\W{1}{4}\km\W{2}{5}\km\W{3}{5}
\\&
+\frac{1}{10} \W{0}{1}\km\W{1}{1}\km\W{2}{3}\km\W{3}{3}\km\W{0}{4}\km\W{0}{5}\km\W{5}{5}
+\W{0}{1}\km\W{1}{1}\km\W{2}{3}\km\W{3}{3}\km\W{0}{4}\km\W{2}{5}\km\W{3}{5}
-\frac{1}{2} \W{0}{1}\km\W{1}{1}\km\W{2}{3}\km\W{3}{3}\km\W{0}{4}\km\W{1}{5}\km\W{4}{5}
\\
\widehat{\mathbf{w}}_{9} &= 
\frac{1}{2} \W{0}{1}\km\W{1}{1}\km\W{0}{3}\km\W{1}{3}\km\W{0}{4}\km\W{4}{5}\km\W{5}{5}
-\W{0}{1}\km\W{1}{1}\km\W{0}{3}\km\W{1}{3}\km\W{1}{4}\km\W{3}{5}\km\W{5}{5}
+\frac{1}{2} \W{0}{1}\km\W{1}{1}\km\W{0}{3}\km\W{1}{3}\km\W{2}{4}\km\W{2}{5}\km\W{5}{5}
+\frac{3}{2} \W{0}{1}\km\W{1}{1}\km\W{0}{3}\km\W{1}{3}\km\W{2}{4}\km\W{3}{5}\km\W{4}{5}
-\W{0}{1}\km\W{1}{1}\km\W{0}{3}\km\W{1}{3}\km\W{3}{4}\km\W{2}{5}\km\W{4}{5}
-\frac{1}{10} \W{0}{1}\km\W{1}{1}\km\W{0}{3}\km\W{1}{3}\km\W{4}{4}\km\W{0}{5}\km\W{5}{5}
\\&
+\frac{1}{2} \W{0}{1}\km\W{1}{1}\km\W{0}{3}\km\W{1}{3}\km\W{4}{4}\km\W{1}{5}\km\W{4}{5}
-\frac{1}{2} \W{0}{1}\km\W{1}{1}\km\W{0}{3}\km\W{1}{3}\km\W{4}{4}\km\W{2}{5}\km\W{3}{5}
-\frac{1}{2} \W{0}{1}\km\W{1}{1}\km\W{0}{3}\km\W{2}{3}\km\W{0}{4}\km\W{3}{5}\km\W{5}{5}
+\W{0}{1}\km\W{1}{1}\km\W{0}{3}\km\W{2}{3}\km\W{1}{4}\km\W{2}{5}\km\W{5}{5}
+\W{0}{1}\km\W{1}{1}\km\W{0}{3}\km\W{2}{3}\km\W{1}{4}\km\W{3}{5}\km\W{4}{5}
-\frac{1}{2} \W{0}{1}\km\W{1}{1}\km\W{0}{3}\km\W{2}{3}\km\W{2}{4}\km\W{1}{5}\km\W{5}{5}
\\&
-2 \W{0}{1}\km\W{1}{1}\km\W{0}{3}\km\W{2}{3}\km\W{2}{4}\km\W{2}{5}\km\W{4}{5}
+\frac{1}{5} \W{0}{1}\km\W{1}{1}\km\W{0}{3}\km\W{2}{3}\km\W{3}{4}\km\W{0}{5}\km\W{5}{5}
+3 \W{0}{1}\km\W{1}{1}\km\W{0}{3}\km\W{2}{3}\km\W{3}{4}\km\W{2}{5}\km\W{3}{5}
-\frac{1}{2} \W{0}{1}\km\W{1}{1}\km\W{0}{3}\km\W{2}{3}\km\W{4}{4}\km\W{1}{5}\km\W{3}{5}
-3 \W{0}{1}\km\W{1}{1}\km\W{0}{3}\km\W{3}{3}\km\W{0}{4}\km\W{3}{5}\km\W{4}{5}
+\frac{7}{6} \W{0}{1}\km\W{1}{1}\km\W{0}{3}\km\W{3}{3}\km\W{0}{4}\km\W{2}{5}\km\W{5}{5}
\\&
-\frac{7}{3} \W{0}{1}\km\W{1}{1}\km\W{0}{3}\km\W{3}{3}\km\W{1}{4}\km\W{1}{5}\km\W{5}{5}
+\frac{11}{3} \W{0}{1}\km\W{1}{1}\km\W{0}{3}\km\W{3}{3}\km\W{1}{4}\km\W{2}{5}\km\W{4}{5}
+\frac{16}{15} \W{0}{1}\km\W{1}{1}\km\W{0}{3}\km\W{3}{3}\km\W{2}{4}\km\W{0}{5}\km\W{5}{5}
+\frac{13}{6} \W{0}{1}\km\W{1}{1}\km\W{0}{3}\km\W{3}{3}\km\W{2}{4}\km\W{1}{5}\km\W{4}{5}
-\frac{53}{6} \W{0}{1}\km\W{1}{1}\km\W{0}{3}\km\W{3}{3}\km\W{2}{4}\km\W{2}{5}\km\W{3}{5}
-\frac{7}{3} \W{0}{1}\km\W{1}{1}\km\W{0}{3}\km\W{3}{3}\km\W{3}{4}\km\W{0}{5}\km\W{4}{5}
\\&
+\frac{11}{3} \W{0}{1}\km\W{1}{1}\km\W{0}{3}\km\W{3}{3}\km\W{3}{4}\km\W{1}{5}\km\W{3}{5}
+\frac{7}{6} \W{0}{1}\km\W{1}{1}\km\W{0}{3}\km\W{3}{3}\km\W{4}{4}\km\W{0}{5}\km\W{3}{5}
-3 \W{0}{1}\km\W{1}{1}\km\W{0}{3}\km\W{3}{3}\km\W{4}{4}\km\W{1}{5}\km\W{2}{5}
-3 \W{0}{1}\km\W{1}{1}\km\W{1}{3}\km\W{2}{3}\km\W{0}{4}\km\W{2}{5}\km\W{5}{5}
+\frac{21}{2} \W{0}{1}\km\W{1}{1}\km\W{1}{3}\km\W{2}{3}\km\W{0}{4}\km\W{3}{5}\km\W{4}{5}
+6 \W{0}{1}\km\W{1}{1}\km\W{1}{3}\km\W{2}{3}\km\W{1}{4}\km\W{1}{5}\km\W{5}{5}
\\&
-15 \W{0}{1}\km\W{1}{1}\km\W{1}{3}\km\W{2}{3}\km\W{1}{4}\km\W{2}{5}\km\W{4}{5}
-\frac{33}{10} \W{0}{1}\km\W{1}{1}\km\W{1}{3}\km\W{2}{3}\km\W{2}{4}\km\W{0}{5}\km\W{5}{5}
+\frac{51}{2} \W{0}{1}\km\W{1}{1}\km\W{1}{3}\km\W{2}{3}\km\W{2}{4}\km\W{2}{5}\km\W{3}{5}
+6 \W{0}{1}\km\W{1}{1}\km\W{1}{3}\km\W{2}{3}\km\W{3}{4}\km\W{0}{5}\km\W{4}{5}
-15 \W{0}{1}\km\W{1}{1}\km\W{1}{3}\km\W{2}{3}\km\W{3}{4}\km\W{1}{5}\km\W{3}{5}
-3 \W{0}{1}\km\W{1}{1}\km\W{1}{3}\km\W{2}{3}\km\W{4}{4}\km\W{0}{5}\km\W{3}{5}
\\&
+\frac{21}{2} \W{0}{1}\km\W{1}{1}\km\W{1}{3}\km\W{2}{3}\km\W{4}{4}\km\W{1}{5}\km\W{2}{5}
-\frac{1}{2} \W{0}{1}\km\W{1}{1}\km\W{1}{3}\km\W{3}{3}\km\W{0}{4}\km\W{2}{5}\km\W{4}{5}
+\frac{1}{5} \W{0}{1}\km\W{1}{1}\km\W{1}{3}\km\W{3}{3}\km\W{1}{4}\km\W{0}{5}\km\W{5}{5}
+3 \W{0}{1}\km\W{1}{1}\km\W{1}{3}\km\W{3}{3}\km\W{1}{4}\km\W{2}{5}\km\W{3}{5}
-2 \W{0}{1}\km\W{1}{1}\km\W{1}{3}\km\W{3}{3}\km\W{2}{4}\km\W{1}{5}\km\W{3}{5}
-\frac{1}{2} \W{0}{1}\km\W{1}{1}\km\W{1}{3}\km\W{3}{3}\km\W{2}{4}\km\W{0}{5}\km\W{4}{5}
\\&
+\W{0}{1}\km\W{1}{1}\km\W{1}{3}\km\W{3}{3}\km\W{3}{4}\km\W{0}{5}\km\W{3}{5}
+\W{0}{1}\km\W{1}{1}\km\W{1}{3}\km\W{3}{3}\km\W{3}{4}\km\W{1}{5}\km\W{2}{5}
-\frac{1}{2} \W{0}{1}\km\W{1}{1}\km\W{1}{3}\km\W{3}{3}\km\W{4}{4}\km\W{0}{5}\km\W{2}{5}
-\frac{1}{10} \W{0}{1}\km\W{1}{1}\km\W{2}{3}\km\W{3}{3}\km\W{0}{4}\km\W{0}{5}\km\W{5}{5}
-\frac{1}{2} \W{0}{1}\km\W{1}{1}\km\W{2}{3}\km\W{3}{3}\km\W{0}{4}\km\W{2}{5}\km\W{3}{5}
+\frac{1}{2} \W{0}{1}\km\W{1}{1}\km\W{2}{3}\km\W{3}{3}\km\W{0}{4}\km\W{1}{5}\km\W{4}{5}
\\&
-\W{0}{1}\km\W{1}{1}\km\W{2}{3}\km\W{3}{3}\km\W{1}{4}\km\W{1}{5}\km\W{3}{5}
+\frac{1}{2} \W{0}{1}\km\W{1}{1}\km\W{2}{3}\km\W{3}{3}\km\W{2}{4}\km\W{0}{5}\km\W{3}{5}
+\frac{3}{2} \W{0}{1}\km\W{1}{1}\km\W{2}{3}\km\W{3}{3}\km\W{2}{4}\km\W{1}{5}\km\W{2}{5}
-\W{0}{1}\km\W{1}{1}\km\W{2}{3}\km\W{3}{3}\km\W{3}{4}\km\W{0}{5}\km\W{2}{5}
+\frac{1}{2} \W{0}{1}\km\W{1}{1}\km\W{2}{3}\km\W{3}{3}\km\W{4}{4}\km\W{0}{5}\km\W{1}{5}
\\
\widehat{\mathbf{w}}_{10} &= 
\frac{1}{9} \W{0}{1}\km\W{1}{1}\km\W{0}{3}\km\W{0}{4}\km\W{1}{4}\km\W{4}{4}\km\W{5}{5}
-\frac{1}{9} \W{0}{1}\km\W{1}{1}\km\W{0}{3}\km\W{0}{4}\km\W{2}{4}\km\W{3}{4}\km\W{5}{5}
-\frac{2}{9} \W{0}{1}\km\W{1}{1}\km\W{0}{3}\km\W{0}{4}\km\W{2}{4}\km\W{4}{4}\km\W{4}{5}
+\frac{2}{9} \W{0}{1}\km\W{1}{1}\km\W{0}{3}\km\W{0}{4}\km\W{3}{4}\km\W{4}{4}\km\W{3}{5}
+\frac{4}{9} \W{0}{1}\km\W{1}{1}\km\W{0}{3}\km\W{1}{4}\km\W{2}{4}\km\W{3}{4}\km\W{4}{5}
+\frac{2}{9} \W{0}{1}\km\W{1}{1}\km\W{0}{3}\km\W{1}{4}\km\W{2}{4}\km\W{4}{4}\km\W{3}{5}
\\&
-\frac{4}{9} \W{0}{1}\km\W{1}{1}\km\W{0}{3}\km\W{1}{4}\km\W{3}{4}\km\W{4}{4}\km\W{2}{5}
+\frac{1}{3} \W{0}{1}\km\W{1}{1}\km\W{0}{3}\km\W{2}{4}\km\W{3}{4}\km\W{4}{4}\km\W{1}{5}
-\frac{2}{9} \W{0}{1}\km\W{1}{1}\km\W{1}{3}\km\W{0}{4}\km\W{1}{4}\km\W{3}{4}\km\W{5}{5}
-\frac{1}{9} \W{0}{1}\km\W{1}{1}\km\W{1}{3}\km\W{0}{4}\km\W{1}{4}\km\W{4}{4}\km\W{4}{5}
+ W{0}{1}\km\W{1}{1}\km\W{1}{3}\km\W{0}{4}\km\W{2}{4}\km\W{3}{4}\km\W{4}{5}
-\frac{2}{9} \W{0}{1}\km\W{1}{1}\km\W{1}{3}\km\W{0}{4}\km\W{3}{4}\km\W{4}{4}\km\W{2}{5}
\\&
-\frac{8}{3} \W{0}{1}\km\W{1}{1}\km\W{1}{3}\km\W{1}{4}\km\W{2}{4}\km\W{3}{4}\km\W{3}{5}
+\frac{2}{3} \W{0}{1}\km\W{1}{1}\km\W{1}{3}\km\W{1}{4}\km\W{2}{4}\km\W{4}{4}\km\W{2}{5}
+\frac{2}{9} \W{0}{1}\km\W{1}{1}\km\W{1}{3}\km\W{1}{4}\km\W{3}{4}\km\W{4}{4}\km\W{1}{5}
-\frac{1}{3} \W{0}{1}\km\W{1}{1}\km\W{1}{3}\km\W{2}{4}\km\W{3}{4}\km\W{4}{4}\km\W{0}{5}
+\frac{1}{3} \W{0}{1}\km\W{1}{1}\km\W{2}{3}\km\W{0}{4}\km\W{1}{4}\km\W{2}{4}\km\W{5}{5}
-\frac{2}{9} \W{0}{1}\km\W{1}{1}\km\W{2}{3}\km\W{0}{4}\km\W{1}{4}\km\W{3}{4}\km\W{4}{5}
\\&
+\frac{2}{9} \W{0}{1}\km\W{1}{1}\km\W{2}{3}\km\W{0}{4}\km\W{1}{4}\km\W{4}{4}\km\W{3}{5}
-\frac{2}{3} \W{0}{1}\km\W{1}{1}\km\W{2}{3}\km\W{0}{4}\km\W{2}{4}\km\W{3}{4}\km\W{3}{5}
+\frac{1}{9} \W{0}{1}\km\W{1}{1}\km\W{2}{3}\km\W{0}{4}\km\W{3}{4}\km\W{4}{4}\km\W{1}{5}
-\W{0}{1}\km\W{1}{1}\km\W{2}{3}\km\W{1}{4}\km\W{2}{4}\km\W{4}{4}\km\W{1}{5}
+\frac{8}{3} \W{0}{1}\km\W{1}{1}\km\W{2}{3}\km\W{1}{4}\km\W{2}{4}\km\W{3}{4}\km\W{2}{5}
+\frac{2}{9} \W{0}{1}\km\W{1}{1}\km\W{2}{3}\km\W{1}{4}\km\W{3}{4}\km\W{4}{4}\km\W{0}{5}
\\&
-\frac{1}{3} \W{0}{1}\km\W{1}{1}\km\W{3}{3}\km\W{0}{4}\km\W{1}{4}\km\W{2}{4}\km\W{4}{5}
+\frac{4}{9} \W{0}{1}\km\W{1}{1}\km\W{3}{3}\km\W{0}{4}\km\W{1}{4}\km\W{3}{4}\km\W{3}{5}
-\frac{2}{9} \W{0}{1}\km\W{1}{1}\km\W{3}{3}\km\W{0}{4}\km\W{1}{4}\km\W{4}{4}\km\W{2}{5}
-\frac{2}{9} \W{0}{1}\km\W{1}{1}\km\W{3}{3}\km\W{0}{4}\km\W{2}{4}\km\W{3}{4}\km\W{2}{5}
-\frac{4}{9} \W{0}{1}\km\W{1}{1}\km\W{3}{3}\km\W{1}{4}\km\W{2}{4}\km\W{3}{4}\km\W{1}{5}
+\frac{2}{9} \W{0}{1}\km\W{1}{1}\km\W{3}{3}\km\W{0}{4}\km\W{2}{4}\km\W{4}{4}\km\W{1}{5}
\\&
-\frac{1}{9} \W{0}{1}\km\W{1}{1}\km\W{3}{3}\km\W{0}{4}\km\W{3}{4}\km\W{4}{4}\km\W{0}{5}
+\frac{1}{9} \W{0}{1}\km\W{1}{1}\km\W{3}{3}\km\W{1}{4}\km\W{2}{4}\km\W{4}{4}\km\W{0}{5} 
\\
\widehat{\mathbf{w}}_{11} &= 
\frac{1}{9} \W{0}{1}\km\W{1}{1}\km\W{0}{3}\km\W{0}{4}\km\W{1}{4}\km\W{4}{4}\km\W{5}{5}
-\frac{4}{9} \W{0}{1}\km\W{1}{1}\km\W{0}{3}\km\W{0}{4}\km\W{2}{4}\km\W{3}{4}\km\W{5}{5}
+\frac{1}{9} \W{0}{1}\km\W{1}{1}\km\W{0}{3}\km\W{0}{4}\km\W{2}{4}\km\W{4}{4}\km\W{4}{5}
-\frac{1}{9} \W{0}{1}\km\W{1}{1}\km\W{0}{3}\km\W{0}{4}\km\W{3}{4}\km\W{4}{4}\km\W{3}{5}
+\frac{16}{9} \W{0}{1}\km\W{1}{1}\km\W{0}{3}\km\W{1}{4}\km\W{2}{4}\km\W{3}{4}\km\W{4}{5}
-\frac{10}{9} \W{0}{1}\km\W{1}{1}\km\W{0}{3}\km\W{1}{4}\km\W{2}{4}\km\W{4}{4}\km\W{3}{5}
\\&
+\frac{8}{9} \W{0}{1}\km\W{1}{1}\km\W{0}{3}\km\W{1}{4}\km\W{3}{4}\km\W{4}{4}\km\W{2}{5}
-\frac{2}{3} \W{0}{1}\km\W{1}{1}\km\W{0}{3}\km\W{2}{4}\km\W{3}{4}\km\W{4}{4}\km\W{1}{5}
+\frac{4}{9} \W{0}{1}\km\W{1}{1}\km\W{1}{3}\km\W{0}{4}\km\W{1}{4}\km\W{3}{4}\km\W{5}{5}
-\frac{7}{9} \W{0}{1}\km\W{1}{1}\km\W{1}{3}\km\W{0}{4}\km\W{1}{4}\km\W{4}{4}\km\W{4}{5}
+\W{0}{1}\km\W{1}{1}\km\W{1}{3}\km\W{0}{4}\km\W{2}{4}\km\W{4}{4}\km\W{3}{5}
-\frac{5}{9} \W{0}{1}\km\W{1}{1}\km\W{1}{3}\km\W{0}{4}\km\W{3}{4}\km\W{4}{4}\km\W{2}{5}
\\&
-\frac{8}{3} \W{0}{1}\km\W{1}{1}\km\W{1}{3}\km\W{1}{4}\km\W{2}{4}\km\W{3}{4}\km\W{3}{5}
+\frac{2}{3} \W{0}{1}\km\W{1}{1}\km\W{1}{3}\km\W{1}{4}\km\W{2}{4}\km\W{4}{4}\km\W{2}{5}
-\frac{4}{9} \W{0}{1}\km\W{1}{1}\km\W{1}{3}\km\W{1}{4}\km\W{3}{4}\km\W{4}{4}\km\W{1}{5}
+\frac{2}{3} \W{0}{1}\km\W{1}{1}\km\W{1}{3}\km\W{2}{4}\km\W{3}{4}\km\W{4}{4}\km\W{0}{5}
-\frac{2}{3} \W{0}{1}\km\W{1}{1}\km\W{2}{3}\km\W{0}{4}\km\W{1}{4}\km\W{2}{4}\km\W{5}{5}
+\frac{4}{9} \W{0}{1}\km\W{1}{1}\km\W{2}{3}\km\W{0}{4}\km\W{1}{4}\km\W{3}{4}\km\W{4}{5}
\\&
+\frac{5}{9} \W{0}{1}\km\W{1}{1}\km\W{2}{3}\km\W{0}{4}\km\W{1}{4}\km\W{4}{4}\km\W{3}{5}
-\frac{2}{3} \W{0}{1}\km\W{1}{1}\km\W{2}{3}\km\W{0}{4}\km\W{2}{4}\km\W{3}{4}\km\W{3}{5}
-\W{0}{1}\km\W{1}{1}\km\W{2}{3}\km\W{0}{4}\km\W{2}{4}\km\W{4}{4}\km\W{2}{5} 
+\frac{7}{9} \W{0}{1}\km\W{1}{1}\km\W{2}{3}\km\W{0}{4}\km\W{3}{4}\km\W{4}{4}\km\W{1}{5}
+\frac{8}{3} \W{0}{1}\km\W{1}{1}\km\W{2}{3}\km\W{1}{4}\km\W{2}{4}\km\W{3}{4}\km\W{2}{5}
-\frac{4}{9} \W{0}{1}\km\W{1}{1}\km\W{2}{3}\km\W{1}{4}\km\W{3}{4}\km\W{4}{4}\km\W{0}{5}
\\&
+\frac{2}{3} \W{0}{1}\km\W{1}{1}\km\W{3}{3}\km\W{0}{4}\km\W{1}{4}\km\W{2}{4}\km\W{4}{5}
-\frac{8}{9} \W{0}{1}\km\W{1}{1}\km\W{3}{3}\km\W{0}{4}\km\W{1}{4}\km\W{3}{4}\km\W{3}{5}
+\frac{1}{9} \W{0}{1}\km\W{1}{1}\km\W{3}{3}\km\W{0}{4}\km\W{1}{4}\km\W{4}{4}\km\W{2}{5}
+\frac{10}{9} \W{0}{1}\km\W{1}{1}\km\W{3}{3}\km\W{0}{4}\km\W{2}{4}\km\W{3}{4}\km\W{2}{5}
-\frac{16}{9} \W{0}{1}\km\W{1}{1}\km\W{3}{3}\km\W{1}{4}\km\W{2}{4}\km\W{3}{4}\km\W{1}{5}
-\frac{1}{9} \W{0}{1}\km\W{1}{1}\km\W{3}{3}\km\W{0}{4}\km\W{2}{4}\km\W{4}{4}\km\W{1}{5}
\\&
-\frac{1}{9} \W{0}{1}\km\W{1}{1}\km\W{3}{3}\km\W{0}{4}\km\W{3}{4}\km\W{4}{4}\km\W{0}{5}
+\frac{4}{9} \W{0}{1}\km\W{1}{1}\km\W{3}{3}\km\W{1}{4}\km\W{2}{4}\km\W{4}{4}\km\W{0}{5}
\\
\widehat{\mathbf{w}}_{12} &= 
\W{0}{1}\km\W{1}{1}\km\W{0}{4}\km\W{1}{4}\km\W{2}{4}\km\W{3}{4}\km \W{4}{4}
\\ 
\widehat{\mathbf{w}}_{13} &= 
\W{1}{1}\km\W{0}{3}\km\W{1}{3}\km\W{2}{3}\km\W{3}{3}\km\W{3}{4}\km\W{1}{5}
+\W{1}{1}\km\W{0}{3}\km\W{1}{3}\km\W{2}{3}\km\W{3}{3}\km\W{1}{4}\km\W{3}{5}
+\frac{1}{4} \W{0}{1}\km\W{0}{3}\km\W{1}{3}\km\W{2}{3}\km\W{3}{3}\km\W{0}{4}\km\W{5}{5}
-\W{0}{1}\km\W{0}{3}\km\W{1}{3}\km\W{2}{3}\km\W{3}{3}\km\W{1}{4}\km\W{4}{5}
+\frac{3}{2} \W{0}{1}\km\W{0}{3}\km\W{1}{3}\km\W{2}{3}\km\W{3}{3}\km\W{2}{4}\km\W{3}{5}
-\W{0}{1}\km\W{0}{3}\km\W{1}{3}\km\W{2}{3}\km\W{3}{3}\km\W{3}{4}\km\W{2}{5}
\\&
+\frac{1}{4} \W{0}{1}\km\W{0}{3}\km\W{1}{3}\km\W{2}{3}\km\W{3}{3}\km\W{4}{4}\km\W{1}{5}
-\frac{1}{4} \W{1}{1}\km\W{0}{3}\km\W{1}{3}\km\W{2}{3}\km\W{3}{3}\km\W{0}{4}\km\W{4}{5}
-\frac{3}{2} \W{1}{1}\km\W{0}{3}\km\W{1}{3}\km\W{2}{3}\km\W{3}{3}\km\W{2}{4}\km\W{2}{5}
-\frac{1}{4} \W{1}{1}\km\W{0}{3}\km\W{1}{3}\km\W{2}{3}\km\W{3}{3}\km\W{4}{4}\km\W{0}{5}
\\
\widehat{\mathbf{w}}_{14} &= 
\W{1}{1}\km\W{0}{3}\km\W{2}{3}\km\W{3}{3}\km\W{0}{4}\km\W{2}{4}\km\W{3}{4}
+\frac{1}{2} \W{0}{1}\km\W{0}{3}\km\W{1}{3}\km\W{3}{3}\km\W{0}{4}\km\W{3}{4}\km\W{4}{4}
-\W{0}{1}\km\W{0}{3}\km\W{1}{3}\km\W{3}{3}\km\W{1}{4}\km\W{2}{4}\km\W{4}{4}
-\frac{1}{2} \W{0}{1}\km\W{0}{3}\km\W{2}{3}\km\W{3}{3}\km\W{0}{4}\km\W{2}{4}\km\W{4}{4}
+4 \W{0}{1}\km\W{0}{3}\km\W{2}{3}\km\W{3}{3}\km\W{1}{4}\km\W{2}{4}\km\W{3}{4}
+\frac{3}{2} \W{0}{1}\km\W{1}{3}\km\W{2}{3}\km\W{3}{3}\km\W{0}{4}\km\W{1}{4}\km\W{4}{4}
\\&
-3 \W{0}{1}\km\W{1}{3}\km\W{2}{3}\km\W{3}{3}\km\W{0}{4}\km\W{2}{4}\km\W{3}{4}
-\frac{3}{2} \W{1}{1}\km\W{0}{3}\km\W{1}{3}\km\W{2}{3}\km\W{0}{4}\km\W{3}{4}\km\W{4}{4}
+3 \W{1}{1}\km\W{0}{3}\km\W{1}{3}\km\W{2}{3}\km\W{1}{4}\km\W{2}{4}\km\W{4}{4}
+\frac{1}{2} \W{1}{1}\km\W{0}{3}\km\W{1}{3}\km\W{3}{3}\km\W{0}{4}\km\W{2}{4}\km\W{4}{4}
-4 \W{1}{1}\km\W{0}{3}\km\W{1}{3}\km\W{3}{3}\km\W{1}{4}\km\W{2}{4}\km\W{3}{4}
-\frac{1}{2} \W{1}{1}\km\W{0}{3}\km\W{2}{3}\km\W{3}{3}\km\W{0}{4}\km\W{1}{4}\km\W{4}{4}
\end{align*}

\subsection{Basis of $\ds\CGF{8}{2}{}{8}$}

\begin{align*}
\widehat{\mathbf{r}_1} =& 
-\frac{1}{5} \W{0}{1}\km\W{1}{1}\km\W{0}{3}\km\W{1}{3}\km\W{2}{3}\km\W{3}{3}\km\W{0}{5}\km\W{5}{5}
+\W{0}{1}\km\W{1}{1}\km\W{0}{3}\km\W{1}{3}\km\W{2}{3}\km\W{3}{3}\km\W{1}{5}\km\W{4}{5} 
-2 \W{0}{1}\km\W{1}{1}\km\W{0}{3}\km\W{1}{3}\km\W{2}{3}\km\W{3}{3}\km\W{2}{5}\km\W{3}{5}
\\
\widehat{\mathbf{r}_2} =& 
-\frac{1}{2} \W{0}{1}\km\W{1}{1}\km\W{0}{3}\km\W{1}{3}\km\W{2}{3}\km\W{0}{4}\km\W{3}{4}\km\W{5}{5}
+\frac{1}{2} \W{0}{1}\km\W{1}{1}\km\W{0}{3}\km\W{1}{3}\km\W{2}{3}\km\W{0}{4}\km\W{4}{4}\km\W{4}{5} 
+3 \W{0}{1}\km\W{1}{1}\km\W{0}{3}\km\W{1}{3}\km\W{2}{3}\km\W{1}{4}\km\W{2}{4}\km\W{5}{5} 
-4 \W{0}{1}\km\W{1}{1}\km\W{0}{3}\km\W{1}{3}\km\W{2}{3}\km\W{1}{4}\km\W{3}{4}\km\W{4}{5} 
+\W{0}{1}\km\W{1}{1}\km\W{0}{3}\km\W{1}{3}\km\W{2}{3}\km\W{1}{4}\km\W{4}{4}\km\W{3}{5} 
\\& 
+6 \W{0}{1}\km\W{1}{1}\km\W{0}{3}\km\W{1}{3}\km\W{2}{3}\km\W{2}{4}\km\W{3}{4}\km\W{3}{5} 
-3 \W{0}{1}\km\W{1}{1}\km\W{0}{3}\km\W{1}{3}\km\W{2}{3}\km\W{2}{4}\km\W{4}{4}\km\W{2}{5}
+\frac{3}{2} \W{0}{1}\km\W{1}{1}\km\W{0}{3}\km\W{1}{3}\km\W{2}{3}\km\W{3}{4}\km\W{4}{4}\km\W{1}{5} 
+\frac{3}{2} \W{0}{1}\km\W{1}{1}\km\W{0}{3}\km\W{1}{3}\km\W{3}{3}\km\W{0}{4}\km\W{3}{4}\km\W{4}{5} 
-\W{0}{1}\km\W{1}{1}\km\W{0}{3}\km\W{1}{3}\km\W{3}{3}\km\W{0}{4}\km\W{4}{4}\km\W{3}{5} 
\\&
-\frac{1}{2} \W{0}{1}\km\W{1}{1}\km\W{0}{3}\km\W{1}{3}\km\W{3}{3}\km\W{0}{4}\km\W{2}{4}\km\W{5}{5}
-\W{0}{1}\km\W{1}{1}\km\W{0}{3}\km\W{1}{3}\km\W{3}{3}\km\W{1}{4}\km\W{2}{4}\km\W{4}{5} 
+\W{0}{1}\km\W{1}{1}\km\W{0}{3}\km\W{1}{3}\km\W{3}{3}\km\W{1}{4}\km\W{4}{4}\km\W{2}{5} 
-2 \W{0}{1}\km\W{1}{1}\km\W{0}{3}\km\W{1}{3}\km\W{3}{3}\km\W{2}{4}\km\W{3}{4}\km\W{2}{5} 
+\frac{1}{2} \W{0}{1}\km\W{1}{1}\km\W{0}{3}\km\W{1}{3}\km\W{3}{3}\km\W{2}{4}\km\W{4}{4}\km\W{1}{5} 
\\&
-\frac{1}{2} \W{0}{1}\km\W{1}{1}\km\W{0}{3}\km\W{1}{3}\km\W{3}{3}\km\W{3}{4}\km\W{4}{4}\km\W{0}{5} 
+\frac{1}{2} \W{0}{1}\km\W{1}{1}\km\W{0}{3}\km\W{2}{3}\km\W{3}{3}\km\W{0}{4}\km\W{1}{4}\km\W{5}{5} 
-\frac{1}{2} \W{0}{1}\km\W{1}{1}\km\W{0}{3}\km\W{2}{3}\km\W{3}{3}\km\W{0}{4}\km\W{2}{4}\km\W{4}{5} 
-\W{0}{1}\km\W{1}{1}\km\W{0}{3}\km\W{2}{3}\km\W{3}{3}\km\W{0}{4}\km\W{3}{4}\km\W{3}{5}
+2 \W{0}{1}\km\W{1}{1}\km\W{0}{3}\km\W{2}{3}\km\W{3}{3}\km\W{1}{4}\km\W{2}{4}\km\W{3}{5} 
\\&
+\W{0}{1}\km\W{1}{1}\km\W{0}{3}\km\W{2}{3}\km\W{3}{3}\km\W{0}{4}\km\W{4}{4}\km\W{2}{5} 
-\frac{3}{2} \W{0}{1}\km\W{1}{1}\km\W{0}{3}\km\W{2}{3}\km\W{3}{3}\km\W{1}{4}\km\W{4}{4}\km\W{1}{5} 
+\W{0}{1}\km\W{1}{1}\km\W{0}{3}\km\W{2}{3}\km\W{3}{3}\km\W{2}{4}\km\W{3}{4}\km\W{1}{5}
+\frac{1}{2} \W{0}{1}\km\W{1}{1}\km\W{0}{3}\km\W{2}{3}\km\W{3}{3}\km\W{2}{4}\km\W{4}{4}\km\W{0}{5}
-\W{0}{1}\km\W{1}{1}\km\W{1}{3}\km\W{2}{3}\km\W{3}{3}\km\W{0}{4}\km\W{3}{4}\km\W{2}{5} 
\\&
-\frac{3}{2} \W{0}{1}\km\W{1}{1}\km\W{1}{3}\km\W{2}{3}\km\W{3}{3}\km\W{0}{4}\km\W{1}{4}\km\W{4}{5} 
+3 \W{0}{1}\km\W{1}{1}\km\W{1}{3}\km\W{2}{3}\km\W{3}{3}\km\W{0}{4}\km\W{2}{4}\km\W{3}{5} 
-6 \W{0}{1}\km\W{1}{1}\km\W{1}{3}\km\W{2}{3}\km\W{3}{3}\km\W{1}{4}\km\W{2}{4}\km\W{2}{5} 
-\frac{1}{2} \W{0}{1}\km\W{1}{1}\km\W{1}{3}\km\W{2}{3}\km\W{3}{3}\km\W{0}{4}\km\W{4}{4}\km\W{1}{5}
+4 \W{0}{1}\km\W{1}{1}\km\W{1}{3}\km\W{2}{3}\km\W{3}{3}\km\W{1}{4}\km\W{3}{4}\km\W{1}{5} 
\\&
+\frac{1}{2} \W{0}{1}\km\W{1}{1}\km\W{1}{3}\km\W{2}{3}\km\W{3}{3}\km\W{1}{4}\km\W{4}{4}\km\W{0}{5}
-3 \W{0}{1}\km\W{1}{1}\km\W{1}{3}\km\W{2}{3}\km\W{3}{3}\km\W{2}{4}\km\W{3}{4}\km\W{0}{5} 
\\
\widehat{\mathbf{r}_3} =& 
-\frac{3}{4} \W{0}{1}\km\W{1}{1}\km\W{0}{3}\km\W{1}{3}\km\W{2}{3}\km\W{0}{4}\km\W{3}{4}\km\W{5}{5}
+\frac{3}{4} \W{0}{1}\km\W{1}{1}\km\W{0}{3}\km\W{1}{3}\km\W{2}{3}\km\W{0}{4}\km\W{4}{4}\km\W{4}{5}
+3 \W{0}{1}\km\W{1}{1}\km\W{0}{3}\km\W{1}{3}\km\W{2}{3}\km\W{1}{4}\km\W{2}{4}\km\W{5}{5} 
-3 \W{0}{1}\km\W{1}{1}\km\W{0}{3}\km\W{1}{3}\km\W{2}{3}\km\W{1}{4}\km\W{3}{4}\km\W{4}{5}
+\frac{9}{2} \W{0}{1}\km\W{1}{1}\km\W{0}{3}\km\W{1}{3}\km\W{2}{3}\km\W{2}{4}\km\W{3}{4}\km\W{3}{5}
\\&
-\frac{3}{2} \W{0}{1}\km\W{1}{1}\km\W{0}{3}\km\W{1}{3}\km\W{2}{3}\km\W{2}{4}\km\W{4}{4}\km\W{2}{5} 
+\frac{3}{4} \W{0}{1}\km\W{1}{1}\km\W{0}{3}\km\W{1}{3}\km\W{2}{3}\km\W{3}{4}\km\W{4}{4}\km\W{1}{5}
+\frac{5}{4} \W{0}{1}\km\W{1}{1}\km\W{0}{3}\km\W{1}{3}\km\W{3}{3}\km\W{0}{4}\km\W{3}{4}\km\W{4}{5}
-\W{0}{1}\km\W{1}{1}\km\W{0}{3}\km\W{1}{3}\km\W{3}{3}\km\W{0}{4}\km\W{4}{4}\km\W{3}{5} 
-\frac{1}{4} \W{0}{1}\km\W{1}{1}\km\W{0}{3}\km\W{1}{3}\km\W{3}{3}\km\W{0}{4}\km\W{2}{4}\km\W{5}{5}
\\&
-2 \W{0}{1}\km\W{1}{1}\km\W{0}{3}\km\W{1}{3}\km\W{3}{3}\km\W{1}{4}\km\W{2}{4}\km\W{4}{5}
+\W{0}{1}\km\W{1}{1}\km\W{0}{3}\km\W{1}{3}\km\W{3}{3}\km\W{1}{4}\km\W{3}{4}\km\W{3}{5}
+\W{0}{1}\km\W{1}{1}\km\W{0}{3}\km\W{1}{3}\km\W{3}{3}\km\W{1}{4}\km\W{4}{4}\km\W{2}{5}
-\frac{5}{2} \W{0}{1}\km\W{1}{1}\km\W{0}{3}\km\W{1}{3}\km\W{3}{3}\km\W{2}{4}\km\W{3}{4}\km\W{2}{5}
+\frac{1}{4} \W{0}{1}\km\W{1}{1}\km\W{0}{3}\km\W{1}{3}\km\W{3}{3}\km\W{2}{4}\km\W{4}{4}\km\W{1}{5} 
\\&
-\frac{1}{4} \W{0}{1}\km\W{1}{1}\km\W{0}{3}\km\W{1}{3}\km\W{3}{3}\km\W{3}{4}\km\W{4}{4}\km\W{0}{5}
+\frac{1}{4} \W{0}{1}\km\W{1}{1}\km\W{0}{3}\km\W{2}{3}\km\W{3}{3}\km\W{0}{4}\km\W{1}{4}\km\W{5}{5}
-\frac{1}{4} \W{0}{1}\km\W{1}{1}\km\W{0}{3}\km\W{2}{3}\km\W{3}{3}\km\W{0}{4}\km\W{2}{4}\km\W{4}{5} 
-\W{0}{1}\km\W{1}{1}\km\W{0}{3}\km\W{2}{3}\km\W{3}{3}\km\W{0}{4}\km\W{3}{4}\km\W{3}{5} 
+\frac{5}{2} \W{0}{1}\km\W{1}{1}\km\W{0}{3}\km\W{2}{3}\km\W{3}{3}\km\W{1}{4}\km\W{2}{4}\km\W{3}{5}
\\&
+\W{0}{1}\km\W{1}{1}\km\W{0}{3}\km\W{2}{3}\km\W{3}{3}\km\W{0}{4}\km\W{4}{4}\km\W{2}{5} 
-\W{0}{1}\km\W{1}{1}\km\W{0}{3}\km\W{2}{3}\km\W{3}{3}\km\W{1}{4}\km\W{3}{4}\km\W{2}{5}
-\frac{5}{4} \W{0}{1}\km\W{1}{1}\km\W{0}{3}\km\W{2}{3}\km\W{3}{3}\km\W{1}{4}\km\W{4}{4}\km\W{1}{5} 
+2 \W{0}{1}\km\W{1}{1}\km\W{0}{3}\km\W{2}{3}\km\W{3}{3}\km\W{2}{4}\km\W{3}{4}\km\W{1}{5} 
+\frac{1}{4} \W{0}{1}\km\W{1}{1}\km\W{0}{3}\km\W{2}{3}\km\W{3}{3}\km\W{2}{4}\km\W{4}{4}\km\W{0}{5}
\\&
-\frac{3}{4} \W{0}{1}\km\W{1}{1}\km\W{1}{3}\km\W{2}{3}\km\W{3}{3}\km\W{0}{4}\km\W{1}{4}\km\W{4}{5}
+\frac{3}{2} \W{0}{1}\km\W{1}{1}\km\W{1}{3}\km\W{2}{3}\km\W{3}{3}\km\W{0}{4}\km\W{2}{4}\km\W{3}{5} 
-\frac{9}{2} \W{0}{1}\km\W{1}{1}\km\W{1}{3}\km\W{2}{3}\km\W{3}{3}\km\W{1}{4}\km\W{2}{4}\km\W{2}{5}
-\frac{3}{4} \W{0}{1}\km\W{1}{1}\km\W{1}{3}\km\W{2}{3}\km\W{3}{3}\km\W{0}{4}\km\W{4}{4}\km\W{1}{5}
+3 \W{0}{1}\km\W{1}{1}\km\W{1}{3}\km\W{2}{3}\km\W{3}{3}\km\W{1}{4}\km\W{3}{4}\km\W{1}{5} 
\\&
+\frac{3}{4} \W{0}{1}\km\W{1}{1}\km\W{1}{3}\km\W{2}{3}\km\W{3}{3}\km\W{1}{4}\km\W{4}{4}\km\W{0}{5} 
-3 \W{0}{1}\km\W{1}{1}\km\W{1}{3}\km\W{2}{3}\km\W{3}{3}\km\W{2}{4}\km\W{3}{4}\km\W{0}{5} 
\\
\widehat{\mathbf{r}_4} =& 
+ \W{0}{1}\km\W{1}{1}\km\W{1}{3}\km\W{2}{3}\km\W{0}{4}\km\W{1}{4}\km\W{3}{4}\km\W{4}{4}
+2 \W{0}{1}\km\W{1}{1}\km\W{0}{3}\km\W{1}{3}\km\W{1}{4}\km\W{2}{4}\km\W{3}{4}\km\W{4}{4}
-\W{0}{1}\km\W{1}{1}\km\W{0}{3}\km\W{2}{3}\km\W{0}{4}\km\W{2}{4}\km\W{3}{4}\km\W{4}{4} 
+\frac{1}{3}\W{0}{1}\km\W{1}{1}\km\W{0}{3}\km\W{3}{3}\km\W{0}{4}\km\W{1}{4}\km\W{3}{4}\km\W{4}{4}
-\W{0}{1}\km\W{1}{1}\km\W{1}{3}\km\W{3}{3}\km\W{0}{4}\km\W{1}{4}\km\W{2}{4}\km\W{4}{4} 
\\& 
+2 \W{0}{1}\km\W{1}{1}\km\W{2}{3}\km\W{3}{3}\km\W{0}{4}\km\W{1}{4}\km\W{2}{4}\km\W{3}{4}
\end{align*}

\end{scriptsize}

For simplicity, we denote $\displaystyle \CGF{\bullet}{2}{}{8}$ by   
$\displaystyle \frakC^{\bullet}$.


For instance, the matrix representations of 
$\displaystyle \mydz : \frakC^{6} \rightarrow \frakC^{7}$ and 
$\displaystyle \mydz : \frakC^{7} \rightarrow \frakC^{8}$ are 
given as follows: 
\begingroup
\renewcommand{\arraystretch}{1.7} 
\begin{small}
\setcounter{MaxMatrixCols}{20}
\begin{equation}
{}^t 
\overline{M} = \begin{bmatrix} 
\dfrac{135}{4} & 0 & -20 & -60 & 15 &\dfrac{5}{4} & -\dfrac{135}{4} & -15 & -15 & 0 & 0 & 0 & 0 & 0 \\ 
-18 & -12 & 0 & 0 & 0 &\dfrac{23}{3} & -45 & 21 & 10 & 0 & 0 & 0 & 0 & 0 \\ 
-\dfrac{27}{4} & 0 & -20 & 12 & 3 & 0 & 0 & 0 & 0 &\dfrac{63}{8} & -18 & 0 & 0 & 0 \\ 
0 & 0 & 2 & -10 &\dfrac{2}{3} & 1 & 2 & -\dfrac{1}{3} & 2 & -7 & 1 & 0 & 0 & 0 \\ 
0 & -\dfrac{11}{2} & -17 & 9 &\dfrac{71}{12} &\dfrac{7}{8} &\dfrac{51}{8} &\dfrac{233}{12} & 10 & -\dfrac{49}{4} & 4 & 0 & 0 & 0 \\ 
0 & -\dfrac{1}{3} & -9 & 5 &\dfrac{49}{18} &\dfrac{5}{12} &\dfrac{55}{4} &\dfrac{47}{9} &\dfrac{19}{3} & -\dfrac{28}{3} & -\dfrac{1}{6} & 0 & 0 & 0 \\ 
0 & -\dfrac{5}{2} & -20 & 0 &\dfrac{65}{12} & -\dfrac{15}{8} & -\dfrac{75}{8} &\dfrac{65}{12} & -5 &\dfrac{35}{4} & -\dfrac{25}{2} & 0 & 0 & 0 \\ 
0 & 0 & 0 & 0 & 0 & 7 & -18 & -9 & 6 & 0 & 0 & 0 & 0 & 0 \\ 
0 & 0 & -2 & -6 & 0 & 0 & 0 & 0 & 0 & -\dfrac{21}{2} & 6 & 70 & 0 & 0 \\ 
-\dfrac{1}{2} & 2 & -2 &\dfrac{14}{3} &\dfrac{1}{3} & 0 & 0 &  0 & 0 & 0 & 0 & 0 & -28 & -\dfrac{14}{3} \\ 
-\dfrac{5}{2} & 2  & 0 &\dfrac{16}{3} & -\dfrac{1}{3} & 0 & 0 & 0 & 0 & 0 & 0 & 0 & -28 & -\dfrac{28}{3} \\ 
0 & 8 & 0 & 0 & 0 & -3 & 0 & 1 & -6 & 0 & 0 & 0 & -112 & 0 \\
0 & 0 & 0 & 0 & 0 & -1 & 15 & -6 & 0 & -3 /2 & -3 & 0 & 48 & -14 \\ 
0 & 0 & 0 & 0 & 0 & -\dfrac{2}{3} & -1 & -1 & -2 &\dfrac{5}{2} & -1 & 0 & 4 & -\dfrac{14}{3} \\ 
0 & 0 & 0 & 0 &\dfrac{4}{3} &\dfrac{11}{3} & -45 &\dfrac{46}{3} & 0 &\dfrac{5}{2} & 11 & 0 & -136 & 14 \\ 
0 & 0 & 2 & 0 & -1 & -\dfrac{3}{2} &\dfrac{123}{8} & -4 & 0 & -\dfrac{21}{16} & -\dfrac{33}{16} & 0 & 42 & -\dfrac{21}{8} \\ 
0 & 0 & 0 & 2 & -\dfrac{1}{3} & -\dfrac{1}{2} &\dfrac{39}{8} & -\dfrac{4}{3} & 0 &\dfrac{7}{16} & -\dfrac{13}{16} & 0 & 14 & -\dfrac{35}{24} \\ 
0 & 0 & 0 & 0 & 0 & 0 & 0 & 0 & 0 & 1 & -1 & -\dfrac{20}{3} & 0 & 2 
\end{bmatrix}
\end{equation}
\end{small}
\endgroup
and 
\begingroup
\renewcommand{\arraystretch}{1.7} 
\setcounter{MaxMatrixCols}{20}
\begin{equation} 
\overline{N} = \begin{bmatrix} 0 & -35 & 0 & 0 & 0 & -30 & 0 & -15
&\dfrac{25}{2} & 0 & 0 & 0 & -\dfrac{5}{2} & 0 \\ 
11 & -9 & -\dfrac{39}{8} & -\dfrac{31}{8} & -\dfrac{61}{2} & -75 & -10 & -10
&\dfrac{85}{2} & -\dfrac{2}{3} & -\dfrac{20}{3} & 0 & -1 & -3 \\ 
-16 & 8 & 9 & 5 & 52 & 20 &\dfrac{8}{3} & -2 & -\dfrac{55}{3} & 0 & 8 & 0 & 1 &
4 \\ 
0 & 0 &\dfrac{63}{2} &\dfrac{21}{2} & 84 & 0 & 0 & 0 & 0 & 0 & -14 & 3 & 0 & 3 
\end{bmatrix} 
\end{equation}
\endgroup

Since $\rank \overline{M} = 9$ and $\rank \overline{N} = 4$, we see the
dimensions of $\mydz ( \frakC^{6})$ and $\ker( \mydz : \frakC^{7} \rightarrow
C^{6})$ , and so on. The precise data of the structures of $\displaystyle
\CGF{\bullet}{2}{}{8}$ and  $\displaystyle \HGF{\bullet}{2}{}{8}$ are in the
table below.

\begin{center}
\tabcolsep=4.5pt
\begin{tabular}{|c|*{14}{c}c|}
\hline
$\displaystyle \frakham^{ }_{2}$, w=8
& $\mathbf{0}$ & $\rightarrow$ &$\displaystyle \frakC^{3}$ & $\rightarrow$ &
$\displaystyle \frakC^{4} $& $ \rightarrow $ &
$\displaystyle \frakC^{5} $& $ \rightarrow $ &
$\displaystyle \frakC^{6} $& $ \rightarrow $ &
$\displaystyle \frakC^{7} $& $ \rightarrow $ &
$\displaystyle \frakC^{8} $& $ \rightarrow $ &
$\mathbf{0}$ \\ \hline
$\dim$  &   && 5 && 13 && 17 && 18 && 14&& 4 && \\
rank        & & 0 && 5 && 8 && 9 && 9 && 4 && 0 &\\
Betti num     &   && 0 && 0 && 0 && 0 && 1 && 0 &&\\
\hline
\end{tabular}
\end{center}
\section{Another proof by Gr\"obner bases} \label{Another:Proof}
Since 
$\displaystyle \HGF{5}{2}{0}{10}$ and 
$\displaystyle \HGF{7}{2}{}{8}$ are both 1-dimensional, 
if $$\displaystyle \omega\wedge : 
\HGF{5}{2}{0}{10} \longrightarrow 
\HGF{7}{2}{}{8}$$ 
is non-zero map, then it is an isomorphism. 

We need to check if 
$\displaystyle \omega\wedge \ker(\mydo) \subset \mydz(
\CGF{6}{2}{}{8})$ or not. For that purpose, choose a basis
$\mathbf{k}_1, \ldots, \mathbf{k}_{8}$ of $\displaystyle \ker(\mydo)$
and linear independent cochains $\mathbf{b}_{1},\ldots,\mathbf{b}_{9}$ in 
$\displaystyle \CGF{6}{2}{}{8}$ such that 
$\mydz(\mathbf{b}_1), \ldots, \mydz(\mathbf{b}_{9})$ is a basis of 
$\displaystyle \mydz( \CGF{6}{2}{}{8})$.  

By taking matrix representation, we see that 
$$\rank( 
\omega\wedge\mathbf{k}_1, \ldots, \omega\wedge\mathbf{k}_{8},
\mydz(\mathbf{b}_1), \ldots, \mydz(\mathbf{b}_{9}) ) =10> 9$$ 
Thus, for an element, say $\mathbf{h}$, which represents the non-trivial
cohomology class, we have to check if $\omega\wedge\mathbf{h}$ is absorbed in
$\displaystyle \mydz( \CGF{6}{2}{}{8})$ or not, namely, if 
$\omega\wedge\mathbf{h}$ realizes the non-trivial cohomology class in 
$\displaystyle \HGF{7}{2}{}{8}$ or not.  
This is our strategy to complete the proof of Theorem. 

\bigskip

Since the both methodologies of using Gr\"obner bases in order to investigate
the cohomology groups $\displaystyle \HGF{5}{2}{0}{10}$ or 
$\displaystyle \HGF{7}{2}{}{8}$ are the same, we discuss in the case 
of $\displaystyle \HGF{5}{2}{0}{10}$ in detail and  write down only the
result for $\displaystyle \HGF{7}{2}{}{8}$.  In particular, we discuss the key
issue where $\omega\wedge$ is involved, carefully.   

Let $\displaystyle \{\mathbf{w}_1, \ldots, \mathbf{w}_{12}\}$ be the  basis
of $\displaystyle C^{5}$ and $\displaystyle \{ \mathbf{q}_1, \ldots,
\mathbf{b}_9\}$ be the basis of $\displaystyle C^{4}$ as before.  
From the matrix representation (\ref{d:one:4and5}) of the coboundary operator 
$\displaystyle \mydo$ of $\displaystyle C^{4} \rightarrow C^{5}$, 
we define the linear functions $$\displaystyle
g_j(y)  = \sum_{k=1}^{12} \lambda_{k\,j} y_k \qquad  (j=1,\ldots,9)$$ where
$\displaystyle ( \lambda_{k\,j}) = M $ and 
$\displaystyle \{ y_1, \ldots, y_{12}\}$ are the auxiliary
variables.

Fixing a monomial order of polynomials induced, say $y_1 \succ  \cdots \succ  y_{12}$, 
we get the Gr\"obner basis $GB_{e}$ of the ideal generated by  
$\displaystyle\{  g_j(y) \mid j=1,\ldots,9\}$.  This corresponds to the
non-zero rows of the elementary matrix of $M$ obtained by the elementary
row operations for $M$.  Thus, 
the cardinality of
$GB_{e}$ is equal to the rank of $M$, namely, 
to $\displaystyle \dim(\mydo( C^{4}))$ and 
$\displaystyle \{ \widehat{g}(\mathbf{w}) \mid \widehat{g} \in GB_{e}\}$
gives a basis of $\displaystyle \mydo( C^{4})$  
(cf.\ Proposition 3.1 in \cite{GB:J:Merker}).  
In our case, 
\begin{align*}  GB_{e} = 
[\; &
21 y_{7}-9 y_{8}-18 y_{9}-15 y_{10}+30 y_{11}-140 y_{12},\\&  
18 y_{6}+9 y_{8}+15 y_{10}-30 y_{11}+140 y_{12},\\& 
1512 y_{5}+75 y_{8}-900 y_{9}-666 y_{10}-1461 y_{11}+3290 y_{12},\\& 
36 y_{4}-3 y_{8}+36 y_{9}-18 y_{10}+57 y_{11}-770 y_{12},\\& 
72 y_{3}+3 y_{8}-36 y_{9}-18 y_{10}+15 y_{11}-70 y_{12},\\& 
63 y_{2}-327 y_{8}+396 y_{9}-258 y_{10}-660 y_{11}+3080 y_{12},\\&  
189 y_{1}-12 y_{8}+144 y_{9}+99 y_{10}+390 y_{11}-1820 y_{12}
\; ]
\end{align*} 
In general, the 
normal form of a given polynomial $g$ with respect to the Gr\"obner basis
is the ``smallest'' remainder of $g$ modulo by the Gr\"obner basis.  

For a linear function $L(y)$ of $\displaystyle y_1,\ldots, y_{12}$,
that 
$\displaystyle L(\mathbf{w})$ belongs to $\displaystyle \mydo( C^{4}) $ is equivalent to the
normal form of $L(y)$ with respect to $GB_{e}$ is zero.  

\bigskip

Let  
$\displaystyle \{\mathbf{r}_1, \mathbf{r}_2, \mathbf{r}_3,\mathbf{r}_{4}\}$ 
be the basis of 
$\displaystyle C^{6}$ as before.     
The kernel space of $\displaystyle \mydo : C^{5} \longrightarrow C^{6}$,
whose element is given by $\displaystyle \sum_{j=1}^{12} c_j \mathbf{w}_j $ satisfying 
$\displaystyle \sum_{j=1}^{12} c_j \mydo(\mathbf{w}_j) =\mathbf{0}$, is characterized by 4 linear functions,  
say $\displaystyle f_1(c), f_2(c), f_3(c), f_4(c)$ of $c_1,\ldots, c_{12}$ given by 
$$  [    
f_1(c), f_2(c), f_3(c), f_4(c)] = [c_1,\ldots, c_{12}]\, {}^t N $$ 
where $N$ is the matrix representing the operator $\displaystyle \mydo:
C^{5}\rightarrow C^{6}$  (This means we deal with the dual map 
$\ds\mydo^{*}:\left(C^{5}\right)^{*} \leftarrow \left(C^{6}\right)^{*}$).     
In our case, 
\begin{align*}
f_{
1} & =140 c_{2}-15 c_{6}+15 c_{7}+30 c_{8}+\frac{5}{2} c_{9} \\
f_{2}
& =-5 c_{1}-4 c_{2}+\frac{1}{4} c_{3}-\frac{11}{2} c_{4}+\frac{31}{12}
c_{5}+\frac{31}{6} c_{6}-3 c_{7}-2 c_{8} +\frac{5}{3} c_{9}-c_{10}+2
c_{11}\\
f_{3}
&=-16 c_{1}+32 c_{2}-2 c_{3}-12 c_{4}+\frac{22}{3} c_{5}+\frac{58}{3} c_{6}
      -18 c_{7}-12 c_{8}-\frac{5}{3} c_{9} +8 c_{11}\\
f_{4}
& =42 c_{4}+7 c_{5}+14 c_{11}+3 c_{12} 
\end{align*}

By taking a monomial order, say   
$c_1 \succ  \cdots \succ  c_{12}$, we get the Gr\"obner basis $GB$ of the ideal  
$\displaystyle\langle f_1(c), f_2(c),
f_3(c), f_4(c)\rangle $. 
In our case, 
\begin{align*} 
GB = [\:  & 42c_{4}+7c_{5}+14c_{11}+3c_{12}, \\&   
42c_{3}+28c_{5}-114c_{6}+198c_{7}+228c_{8}+117c_{9}-48c_{10}+4c_{11}+6c_{12},
\\& 
56c_{2}-6c_{6}+6c_{7}+12c_{8}+c_{9},\\&  
168c_{1}-112c_{5}-182c_{6}+126c_{7}+84c_{8}-35c_{9}+24c_{10}-128c_{11}-12c_{12}
\: ] 
\end{align*}
The $GB$ gives a basis of the subspace 
$\left(\ds\mydo^{*}:\left(C^{5}\right)^{*} \leftarrow \left(C^{6}\right)^{*}
\right)
\left( (C^{6})^{*} \right)$.

Consider the polynomial $\displaystyle h = \sum_{j=1}^{12} c_j y_j$
where $\displaystyle \{ y_1, \ldots, y_{12}\}$ are the other   
auxiliary variables. 

Proposition 3.3 in \cite{GB:J:Merker}
says that the normal form of $h$ with respect to the
Gr\"obner basis  $GB$ is written as $\displaystyle \sum_{j\in J} c_j
\tilde{f}_j (y)$ where $J$ is a subset of $\{1,2,\ldots, 12\}$,
$\displaystyle  \tilde{f}_j (y)$ is linear in $\displaystyle \{ y_1,
\ldots, y_{12}\}$, the cardinality of $J$ is $\displaystyle \dim \ker(
\mydo)$, and $\displaystyle \{ \tilde{f}_j (\mathbf{w}) \mid j\in J\}$
gives a basis of $\displaystyle \ker(\mydo)$.  We
continue the discussion in our case, then we have
\begin{alignat*}{2}
\tilde{f}_1 = & 0,  \quad \tilde{f}_2 =  0, \quad 
\tilde{f}_3 =   0,  \quad \tilde{f}_4 =  0, \\
\tilde{f}_5 = & 
 \frac{2}{3} y_{1}- \frac{2}{3} y_{3}- \frac{1}{6} y_{4}+y_{5}, &  
\tilde{f}_6 = & 
\frac {13}{12} y_{1}+\frac {3}{28} y_{2}+\frac {19}{7} y_{3} +y_{6},\\
\tilde{f}_7 = & 
-\frac{3}{4} y_{1}-\frac {3}{28} y_{2}-\frac {33}{7} y_{3}+y_{7}, & 
\tilde{f}_8 = & 
 -\frac{1}{2} y_{1}-\frac{3}{14} y_{2}-\frac{38}{7} y_{3}+y_{8}, \\
\tilde{f}_9 = & 
 \frac {5}{24} y_{1}-\frac {1}{56} y_{2} -\frac {39}{14} y_{3}+y_{9}, & 
\tilde{f}_{10} = & 
-\frac{1}{7} y_{1}+\frac {8}{7} y_{3}+y_{10}, \\
\tilde{f}_{11} = & 
 \frac {16}{21} y_{1}-\frac{2}{21} y_{3}-\frac{1}{3} y_{4}+y_{11}, & 
\tilde{f}_{12} = & 
 \frac{1}{14} y_{1}- \frac{1}{7} y_{3}- \frac{1}{14} y_{4}+y_{12}
\end{alignat*}
Again, fixing the monomial order of $\{y_j\}$, 
we get the Gr\"obner basis  $\displaystyle GB_k$ 
of the ideal generated by  
$\displaystyle 
\tilde{f}_{j}$ ($j\in J$) as  
\begin{align*}
     GB_{k} 
= [\; & 3y_{{8}}-36y_{{9}}-72y_{{10}}-3y_{{11}}+14y_{{12}}, \quad 
3y_{{7}}-18y_{{9}}-33y_{{10}}+3y_{{11}}-14y_{{12}}, \\& 
18y_{{6}}+108y_{{9}}+231y_{{10}}-21y_{{11}}+98y_{{12}}, \quad  
36y_{{5}}+27y_{{10}}-33y_{{11}}+70y_{{12}}, \\&
2y_{{4}}-5y_{{10}}+3y_{{11}}-42y_{{12}}, \quad  
12y_{{3}}+9y_{{10}}+3y_{{11}}-14y_{{12}}, \\&  
9y_{{2}}-504y_{{9}}-1158y_{{10}}-141y_{{11}}+658y_{{12}}, \quad 
3y_{{1}}-3y_{{10}}+6y_{{11}}-28y_{{12}}\; ] 
\end{align*} 
in our case.  

The cohomology $\displaystyle \HGF{5}{2}{0}{10}$ 
corresponds to the Gr\"obner basis  
$\displaystyle GB_{k/e}$  
of the ideal generated by the normal form of $\displaystyle \widehat{g} \in GB_{k}$ with respect to 
$\displaystyle GB_{e}$.  
In our case, this is given by 
\begin{equation} \displaystyle GB_{k/e} = 
[\; 3 y_{8}-36 y_{9}-72 y_{10}-3 y_{11}+14 y_{12}\; ] \label{gb:kmode}
\end{equation}   
 
\bigskip

\paragraph{$\displaystyle \HGF{7}{2}{}{8}$ case:} 
In the case of 
$\displaystyle \HGF{7}{2}{}{8}$, we use the notations 
$\displaystyle \overline{GB}_{k}$,  
$\displaystyle \overline{GB}_{e}$ and   
$\displaystyle \overline{GB}_{k/e}$ for the Gr\"obner bases
corresponding to the kernel, $\mydz$-image and     
$\displaystyle \HGF{7}{2}{}{8}$ respectively.    
The space $\displaystyle \mydz( \frakC^{6} )$ is characterized by the
following Gr\"obner basis:
\begin{align*}
\overline{GB}_{e}= [\; & 
3 y_{{10}}-3 y_{{11}}-20 y_{{12}}+6 y_{{14}},
\\&
100 y_{{8}}+36 y_{{9}}-15 y_{{11}}-420 y_{{12}}-420 y_{{13}}+350
y_{{14}},\\&
300 y_{{7}}+84 y_{{9}}-135 y_{{11}}-980 y_{{12}}+420 y_{{13}}+350
y_{{14}},\\&
100 y_{{6}}+204 y_{{9}}-135 y_{{11}}-1380 y_{{12}}-180 y_{{13}}+750 y_{{14}},
\\& 
40 y_{{5}}-12 y_{{9}}+15 y_{{11}}-460 y_{{12}}-60 y_{{13}}-590
y_{{14}},\\&
4800 y_{{4}}+84 y_{{9}}+2565 y_{{11}}+6020 y_{{12}}+420 y_{{13}}-10850
y_{{14}},\\&
1600 y_{{3}}-84 y_{{9}}+1035 y_{{11}}-6020 y_{{12}}-420 y_{{13}}-5950
y_{{14}},\\&
400 y_{{2}}-12 y_{{9}}-195 y_{{11}}-1860 y_{{12}}-5660 y_{{13}}+950
y_{{14}},\\&
450 y_{{1}}+24 y_{{9}}+315 y_{{11}}+220 y_{{12}}+120 y_{{13}}+1250
y_{{14}}\; ]
\end{align*}
The kernel space of $\displaystyle \mydz : \frakC^{7} \rightarrow
\frakC^{8}$ is generated by 
\begin{align*}
f_1 =&-35 c_{2}-30 c_{6}-15 c_{8}+\frac{25}{2} c_{9}-\frac{5}{2}
c_{13}\\ 
f_2 =&11 c_{1}-9 c_{2}-\frac{39}{8} c_{3}-\frac{31}{8}
c_{4}-\frac{61}{2} c_{5}-75 c_{6} -10 c_{7} \\& -10 c_{8} 
+\frac{85}{2}
c_{9}-\frac{2}{3} c_{10}
-\frac{20}{3} c_{11}-c_{13}-3 c_{14}\\ 
f_3 =&-16 c_{1}+8 c_{2}+9 c_{3}+5 c_{4}+52 c_{5}+20 c_{6}+\frac{8}{3}
c_{7}\\& -2 c_{8} -\frac{55}{3} c_{9}
+8 c_{11}+c_{13}+4 c_{14}\\ 
f_4 =&\frac{63}{2} c_{3}+\frac{21}{2} c_{4}+84 c_{5}-14 c_{11}+3 c_{12}+3 c_{14}
\end{align*}
and the 
kernel space of $\displaystyle \mydz : \frakC^{7} \rightarrow
\frakC^{8}$ is characterized by the following Gr\"obner basis.  
\begin{align*} 
\overline{GB}_{k} = [\; & 
3 y_{{10}}-3 y_{{11}}-20 y_{{12}}+6 y_{{14}}, \quad  
12 y_{{9}}+495 y_{{11}}+3260 y_{{12}} +60 y_{{13}}-950 y_{{14}}, \\&  
y_{{8}}-15 y_{{11}}-102 y_{{12}}-6 y_{{13}}+32 y_{{14}}, \quad 
3 y_{{7}}-36 y_{{11}}-238 y_{{12}}+70 y_{{14}}, \\&
2 y_{{6}}-171 y_{{11}}-1136 y_{{12}}-24 y_{{13}}+338 y_{{14}}, \quad 
4 y_{{5}}+51 y_{{11}}+280 y_{{12}}-154 y_{{14}}, \\&
16 y_{{4}}-3 y_{{11}}-56 y_{{12}}-14 y_{{14}}, \quad  
16 y_{{3}}+45 y_{{11}}+168 y_{{12}}-126 y_{{14}}, \\&
4 y_{{2}}+3 y_{{11}}+14 y_{{12}}-56 y_{{13}}, \quad
2 y_{{1}}-3 y_{{11}}-28 y_{{12}}+14 y_{{14}}\; ] \\[2mm]
\noalign{thus, $\displaystyle \HGF{7}{2}{}{8}$ is characterized by }
\overline{GB}_{k/e} =[\; & 12 y_{{9}}+495 y_{{11}}+3260 y_{{12}}+60
y_{{13}}-950 y_{{14}}\; ]
\end{align*}

\medskip
Take $ h(y) = 
3 y_{8}-36 y_{9}-72 y_{10}-3 y_{11}+14 y_{12}$ from 
$\displaystyle GB_{k/e}$ of (\ref{gb:kmode}).  Now   
$ h(\mathbf{w}) $ is in $\displaystyle \ker( \mydo : C^{5} \rightarrow
C^{6} ) \setminus \mydo( C^{4})$. 
We express the following element $$\displaystyle \omega \wedge h(\mathbf{w}) = \nz{0}{1} \wedge \nz{1}{1}
\wedge  h(\mathbf{w})$$ by the basis of $\displaystyle \frakC^{7}$.   
We see that 
$$\omega \wedge h(\mathbf{w}) = \nz{0}{1} \wedge \nz{1}{1}
\wedge  h(\mathbf{w}) = 
-9 \overline{\mathbf{w}}_{7} + 
105 \overline{\mathbf{w}}_{10} + 
3 \overline{\mathbf{w}}_{11} + 
14 \overline{\mathbf{w}}_{12}  = \overline{ h} ( \overline{\mathbf{w}}) $$ 
where  
$\displaystyle \overline{h} = 
-9 y_{7} + 
105 y_{10} + 
3 y_{11} + 
14 y_{12} $.   
The normal form of $\displaystyle \overline{h} $ 
with respect to 
$\displaystyle \overline{GB}_{e}$ is   
$$
 \frac{63}{25}y_{9}+\frac{ 2079}{20}y_{11}+\frac{ 3423}{5}y_{12}+\frac{ 63}{5}y_{13}-\frac{399}{2}y_{14}
$$ and is not zero. This finishes the proof of the Theorem.  \kmqed

\bigskip

\begin{kmRemark} \label{rem:stress} 
We emphasize that 
everything starts from 
the concrete bases of cochain complexes  
$\ds\CGF{4}{2}{0}{10}$, $\ds\CGF{5}{2}{0}{10}$, 
$\ds\CGF{6}{2}{0}{10}$, $\ds\CGF{6}{2}{}{8}$ and \\
$\ds\CGF{7}{2}{}{8}$.     

Even though we make use of 
Gr\"obner Base theory or use of classical linear algebra argument, we are
based on 
some concrete matrix representations. 
\end{kmRemark}



\def\cprime{$'$} \def\cprime{$'$}

\vspace{1cm}

{Akita University\\
1-1 Tegata Akita City, Japan\\
mikami@math.akita-u.ac.jp}


\begin{appendices} 
        \section{} \label{t1:w10:betti}
In this Appendix, we make use of Risa/Asir, which is another
Symbol Calculus Software, and show the results we got by Maple and
Risa/Asir are the same up to non-zero scalar multiples.  

We remark that we added some line breaks so that we get better look.

\paragraph{Basis of $\ds\mydo( \CGF{4}{2}{0}{10}) \subset
\CGF{5}{2}{0}{10}$:}\ \\ 
Our source file for Risa/Asir is this:
\begin{rmfamily}
\begin{verbatim}
/* #####  On C^{4} -> C^{5} #### */ 
G1 = -135/4*y1-60*y3+15/2*y4-45*y5-15*y6+5/4*y7-45/4*y8+75/2*y9$
G2 = 108/11*y1+18/11*y2+60/11*y6+46/11*y7-90/11*y8+156/11*y9$
G3 = 27/4*y1+12*y3-9/2*y4-9*y5+27/4*y10+18*y11$
G4 = -10*y3+2/3*y4-2*y5+2*y6+y7+6*y9+4*y10-y11$
G5 = 5/2*y2+29*y3+47/3*y4-23*y5+43*y6+13/2*y7+9/2*y8+25*y9+16*y10-71/2*y11$
G6 = 5*y2+45*y3+155/6*y4-40*y5+65*y6+10*y7+50*y9+20*y10-115/2 *y11$
G7 = 3/2*y2+18*y3+23/2*y4-3*y5+30*y6+11/2*y7+9/2*y8+9*y9+6* y10-33*y11$
G8 = 6*y6+7*y7-6*y9$        G9 = -6*y3-3*y4+3*y10-6*y11+70*y12$ 
GBe = gr([G1,G2,G3,G4,G5,G6,G7,G8,G9],[y1,y2,y3,y4,y5,y6,y7,y8,y9,y10,y11,y12],1) ;
\end{verbatim}
\end{rmfamily}
\paragraph{The output of Groebner Basis is the next:} 
\begin{verbatim}
[-140*y12+30*y11-15*y10-18*y9-9*y8+21*y7,
 140*y12-30*y11+15*y10+9*y8+18*y6,
 -3290*y12+1461*y11+666*y10+900*y9-75*y8-1512*y5,
 -770*y12+57*y11-18*y10+36*y9-3*y8+36*y4,
 70*y12-15*y11+18*y10+36*y9-3*y8-72*y3,
 3080*y12-660*y11-258*y10+396*y9-327*y8+63*y2,
 -1820*y12+390*y11+99*y10+144*y9-12*y8+189*y1]
\end{verbatim}

\paragraph{Kernel space of $\ds  \mydo : \CGF{5}{2}{0}{10}\rightarrow \CGF{6}{2}{0}{10}$:}
\ \\
Our source file for Risa/Asir is this:
\begin{verbatim} 
/* #####  On C^{5} -> C^{6} #### */ 
F1 = -5*w2-16*w3$               F2 = 140*w1-4*w2+32*w3$
F3 = 1/4*w2-2*w3$               F4 = -11/2*w2-12*w3+42*w4$
F5 = 31/12*w2+22/3*w3+7*w4$     F6 = -15*w1+31/6*w2+58/3*w3$
F7 = 15*w1-3*w2-18*w3$          F8 = 30*w1-2*w2-12*w3$
F9 = 5/2*w1+5/3*w2-5/3*w3$      F10 = -w2$
F11 = 2*w2+8*w3+14*w4$          F12 = 3*w4$ 
FList = [F1,F2,F3,F4,F5,F6,F7,F8,F9,F10,F11,F12]$   WList = [w1,w2,w3,w4]$ 
/* ########################################################## */ 
NagawaW = length(WList)$  NagasaF = length(FList)$ 
CC = [c1,c2,c3,c4,c5,c6,c7,c8,c9,c10,c11,c12]$ 
YY = [y1,y2,y3,y4,y5,y6,y7,y8,y9,y10,y11,y12]$ 
for ( Uke = [], J=1; J <= NagawaW; J++ ) { MyA = WList[J-1];  Atai = 0;
        for (K=1 ; K <= NagasaF; K++ ){ MyB = FList[K-1]; 
            Atai += diff( MyB, MyA)* CC[K-1];}; Uke = cons(Atai, Uke ); 
}
print("mark A")$    Uke = reverse( Uke ); 
print("mark B")$    GBadj = gr( Uke, CC, 0); /* Groebner Basis */
for (H=0, I=1; I <= NagasaF; I++){ H += CC[I-1]* YList[I-1]; } 
Hnf = p_nf(H, GBadj, CC, 0)$  /* Normal Form */ 
for( MyUkez = [], T=CC; T != []; T = cdr(T)){ 
    MyA = car(T);   MyV = diff( Hnf, MyA); MyUkez = cons( MyV, MyUkez);} 
print("mark C")$    MyUkez = reverse(MyUkez); 
print("mark D")$    GBk = gr( MyUkez, YY, 0); /* Groebner Basis */
end$
\end{verbatim}

\paragraph{The outputs are the follows:}
\begin{verbatim}
mark A
[140*c2-15*c6+15*c7+30*c8+5/2*c9,
-5*c1-4*c2+1/4*c3-11/2*c4+31/12*c5+31/6*c6-3*c7-2*c8+5/3*c9-c10+2*c11,
-16*c1+32*c2-2*c3-12*c4+22/3*c5+58/3*c6-18*c7-12*c8-5/3*c9+8*c11,
42*c4+7*c5+14*c11+3*c12] 
mark B
[42*c4+7*c5+14*c11+3*c12,
-42*c3-28*c5+114*c6-198*c7-228*c8-117*c9+48*c10-4*c11-6*c12,
56*c2-6*c6+6*c7+12*c8+c9,
168*c1-112*c5-182*c6+126*c7+84*c8-35*c9+24*c10-128*c11-12*c12] 
mark C
[0, 0, 0, 0,
-168*y5+28*y4+112*y3-112*y1,    -168*y6-456*y3-18*y2-182*y1,
-168*y7+792*y3+18*y2+126*y1,    -168*y8+912*y3+36*y2+84*y1,
-168*y9+468*y3+3*y2-35*y1,      -168*y10-192*y3+24*y1,
-168*y11+56*y4+16*y3-128*y1,    -168*y12+12*y4+24*y3-12*y1] 
mark D
[14*y12-3*y11-72*y10-36*y9+3*y8,        14*y12-3*y11+33*y10+18*y9-3*y7,
98*y12-21*y11+231*y10+108*y9+18*y6,     70*y12-33*y11+27*y10+36*y5,
42*y12-3*y11+5*y10-2*y4,                14*y12-3*y11-9*y10-12*y3,
-658*y12+141*y11+1158*y10+504*y9-9*y2,  28*y12-6*y11+3*y10-3*y1]
\end{verbatim}

\paragraph{Basis of $\ds\HGF{5}{2}{0}{10}$}\ \\   
The next is a source file for Risa/Asir. GBe and GBk are data gotten above. 

\begin{verbatim}
GBe = [-140*y12+30*y11-15*y10-18*y9-9*y8+21*y7, 
        140*y12-30*y11+15*y10+9*y8+18*y6, 
       -3290*y12+1461*y11+666*y10+900*y9-75*y8-1512*y5,
       -770*y12+57*y11-18*y10+36*y9-3*y8+36*y4, 
        70*y12-15*y11+18*y10+36*y9-3*y8-72*y3,
        3080*y12-660*y11-258*y10+396*y9-327*y8+63*y2,
       -1820*y12+390*y11+99*y10+144*y9-12*y8+189*y1]$ 
GBk = [ 14*y12-3*y11-72*y10-36*y9+3*y8,         14*y12-3*y11+33*y10+18*y9-3*y7,
        98*y12-21*y11+231*y10+108*y9+18*y6,     70*y12-33*y11+27*y10+36*y5,
        42*y12-3*y11+5*y10-2*y4,                14*y12-3*y11-9*y10-12*y3,
       -658*y12+141*y11+1158*y10+504*y9-9*y2,   28*y12-6*y11+3*y10-3*y1]$ 
YY = [y1,y2,y3,y4,y5,y6,y7,y8,y9,y10,y11,y12]$

for(Uke=[], T= GBk; T != []; T = cdr(T)) {
    MyA = car(T); Atai = p_nf( MyA, GBe, YY , 0) ;  /* NormalForm */
    Uke = cons(Atai, Uke); }
Uke = reverse(Uke)$ 
GBh = gr( Uke, YY , 0);  /* Groebner Basis */
\end{verbatim}
\paragraph{A basis of $\ds\HGF{5}{2}{0}{10}$ is given by the output}
\begin{verbatim}
[-14*y12+3*y11+72*y10+36*y9-3*y8] 
\end{verbatim}

\paragraph{Check $\ds 
  \mydo \circ \mydo : \CGF{4}{2}{0}{10}\rightarrow \CGF{6}{2}{0}{10}$ is zero
  identically:}\ \\   
The next is a source file for Risa/Asir.  
\begin{verbatim}
/* Feb 07, 2014 n=1, type 1, weight 10, C^{4} --> C^{5} --> C^{6} */ 
G1 = -135/4*y1-60*y3+15/2*y4-45*y5-15*y6+5/4*y7-45/4*y8+75/2*y9$
G2 = 108/11*y1+18/11*y2+60/11*y6+46/11*y7-90/11*y8+156/11*y9$
G3 = 27/4*y1+12*y3-9/2*y4-9*y5+27/4*y10+18*y11$
G4 = -10*y3+2/3*y4-2*y5+2*y6+y7+6*y9+4*y10-y11$
G5 = 5/2*y2+29*y3+47/3*y4-23*y5+43*y6+13/2*y7+9/2*y8+25*y9+16 *y10-71/2*y11$
G6 = 5*y2+45*y3+155/6*y4-40*y5+65*y6+10*y7+50*y9+20*y10-115/2 *y11$
G7 = 3/2*y2+18*y3+23/2*y4-3*y5+30*y6+11/2*y7+9/2*y8+9*y9+6* y10-33*y11$
G8 = 6*y6+7*y7-6*y9$
G9 = -6*y3-3*y4+3*y10-6*y11+70*y12$ 
/* The next data are gotten by replacing  y to F and G to GG  in the above.  */ 
GG1 = -135/4*F1-60*F3+15/2*F4-45*F5-15*F6+5/4*F7-45/4*F8+75/2*F9$
GG2 = 108/11*F1+18/11*F2+60/11*F6+46/11*F7-90/11*F8+156/11*F9$
GG3 = 27/4*F1+12*F3-9/2*F4-9*F5+27/4*F10+18*F11$
GG4 = -10*F3+2/3*F4-2*F5+2*F6+F7+6*F9+4*F10-F11$
GG5 = 5/2*F2+29*F3+47/3*F4-23*F5+43*F6+13/2*F7+9/2*F8+25*F9+16 *F10-71/2*F11$
GG6 = 5*F2+45*F3+155/6*F4-40*F5+65*F6+10*F7+50*F9+20*F10-115/2 *F11$
GG7 = 3/2*F2+18*F3+23/2*F4-3*F5+30*F6+11/2*F7+9/2*F8+9*F9+6* F10-33*F11$
GG8 = 6*F6+7*F7-6*F9$
GG9 = -6*F3-3*F4+3*F10-6*F11+70*F12$ 
/* ### On C^{5} -> C^{6}: */
F1 = -5*w2-16*w3$               F2 = 140*w1-4*w2+32*w3$
F3 = 1/4*w2-2*w3$               F4 = -11/2*w2-12*w3+42*w4$
F5 = 31/12*w2+22/3*w3+7*w4$     F6 = -15*w1+31/6*w2+58/3*w3$
F7 = 15*w1-3*w2-18*w3$          F8 = 30*w1-2*w2-12*w3$
F9 = 5/2*w1+5/3*w2-5/3*w3$      F10 = -w2$
F11 = 2*w2+8*w3+14*w4$          F12 = 3*w4$
/* ########################################## */ 
L = [GG1,GG2,GG3,GG4,GG5,GG6,GG7,GG8,GG9]$      print(L)$   end$ 
\end{verbatim}

\paragraph{The output is as expectd.} 
\begin{verbatim}
[0,0,0,0,0,0,0,0,0]
\end{verbatim}

\section{} \label{t0:w8:betti}
In section  \ref{Another:Proof}, we have shown 
a basis of $\ds \HGF{7}{2}{ }{8})$ concretely by Groebner Basis Package of
Maple.   

In this Appendix, we do the same jobi by 
using Risa/Asir, which is another Symbol Calculus Software, and show that  
the results we got by Maple and Risa/Asir are the same up to non-zero scalar multiples.  

We remark that 
in this note  we added some line breaks so that we get better look and 
we use \texttt{nd\_gr()} instead of \texttt{gr()}.

We stock two matrix representations of $\ds\mydz$ in the two files: 

\medskip

\textbf{Mat\_w8\_6and7\_type0.rr}  
\begin{small}
\verbatiminput{Mat_w8_6and7_type0.rr}
\end{small}

\textbf{Mat\_w8\_7and8\_type0.rr}.    

\begin{small}
\verbatiminput{Mat_w8_7and8_type0.rr}
\end{small}

\paragraph{Basis of $\ds\mydz( \CGF{6}{2}{ }{8}) \subset
\CGF{7}{2}{ }{8}$:}\ \\ 
Our source file for Risa/Asir is this:
\begin{verbatim} 
load("./Mat_w8_6and7_type0.rr")$ /* GB1 = gr( [ ], [ ], 0) $ */ 
ord( YList ) $ 
GB1 = reverse( 
nd_gr(GList , YList ,0,0)) $ 
print(["GBe",GB1])$ 
end$
\end{verbatim}

\paragraph{The output of Groebner Basis is the next:} 
\begin{verbatim} 
[GBe, [3*y10-3*y11-20*y12+6*y14, 100*y8+36*y9-15*y11-420*y12-420*y13+350*y14,
      -300*y7-84*y9+135*y11+980*y12-420*y13-350*y14,
       100*y6+204*y9-135*y11-1380*y12-180*y13+750*y14,
       40*y5-12*y9+15*y11-460*y12-60*y13-590*y14,
       4800*y4+84*y9+2565*y11+6020*y12+420*y13-10850*y14,
       1600*y3-84*y9+1035*y11-6020*y12-420*y13-5950*y14,
       400*y2-12*y9-195*y11-1860*y12-5660*y13+950*y14,
      -450*y1-24*y9-315*y11-220*y12-120*y13-1250*y14]] 
\end{verbatim}

\paragraph{Kernel space of $\ds  \mydz : \CGF{7}{2}{ }{8}\rightarrow \CGF{8}{2}{ }{8}$:}
\ \\
Our source file for Risa/Asir is this:
\begin{verbatim} 
load("Mat_w8_7and8_type0.rr")$  NagasaW = length(WList)$  NagasaF = length(FList)$ 
for ( Uke = [], J=1; J <= NagasaW; J++ ) { MyA = WList[J-1];  Atai = 0;
    for (K=1 ; K <= NagasaF; K++ ){ MyB = FList[K-1]; 
    Atai += diff( MyB, MyA)* CList[K-1];}
    Uke = cons(Atai, Uke ); }
print("mark A")$    Uke = reverse( Uke ); 
print("mark B")$    GBadj = nd_gr( Uke, CList, 0, 0 ); 
for (H=0, I=1; I <= NagasaF; I++){ H += CList[I-1]* YList[I-1]; } 
Hnf = p_nf(H, GBadj, CList, 0)$ 
for( MyUkez = [], T=CList; T != []; T = cdr(T)){ 
    MyA = car(T); MyV = diff( Hnf, MyA);
    MyUkez = cons( MyV, MyUkez);} 
print("mark C")$    MyUkez = reverse(MyUkez); 
ord(YList)$ 
print("mark D")$ GBk = reverse( nd_gr( MyUkez, YList, 0, 0 )); 
end$ 
\end{verbatim}

\paragraph{The outputs are the follows:}
\begin{verbatim} mark A
[-35*c2-30*c6-15*c8+25/2*c9-5/2*c13,
11*c1-9*c2-39/8*c3-31/8*c4-61/2*c5-75*c6-10*c7
    -10*c8+85/2*c9-2/3*c10-20/3*c11-c13-3*c14,
-16*c1+8*c2+9*c3+5*c4+52*c5+20*c6+8/3*c7-2*c8-55/3*c9+8*c11+c13+4*c14,
63/2*c3+21/2*c4+84*c5-14*c11+3*c12+3*c14]
mark B
[-168*c1+336*c5-1260*c6-168*c7-294*c8+525*c9-16*c10+112*c11-12*c12+3*c13+24*c14,
-14*c2-12*c6-6*c8+5*c9-c13,
-126*c3-420*c5+2796*c6+392*c7+474*c8-1375*c9+32*c10+84*c11-6*c12+3*c13+6*c14,
-42*c4+84*c5-2796*c6-392*c7-474*c8+1375*c9-32*c10-28*c11-6*c12-3*c13-18*c14]
mark C
[0, 0, 0, 0,
1008*y1-1680*y3+1008*y4+504*y5,
-3780*y1-432*y2+11184*y3-33552*y4+504*y6,
-504*y1+1568*y3-4704*y4+504*y7,
-882*y1-216*y2+1896*y3-5688*y4+504*y8,
1575*y1+180*y2-5500*y3+16500*y4+504*y9,
-48*y1+128*y3-384*y4+504*y10,
336*y1+336*y3-336*y4+504*y11,
-36*y1-24*y3-72*y4+504*y12,
9*y1-36*y2+12*y3-36*y4+504*y13,
72*y1+24*y3-216*y4+504*y14]
mark D
[-3*y10+3*y11+20*y12-6*y14, 
-12*y9-495*y11-3260*y12-60*y13+950*y14,
y8-15*y11-102*y12-6*y13+32*y14,
3*y7-36*y11-238*y12+70*y14,
-2*y6+171*y11+1136*y12+24*y13-338*y14,
4*y5+51*y11+280*y12-154*y14,
16*y4-3*y11-56*y12-14*y14,
16*y3+45*y11+168*y12-126*y14,
4*y2+3*y11+14*y12-56*y13,
-2*y1+3*y11+28*y12-14*y14] 
\end{verbatim}

\paragraph{Basis of $\ds\HGF{7}{2}{ }{8}$}\ \\   
The next is a source file for Risa/Asir. GBe and GBk are data gotten above. 

\begin{verbatim} 
GBe = [ 3*y10-3*y11-20*y12+6*y14,
        100*y8+36*y9-15*y11-420*y12-420*y13+350*y14,
        -300*y7-84*y9+135*y11+980*y12-420*y13-350*y14,
        100*y6+204*y9-135*y11-1380*y12-180*y13+750*y14,
        40*y5-12*y9+15*y11-460*y12-60*y13-590*y14,
        4800*y4+84*y9+2565*y11+6020*y12+420*y13-10850*y14,
        1600*y3-84*y9+1035*y11-6020*y12-420*y13-5950*y14,
        400*y2-12*y9-195*y11-1860*y12-5660*y13+950*y14,
        -450*y1-24*y9-315*y11-220*y12-120*y13-1250*y14]$ 
GBk = [ -3*y10+3*y11+20*y12-6*y14,
        -12*y9-495*y11-3260*y12-60*y13+950*y14,
        y8-15*y11-102*y12-6*y13+32*y14,
        3*y7-36*y11-238*y12+70*y14,
        -2*y6+171*y11+1136*y12+24*y13-338*y14,
        4*y5+51*y11+280*y12-154*y14,
        16*y4-3*y11-56*y12-14*y14,
        16*y3+45*y11+168*y12-126*y14,
        4*y2+3*y11+14*y12-56*y13,
        -2*y1+3*y11+28*y12-14*y14]$ 
YList = [y1,y2,y3,y4,y5,y6,y7,y8,y9,y10,y11,y12,y13,y14]$ 
for(Uke=[], T= GBk; T != []; T = cdr(T)) {
    MyA = car(T); /* print(MyA); */
    Atai = p_nf( MyA, GBe, YList , 0) ; Uke = cons(Atai, Uke); } 
Uke = reverse(Uke)$ 
ord(YList)$
GBb = reverse( nd_gr( Uke, YList , 0, 0) );
end$ 
\end{verbatim}
\paragraph{A basis of $\ds\HGF{7}{2}{ }{8}$ is given by the output} 
\begin{verbatim} 
[-12*y9-495*y11-3260*y12-60*y13+950*y14] 
\end{verbatim}

\bigskip

We may omit the job of 
\textbf {Check $\ds 
  \mydz \circ \mydz : \CGF{6}{2}{ }{8}\rightarrow \CGF{8}{2}{ }{8}$ is 
  identically zero}.


        \section{Final stage by Risa/Asir} \label{t0:w8:final}
We have already studied of 
$\ds \HGF{5}{2}{0}{10}$, $\ds GB_{e}$, $\ds GB_{k}$ and  $GB_{k/e}$, 
also of   
$\ds \HGF{7}{2}{ }{8}$, $\ds \overline{GB}_{e}$, $\ds \overline{GB}_{k}$ and
$\overline{GB}_{k/e}$.  

We calculate $\ds \omega \wedge h( \mathbf{w}_j ) $ and we have 
$\ds \overline{h} = -9 y_7 + 105 y_{10} + 3 y_{11} + 14 y_{12}$.  

We will check $\text{NormalForm}( \overline{h}, \overline{GB}_{e},
\text{Ord}_{y})$ does not vanish.  Then a proof to a Theorem in {The
{G}el'fand-{K}alinin-{F}uks class and characteristic classes of transversely
symplectic foliations}, {\em arXiv:0910.3414}, October 2009 by D.~Kotschick
and S.~Morita will be done.  

We remark that 
in this note  we added some line breaks so that we get better look and 
we use \texttt{nd\_gr()} instead of \texttt{gr()}.

\medskip

\paragraph{Final stage:}\ \\ 
Our source file for Risa/Asir is this:
\begin{verbatim} 
YList = [y1,y2,y3,y4,y5,y6,y7,y8,y9,y10,y11,y12,y13,y14]$
ord(YList)$

GBe =[3*y10-3*y11-20*y12+6*y14,
        100*y8+36*y9-15*y11-420*y12-420*y13+350*y14,
        -300*y7-84*y9+135*y11+980*y12-420*y13-350*y14,
        100*y6+204*y9-135*y11-1380*y12-180*y13+750*y14,
        40*y5-12*y9+15*y11-460*y12-60*y13-590*y14,
        4800*y4+84*y9+2565*y11+6020*y12+420*y13-10850*y14,
        1600*y3-84*y9+1035*y11-6020*y12-420*y13-5950*y14,
        400*y2-12*y9-195*y11-1860*y12-5660*y13+950*y14,
        -450*y1-24*y9-315*y11-220*y12-120*y13-1250*y14]$

H = -9 * y7 + 105 * y10 + 3 * y11 + 14* y12$

p_nf(H, GBe, YList, 0);
end$
\end{verbatim}

\paragraph{The output of Groebner Basis is the next:} 
\begin{verbatim} 
-252*y9-10395*y11-68460*y12-1260*y13+19950*y14
\end{verbatim}

\end{appendices}
\end{document}